\documentclass[12pt,amstex]{book}
\usepackage{amssymb,amsmath}
\usepackage{stmaryrd}
\usepackage{amsfonts}
\usepackage{graphics}
\usepackage{rotating}
\usepackage{pspicture}
\usepackage{epic,eepic}
\begin{document}
\title{The Triangle}
\author{Richard K.~Guy}
\date{}
\maketitle

\subsection*{Abstract}
If we label the vertices of a triangle with 1, 2 and 4, and the orthocentre with 7, then any of the
four numbers 1, 2, 4, 7 is the nim-sum of the other three and is their orthocentre. Regard the
triangle as an orthocentric quadrangle. Steiner's theorem states that the reflexions of a point on
a circumcircle in each of the three edges of the corresponding triangle are collinear and collinear
with its orthontre. This line intersects the circumcircles in new points to which the theorem may
be applied. Iteration of this process with the triangle and the points rational leads to a
``trisequence'' whose properties merit study.

\chapter{Preamble}  This might have been titled
{\bf The Triangle Book}, except that John Conway
already has a project in hand for such a book.
Indeed, Conway's book might well have been completed but
for the tragically early death of Steve Sigur.
It might also have been finished, had I been in
closer proximity to John.

                    It is a very badly edited version
of a paper [roughly \S2.1 onwards] that was rejected
by the {\sc Monthly}.  It displays confusion about what
is likely to be the best notation.  It may be thought
of more as a film script than a book --- who's going
to make the film\,?  If John likes to use any of what
follows I would be very flattered.

LATER:  For more completeness I've added Chapters 6 \& 7.
The former is an almost verbatim copy of the Lighthouse
Theorem paper \cite{G}; the latter is an unfinished paper
with John Conway, which corrects some errors in the
former.  They don't always mix very well. There are
repetitions and probably some contradictions!

LATER STILL (2015-08-12) I've added a Chapter 8, which is
a version of the presentation made at Alberta Math Dialog,
Lethbridge. 2015-05-08; MOVES Conference, New York,
2015-08-03; and MathFest, Washington, DC, 2015-08-08.
This further duplicates some things, but enables me to
repair and replace much of the unsatisfactory and incomplete
Chapter 7.

2016-03-12.  Just woken up with the realization that I was
right to title this {\bf The Triangle}, because, instead
of {\bf The Triangle Book}, this should be {\bf The Triangle
Movie}.  Only that way can we have the infinite number
of pictures that we need.

\begin{itemize}
\item  The Wallace line enveloping theS deltoid.
\item  The Droz-Farny line enveloping a conic.
\item  The Lighthouse Theorem generating Morley triangles.
\item  Feuerbach's Theorem: the 11-point conic.
\item  The dual of Feuerbach's Theorem and the Droz-Farny
Theorem.

\end{itemize}

2016-06-08.  See notes on Notation at \S7.6

\section{What is a triangle?}

Two women pursuing the same man? Or, dually, two men
chasing the same woman.
Usually thought of as three points whose pairs
determine three straight lines, or, dually, as
three straight lines whose pairs intersect in three
points.  We shall call the points vertices, and
the lines edges.  More numerically, it can be thought
of as three positive real numbers, $a$, $b$, $c$,
which satisfy the {\bf triangle inequality}
$$b+c>a, \quad c+a> b, \quad a+b > c.$$

\section{What's new?}
What's new(?) here is the
introduction of {\bf Quadration} and {\bf Twinning},
which, together with Conway's
{\bf Extraversion} yield a remarkably general view
of the {\bf Triangle}, sometimes with as many as 32
items for the price of one.  {\bf Eight
vertices}, which are also {\bf circumcentres}, as
well as {\bf orthocentres}.  Note that an orthocentre
is the perspector of the 9-point circle with the relevant
circumcircle.

The triangle now has six pairs of
parallel edges which form three rectangles
whose twelve vertices are the ends of six diameters of
the {\bf Central Circle}.  Also known as the fifty-point circle
The other 38 points are the points of contact with the 32
{\bf touch-circles} and six points of contact with the
{\bf double-deltoid}.  Six of the points are associated with
the name of Euler, 9+32 with that of Feuerbach, and six with
Steiner; the double-deltoid having the appearance of a
Star-of-David, and is the envelope of pairs of parallel
Wallace (Simson) lines, and homothetic to all of the 144 Morley triangles.
There are 32 Gergonne points, 32 Nagel points and 256
radpoints.  But I'm getting ahead of myself.

On the other hand, I shall probably never catch up with myself,
so let me list here some of the topics that I'd like to cover.

\section{What's needed?}

A nice notation!!  The vertices might be $V_1$, $V_2$, $V_4$
with orthocentre $V_7$, circumcentre $V_{\bar7}$, etc.
But see \S7.6.

\section{Clover-leaf theorems}

These usually take the shape of the concurrence of the three
radical axes of pairs of circles, chosen from three, in the radical
centre.

The original Clover-leaf Theorem arose from what I first called
{\bf five-point circles} but which turned out to be
nine-point circles (no, not the traditional nine-point circles)
and which I will re-christen {\bf medial circles}.  The radical
axis of two medial circles of the same triangle is an altitude
of that triangle, so that a medial circle contains a vertex,
a midpoint, a diagonal point, two mid-foot points and four
altitude points.  See Figure \ref{clover} in Chapter 2.

Other examples of what might be called clover-leaf theorems
are proofs (2.) and (3.) of the orthocentre.  Again, see
Chapter 2.  Three edge-circles
through the same vertex, which is then the radical centre,
and orthocentre or fourth vertex.  Three edge-circles of a
triangle as in Figure \ref{proof2}, so that the radical
axes are altitudes again, concurring in the orthocentre.

Clover-leaf theorems may extend into four-leaf clover theorems.
For example, add the circumcircle to Figure \ref{proof2} and
the 6 radical axes are the 6 edges of the (generalized, quadrated)
triangle and the 4 radical centres are the 4 vertices-orthocentres.
Another example is adding the 50-point circle to the medial circles
as in Figure \ref{fourleaf} in Chapter 2, \S2.2.

\section{Radical axes of circumcircles}

There are $\binom{8}{2}=28$ of these.  Twelve are the
12 edges of the (generalized) triangle.  As the circumcircles
are all congruent, the other 16 will be the perpendicular
bisectors of the segments $V_iV_j$ where $i\in\{1,2,4,8\}$,
$j\in\{14,13,11,7\}$, where we may write the latter in
hexadecimal as $j\in\{e,d,b,7\}$.  [Perhaps better is
$j\in\{\bar1,\bar2,\bar4,\bar7\}$.]  Four of these 16 will
be perpendicular bisectors of diametral segments joining
a twin pair of vertices.  These diametral segments are
{\bf Euler lines}.

\section{Touch-circles}

Conway's {\bf extraversion} generalizes the incircle to a set
of 4 {\bf touch-circles} which touch the 3 edges of a triangle.
And {\bf quadration} turns a triangle into an orthocentric
quadrangle of 4 triangles, giving 16 touch-circles.  Finally,
{\bf twinning} doubles this number to 32 touch-circles.
Feuerbach tells us that these all touch the 50-point circle
(Central Circle).

\subsection{Interlude: Proof of Feuerbach's theorem}

(See also Altshiller-Court, \cite[pp.105,273]{C3}; Roger A. Johnson,
\cite[pp.200,244]{J}).

Area $\Delta = \frac{1}{2}ab\sin C = 2R^2\sin A\sin B\sin C$ \\
Semiperimeter $s = \frac{1}{2}(a+b+c) = R(\sin A + \sin B + \sin C)$ \\
Inradius = $$r=\frac{\Delta}{s}=
\frac{2R\sin A\sin B\sin C}{\sin A + \sin B + \sin C}$$

50-point radius = $R/2$.
Square of distance of incentre, $I$, from 50-pt centre is
$$\left[\frac{1}{2}\left(b\sin C-R\cos A\right)-r\right]^2+
\left[\frac{1}{2}\left(\frac{1}{2}a+
b\cos C\right)-(s-c)\right]^2$$
We need this to simplify to $(\frac{1}{2}R-r)^2$.  Not the
right way to go!!

I think we can find a proof which uses the spirit of
quadration (1.6 below), edge-circles and 5 or 6 (now 9)
proofs that a triangle has an orthocentre (Chap.2 below), etc.

In fact, see \$\,8.5 Hexaflexing, below.

[resume {\bf 1.5 Touch-circles}]

These 32 touch-circles touch the 12 edges in 96 touch-points.
Their joins to appropriate vertices concur in triples at
32 {\bf Gergonne points} and 32 {\bf Nagel points}. [No!!
We are not so lucky with Nagel points!]
Gerg{\bf onne} points are obtained by taking the three
touch-points of {\bf onne} (one and only one) touch-circle
and joining them to the respective vertices of the
relevant triangle.  By contrast, the Nagel points are formed
using the touch-points of three (all but one) touch-circles
out of a set of four, [Here's the snag!  We can't use any
old three.  They have to be the three EX-circles.]

\begin{figure}[h]
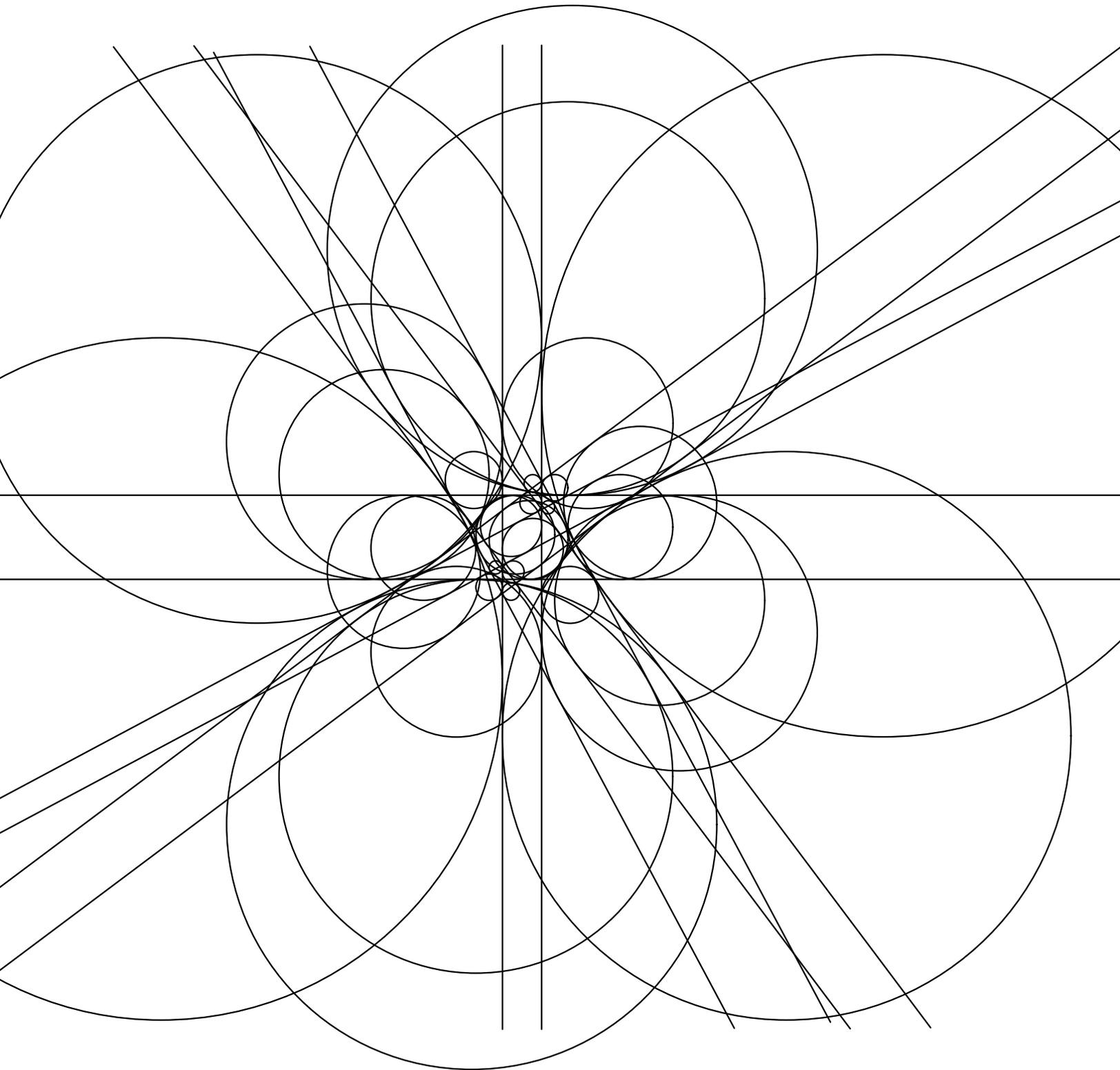
 

\caption{32 touch-circles each touch 3 of the 12 edges at
one of 8 points, and each touches the 50-point circle
(which isn't drawn, but you can ``see'' it!)}
\label{touch}
\end{figure}

\clearpage

\begin{figure}[h]
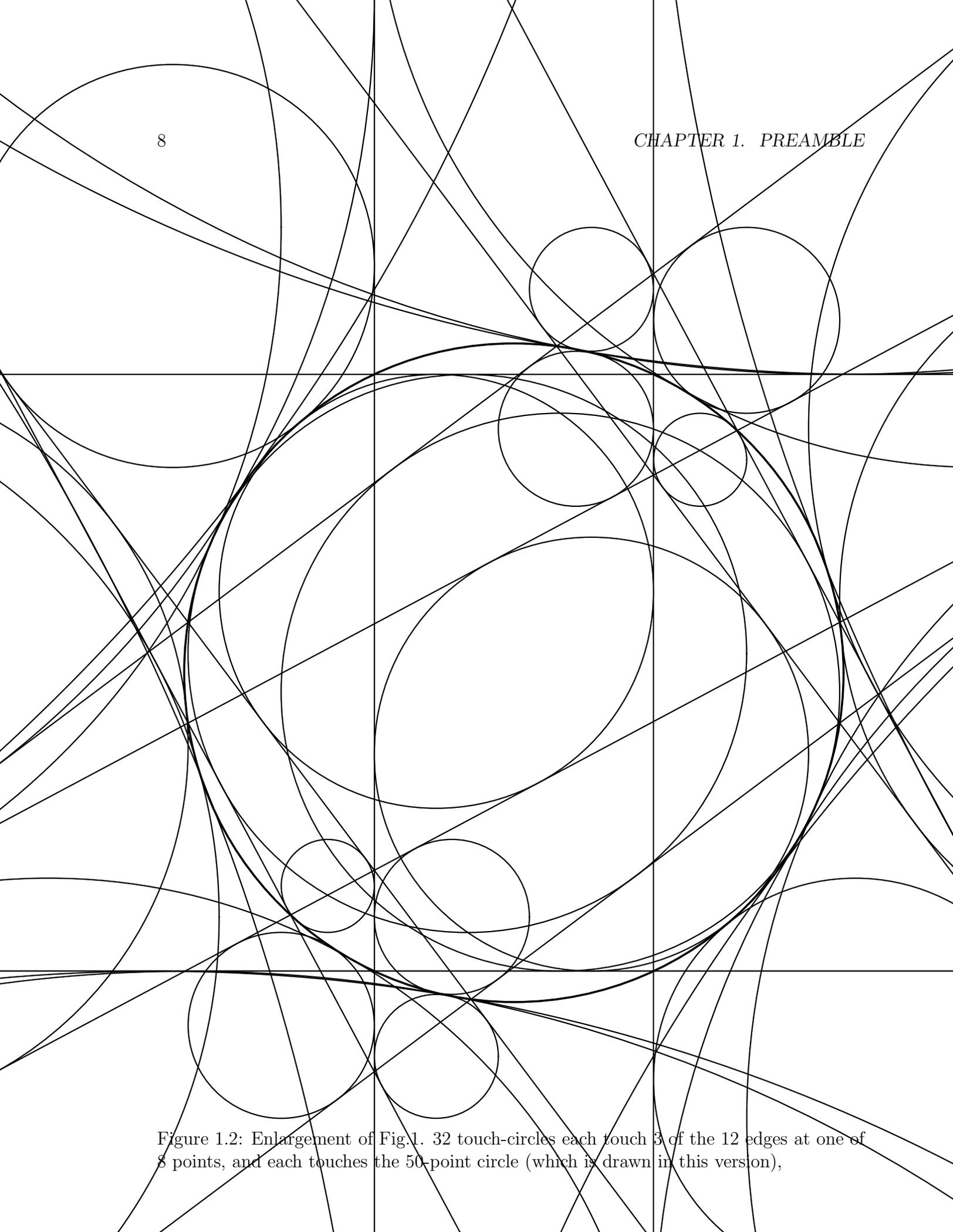
 
%
\caption{Enlargement of Fig.1.  32 touch-circles each
touch 3 of the 12 edges at one of 8 points, and each
touches the 50-point circle (which is drawn in this version),}
\label{touch2}
\end{figure}

\clearpage

\section{Quadration}

This is to grant the same status to the orthocentre as
to the vertices, so that each of the four points is the
orthocentre of the triangle formed by the other three.
That is, to regard the triangle as an {\bf orthocentric
quadrangle}.  It now has 4 vertices, 6 edges and 3
{\bf diagonal points}, $D_6$, $D_5$, $D_3$.

This may be a good place to give some formulas.  We will take
THE (was 50-point) Centre of the triangle as {\bf origin}.
If the vertices of a triangle are given by vectors
{\bf a, b, c}, then the centroid (which we will soon forget!)
is $\frac{1}{3}$({\bf a} + {\bf b} + {\bf c}), the
circumcentre is {\bf a} + {\bf b} + {\bf c}, and the
orthocentre is --({\bf a} + {\bf b} + {\bf c}).  That is,
if {\bf d} is the orthocentre, then
{\bf d} = --({\bf a} + {\bf b} + {\bf c}), and
\centerline{{\bf a} + {\bf b} + {\bf c} + {\bf d} = 0}
reminding us that each of the four points is the orthocentre
of the other three.

Note that ${\bf a}\centerdot {\bf b} + {\bf c}\centerdot {\bf d}
= {\bf a}\centerdot {\bf c} + {\bf b}\centerdot {\bf d}
= {\bf a}\centerdot {\bf d} + {\bf b}\centerdot {\bf c}.$

Notice that in quadration (and in twinning; see below)
the circumradius, $R$, of all four (all eight) triangles is
the same.  If the angles of a triangle are $A$, $B$, $C$,
then, by the sine formula, the edges are
$$2R\sin A, \quad 2R\sin B, \quad 2R\sin C.$$

and the angles and edges of the quadrations are
$$\pi - A, \quad \tfrac{\pi}{2} - B, \quad \tfrac{\pi}{2} - C,
\mbox{\quad and \quad} 2R\sin A, \quad 2R\cos B, \quad 2R\cos C, $$
$$\tfrac{\pi}{2} - A, \quad \pi - B, \quad \tfrac{\pi}{2} - C,
\mbox{\quad and \quad} 2R\cos A, \quad 2R\sin B, \quad 2R\cos C, $$
$$\tfrac{\pi}{2} - A, \quad \tfrac{\pi}{2} - B, \quad \pi - C,
\mbox{\quad and \quad} 2R\cos A, \quad 2R\cos B, \quad 2R\sin C. $$

Note that the triangle inequality appears as
$$\sin B + \sin C > \sin B\cos C + \sin C\cos B = \sin(B+C)=
\sin(\pi-A)=\sin A$$
and, for the quadrations
$$\cos B + \cos C > \sin B\cos C + \sin C\cos B = \sin(B+C)=
\sin(\pi-A)=\sin A$$

We can contrast and combine this with Conway's extraversion,
in which the extraverted triangles have angles and edges
$$-A,\quad \pi - B, \quad \pi - C, \mbox{\quad and \quad}
-2R\sin A, \quad 2R\sin B, \quad 2R\sin C,$$
$$\pi - A,\quad - B, \quad \pi - C, \mbox{\quad and \quad}
2R\sin A, \quad -2R\sin B, \quad 2R\sin C,$$
$$\pi - A,\quad \pi - B, \quad - C, \mbox{\quad and \quad}
2R\sin A, \quad 2R\sin B, \quad -2R\sin C.$$

\section{Twinning}

If you draw the perpendicular bisectors of each of the six
edges of an orthocentric quadrangle, you produce a quadrangle
of circumcentres, congruent to the original one; in fact the
two {\bf twins} are obtained from each other by reflexion in,
or rotation through $180^{\circ}$ about, a common {\bf centre},
$O$.  In fact, if one wishes to select just one out of Clark
Kimberling's six thousand, six hundred-odd triangle centres,
then a good case can be made for THE centre, was 50-point
centre!  Indeed, any other candidate has its twin to compete with.

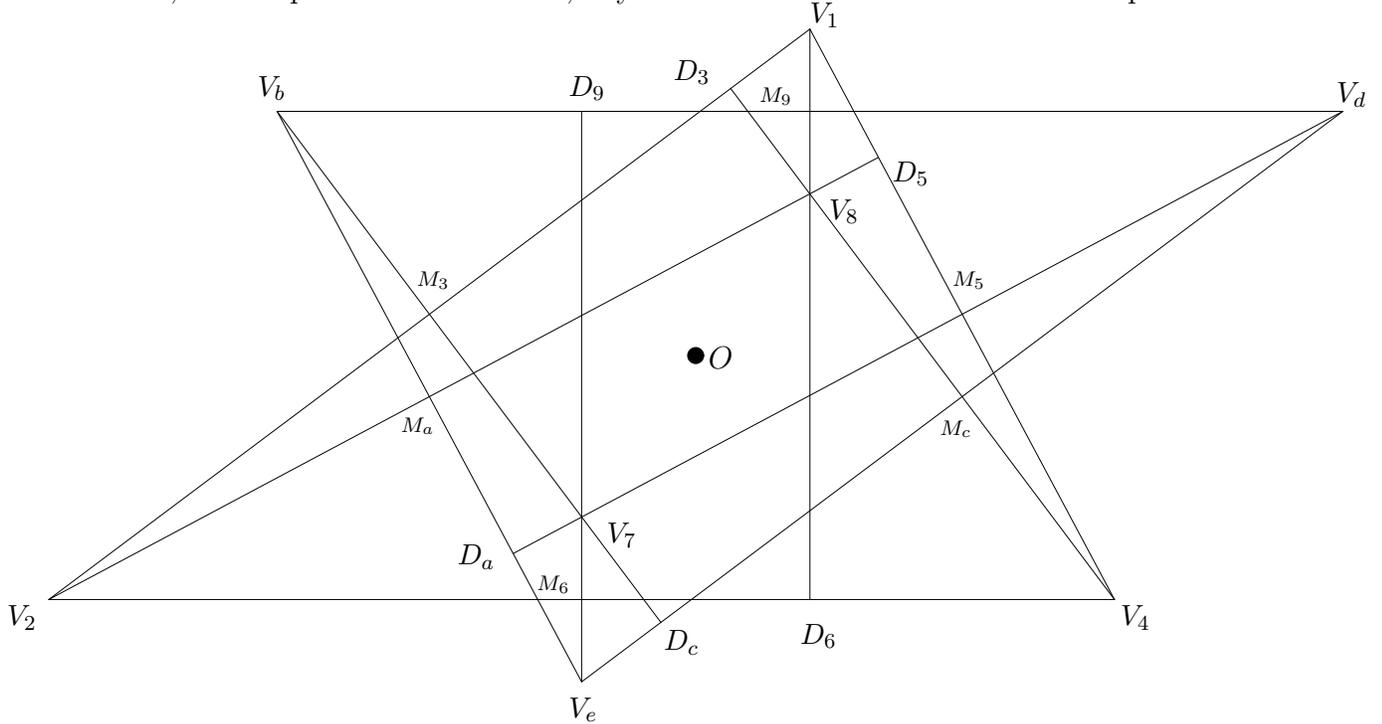
\begin{figure}[h] 
\begin{picture}(440,250)(-230,-130)
\setlength{\unitlength}{1.2pt}
\drawline(-204,-77)(132,-77)(36,103)(-204,-77)  
\drawline(-204,-77)(57.59,62.51)                
\drawline(132,-77)(11.04,84.28)                 
\drawline(36,-77)(36,103)                       
\drawline(204,77)(-132,77)(-36,-103)(204,77)    
\drawline(204,77)(-57.59,-62.51)                
\drawline(-132,77)(-11.04,-84.28)               
\drawline(-36,77)(-36,-103)                     
\put(36,105){\small$V_1$}
\put(-217,-85){\small$V_2$}
\put(134,-85){\small$V_4$}
\put(42,43){\small$V_8$}
\put(-40,-114){\small$V_e$}
\put(202,80){\small$V_d$}
\put(-138,82){\small$V_b$}
\put(-28,-59){\small$V_7$}
\put(33,-91){\small$D_6$}
\put(-7,87){\small$D_3$}
\put(62,55){\small$D_5$}
\put(-40,82){\small$D_9$}
\put(-10,-93){\small$D_c$}
\put(-75,-66){\small$D_a$}
\put(0,0){\circle*{5}}
\put(4,-4){$O$}
\put(-88,22){\scriptsize$M_3$}
\put(81,22){\scriptsize$M_5$}
\put(-50,-74){\scriptsize$M_6$}
\put(77,-25){\scriptsize$M_c$}
\put(-93,-24){\scriptsize$M_a$}
\put(20,80){\scriptsize$M_9$}
\end{picture}
\caption{Hello, twins!  Hexadecimal subscripts:
$a=10$, $b=11$, $c=12$, $d=13$, $e=14$.}
\label{twins}
\end{figure}

\bigskip

~\hspace{20pt} To the tune of ``My Bonny Lies Over the Ocean''

\medskip

\begin{verse}

Now that we've heard of Quadration, \\
And Twinning gives Two -- symmetry, \\
And Conway has found Extraversion. \\
We can bring back that gee-om-met-tree!

Bring back! Oh, bring back! \\
Oh bring back that gee-om-met-tree, to me! \\
Bring back! Oh, bring back! \\
Oh, bring back that gee-om-met-tree !!
\end{verse}

More verses about various geometrical objects??

Since then I've been invited to preach at Alberta Math
Dialog in Lethbridge in May, at the MOVES conference in
NYC, and at Mathfest Washington DC in August.  So here's
a shot at an

\bigskip

\centerline{\LARGE ABSTRACT}

\bigskip

\centerline{\Large A triangle has eight vertices, but only one
centre}

\medskip

Quadration regards a triangle as an orthocentric
quadrangle.  Twinning is an involution between
orthocentres and circumcentres.  Together with
variations of Conway's Extraversion, these give rise
to symmetric sets of points, lines and circles.
There are eight vertices, which are also both
orthocentres and circumcentres.  Twelve edges share
six midpoints, which with six diagonal points,
lie on the was 50-point circle, better known as the 9-point
circle.  There are 32 circles which touch three edges and
also touch the 50-point circle.  32 Gergonne points, when
joined to their respective touch-centres, give sets of four
segments which concur in eight deLongchamps points, which,
with the eight centroids, form two harmonic ranges with the
ortho- and circum-centres on each of the four Euler lines.
Corresponding points on the eight circumcircles generate
pairs of parallel Simson-Wallace lines, each containing six
feet of perpendiculars.  In three symmetrical positions
these coincide, with twelve feet on one line.  In the three
orthogonal positions they are pairs of parallel tangents to
the 50-point circle, forming the Steiner Star of David.
This three-symmetry is shared with the 144 Morley triangles
which are all homothetic.  Time does not allow investigation
of the 256 Malfatti configurations, whose 256 radpoints
probably lie in fours on 64 guylines, eight through each
of the eight vertices.

\newpage

\section{Radical axes of touch-circles}

On 2015-04-28 Peter Moses told me that
The radical centre of the ex-circles is X(10) in \cite{K}.

[The Spieker centre, barycentric $(b+c,c+a,a+b)$.  So, with quadration
\& twinning, there are 8 Spieker centres. Lies on IG. And on HM (orthoc
\& mittenpunkt).]

The radical center of the B- \& C-excircles with the incircle is {b + c, c -
a, b - a} (barycentric).

Call the triangle formed by the cyclic permutation of {b + c, c - a, b - a}
A'B'C'.
The excentral triangle, IaIbIc, is perspective to A'B'C' at X(2), G.
Indeed |Ia G| / |G A'| = 2.

\section{Apollonian interlude}

This section was called ``Cartesian points and axes'' until
Peter Moses reminded me that I'd already dealt with many parts
of the subject in  below, where it can be seen that
these points and line are the Isoperimetric Point, the Equal
Detour Point, and the Soddy Line!  In fact, let
the two points be $X_{175}$ and $X_{176}$ and the radii of the
two circles be $r_{175}$ and $r_{176}$, then the perimeters of the
three triangles $X_{175}BC$, $X_{175}CA$, $X_{175}AB$ are
$(r_{175} - (s-b)) + a + (r_{175} - (s-c)) = 2\,r_{175}$ and two
equal expressions, while the extra detour in travelling from
$B$ to $C$ via $X_{176}$ instead of travelling directly is
$((s-b)+r_{176})+(r_{176}+(s-c))-a=2\,r_{176}$, and similarly
in detouring through $X_{176}$ when travelling from $C$ to
$A$, or from $A$ to $B$.

Then I heard Kate Stange's lecture at the Alberta Number
Theory Days in Banff, and realised that while much has been
said recently (\cite{Lag,Gra1,Gra2,Gra3,Gra4} and about
100 other papers) about Apollonian packings, there is
plenty more to be explored.

\newpage

If we draw circles with centres at the vertices (for which, pro tem,
we'll use the traditional notation $A$, $B$, $C$) and respective radii
$s-a$, $s-b$, $s-c$, they will touch one another at the touch points
of the edges with the incircle, with which they are orthogonal.
[I proposed to call them ``vertical circles'', but Peter
Moses reminds me that I have already called
them ``tangent circles'' in Chapter 5 below, written some time
before the present sections!] According to Descartes, following
Apollonius, there are two circles which touch all three vertical
circles.  I called their centres the Cartesian points
and their join the Cartesian axis, but later
(earlier!) found that these are the so-called Soddy circles,
centres the isoperimetric point and the equal detour point,
X(175) and X(176) in \cite{K,V}, discussed in some detail
in Chapter 5 below!!]

[I can't use ``Apollonian point'' since
``Apollonius point'' is already used for something else, X(181)
in Clark Kimberling's Encyclopedia \cite{K}.  Suppose
that the circle circumscribing the three excircles touches them at
$A'$, $B'$, $C'$.  Then $AA'$, $BB'$, $CC'$ concur at X(181).  Its
barycentric coordinates are $\{\ldots, b^3\cos^2[(C-A)/2],\ldots\}$.
See \cite{KIF}.  This is a pity, since this is {\bf not} an example
of what has recently received a good deal of attention, namely
``Apollonian packings''.]

Indeed, look at an Apollonian packing.  Draw the segment joining
the centres of each pair of circles which touch.  We have a
triangulation which Kate Stange called an {\bf Apollonian city}.
Now apply much of the enormous amount that is known about
triangle geometry to each of the triangles in the triangulation.

By Conway's extraversion we will have four isoperimetric points,
four equal detour points and four Soddy lines.  Quadration and
twinning give 32 of each of these items.

What coincidences, collinearities and concurrences are there??
Skip ahead to just after Figure \ref{4quads}.

\begin{figure}[h]
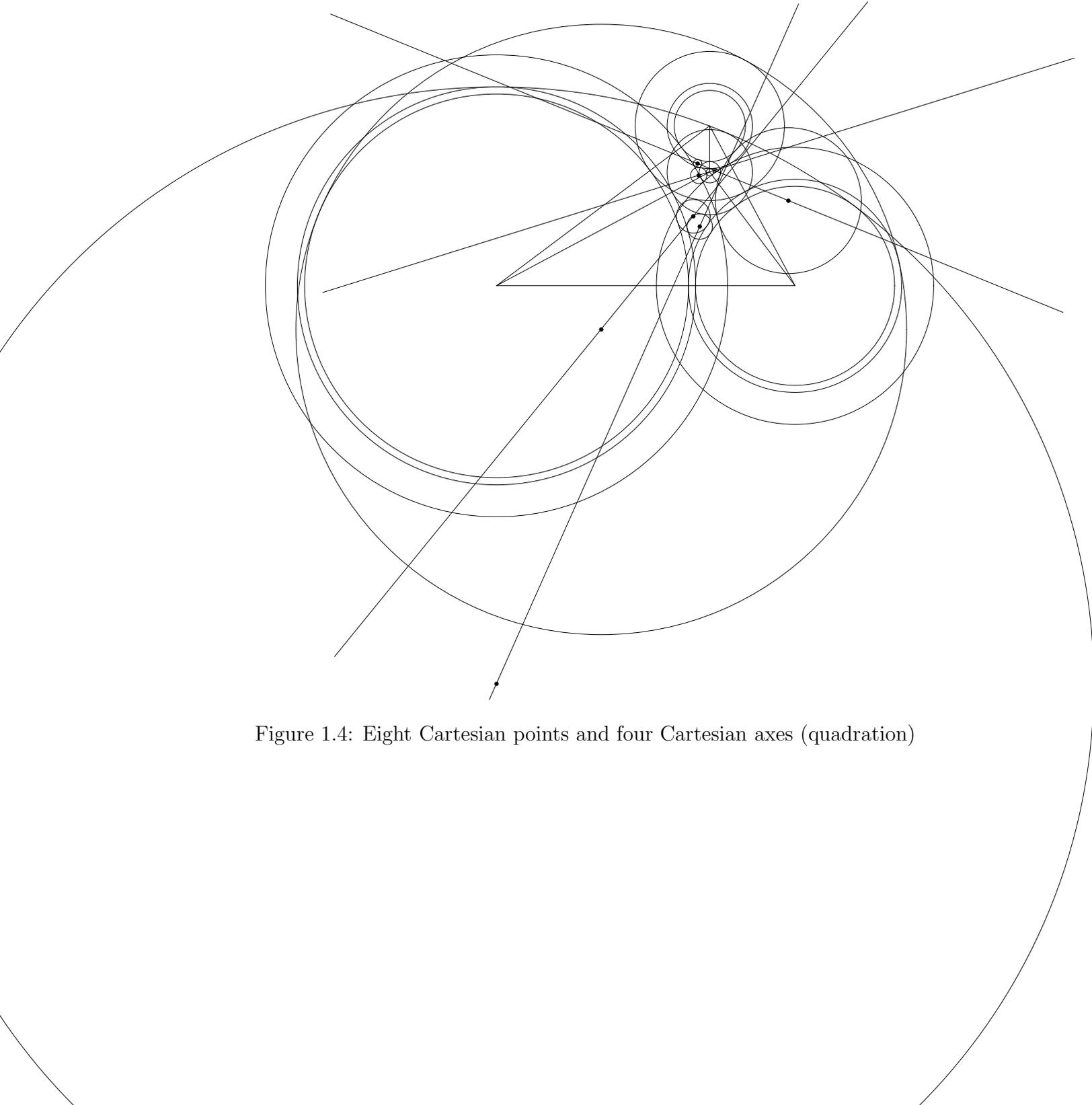
 

\caption{Eight Cartesian points and four Cartesian axes (quadration)}
\label{quadapollo}
\end{figure}

\clearpage

That was the result of quadration.  Now let's try extraversion:

\begin{figure}[h]
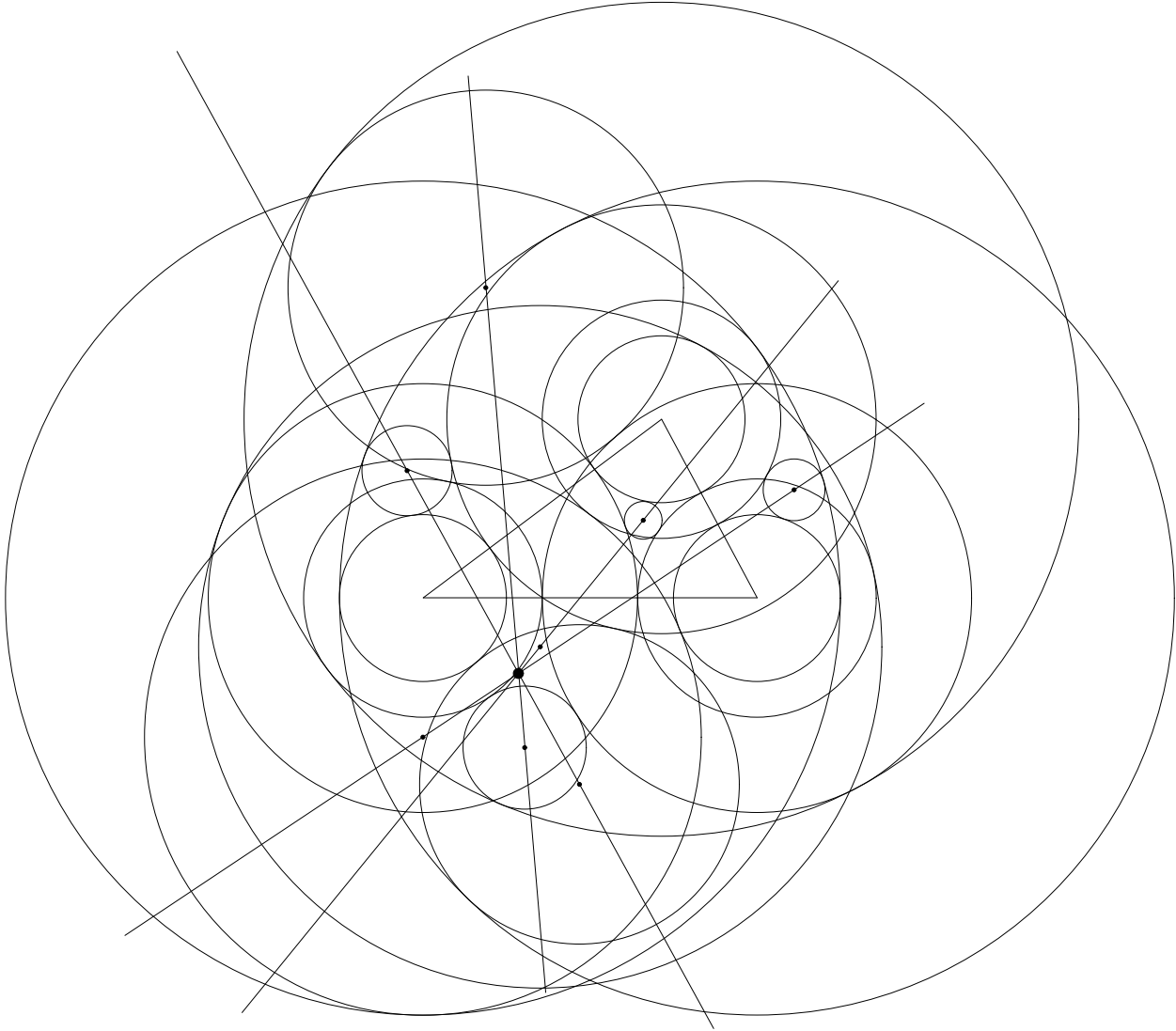
 
%
\caption{Eight Cartesian points and four Cartesian axes (extraversion of $V_8$)}
\label{extrapollo8}
\end{figure}

Surprise! Surprise! The four Cartesian axes in the extraversion concur!!
Will try extraverting the other 3 triangles in the orthocentric
quadrangle.  NOT A SURPRISE TO THOSE WHO WELL KNOW IT!!!(*)  This is the
De Longchamps point [X(20) in \cite{K}].  In the notation of \S1.6,
The orthocentre is {\bf d}, the circumcentre is {\bf --d}, the centroid
is \\
$\frac{1}{3}$({\bf a} + {\bf b} + {\bf c}) = $-\frac{1}{3}${\bf d},
and the De Longchamps point is --3{\bf d}.  The four points form a
harmonic range.

(*) IN FACT, SEE FIRST PAGE OF CHAPTER 5 BELOW!!  No fool like an old
fool!

\clearpage

\begin{figure}[h]
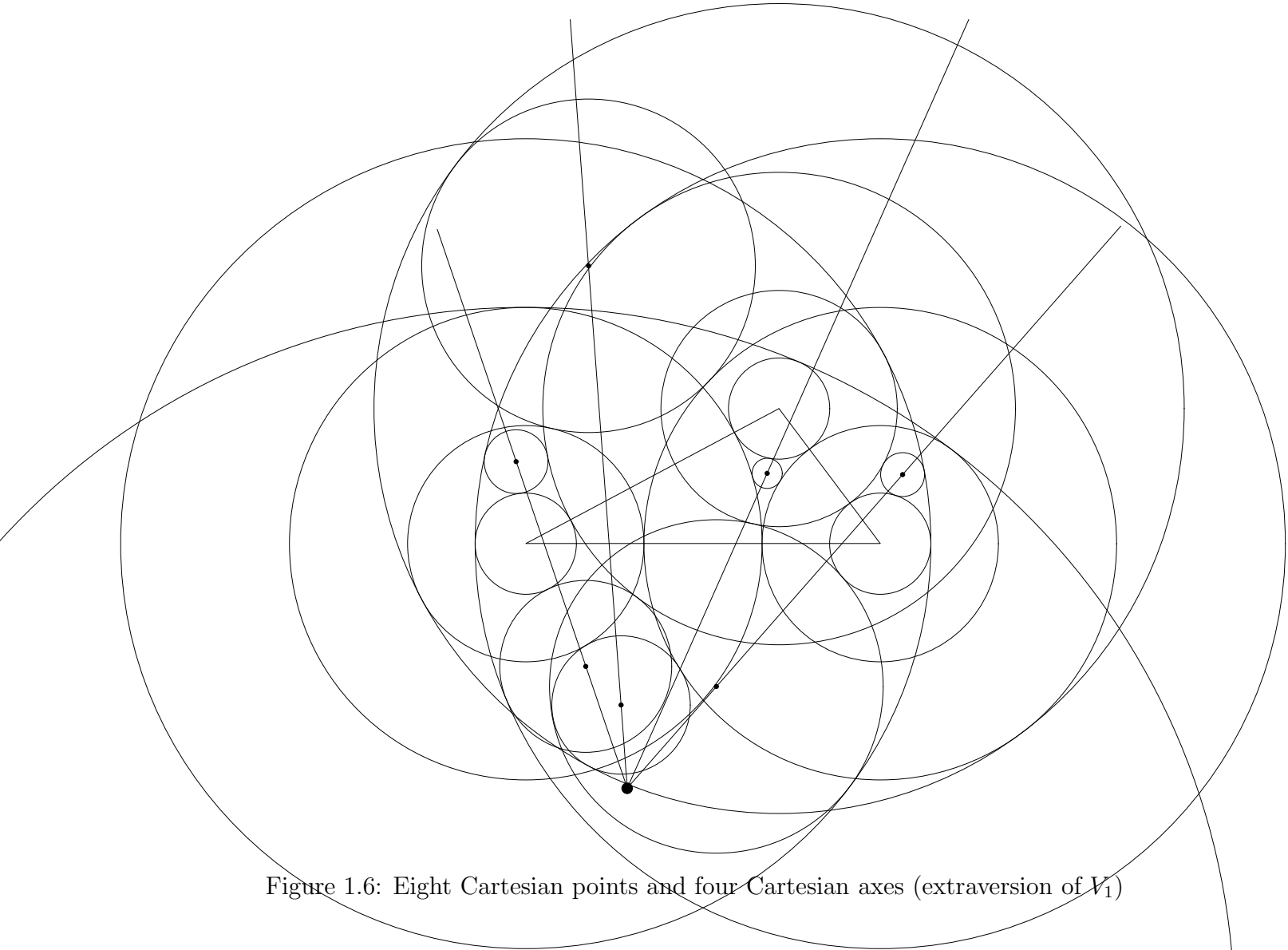
 

\caption{Eight Cartesian points and four Cartesian axes (extraversion of $V_1$)}
\label{extrapollo1}
\end{figure}

\clearpage

\begin{figure}[h] 
%
\caption{Eight Cartesian points and four Cartesian axes (extraversion of $V_2$)}
\label{extrapollo2}
\end{figure}

\clearpage

\begin{figure}[h] 
%
\caption{Eight Cartesian points and four Cartesian axes (extraversion of $V_4$)}
\label{extrapollo4}
\end{figure}

\clearpage

Let us collect together the 16 Cartesian axes:

\begin{figure}[h]
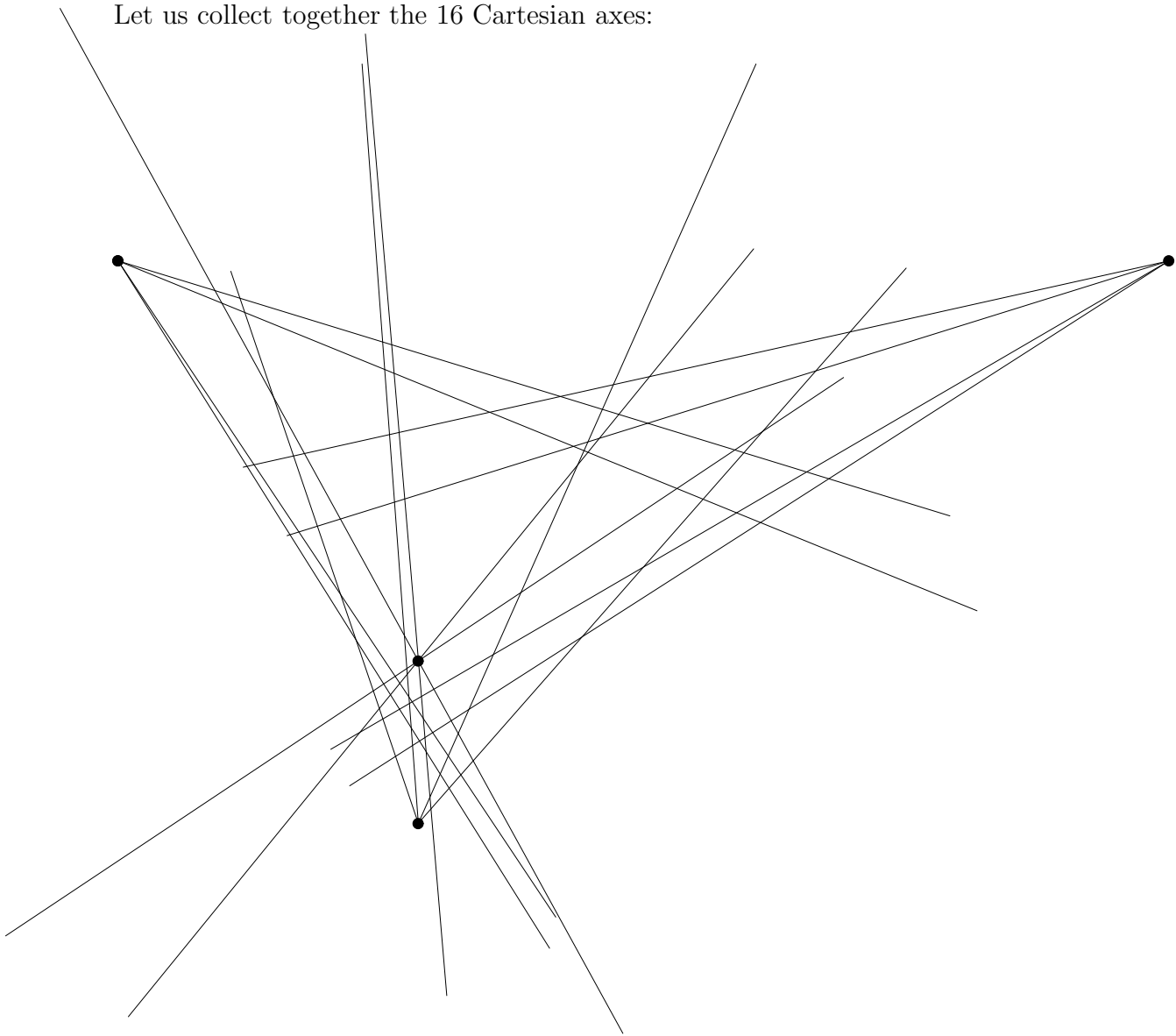
 
%
\caption{Sixteen Cartesian axes}
\label{16axes}
\end{figure}

The four points of concurrence form an orthocentric quadrangle
which is homothetic to the original orthocentric quadrangle
with ratio --3 and centre of similitude the 50-point centre
(which is, of course, the 50-point centre of the new
(Cartesian?) orthocentric quadrangle).

\clearpage

In Figure \ref{4quads} we have an orthocentric quadrangle,
together with its {\bf twin}, and their {\bf Cartesian
quadrangles} (dashed lines).  There are four sets of four
equally spaced vertices (dotted lines) which concur and
bisect each other at THE centre, $O$.

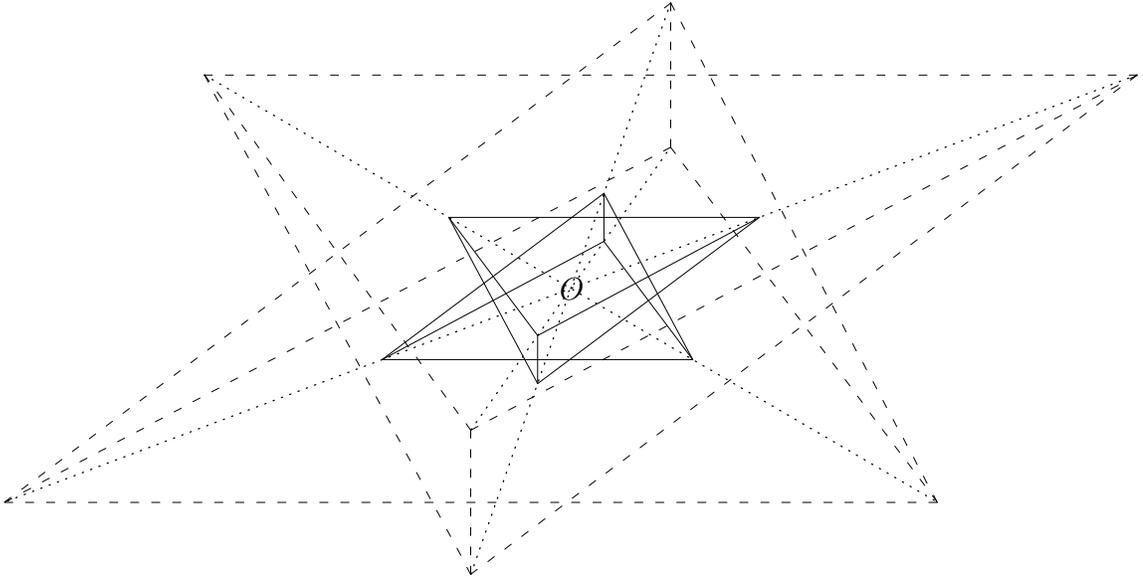
\begin{figure}[h] 
\begin{picture}(440,220)(-230,-100)
\setlength{\unitlength}{0.35pt}
\drawline(36,103)(36,51)(-204,-77)(132,-77)
\drawline(36,51)(132,-77)(36,103)(-204,-77)
\dashline{10}(-108,-309)(-108,-153)(612,231)(-396,231)
\dashline{10}(-108,-153)(-396,231)(-108,-309)(612,231)
\drawline(-36,-103)(-36,-51)(204,77)(-132,77)
\drawline(-36,-51)(-132,77)(-36,-103)(204,77)
\dashline{10}(108,309)(108,153)(-612,-231)(396,-231)
\dashline{10}(108,153)(396,-231)(108,309)(-612,-231)
\dottedline{10}(-108,-309)(108,309)
\dottedline{10}(-108,-153)(108,153)
\dottedline{10}(-396,231)(396,-231)
\dottedline{10}(612,231)(-612,-231)
\put(-12,-12){$O$}
\end{picture}
\caption{Four quadrangles in perspective}
\label{4quads}
\end{figure}

NOT A SURPRISE TO THOSE WHO WELL KNOW IT!!!  THE DOTTED LINES
(through $O$) ARE EULER LINES!!  The 8 points are De Longchamps
points of the original 8 triangles.  In the notation of
\S1.6, if we take THE (50-point) centre as origin and the
orthocentre as {\bf d}, then the circumcentre is {\bf --d},
the centroid is
$\frac{1}{3}$({\bf a} + {\bf b} + {\bf c}) = $-\frac{1}{3}${\bf d},
and the De Longchamps point is --3{\bf d}, or, as Clark Kimberling
(\cite{K}) puts it, the orthocentre of the anticomplementary
triangle!  A harmonic range:-- \\[-9pt]

\centerline{\{Orthocentre , Circumcentre ; Centroid , De Longchamps point\}
\ = \ --1.}

Reverting to the Apollonian city: for each triangle, the circles
of the packing are what I variously called the ``vertical circles''
or the ``tangent circles''.  The positions of their centres are
well known in terms of Gaussian rational numbers.  The incircles
form the orthogonal Apollonian packing, so the coordinates of the
incentres are similarly well known.  But what about the excircles?
And more generally, what does Conway's extraversion bring??

Also what does quadration bring??  Where are the orthocentres?
They must be rational points also.  Similarly for the
circumcentres.  And for the 50-point (9-point) centres.
So little done!  So much to do!

\clearpage

\section{Nagel points}

The joins of the vertices to the touch-points of the excircles
with the opposite edges, concur in the Nagel point.  I don't
see how to extravert this.  But let's first quadrate it.

\begin{figure}[h] 
\begin{picture}(440,220)(-230,-100)
\setlength{\unitlength}{1pt}
\drawline(36,103)(36,51)(-204,-77)(132,-77)
\drawline(36,51)(132,-77)(36,103)(-204,-77)
\put(-60,-41){\circle*{3}}       
\put(12,-5){\circle*{3}}         
\put(-12,-17){\circle*{3}}       
\drawline(12,-5)(-60,-41)
\put(-76,-61){\circle*{3}}       
\put(20,-21){\circle*{3}}        
\put(-12,-34.333){\circle*{3}}   
\drawline(-76,-61)(20,-21)
\put(-164,-49){\circle*{3}}      
\put(16,63){\circle*{3}}         
\put(-44,25.667){\circle*{3}}    
\drawline(-164,-49)(16,63)
\put(108,-33){\circle*{3}}       
\put(48,55){\circle*{3}}         
\put(68,25.667){\circle*{3}}     
\drawline(108,-33)(48,55)
\dashline{2}(36,103)(-84,-77)
\dashline{2}(-204,-77)(92.47,-2.882)
\dashline{2}(132,-77)(-136.8,-26.6)

\dashline{2}(36,51)(-92,-77)
\dashline{2}(-204,-77)(103.2,-38.6)
\dashline{2}(132,-77)(-161.647,-54.412)

\dashline{2}(-204,-77)(36,91)
\dashline{2}(36,51)(-172,-53)
\dashline{2}(36,103)(-193.412,-71.353)

\dashline{2}(36,103)(129.6,-73.8)
\dashline{2}(132,-77)(36,99)
\dashline{2}(36,51)(109.412,-34.647)

\end{picture}
\caption{Quadration of Nagel points.  The four segments
illustrate the collinearity of the Nagel point, incentre
and centroid.}
\label{quadnagel}
\end{figure}
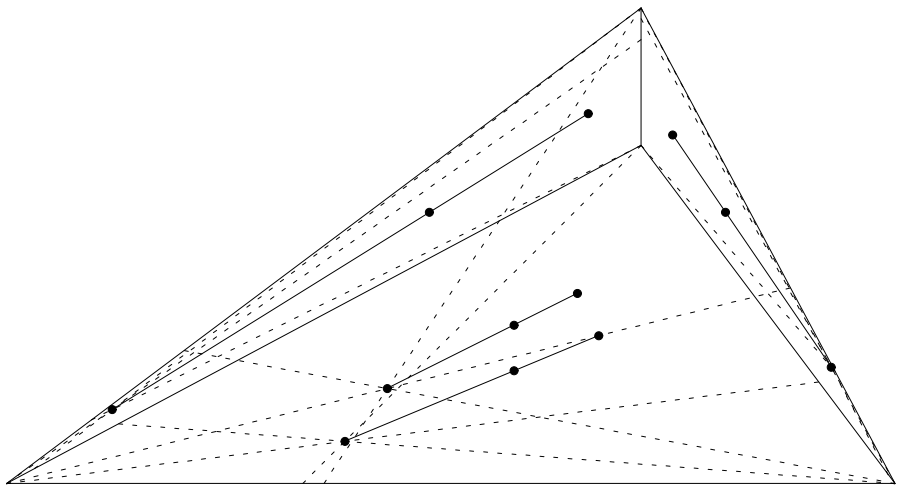

The results are disappointing.  Slopes of joins of
Nagel points are:

\medskip

\centerline{$N_1N_8: \frac{5}{4}$ \quad
$N_2N_4: -\frac{1}{17}$}

\medskip

\centerline{$N_1N_2: \frac{7}{46}$ \quad
$N_4N_8: \frac{1}{13}$}

\medskip

\centerline{$N_1N_4: -\frac{3}{22}$ \quad
$N_2N_8: \frac{1}{21}$}

\medskip

while the slopes of $N_1I_1$,  $N_2I_2$, $N_4I_4$, $N_8I_8$
are $$\frac{5}{12}, \quad -\frac{22}{15}, \quad \frac{28}{45},
\quad \frac{1}{2},$$
--- not particularly inspiring!  Let's hope for better luck with the \ldots

\clearpage

\section{Gergonne points}

\begin{figure}[h]
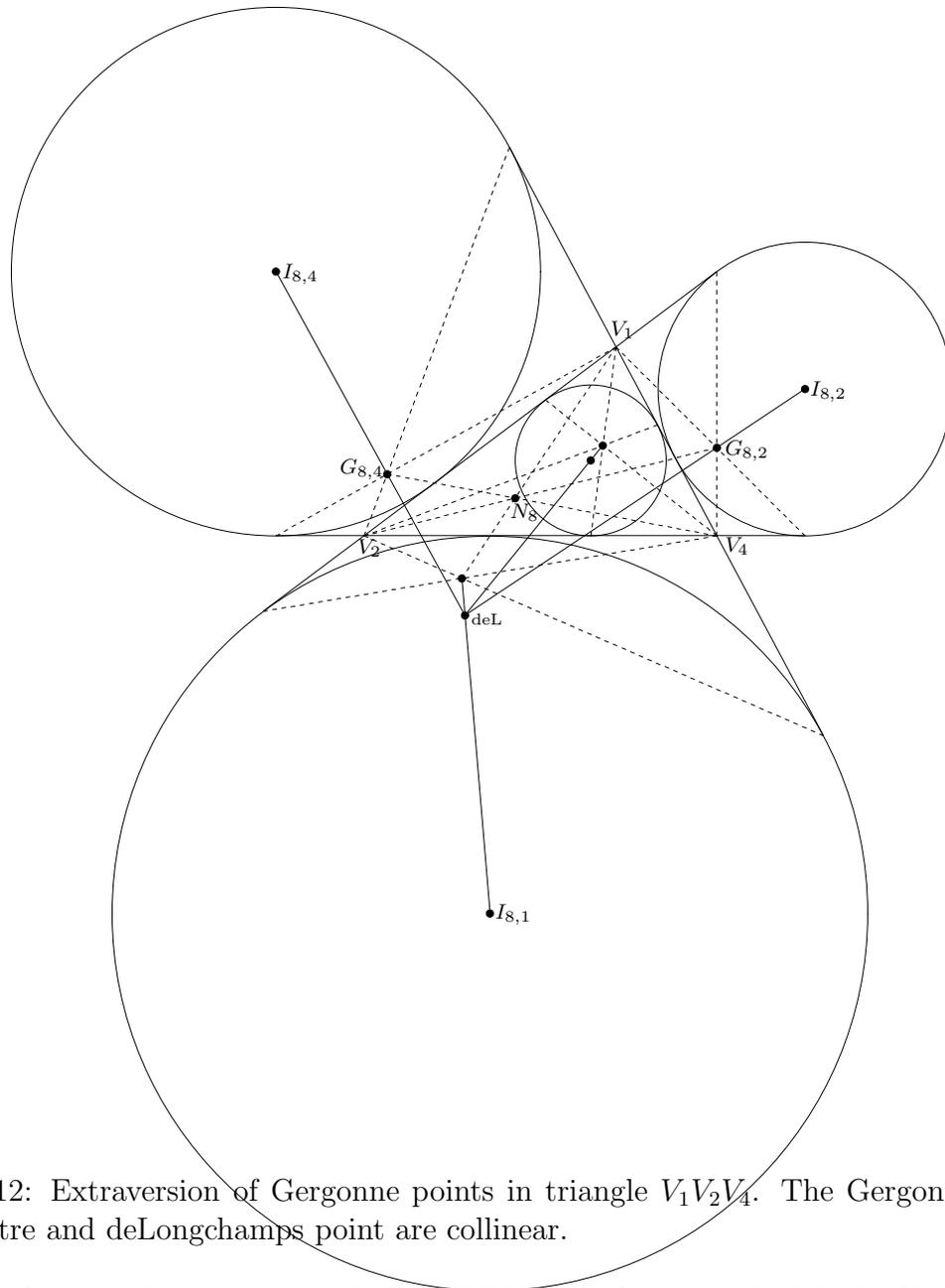
 

\caption{Extraversion of Gergonne points in triangle
$V_1V_2V_4$.  The Gergonne point,
touch-centre and deLongchamps point are collinear.}
\label{extragerg8}
\end{figure}

Note that the joins $V_1G_{8,1}$, $V_2G_{8,2}$, $V_4G_{8,4}$
concur.  In what point??  It's the Nagel point, $N_8$\,!!!

\clearpage

\begin{figure}[h]
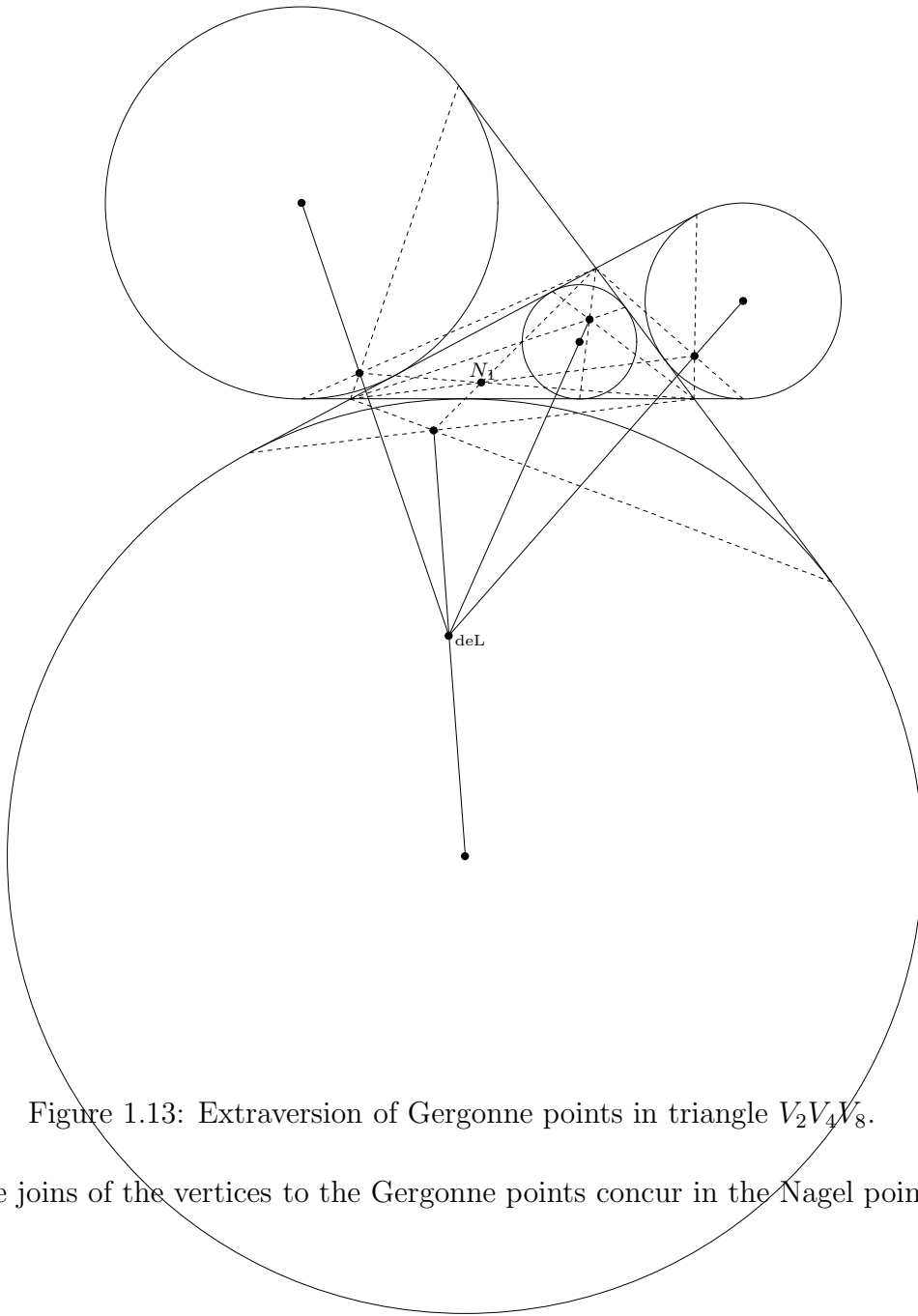
 
%
\caption{Extraversion of Gergonne points in triangle $V_2V_4V_8$.}
\label{extragerg1}
\end{figure}

Again, the joins of the vertices to the Gergonne points
concur in the Nagel point, $N_1$\,!!

\clearpage

\begin{figure}[h]
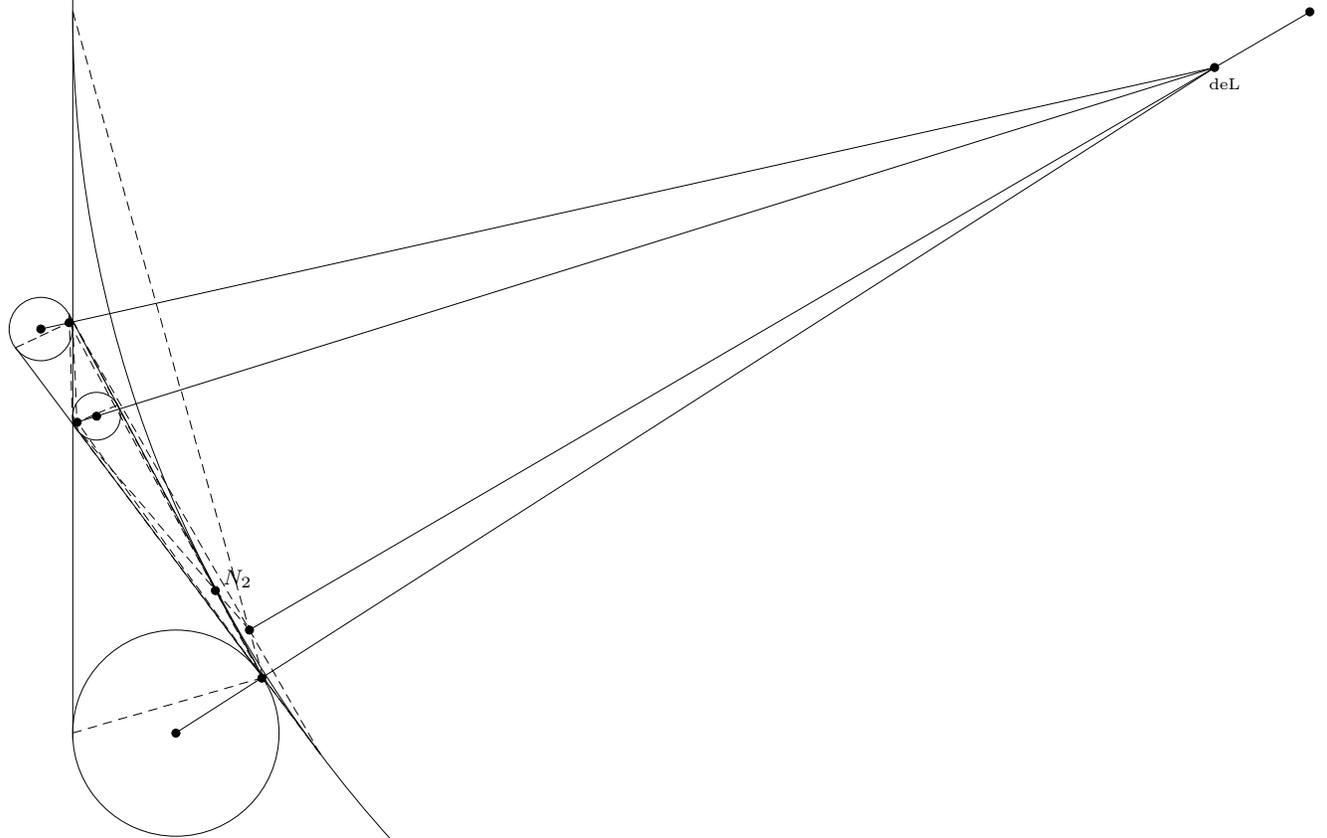
 
%
\caption{Extraversion of Gergonne points in triangle $V_1V_4V_8$.}
\label{extragerg2}
\end{figure}

\clearpage

\begin{figure}[h]
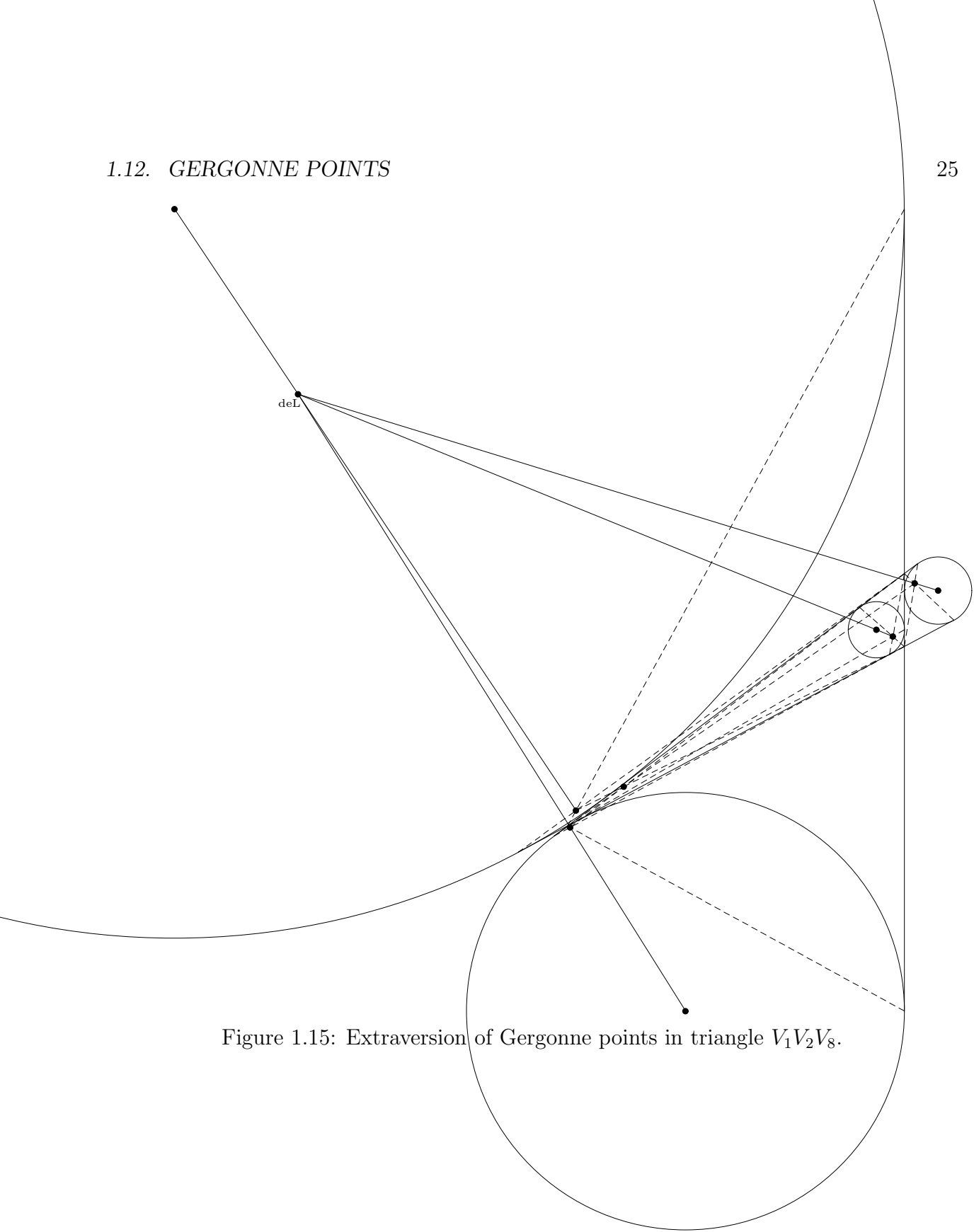
 
%
\caption{Extraversion of Gergonne points in triangle $V_1V_2V_8$.}
\label{extragerg4}
\end{figure}

\clearpage

Let's collect the 16 (32??) Gergonne points into one diagram.
The lines connect them to the corresponding touch-centres
(without dots) and deLongchamps points (with dots and labels),
but the result is not very inspiring!

\begin{figure}[h]
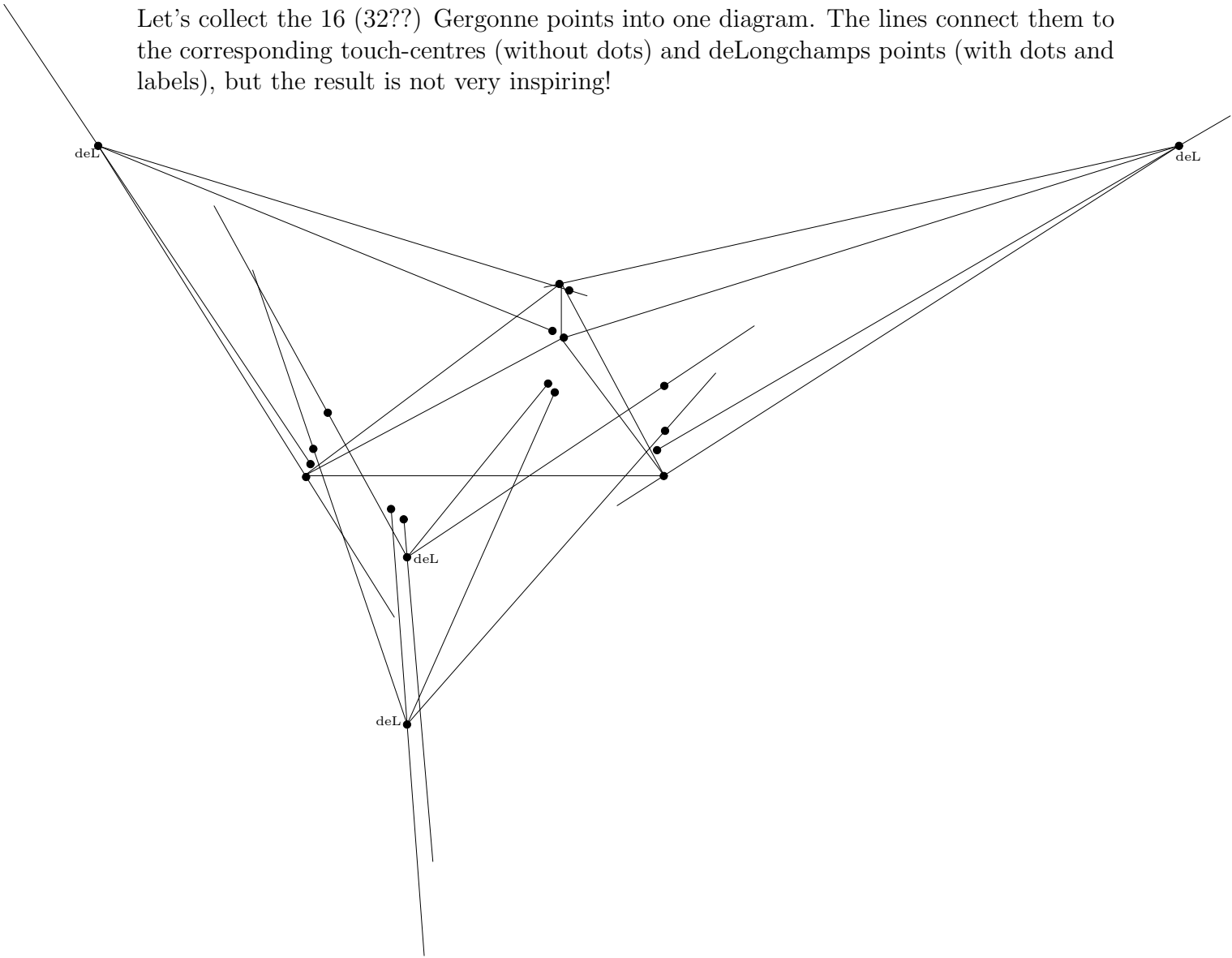
 
%
\caption{Sixteen Gergonne points.  There are 32 if you include their twins.}
\label{16gerg}
\end{figure}

\clearpage

\section{Simson-Wallace lines and the (double-) deltoid}

A theorem that used to be well known is that the feet of the
perpendiculars from a point on the circumcircle to the edges
of a triangle are collinear. This is usually known as the
{\bf Simson line} of the point, though it was discovered
by William {\bf Wallace} and doesn't occur in Simson's work
\cite[p.~16]{Cox}, \cite[p.~41]{C+G}
\footnote{If you google ``Wallace line'' you'll find
something closely related to Darwin's {\it Origin of Species},
whereas if you google ``Wallace Sim(p)son'' you'll find
the wife-to-be of the uncrowned King of England.}.
Steiner showed that as the point moves round the circumcircle,
the Simson-Wallace line envelops a ({\bf Steiner}) {\bf deltoid}
or three-cusped hypocycloid \cite[p.~115]{Cox},
\cite[p.~44]{C+G}. This is usually described as the locus
of a point on a circle rolling inside a fixed circle of
three times its radius.  But here we want the curve as a
{\bf roulette}, and it is easier to visualize and to
construct if you regard it as the envelope of a diameter
of a circle of radius $R$ rolling inside a circle, centre
$O$ and radius $3R/2$.

\begin{figure}[h] 
\begin{picture}(440,240)(-245,-130)
\setlength{\unitlength}{0.6pt}
\put(-36,-51){\circle{340}}  
\drawline(-4,178)(132,-77)(-204,-77)(36,103) 
\put(-62,117){\circle*{3}}
\drawline(36,103)(-62,117)(-62,-77)  
\drawline(-204,-77)(-62,117)(8.4844,154.5917)  
\drawline(132,-77)(-62,117)(-20,61)  
\put(37,104){\scriptsize$V_1$}
\put(-218,-84){\scriptsize$V_2$}
\put(134,-84){\scriptsize$V_4$}
\put(-76,121){\scriptsize$S_8$}
\put(-61,-87){\scriptsize$F_6$}
\put(-20,52){\scriptsize$F_3$}
\put(10,155){\scriptsize$F_5$}
\dashline[+50]{3}(-104,-215)(29,222)  

\end{picture}
\caption{The Simson-Wallace line theorem.}
\label{simson}
\end{figure}
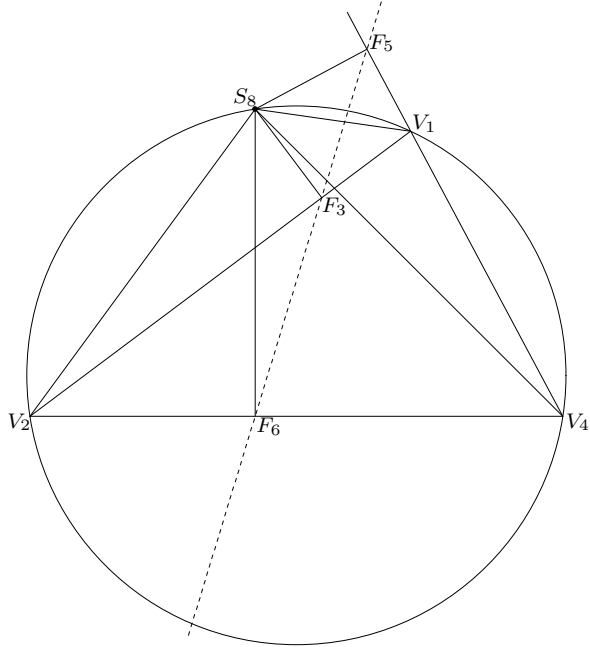

\begin{figure}[h]
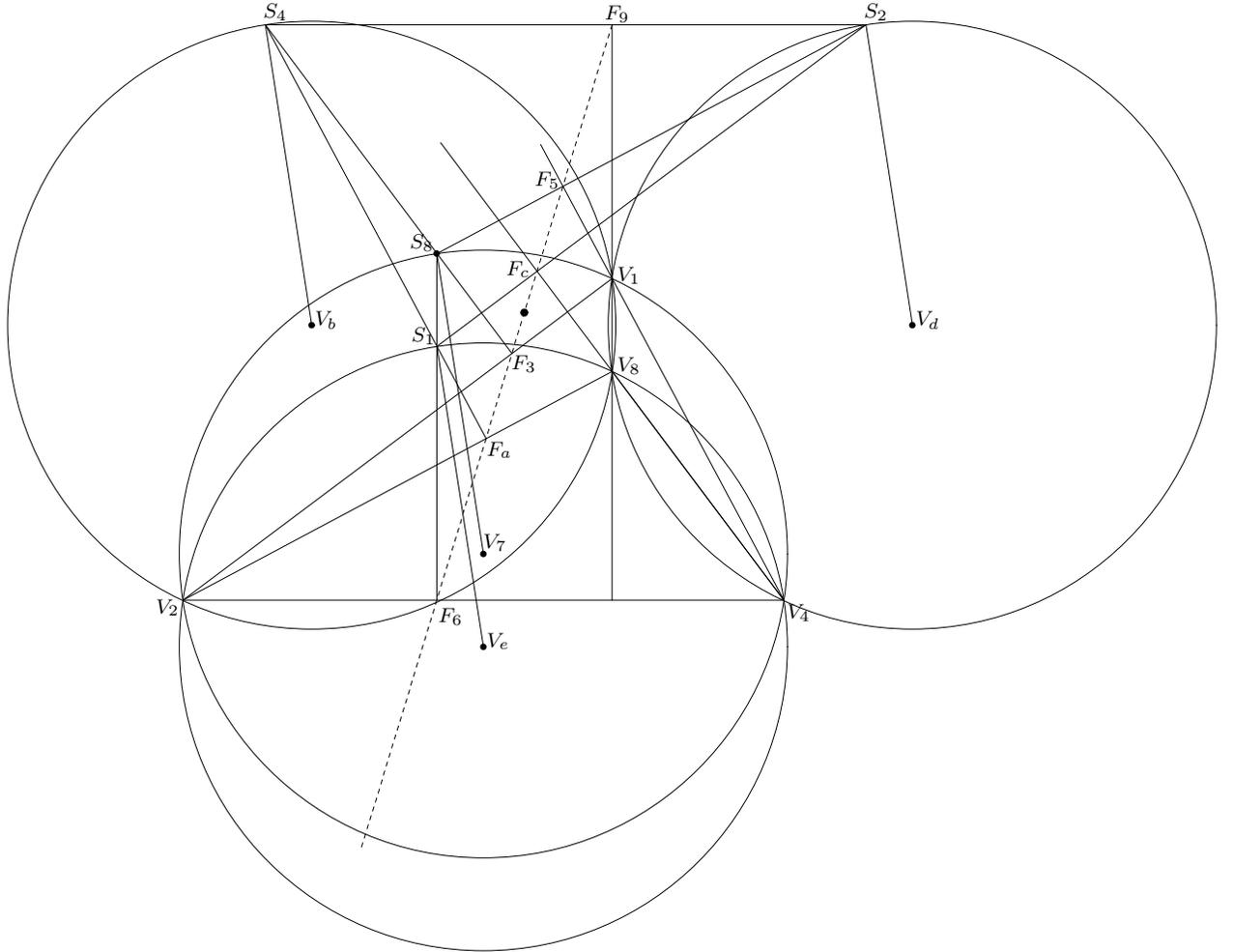
 

\caption{Quadration of the Simson-Wallace line theorem.}
\label{quadsimson}
\end{figure}

We first give a proof of the Simson-Wallace line theorem
(Figure \ref{simson}) and notice the generalization that comes
from quadration (Figure \ref{quadsimson})

Notice that the quadrilaterals $F_6F_3S_8V_2$, $F_5S_8F_6V_4$,
$V_1V_4V_2S_8$, $F_5S_8F_3V_1$ are all cyclic, so that
$\angle F_6F_3V_2=\angle F_6S_8V_2=\angle F_5S_8V_2-\angle F_5S_8F_6=
\angle F_5S_8V_2-(\pi-\angle F_5V_4F_6)=
\angle F_5S_8V_2-(\pi-\angle V_1V_4V_2)=
\angle F_5S_8V_2-\angle V_1S_8V_2=\angle F_5S_8V_1=
\angle F_5F_3V_1$, and the vertically opposite $\angle F_6F_3V_2$
and $\angle F_5F_3V_1$ show that $F_6$, $F_3$ and $F_5$
are collinear.

\clearpage

\section{Quadration and Twinning of the Simson-Wallace line}

Draw radii $V_eS_1$, $V_dS_2$, $V_bS_4$, of the circumcircles
of triangles $V_2V_4V_8$, $V_4V_8V_1$, and $V_8V_1V_2$,
parallel to the radius $V_7S_8$ of the circumcircle of
triangle $V_1V_2V_4$, yielding points $S_1$, $S_2$, $S_4$
on the respective circumcircles.  Then the feet of the
perpendiculars from these points onto appropriate
edges of the respective triangles give triples of points \\
\centerline{\{$F_c$, $F_a$, $F_6$\}, \{$F_9$, $F_5$, $F_c$\},
\{$F_3$, $F_a$, $F_9$\},}
all lying on the same, original, Simson-Wallace line, \{$F_6$, $F_5$, $F_3$\}.

We know that $V_eV_dV_bV_7$ form a quadrangle congruent to
the original $V_1V_2V_4V_8$; note that $S_1S_2S_4S_8$ is
also such a congruent quadrangle.  In fact
$S_1V_1$, $S_2V_2$, $S_4V_4$, $S_8V_8$
bisect each other at a point on the Simson-Wallace line.  See the
second theorem ahead.

Notice that twinning now gives a pair of parallel Simson-Wallace
lines.

On three, equally spaced, occasions, this pair will coincide,
giving twelve collinear points.  These are the feet of the
perpendiculars from points on eight circumcircles onto the three
edges of the corresponding triangles; note that an edge belongs
to just two of the eight triangles.

On three other equally spaced occasions, the pair will be
perpendicular to one of the twelve-point lines and are
tangents to the 50-point circle at its points of intersection
with these twelve-point lines, forming the

\newpage

\section{Steiner Star of David}


\begin{figure}[h]
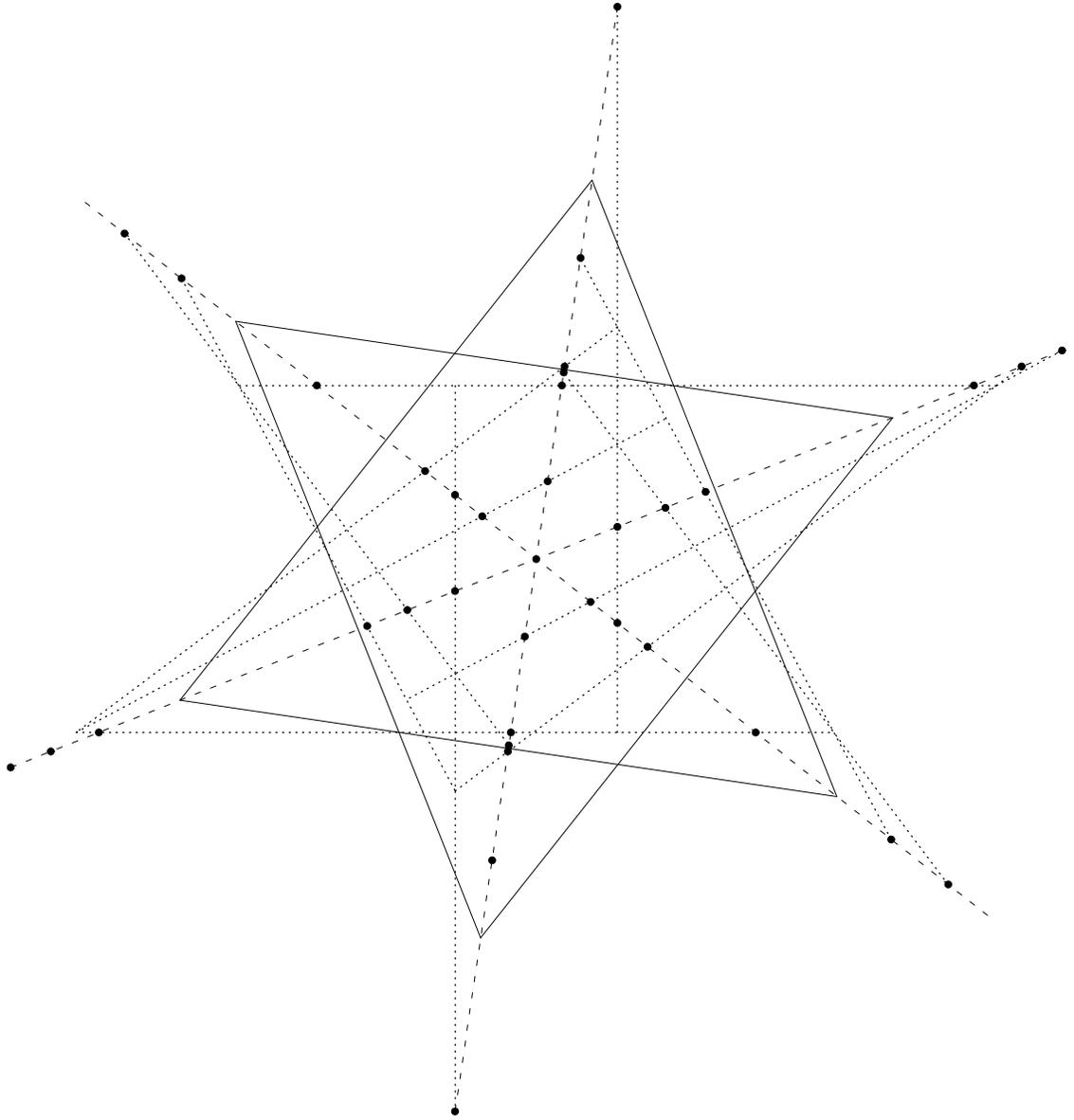
 

\caption{The Steiner Star of David}
\label{starofd}
\end{figure}

\clearpage

We next notice \cite[pp.~44-45]{C+G}:

{\bf Theorem.}  {\it The angle between the Simson-Wallace lines
of points $S$ and $S'$ on the circumcircle is half the
angular measure of the arc $S'S$.}

In Figure 1.20 the perpendiculars from
$S$ and $S'$ to $V_2V_4$ meet the circumcircle again at
$U$ and $U'$.  The Simson-Wallace lines of $S$ and $S'$
(shown dashed in the figure) are respectively parallel
to $V_1U$ and $V_1U'$.  [To see this, note that the
quadrilaterals $SV_1UV_2$ and $SF_3F_6V_2$ are both cyclic,
so that $\angle SUV_1=\angle SV_2V_1=\angle SV_2F_3=
\angle SF_6F_3$, and $UV_1$ is parallel to $F_6F_3$.]
The Simson-Wallace line turns at half the angular velocity with
which its point is describing the circumcircle, and in
the opposite sense.

\begin{figure}[ht]  
\begin{picture}(350,260)(-220,-120)
\setlength{\unitlength}{6pt}
\put(0,0){\circle{39}}
\drawline(15.12,15.54)(-15.6,11.7)(-15.6,-11.7)
(15.6,11.7)(-18,-7.5)(18,-7.5)(14,24.5)
\drawline(14.4,21.3)(0,19.5)(0,-19.5)(15.6,11.7)
\drawline(0,19.5)(7.2,6.9)
\drawline(-15.6,11.7)(-7.92,-1.74)
\dashline{0.5}(-20,-10.8)(19,18.45)
\dashline{0.5}(-2,-11.5)(16,24.5)
\put(16,11.5){\small$V_1$}
\put(-20,-8){\small$V_2$}
\put(18.4,-8){\small$V_4$}
\put(-17.3,11.3){\small$S$}
\put(-0.5,19.9){\small$S'$}
\put(-17.4,-12.2){\small$U$}
\put(-0.5,-21){\small$U'$}
\put(-15.2,-9){\small$F_6$}
\put(-10.3,-1.5){\small$F_3$}
\put(15.4
,14.4){\small$F_5$}
\end{picture}
\caption{Angle between Simson-Wallace lines of $S$ and $S'$ = $\angle UV_1U'$ =
${1\over2}\mbox{arc}UU'={1\over2}\mbox{arc}S'S$.}
\label{angle}
\end{figure}

{\bf Theorem.} {\it The Simson-Wallace line of a point $S$ on the
circumcircle of triangle $V_1V_2V_4$ bisects the segment $SV_8$,
where $V_8$ is the orthocentre of triangle $V_1V_2V_4$.}

To see this, look at Figure 1.21.  Quadrilaterals
$S_8V_1V_4U_6$ and $S_8F_5V_4F_6$ are cyclic, so that
$$\angle S_8U_6V_1=\angle S_8V_4V_1=\angle S_8V_4F_5=\angle S_8F_6F_5$$
and $V_1U_6$ is parallel to the Simson-Wallace line $F_5F_6$.

[In fact we've already proved this further up the page.  If the
perpendiculars $S_8F_3$ and $S_8F_5$ from $S_8$ to the edges
$V_1V_2$ and $V_2V_4$ meet the circumcircle $V_1V_2V_4$
again at $U_3$, $U_5$ respectively, then $V_2U_5$ and $V_4U_3$
are also parallel to the Simson-Wallace line.]

\clearpage

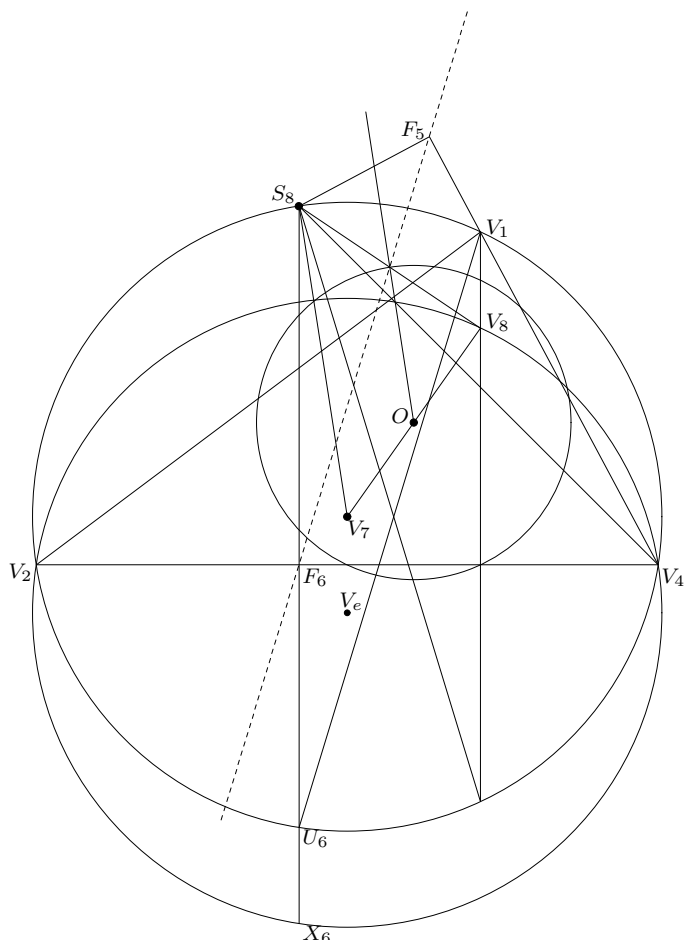
\begin{figure}[ht]  
\begin{picture}(440,350)(-250,-185)
\setlength{\unitlength}{0.7pt}
\put(-36,-51){\circle{340}}   
\put(-36,-103){\circle{340}}  
\put(0,0){\circle{170}}   
\put(-36,-51){\circle*{4}}    
\put(-36,-103){\circle*{3}}   
\put(0,0){\circle*{4}}        
\drawline(36,-205)(36,103)(-204,-77)(132,-77)  
(8.4844,154.5917)(-62,117)
\drawline(0,0)(-26,168)  
\put(-62,117){\circle*{4}}          
\drawline(36,103)(-62,-219) 
\drawline(132,-77)(-62,117)(-62,-271)   
\drawline(36,-205)(-62,117)(36,51)(-36,-51)   
\put(39,102){\scriptsize$V_1$}
\put(-219,-84){\scriptsize$V_2$}
\put(134,-87){\scriptsize$V_4$}
\put(39,51){\scriptsize$V_8$}
\put(-77,120){\scriptsize$S_8$}
\put(-60,-87){\scriptsize$F_6$}
\put(-60,-228){\scriptsize$U_6$}
\put(-60,-280){\scriptsize$X_6$}
\put(-7,155){\scriptsize$F_5$}
\put(-36,-61){\scriptsize$V_7$}
\put(-40,-100){\scriptsize$V_e$}
\put(-12,0){\scriptsize$O$}
\dashline[+50]{3}(-104,-215)(29,222) 
\drawline(-36,-51)(-62,117)   
\end{picture}
\caption{The Simson-Wallace line bisects the segment $S_8V_8$}
\label{midpoint}
\end{figure}

Moreover, if the perpendicular $S_8F_6U_6$ to edge $V_2V_4$ meets
the circumcircle $V_2V_4V_8$ at $X_6$, then the segments
$V_1V_8$, $V_7V_e$, $U_6X_6$ are equal and parallel so that,
in the triangle $S_8X_6V_8$ the Simson-Wallace line passes through the
midpoint $F_6$ of edge $S_8X_6$ and is parallel to the edge
$X_6V_8$, and so passes through the midpoint, call it $T_8$
for the time being, of the third edge $S_8V_8$.

Now look at triangle $S_8V_7V_8$.  The midpoint of edge
$V_7V_8$ [the Euler line of triangle $V_1V_2V_4$] where
$V_7$, $V_8$ are the circumcentre and orthocentre of
triangle $V_1V_2V_4$, is the 50-point centre, $O$.
We have just seen that the midpoint of edge $S_8V_8$ is
$T_8$, so that $OT_8$ is parallel to the third edge $V_7S_8$
and has length equal to half of $V_7S_8$, that is
$\frac{1}{2}R$, and $T_8$ lies on the 50-point circle.
This circle is fixed, and serves for all four (quadrated)
triangles.

\newpage

In Figure 1.22, as corresponding points $S_i$ move round
the four circumcircles their common Simson-Wallace line passes through
the point $T_8$ on the 50-point circle, where $OT_8$ is parallel
to the radii $V_{15-i}S_i$ ($i=8,4,2,1$).  If the angular
velocity of the $S_i$ is $\omega$, then that of the Simson-Wallace line
is $-\frac{1}{2}\omega$.

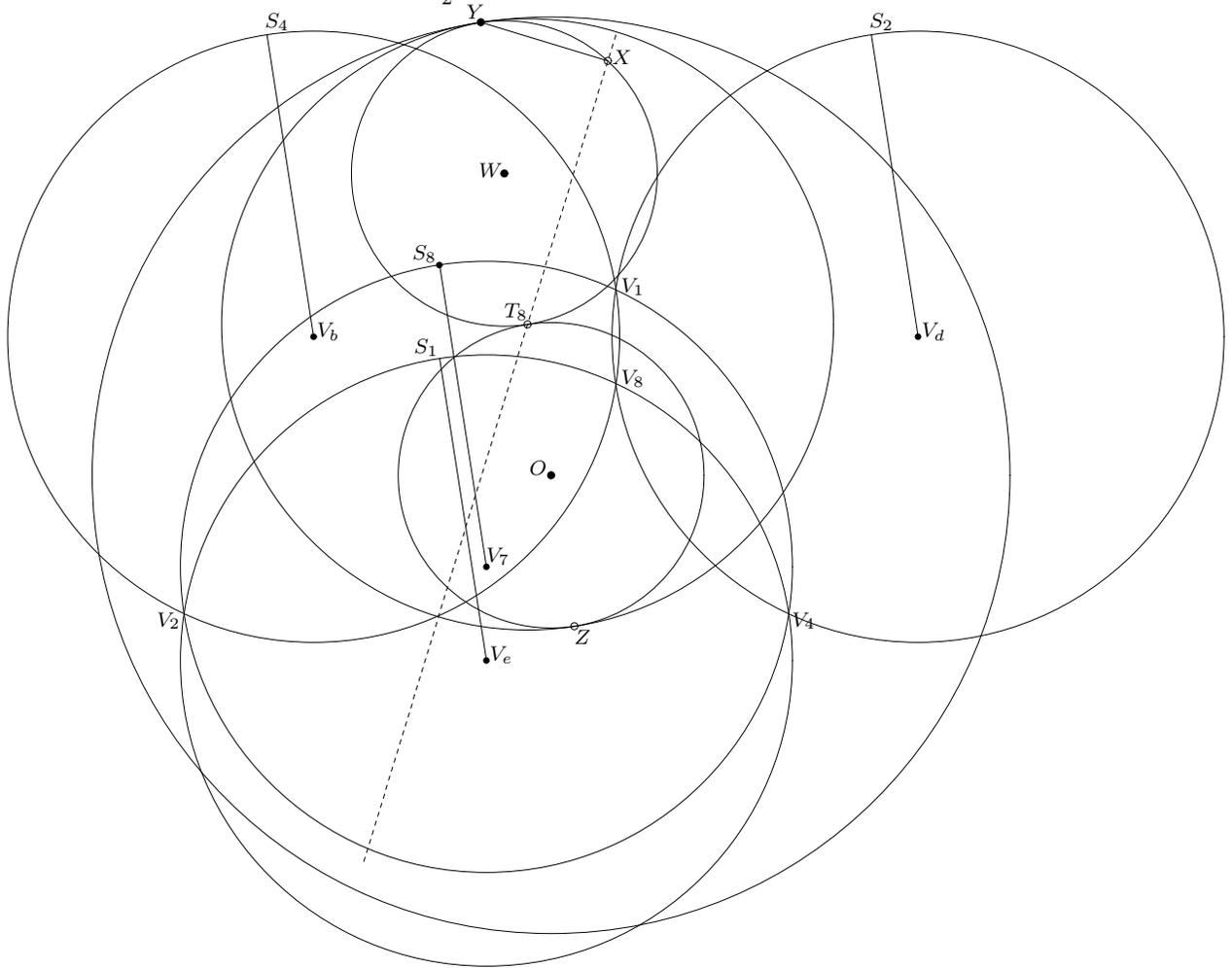
\begin{figure}[ht]  
\begin{picture}(440,350)(-205,-180)
\setlength{\unitlength}{0.7pt}
\put(-36,-51){\circle{340}}   
\put(-36,-103){\circle{340}}  
\put(-132,77){\circle{340}}   
\put(204,77){\circle{340}}    
\put(-36,-51){\circle*{3}}    
\put(-36,-103){\circle*{3}}   
\put(-132,77){\circle*{3}}    
\put(204,77){\circle*{3}}     
\put(0,0){\circle{170}}       
\put(0,0){\circle*{4}}        
\put(-62,117){\circle*{3}}          
\put(39,102){\scriptsize$V_1$}
\put(-219,-84){\scriptsize$V_2$}
\put(134,-84){\scriptsize$V_4$}
\put(39,51){\scriptsize$V_8$}
\put(-77,120){\scriptsize$S_8$}
\put(-26,88){\scriptsize$T_8$}
\put(-40,166){\scriptsize$W$}
\put(-47,254){\scriptsize$Y$}
\put(-76,68){\scriptsize$S_1$}
\put(177,249){\scriptsize$S_2$}
\put(-159,249){\scriptsize$S_4$}
\put(-36,-48){\scriptsize$V_7$}
\put(-34,-103){\scriptsize$V_e$}
\put(-130,77){\scriptsize$V_b$}
\put(206,77){\scriptsize$V_d$}
\put(-12,0){\scriptsize$O$}
\dashline[+50]{3}(-104,-215)(36,245) 
\drawline(-36,-51)(-62,117)   
\drawline(-36,-103)(-62,65)   
\drawline(-132,77)(-158,245)  
\drawline(204,77)(178,245)    
\put(-13,84){\circle{340}}   
\put(-13,84){\circle{4}}   
\put(-26,168){\circle{170}}   
\put(-26,168){\circle*{4}}   
\put(-39,252){\circle*{4}}   
\put(0,0){\circle{510}}   
\put(13,-84){\circle{4}}   
\put(13,-94){\scriptsize$Z$}
\put(34,229){\scriptsize$X$}
\put(31.5917,230.51557){\circle{4}}  
\drawline(31.5917,230.51557)(-39,252)   
\end{picture}
\caption{The deltoid is a locus and an envelope.}
\label{delt}
\end{figure}

Figure 1.22 shows the 50-point (9-point) circle,
centre $O$ and radius $\frac{1}{2}R$, together with a
concentric circle of radius $\frac{3}{2}R$ (the {\bf Steiner
circle}).  As the points
$S_i$ ($i=8$, 4, 2, 1) move round their respective circumcircles
these concentric circles remain fixed.  But the circles
with centres $T_8$ and $W$ and respective radii $R$ and $\frac{1}{2}R$
roll round in contact (at $T_8$, $Y$ and $Z$) with the
fixed circles.  The Simson-Wallace line
is a diameter of the larger of the two rolling circles. The
locus of the ends of the diameter is the deltoid, which is
also the envelope of the Simson-Wallace line, since its instantaneous
centre of rotation is the contact point $Y$, and the foot, $X$,
of the perpendicular from $Y$ to the Simson-Wallace line, lies
on the smaller of the two rolling circles.  Note that
$\angle T_8XY$ is a right angle in a semicircle.

To see this analytically, take $R=2$ and the origin at $O$,
the 50-point centre.  Then a point on the deltoid is
$(2\sin\theta-\sin2\theta,2\cos\theta+\cos2\theta)$, or,
with $t=\tan\frac{1}{2}\theta$,
$$\left(\frac{8t^3}{(1+t^2)^2},\frac{3-6t^2-t^4}{(1+t^2)^2}\right)$$
$$\frac{dx}{dt}=\frac{8t^2(3-t^2)}{(1+t^2)^3}, \quad
\frac{dy}{dt}=\frac{8t(t^2-3)}{(1+t^2)^3}, \quad
\frac{dy}{dx}=-\frac{1}{t}$$
and the equation to the Simson-Wallace line is
$$y=-\frac{x}{t}+\frac{3-t^2}{1+t^2}.$$

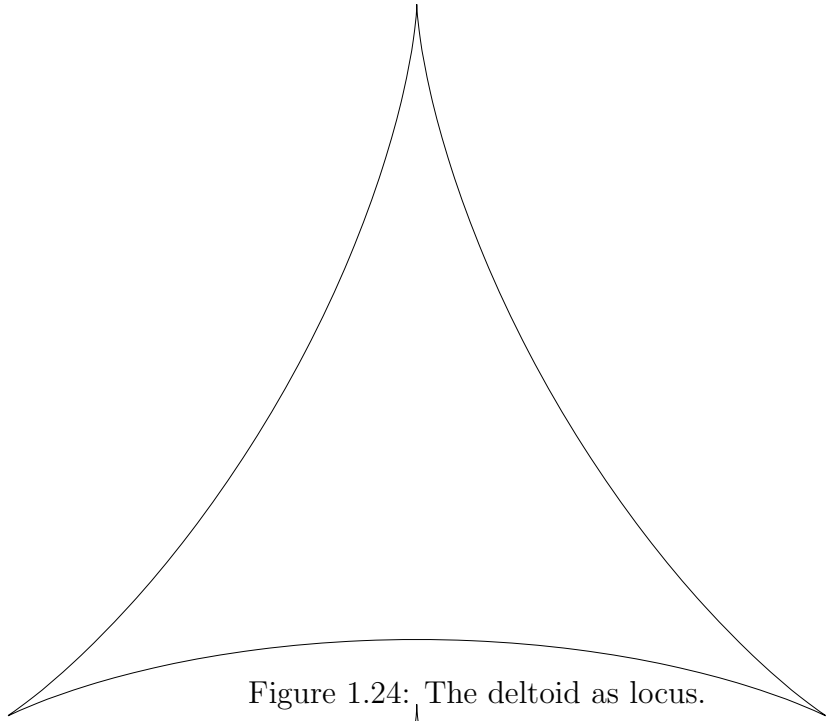
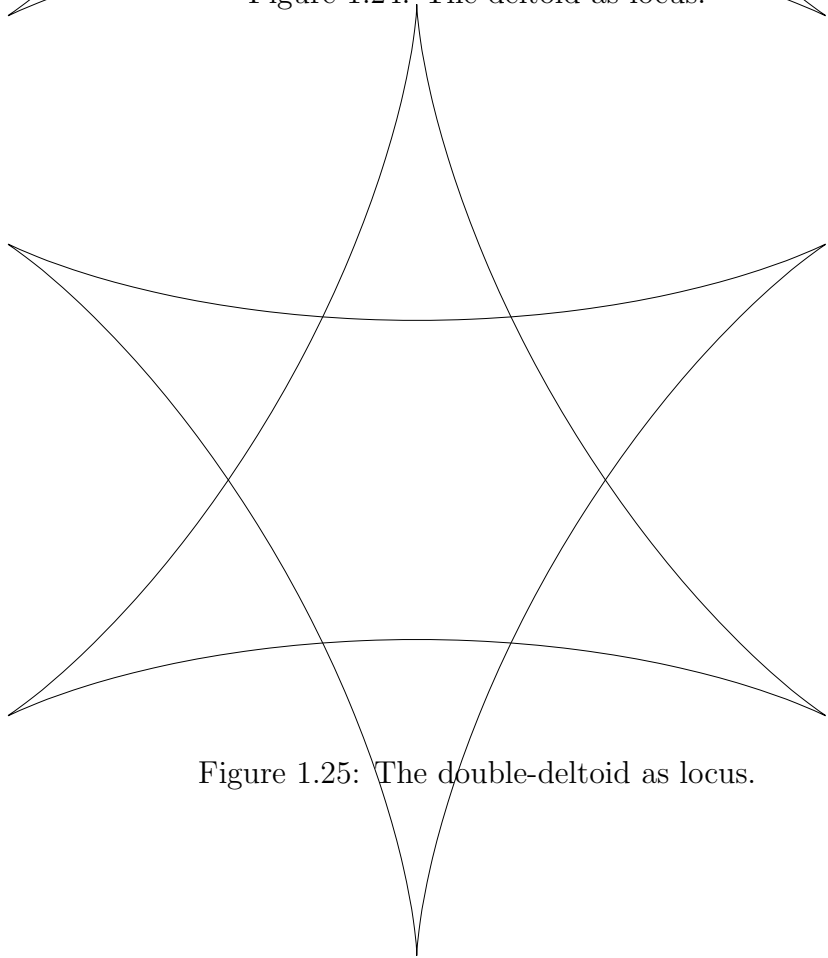
\begin{figure}[ht]  

\caption{The deltoid as envelope.}
\label{delt2}
\end{figure}

\begin{figure}[ht]  
%
\caption{The deltoid as locus.}
\label{delt3}
\end{figure}

\begin{figure}[ht]  
%
\caption{The double-deltoid as locus.}
\label{doubdelt2}
\end{figure}

\begin{figure}[ht]
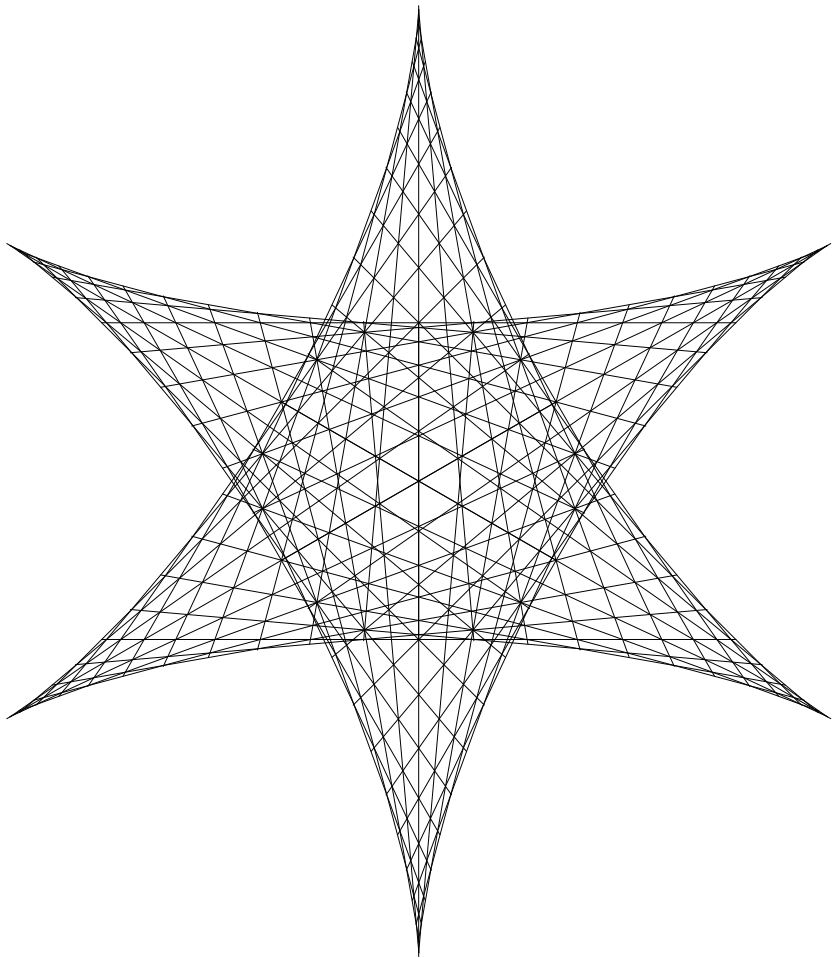
  
%
\caption{The double-deltoid as envelope.}
\label{doubdelt}
\end{figure}

\clearpage

\section{Steiner's Star of David}

\begin{figure}[h]
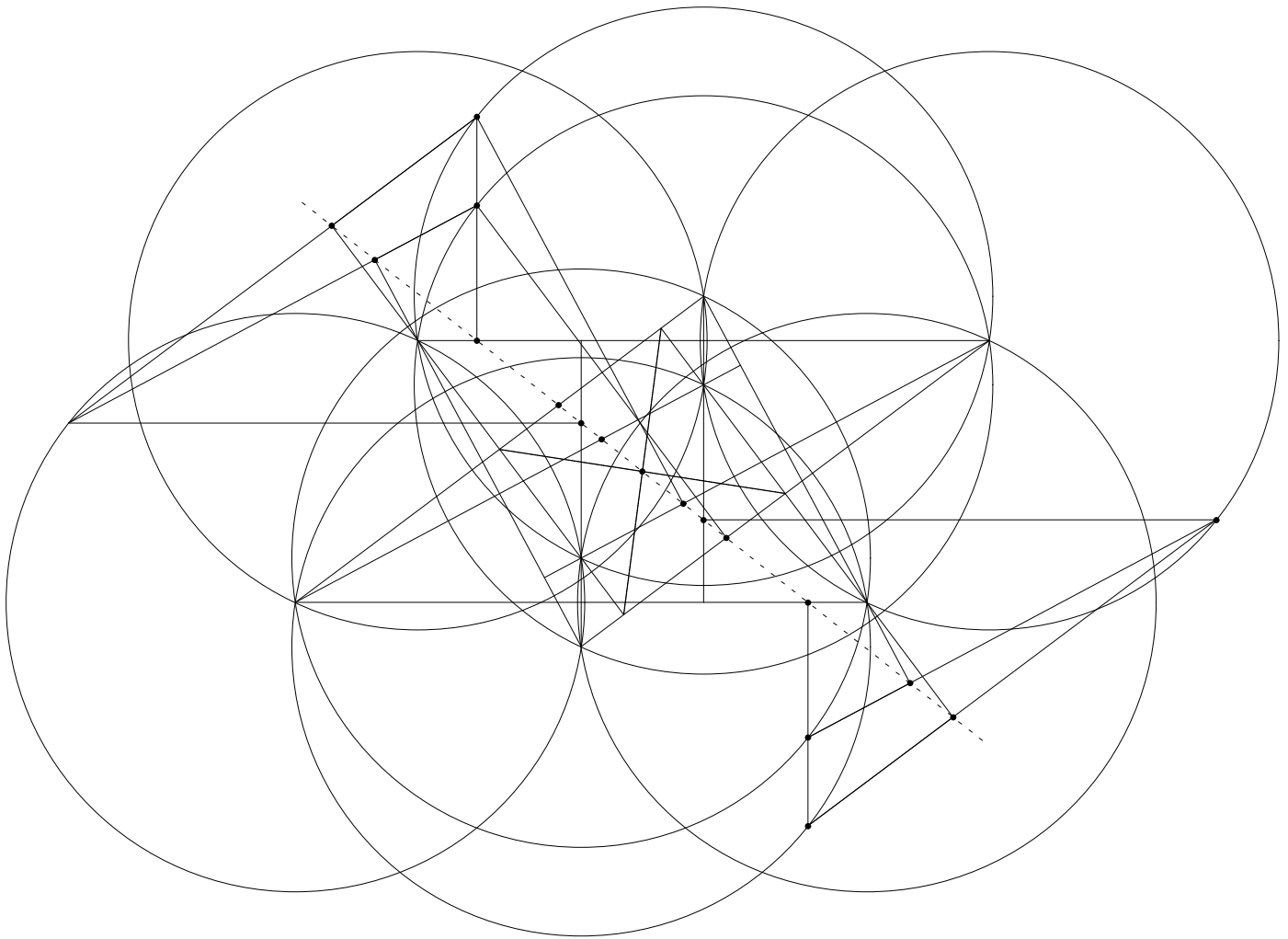
 

\caption{The $S_3$ symmetry of the triangle}
\label{threesym}
\end{figure}

It is easy to see, from the way it was generated, that the
double deltoid has the symmetry, $S_3$, of an equilateral
triangle.  There are three special positions of the
Simson-Wallace line, and three more, perpendicular to the
first three.  One of these first three is shown dashed in
Figure \ref{threesym}.  Each is a trisector of the angle
between the diagonals (pairs of diameters of the 50-point
circle, dotted lines in Figure \ref{threesym}) of one of the
three rectangles formed by the twelve edges of the triangle.

\bigskip

\newpage

\section{The BEAT}

That is: the Best Equilateral Approximating Triangle, or, together
with its twin, the Best Equilateral Approximating Twin, which is an
alternative name for the Steiner Star of David.

\newpage

\clearpage

\section{The Simson-Wallace line and the cardioid}

In \cite{Gu} we read:

\begin{quote}

Given a circle, a point $P$ on it, and a line: is it
possible to inscribe a triangle in the circle so that
the given line is the Simson-Wallace line of $P$
w.r.t.\ the triangle?  The answer is that there is an
infinity of solutions, provided that the line is
``not too far away''.

Take any point $A$ on the circle.  Let the circle on
$PA$ as diameter cut the given line in $M$, $N$.  If
$AM$, $AN$ meet the circle again in $C$ and $B$, then
it is easily shown that $ABC$ is a satisfactory triangle,
and that if $BC$ meets the given line at $L$, then
$PL$ is perpendicular to $BC$.

The condition for a solution is that it is possible to
find $A$ so that the circle on $PA$ intersects the given
line.  The envelope of such circles, for various $A$, is
the cardioid having its cusp at $P$ and vertex [$V$] diametrically
across the given circle.  The condition is that the given
line must intersect this cardioid.

\begin{figure}[ht]  
\begin{picture}(400,220)(-140,-95)
\setlength{\unitlength}{1pt}
\put(100,0){\circle{200}}
\put(64,48){\circle{160}}
\drawline(0,0)(128,96)(180,-60)(12,-84)
\drawline(128,96)(40,-80)
\drawline(224,128)(-4,-100)
\dashline[+30]{3}(144,48)(0,0)(64,-32)
\dashline[+30]{3}(0,0)(12,-84)
\put(-11,-5){$P$}
\put(129,98){$A$}
\put(32,-91){$C$}
\put(63,-41){$M$}
\put(181,-68){$B$}
\put(147,45){$N$}
\put(1,-85){$L$}
\put(202,-5){$V$}

\end{picture}
\caption{Finding triangles to fit a given Simson-Wallace line.}
\label{find}
\end{figure}
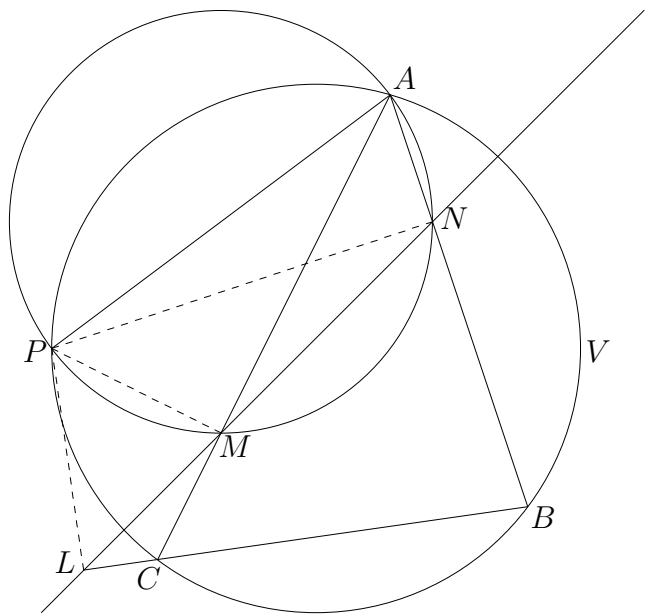

It is instructive to invert the figure w.r.t.\ $P$.  The
result is a like configuration with the roles of the line
$LMN$ and the circle $ABC$ interchanged.  The inverse of
the cardioid is the parabola, which is the envelope of lines
through A, perpendicular to $PA$, as $A$ moves along a given
line.

\end{quote}

\clearpage

\begin{figure}[ht]  
\begin{picture}(400,190)(-160,-95)
\setlength{\unitlength}{0.8pt}
\put(100,0){\circle{200}}
\put(93.3,25){\circle{193.2}}
\put(75,43.3){\circle{173.2}}
\put(50,50){\circle{144.4}}
\put(25,43.3){\circle{100}}
\put(6.7,25){\circle{51.8}}
\put(6.7,-25){\circle{51.8}}
\put(25,-43.3){\circle{100}}
\put(50,-50){\circle{144.4}}
\put(75,-43.3){\circle{173.2}}
\put(93.3,-25){\circle{193.2}}
\put(-18,-5){$P$}
\put(202,-5){$V$}
\end{picture}
\caption{The cardioid enveloped by the circles on $PA$ as diameter.}
\label{card}
\end{figure}
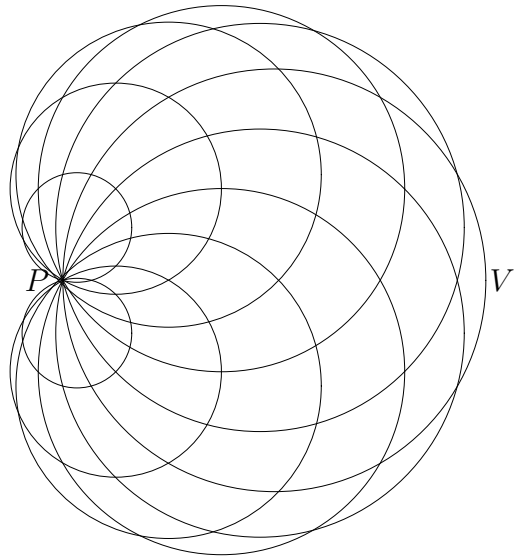

In Figure \ref{find} $\angle PNA$ and $\angle PMA$ are angles
in a semicircle, showing that $PN$, $PM$ are perpendicular to
$AB$, $AC$ respectively, and it follows from the main theorem that
$PL$ will be perpendicular to $BC$.  Figure \ref{card} shows
the cardioid enveloped by the circles on $PA$ as diameter, as
the point $A$ moves round the given circle.

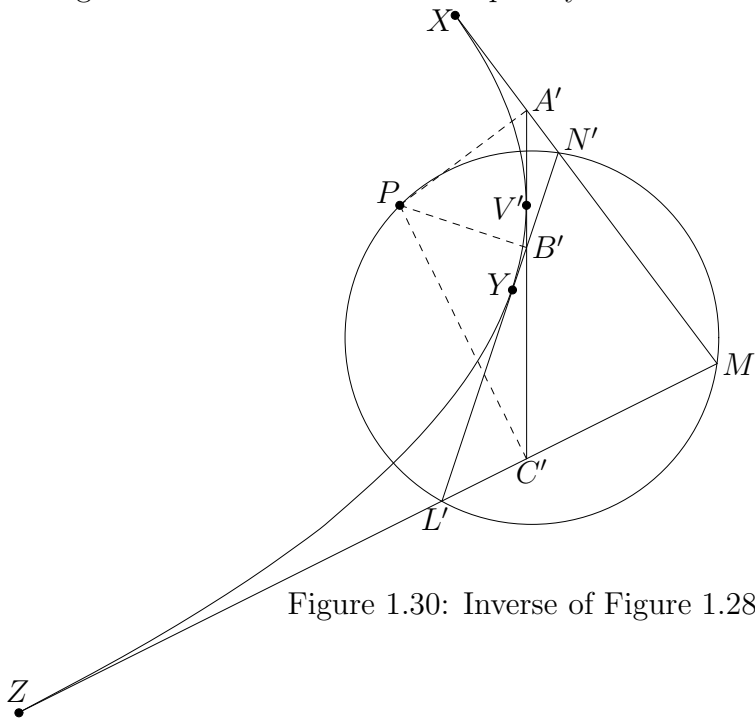
\begin{figure}[ht]  
\begin{picture}(300,185)(-180,-130)
\setlength{\unitlength}{1pt}
\put(0,0){\circle*{3}}                    
\put(48,0){\circle*{3}}                   
\put(21,72){\circle*{3}}                  
\put(42.67,-32){\circle*{3}}              
\put(-144,-192){\circle*{3}}              
\put(50,-50){\circle{141.42}}             
\drawline(48,36)(48,-96)                  
\drawline(60,20)(16,-112)(120,-60)(21,72) 
\drawline(-144,-192)(16,-112)             
\dashline[+30]{3}(0,0)(48,36)             
\dashline[+30]{3}(48,-96)(0,0)(48,-16)    
\spline
(21,72)(32,55.42)
(39.77,39.77)(44.55,25.72)(47.17,12.64)
(48,0)(47.17,-12.64)(44.55,-25.72)
(39.77,-39.77)(32,-55.42)(19.737,-73.67)
(0,-96)(-27,-120)(-33.52,-125.116)(-60,-144)
(-99,-168)(-144,-192)
\put(-9,1){$P$}
\put(50,36){$A'$}
\put(50,-20){$B'$}
\put(44,-105){$C'$}
\put(8,-123){$L'$}
\put(122,-64){$M'$}
\put(62,21){$N'$}
\put(35,-5){$V'$}
\put(10,66){$X$}
\put(33,-34){$Y$}
\put(-149,-188){$Z$}
\end{picture}
\caption{Inverse of Figure \ref{find}.}
\label{invert}
\end{figure}

\clearpage

If we invert Figure \ref{find} with respect to $P$, we obtain
Figure \ref{invert}.  The original given circle inverts into
the line $A'B'C'$ and the given line into the circle $L'M'N'$.
Notice that the four lines $A'B'C'$, $A'M'N'$, $B'N'L'$, $C'L'M'$
are the respective inverses of the circles $ABC$, $AMN$, $BNL$,
$CLM$, which all pass through $P$, and are in fact the circles
having $PV$, $PA$, $PB$, $PC$ as diameters, all of which touch
the cardioid, so that the four lines touch the parabola which
is the inverse of the cardioid, at $V'$, $X$, $Y$, $Z$.

\section{A converse theorem?}

Given a point $P$ and a line $LMN$, do the lines through
$L$, $M$, $N$ perpendicular to $PL$, $PM$, $PN$ respectively,
form a triangle whose circumcircle passes through $P$\,?

\begin{figure}[ht]  
\begin{picture}(300,235)(-210,-170)
\setlength{\unitlength}{1pt}
\put(0,0){\circle*{3}}                      
\put(48,0){\circle*{3}}                     
\put(21,72){\circle*{3}}                    
\put(42.67,-32){\circle*{3}}                
\put(-144,-192){\circle*{3}}                
\put(8,-56){\circle{113.14}}                
\put(60,-30){\circle{134.164}}              
\put(30,10){\circle{63.2455}}               

\drawline(48,36)(48,-96)                    
\drawline(60,20)(16,-112)(120,-60)(21,72)   
\drawline(-144,-192)(16,-112)               
\dashline[+30]{3}(0,0)(48,36)               
\dashline[+30]{3}(48,-96)(0,0)(48,-16)      
\spline
(21,72)(32,55.42)
(39.77,39.77)(44.55,25.72)(47.17,12.64)
(48,0)(47.17,-12.64)(44.55,-25.72)
(39.77,-39.77)(32,-55.42)(19.737,-73.67)
(0,-96)(-27,-120)(-33.52,-125.116)(-60,-144)
(-99,-168)(-144,-192)
\put(-9,1){$P$}
\put(50,36){$L$}
\put(50,-20){$M$}
\put(44,-105){$N$}
\put(8,-123){$A$}
\put(122,-64){$B$}
\put(62,21){$C$}
\put(35,-5){$V$}
\put(10,66){$X$}
\put(33,-34){$Y$}
\put(-149,-188){$Z$}
\end{picture}
\caption{Converse of Simson-Wallace theorem.}
\label{converse}
\end{figure}
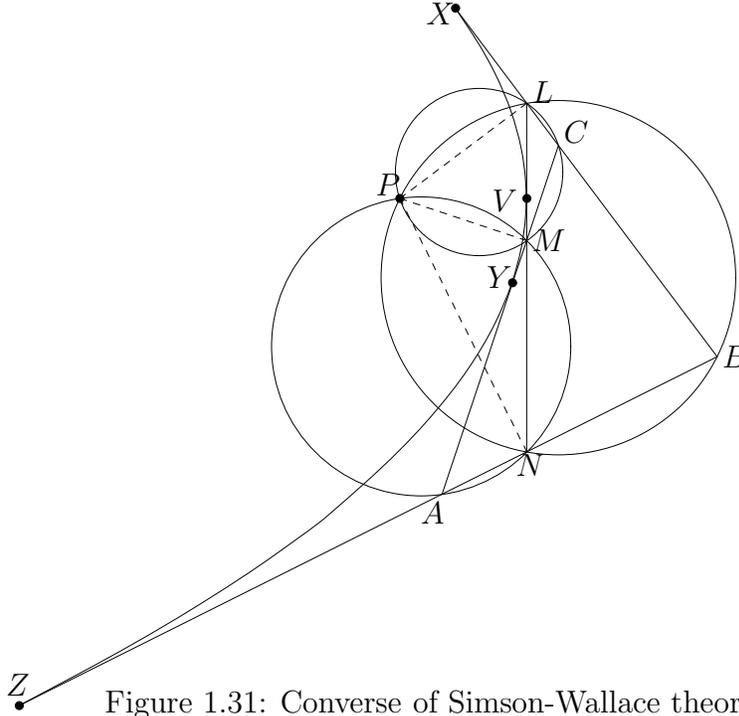

From the right angles at $L$, $M$ and $N$ we see
that $PANM$, $PNBL$, $PLCM$ are cyclic quadrilaterals,
inscribed in circles on $PA$, $PB$, $PC$ as diameters.
Hence

\centerline{$\angle LPC=\angle LMC=\angle NMA=\angle NPA$, so that}

\centerline{$\angle APC=\angle NPL=\pi-\angle NBL=\pi-\angle ABC$}

and $P$ lies on the circumcircle of $ABC$.

Notice that the edges of triangle $ABC$ are tangents to
the parabola having $P$ as focus and $LMN$ as tangent
at the vertex $V$.  The following theorems are well known
to those who well know them (see \cite[pp.26--39]{Du}
for example).

{\bf Theorem.}  The foot of the perpendicular from the
focus onto a tangent of a parabola lies on the tangent
at the vertex.

That is, given a point $P$ and a line, then, as $L$
moves along the line, the perpendicular to $PL$
through $L$ envelopes a parabola.

{\bf Theorem.}  Three tangents to a parabola form a triangle
whose circumcircle passes through the focus of the parabola.

{\bf Steiner's Theorem.}  The orthocentre of a triangle
circumscribing a parabola lies on the directrix

These form a sort of dual of the theorem (see the fifth
proof in Chapter 2 below) concerning the orthocentre
and the circumcircle of a triangle inscribed in a
rectangular hyperbola.

Compare the dual of Feuerbach's theorem.  If we consider
a family of four-line conics, that is the system of conics
touching four lines, in the case where one of the lines is
the line at infinity, then the conics will be parabolas.

\section{The Steiner Line Theorem.}

If, instead of dropping perpendiculars onto the edges from
a point on the circumcircle, we reflect the point in the edges
then we get (sets of) three points lying on a line parallel
to the Wallace line and twice as far from from the point(s)
on the circumcircle(s).  This theorem is due to Steiner.
Note that the circumcircles, in general, are reflexions of
one another in the edges of the triangle (i.e., the orthocentric
quadrangle).

Figure 1.27 should be completed by adding the other two special
positions of the Simson-Wallace line at angles of $60^{\circ}$
with the dashed line already exhibited. These each contain 12 points
as well as the Centre of the triangle.  Then there will be 8
parallels to each of these three lines, passing one through
each of the eight vertices.  Compare Figure 1.19, which also
represents the BEAT twins.

\newpage

\section{The Steiner lines and the Droz-Farny theorem.}

[The reference to the Droz-Farny theorem points to the fact that
the doubled Wallace lines of diametrically opposite points on the
circumcircle are perpendicular lines through the orthocentre.
See next section, 1.22.]

Take, as before, an arbitrary point, 7B, on the circumcircle of the
triangle 124 and reflect it, and the circumcircle, in the
edge 24.  This gives a point 1A on the circumcircle of
the triangle 247.  Reflect also in the edges 41 and
12 giving points 2A and 4A on the respective
circumcircles of triangles 147 and 271.

Then the points 1A, 2A, 4A lie on a Steiner line, line
A say, parallel to the Wallace line of the point 7B
and twice its distance from 7B.  As we know that the Wallace
line bisects the segment joining 7 to 7B the Steiner line A will
pass through the vertex (orthocentre) 7.

Now consider the Steiner line of the point 1A
which lies on the circumcircle of triangle 247.  The reflexions
of 1A in the edges 47, 72, and 24 are
the points 2B, 4B, 7B, which lie on a line,
B say, through vertex 1.

Similarly, for the Steiner line of the point 2A,
which lies on the circumcircle of triangle 471.  The
reflexions of 2A in the edges 71, 14, and 47
are the points 4C, 7B, and  1C.  These three
points lie on the line C which passes through vertex 2.

Finally, for the Steiner line of the point 4A,
which lies on the circumcircle of triangle 712.  The
reflexions of 4A in the edges 12, 27, and 71
are the points 7B, 1D, and 2D, which lie on the
line D which passes through vertex 4.

The previous paragraph should not have started ``Finally'' since
there are still some (how many?) points on circumcircles which
reflect into other points on other circumcircles, some of which
are new \ldots ?

I suspect that, in general, we have an infinite ``trisequence".
Each point has 3 followers, of which one is also a precedessor.

In the following table it appears that 7H = 7E and 7I = 7F.
Also 7G is diametrically across circle 124 from 7B.  So we
have some finite cycles.  Line O will coincide with line A'
(perpendicular to line A; see Section 1.22 and Fig. 1.34).

This is because, in the notation given below, I chose 7B with
$\theta = \pi/2$ !  I'll try again.

\newpage

Here's a new table and figure based on a new 7B which is not
a rational multiple of $\pi$ from any of the vertices 1, 2, or 4.
Lack-a-day!! That's not true!! The new $\theta=\arctan(161/240)
=\pi/2 - \beta$ !!

For the old table, see the Apocrypha below.

\begin{center}

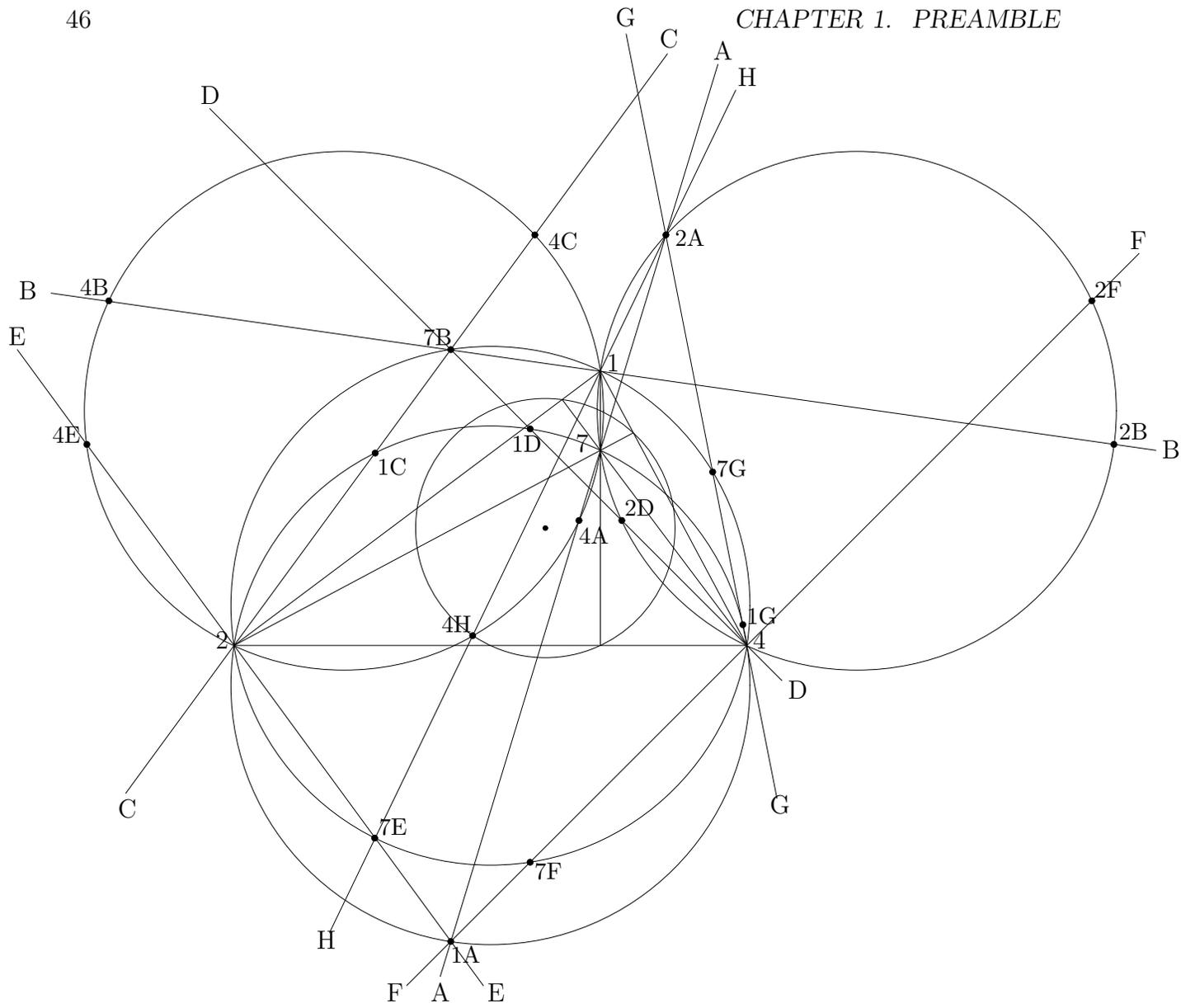
\caption{Start of new ``trisequence''.}
\label{olddoubwall}
\end{figure}

\newpage

Note that if a point lies on a circumcircle, then the
midpoint of its join to the corresponding orthocentre
lies on the Central Circle (9-pt circle) of the triangle
(orthocentric quadrangle).  If the figure is scaled so
that the radius of the Central Circle is 1, then the
midpoints in the following list have coordinates
\vspace{-15pt}

$$\left(\frac{r^2-s^2}{r^2+s^s},\frac{2rs}{r^2+s^2}\right)$$
\vspace{-20pt}

where $r$, $s$ are given in the right hand column.

\small
\begin{center}
\begin{tabular}{ccccc}
point & coords & orthocentre & midpoint & (r,s) \\
7B  & (--62,117) & 7 (36,51) & (--13,84) &  (6,7)\\
1A  & (--62,--271) & 1 (36,103) & (--13,--84) & (6.-7) \\
2A  & (22822,55541)/$17^2$ & 2 (--204,--77) & (--18067,-16644)/$17^2$ & (57,146) \\
4A  & (22,5) & 4 (132,--77) & (77,--36) & (9,--2) \\
2B  & (9314,1373)/25 & 2 (--204,--77) & (2107,276)/25 & (46,3)\\
4B  & (--286,149)/289 & 4 (132,--77) & (--74,36) & (2,9) \\
4C  & (--2014,55541)/$17^2$ & 4 (132,--77) & (18067,16644) & (73,--57) \\
1C  & (--806074,356057)/$85^2$  & 1 (36,103) & (--272987,550116)/$85^2$ & (413,666) \\
1D  & (--10,65) & 1 (36,103) & (13,84) & (7,6) \\
2D  & (50,5) & 2 (--204,--77) & (--77,--36) & (2,--9) \\
4E  & (--7514,1573)/25  & 4 (132,--77) & (--2107,--276)/25 & (3,46) \\
7E  & (--806074,--1468707)/$85^2$ & 7 (36,51) & (--272987,--550116)/$85^2$ & (413,--46) \\
7F  & (--10,--219)/289 & 7 (36,51) & (13,--84) & (7,--6) \\
2F  & (359,149)/289 & 2 (--204,--77) & (79,36) & (9,2) \\
7G  & (31670,10629)/$17^2$ & 7 (36,51) & (21037,12684) & (151,42)  \\
1G  & (10799938,--5289871)/$17^4$ & 1 (36,103) & (6903347,--1656396) & (2646,313) \\
4H  & (--99465434,--147393187)/$1445^2$ & 4 (132,--77) & (88976933,--154085556)/$1445^2$ & (11523,6686) \\
2I  &   & 2 (--204,--77) &  \\
7I  &   & 7 (36,51)   &    &    \\
1J  &   & 1 (36,103) &  \\
7K  &  & 7 (36,51) &  &  \\
1K  &  &  1 (36,103) &  &  \\

\end{tabular}
\end{center}
\normalsize

\newpage

Let the angles of the triangle 124 be $\alpha$, $\beta$ and
$\gamma$, so that $\alpha+\beta+\gamma=\pi$.

We will measure arc-lengths on the four circles counter-clockwise
starting at 1 on circle 124, so that vertex 2 is at $2\gamma$ and
vertex 4 is at $2\alpha+2\gamma$.  Start at 7 for each of the circles
724, 741 and 712, so that vertices 7, 2, and 4 on circle 724 are at
0, $\pi-2\beta$, and $\pi+2\gamma$; vertices 7, 4, and 1 on circle
741 are at 0, $\pi-2\gamma$, and $\pi+2\alpha$; and vertices
7, 1, and 2 on circle 712 are at 0, $\pi-2\alpha$, and $\pi+2\beta$.

Let the point 7B (= P7) be at arc-length $\theta$ from vertex 1 on
circle 124, so that arc 7 1A on circle 724 is $\pi-2\beta+2\gamma-\theta$
and arc 7 2A on circle 741 is
$2\pi-(\pi-2\alpha+\theta)=\pi+2\alpha-\theta$.  Also arc 7 4A
on circle 712 is $2\pi-(\theta-(\pi-2\alpha))=3\pi-2\alpha-\theta$.

Next we reflect point 1A (= P71) in edges 24, 47 and 72,
giving points P7, P712 (= 2B) and P714 (= 4B).  These last two
points are at arc-lengths $\pi-2\gamma+2\beta+\theta$ from
7 on circle 741 and $2\pi-(2\gamma-\theta)-(\pi-2\beta)=
\pi+2\beta-2\gamma+\theta$ from 7 on circle 712.

Now reflect point 2A (= P72) in edges 14, 47 and 71,
giving points P7, P721 (= 1C) and P724 (= 4C).  These last
two are at arc-lengths $\theta+\pi-2\alpha$ from 7 on
circle 724 and $\pi-2\alpha+\theta$ from 7 on circle 714.

Then 4A (= P74), reflected in edges 12, 27 and 71, gives
points P7, P741 (= 1D) and P742 (= 2D)

Note:  4H appears to be on the Central Circle, but its distance
from the (9-pt) Circle is 85.159353528247..., not 85.

{\bf Thanks to Andrew Bremner,}

for discovering that there are four 3-cycles;
i.e. 12 points which each reflect into
two others.  They are 3 points on each of the four circumcircles,
diametrically opposed to the four vertices 1, 2, 4, and 7 of
any orthocentric quadrangle.  In our particular numerical example,

the point (--108,--205) is diametrically opposed to the vertex
1 (36,103) on the 7-circle (124-circle), centre (--36,--51).
Its reflexions in edges 24, 41, 12 are (--108,51), (372,51),
(--300,51), which lie respectively on the 1-circle, 2-circle,
and 4-circle.  The first of these is diametrically opposite to
the point (36,--257) on the 1-circle; the other two are both
vertex 1. They all lie on the line $y=51$ which passes through
the vertex 7 (36,51) and is parallel to the edge 24.

\newpage

{\bf Thanks also to Alex Fink,}

who observes that {\bf every} point leads to three 6-cycles.
For example, start from a point on the 7-circle and reflect it
in the edges 14, 24, 47, 14, 24, 47; a hexagonal path whose
edges are parallel to those of a 3-cycle, traversed twice.

He also noted that the six collinearities 12,21,12$'$,21$'$; \
14,41,14$'$,41$'$, \\ 17,71,17$'$,71$'$; \ 24,42,24$'$,42$'$; \
27,72,27$'$,72$'$; \ 47,74,47$'$,74$'$; concur in threes at
the vertices (--108,--30) \ (612,231) \ (-396,231) \
(--108,--153) of an orthocentric quadrangle, homothetic to,
and three times the size of, the original quadrangle 1247.
The centre of perspective is, of course, the Centre of the
triangle.  Moreover

{\bf The Central Circle of the Trebled Triangle is the Steiner
Circle of the original triangle.}  And the centroids of the
Trebled Triangle are the de Longchamp points of the original
triangle (orthocentric quadrangle).


\begin{figure}[h]
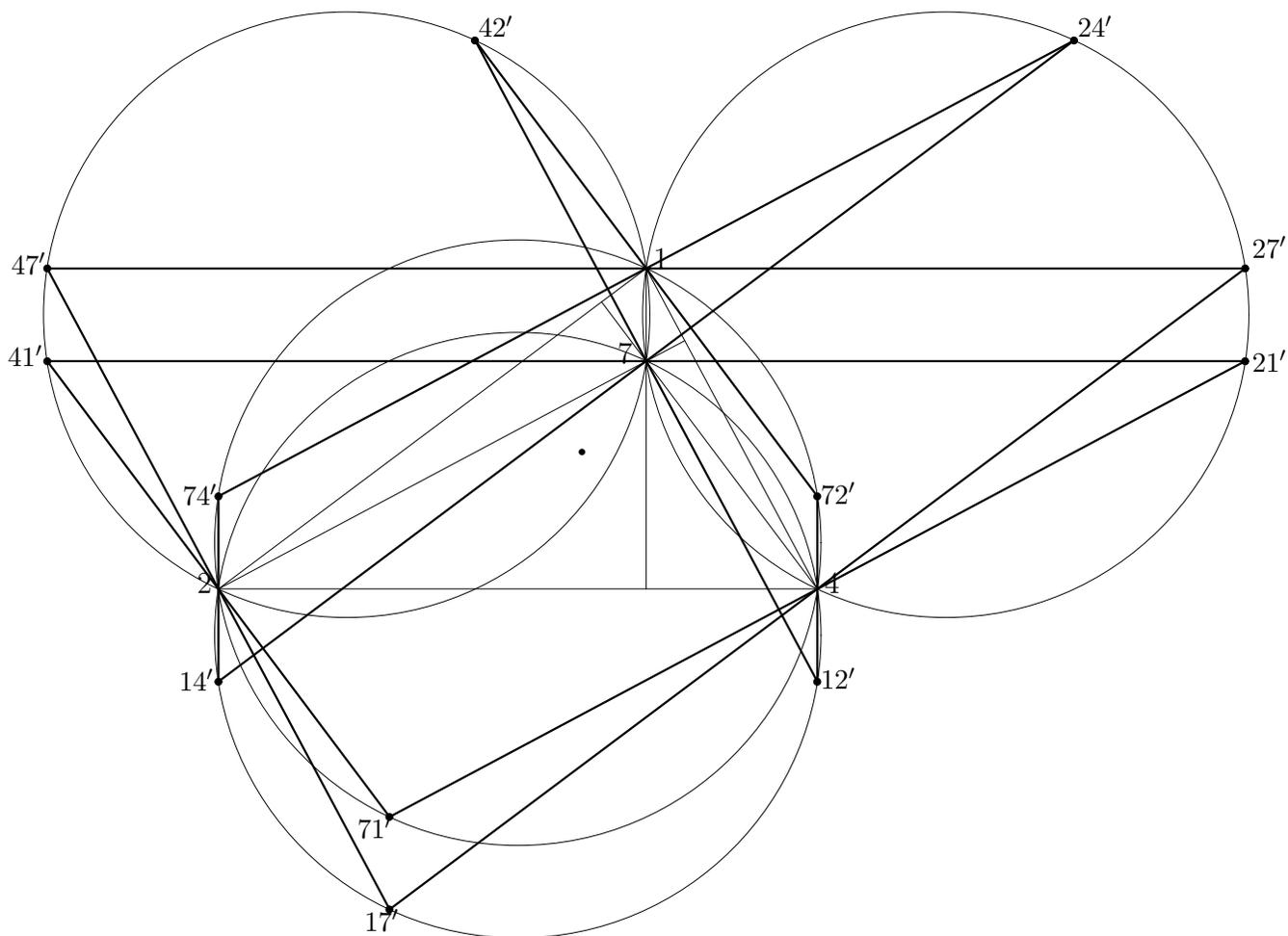
 

\caption{Four three-cycles.}
\label{4 3-cycles}
\end{figure}

\newpage

In the two-digit labels of the twelve points the first digit
refers to the circle on which the point is situated and the
second is the vertex to which it is opposite.  For example,
the point 47$'$ is on the 4-circle, or 127-circle, and is
diametrically opposite to the vertex 7.

\begin{figure}[h]
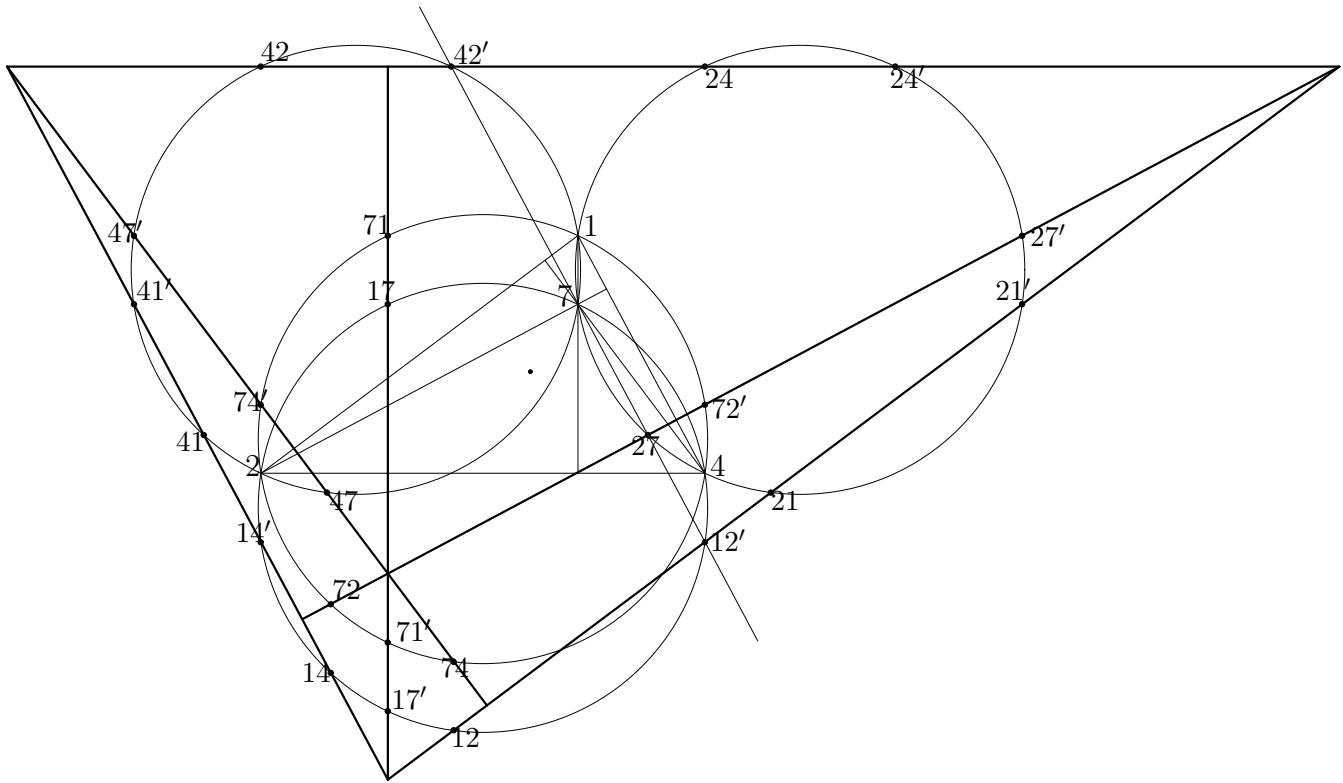
 
%
\caption{A trebled triangle.}
\label{trebled}
\end{figure}

\newpage

\begin{figure}[h] 
%
\caption{Several lines.}
\label{several}
\end{figure}

\newpage

{\bf{\Large Apocrypha}}

This is the first numerical example, where $\theta=\pi/2$,
leading to non-typical coincidences.

\begin{center}
%
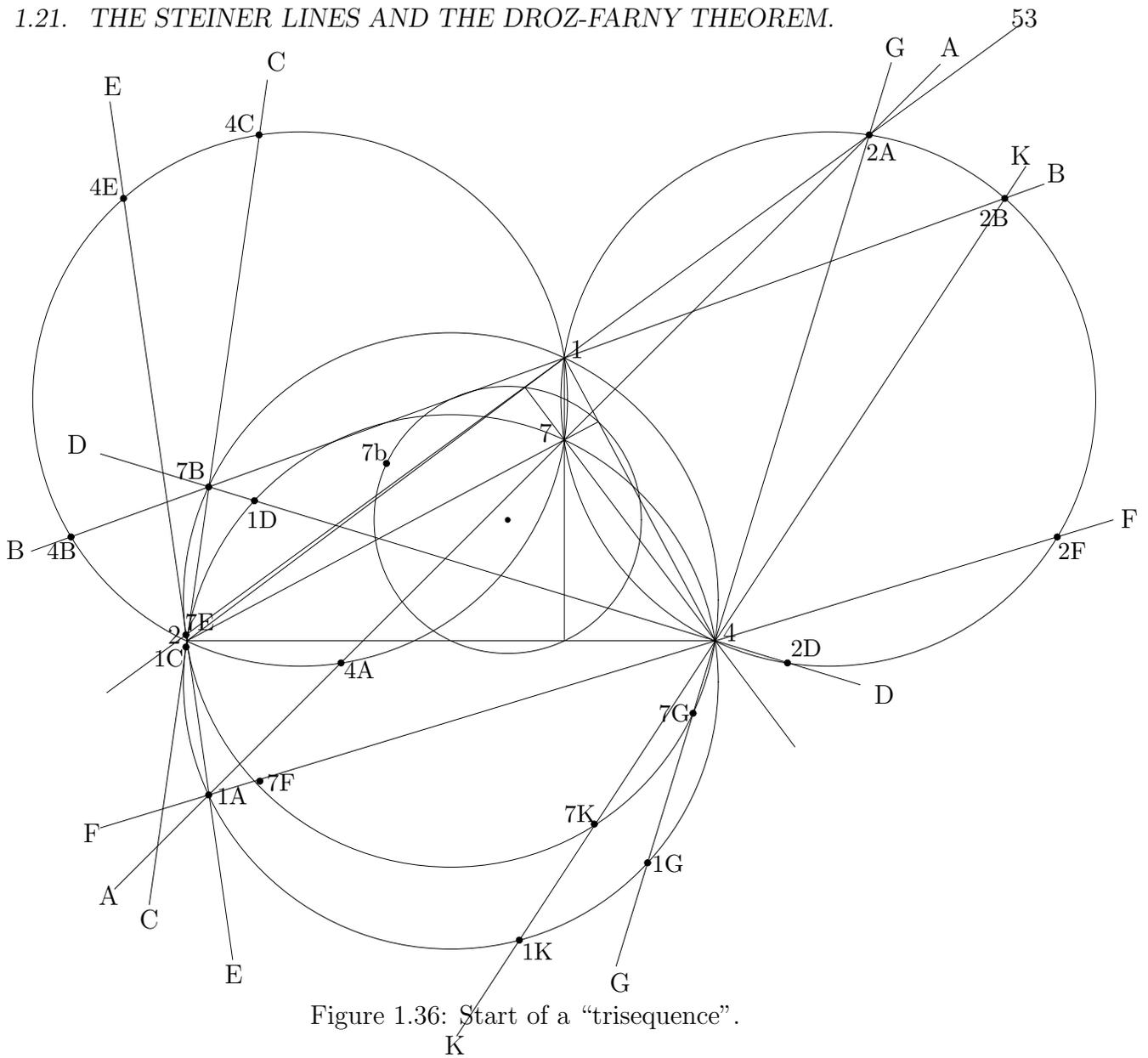
\caption{Start of a ``trisequence''.}
\label{doubwall}
\end{figure}

\newpage

Note that if a point lies on a circumcircle, then the
midpoint of its join to the corresponding orthocentre
lies on the Central Circle (9-pt circle) of the triangle
(orthocentric quadrangle).  If the figure is scaled so
that the radius of the Central Circle is 1, then the
midpoints in the following list have coordinates
\vspace{-15pt}

$$\left(\frac{r^2-s^2}{r^2+s^s},\frac{2rs}{r^2+s^2}\right)$$
\vspace{-20pt}

where $r$, $s$ are given in the right hand column.

\vspace{-18pt}
\small
\begin{center}
%
\caption{The neighborhood of vertex {\bf2} and point {\bf7B}}
\label{enlarge}
\end{figure}

\newpage

\section{Doubling Wallace(-Simson) and Droz-Farni.}

It is known (see elsewhere) that as a point moves round the
circumcircle, the associated Wallace(-Simson) line rotates
with half the angular velocity and in the opposite sense.
So also does the doubled Wallace(-Simson) line.  If we draw
such lines for two points on the circumcircle which are
diametrically opposed they will form a pair of perpendicular
lines through the orthocentre, to which we may apply the
Droz-Farni theorem.

\begin{figure}[h]
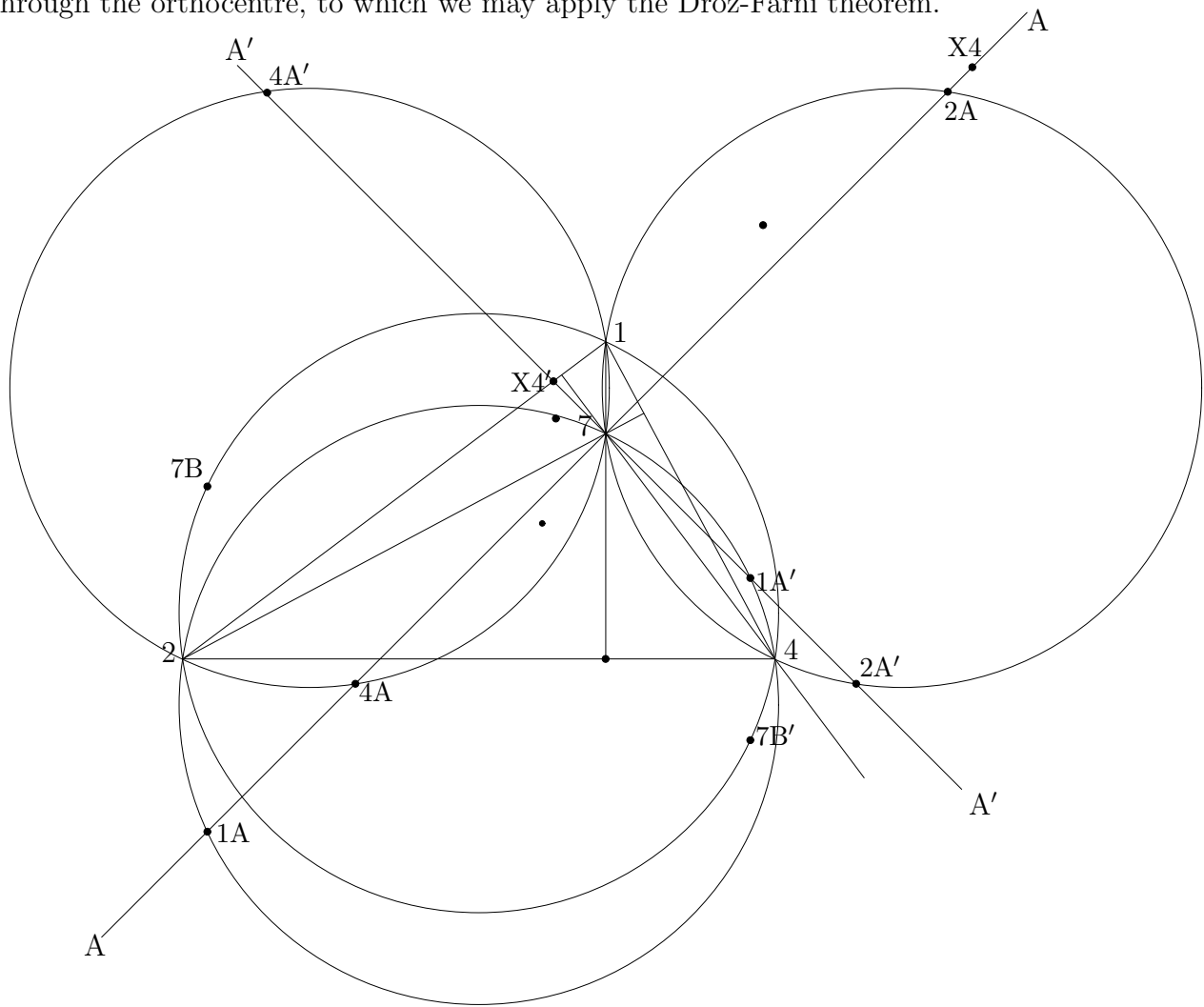
 

\caption{Two perpendicular doubwall lines, A and A$'$.}
\label{wall-droz}
\end{figure}

In Fig.\ref{wall-droz} the point 7B generates the doubwall
line A, and 7B$'$ generates the perpendicular line A$'$.
It appears that line A$'$ is the same as both line O,
and line P, passing through 1O, 2D, and  4C on the one hand
and 2D, 4C and 7P on the other.  Is this true\,?  If so,
is it just a coincidence\,?  Evidently 7B$'$ is the same as 7G.

\newpage

\chapter{Five proofs that a triangle has an orthocentre.}

{\bf First} we'll start with Euclid, who showed that the
perpendicular bisectors of the edges of a triangle
concurred at a point, the circumcentre, equidistant
from the vertices of the triangle.  On drawing, through
each vertex of a triangle, a parallel to the opposite edge,
we produce three parallelograms which form a twice-sized
triangle.  Since opposite edges of a parallelogram are
equal in length, the larger triangle has the vertices of
the original triangle for the midpoints of its edges,
and the altitudes of the original triangle are the
perpendicular bisectors of the edges of the larger
triangle, and they consequently concur.

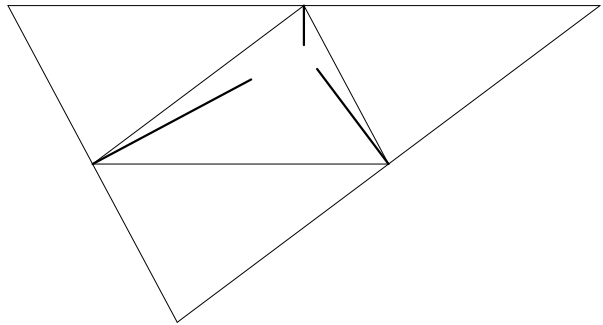
\begin{figure}[h] 
\begin{picture}(440,120)(-230,-55)
\setlength{\unitlength}{1pt}
\drawline(-80,0)(-112,60)(112,60)(-48,-60)
(-80,0)(0,60)(32,0)(-80,0)
\thicklines
\drawline(32,0)(5,36)
\drawline(-80,0)(-20,32)
\drawline(0,60)(0,45)
\end{picture}
\caption{Perpendicular bisectors of edges of large
triangle are altitudes of small triangle}
\label{proof1}
\end{figure}

\clearpage

{\bf Second} we draw the {\bf edge-circles} of the triangle,
that is, circles having an edge of the triangle for
diameter.  In Figure \ref{proof2} the circles with
diameters $V_2V_4$, $V_4V_1$, $V_1V_2$ cut the other
other edges in $D_6$, $D_5$, $D_3$.  These points are
unique, since the angle in a semi-circle is a right-angle,
and $V_1D_6$, $V_2D_5$, $V_4D_3$ are altitudes of the
triangle.  But these are also radical axes of pairs
of the edge-circles, and therefore concur in their
radical centre.

\begin{figure}[h]
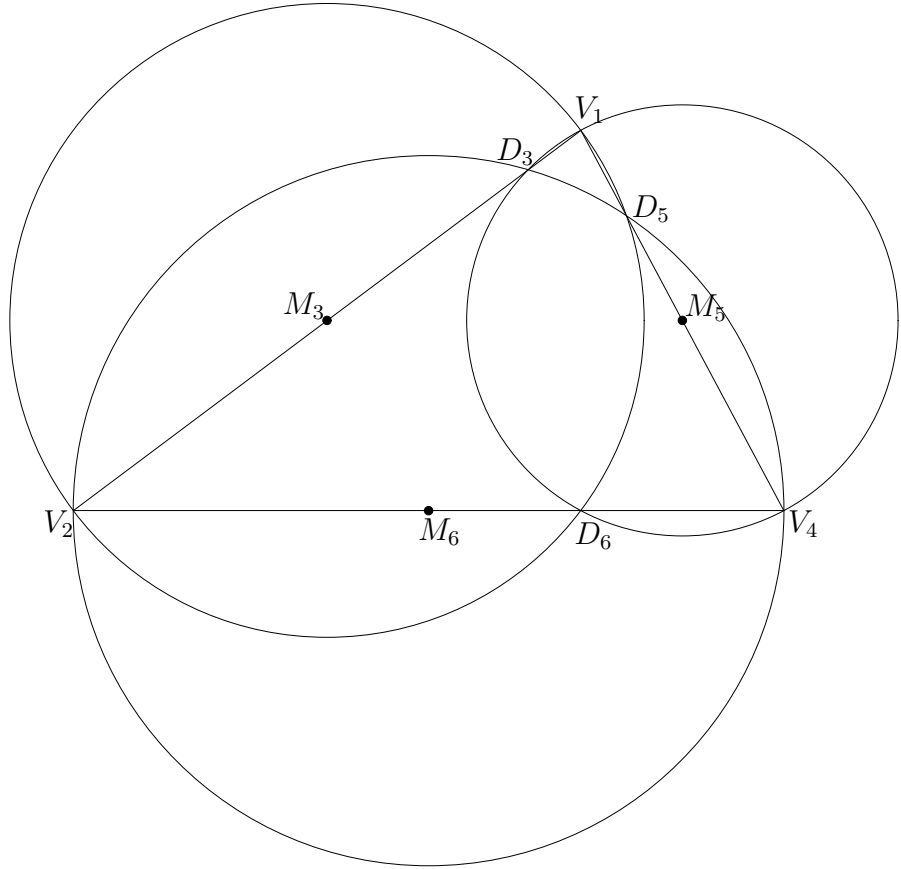
 

\caption{Edge-circles}
\label{proof2}
\end{figure}

\clearpage

{\bf Third}: this proof also uses edge-circles, though we have
to extend our definition of an edge of the triangle
$V_1V_2V_4$ to include $V_1V_8$, $V_2V_8$ and $V_4V_8$.
First, let the altitudes $V_2D_5$ and $V_4D_3$ meet at $V_8$,
so that $V_1D_3V_8D_5$ is cyclic with $V_1V_8$ as a diameter.
This new edge-circle, taken with new edge-circles on
$V_2V_8$ and $V_4V_8$ as diameters, again have altitudes
for radical axes and orthocentre for radical centre.
See Figure \ref{proof3}.

\begin{figure}[h]
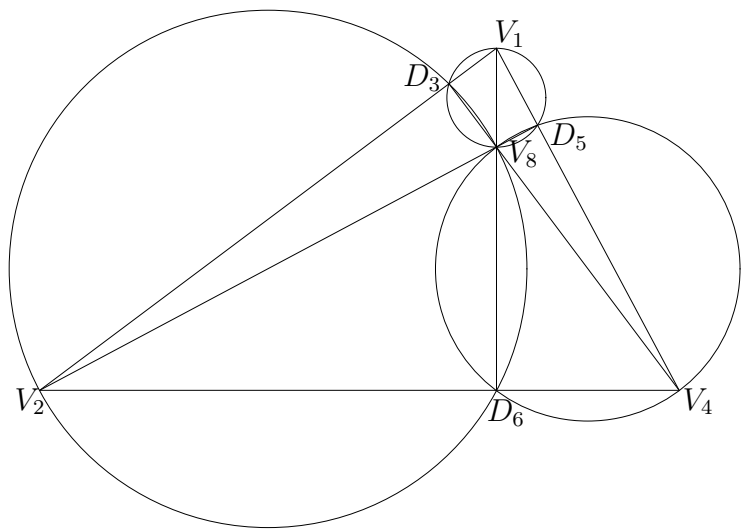
 
%
\caption{More edge-circles}
\label{proof3}
\end{figure}

\vspace*{-10pt}

{\bf Fourth}: the (original) {\bf Clover-leaf Theorem}.

\begin{figure}[h]
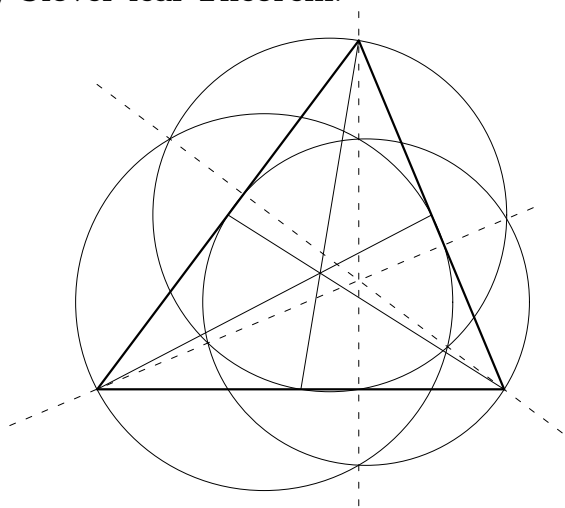
  
%
\caption{The Clover-leaf theorem.}
\label{clover}
\end{figure}

In Figure \ref{clover} the radical centre of the three
circles drawn with the medians of a triangle as diameters is the
orthocentre of the triangle.  The radical axes of pairs of circles
are the altitudes, shown dashed.

\begin{figure}[h]  
%
\caption{Proof of the clover-leaf theorem.  For future
reference, draw in the perpendiculars from $M_6$ onto
the other edges, $V_1V_2$, $V_1V_4$, of the triangle;
where do they land?}
\label{proof}
\end{figure}

To prove the Clover-leaf theorem, we show that the radical axis
of two of the circles, for example, in Figure \ref{proof}, those
with diameters $V_1M_6$ and $V_4M_3$, is an altitude of the triangle.
The line of centres (dashed in Figure \ref{proof}) is the join of the
midpoints of $V_1M_6$ and $V_4M_3$ and so is midway between
$M_3M_6$ and $V_1V_4$, and parallel to them.  So the radical
axis is perpendicular to $V_1V_4$, and it suffices to show
that $V_2$ is on the radical axis.  First note that the
$V_4M_3$-circle passes through $K$, the midpoint of $V_2D_6$
(we shall call such a point a {\bf midfoot point}), since $KM_3$
is parallel to $D_6V_1$, making $V_4KM_3$ a right angle.
Now the square on the tangent from $V_2$ to the $V_1M_6$-circle
equals $V_2M_6\cdot V_2D_6=(\frac{1}{2}V_2V_4)(2V_2K)=V_2V_4\cdot V_2K$
which equals the square on the tangent from $V_2$ to the
$V_4M_3$-circle, and indeed $V_2$ is on the radical axis
of the two circles.

\newpage

{\bf Fifth}, a pair of perpendicular lines $V_1V_8$, $V_2V_4$
intersects another pair $V_2V_8$, $V_1V_4$ in the four
points $V_1$, $V_2$, $V_4$, $V_8$.  These two pairs are
rectangular hyperbolas.  Any conic passing through these four
points is a rectangular hyperbola.  In particular, the
pair of lines $V_1V_2$, $V_4V_8$ is a rectangular hyperbola,
so that this pair of lines is also perpendicular.

Somewhat more generally, choose axes and scale so that the
three vertices of the triangle lie on a rectangular hyperbola,
say at $(t_1,1/t_1)$, $(t_2,1/t_2)$ and $(t_3,1/t_3)$.

\centerline{[This can be done in infinitely many ways;
see next paragraph.]}
\vspace*{-7pt}
Then the slope of the join of the last two is
$\{1/t_3-1/t_2\}/\{t_3-t_2\}=-1/t_2t_3$.  The
equation to the line through the first of the three vertices
and having a perpendicular slope is $y-1/t_1=t_2t_3(x-t_1)$.
This cuts the hyperbola again in the point
$(-1/t_1t_2t_3,-t_1t_2t_3)$.  The symmetry shows that the
other two perpendiculars also pass through the same point.

Moreover, if the circumcircle of the triangle cuts the hyperbola
again in the point $(t_4,1/t_4)$ and has equation
$x^2+y^2-2gx-2fy+c=0$, then $t_1$, $t_2$, $t_3$, $t_4$
are the roots of
$$t^2+1/t^2-2gt-2f/t+c \mbox{ \quad or \quad } t^4-2gt^3+ct^2-2ft+1.$$
The product of the roots is 1,
so $t_4=1/t_1t_2t_3$ and the 4th point of intersection
of the circumcircle is diametrically across the hyperbola
from the orthocentre.  The centre of the hyperbola is the
midpoint of the join of the orthocentre to a point on
the circumcircle.  The locus of the centres of rectangular
hyperbolas drawn through the three vertices of a triangle
(i.e., through the four vertices of an orthocentric
quadrangle) is midway between the orthocentre and the
circumcircle; it is the 50-point circle!  This is a
special case of Feuerbach's theorem, which tells us that
the locus of the centres of conics passing through 4
points is a conic. In this case, six of the points on
the conic are the 3 midpoints and the 3 diagonal points,
so that here the conic is indeed a circle.

\bigskip

{\bf Sixth}, in case five aren't enough, use the case $n=2$ of the
Lighthouse Theorem (\S 6.9 in Chapter 6 below).

\newpage

\section{What is a 5-point (medial) circle?}

The circles of the clover-leaf theorem were discovered
as 5-point circles, but we will call them {\bf medial circles}.
In the notation of Figure \ref{proof},
the points $V_1$, $M_6$, $D_6$, and the midpoints of $V_2D_3$ and
$V_4D_5$, are concyclic.  The last two points are the feet of the
perpendiculars from $M_6$ onto $V_1V_2$ and $V_1V_4$; call them
{\bf midfoot points}.

{\bf A medial circle is a ``9-point circle''.}

In Figure \ref{proof} the $V_1M_6$-circle and the $V_4M_3$-circle
intersect in two points on the $A_2D_5$ altitude, say $B_i$
for the interior point and $B_e$ for the exterior point.
[Better notation needed for the altitude points.]
Similarly, the $V_2M_5$-circle and the $V_4M_3$-circle intersect
$V_1D_6$ in $A_i$ and $A_e$ and the $V_1M_6$-circle and the
$V_2M_5$-circle concur at $C_i$ and $C_e$ on $A_4D_3$.  How does
one characterize these points?

We can define them as follows.  E.g., $A_i$ and $A_e$
are given by the product and the sum of their distances from $V_1$:
$$V_1A_i\cdot V_1A_e=(\tfrac{1}{2}V_1D_5)\cdot V_1V_4
=V_1D_3\cdot V_1M_3=
b\cos A\cdot\tfrac{1}{2}c=\tfrac{1}{2}bc\cos A=
2R^2\cos A\sin B\sin C$$
$$\mbox{and \qquad } V_1A_i + V_1A_e=\tfrac{3}{2}V_1D_6=
\tfrac{3}{2}b\sin C=\tfrac{3}{2}c\sin B=3R\sin B\sin C$$
where $R$ is the circumradius of $V_1V_2V_4$.

Each medial circle passes through two such points
on each of two altitudes, making it a 9-point circle,
though not a traditional one.
These points deserve a name.  Call them {\bf altitude
points}.  The interior altitude points and the exterior
altitude points each form triangles in perspective
with $V_1V_2V_4$, and with the pedal triangle $D_6D_5D_3$,
having the orthocentre of $V_1V_2V_4$ as perspector.
At first glance they appear to be respectively
homothetic to $V_1V_2V_4$ and to the pedal triangle
$D_6D_5D_3$, but that is not so.  What are the various
perspectrices?

Notice that the formulas for the product and the sum of the
distances of a pair of altitude points from their vertex,
namely
$$2R^2\cos A\sin B\sin C \mbox{\quad and \quad}3R\sin B\sin C$$
$$2R^2\sin A\cos B\sin C \mbox{\quad and \quad}3R\sin C\sin A$$
$$2R^2\sin A\sin B\cos C \mbox{\quad and \quad}3R\sin A\sin B$$
extravert into

[left to reader for the time being!]

and quadrate into

[left to reader for the time being!]

{\bf The degenerate case}

If our triangle is right-angled at $V_1$, then the
midpoints of $V_1V_2$ and $V_1V_4$ coincide with the
midpoints of the joins of $V_2$ and $V_4$ to the
orthocentre, $V_8$, which is also the foot of
each of the perpendiculars from $V_2$ and $V_4$
onto $V_1V_4$ and $V_1V_2$ respectively, and the
midpoint of the join of vertex $V_1$ to itself as
orthocentre, so

{\bf The 9-point circle becomes a 5-point (medial) circle}

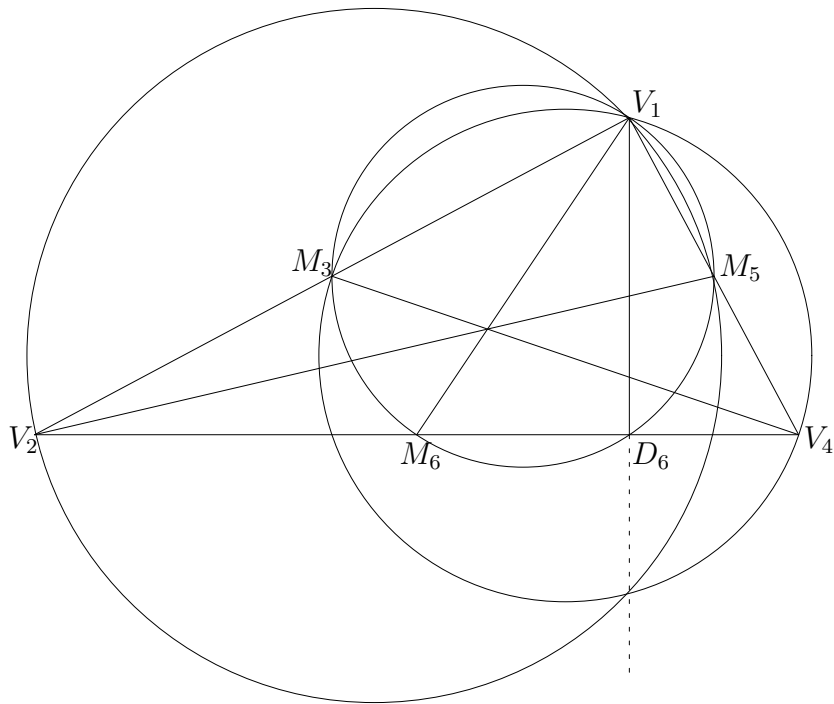
\begin{figure}[h]  
\begin{picture}(400,280)(-310,-100)
\setlength{\unitlength}{1pt}
\drawline(-112.5,60)(64,0)(-225,0)(0,120)(64,0)  
\drawline(-80.5,0)(0,120)(0,0)     
\drawline(-225,0)(32,60)        
\dashline{2}(0,0)(0,-90)         
\put(-40.25,60){\circle{144.5}}
\put(-96.5,30){\circle{262.77}}
\put(-24.25,30){\circle{186.42}}
\put(1,122){$V_1$}
\put(-235,-5){$V_2$}
\put(66,-5){$V_4$}
\put(-87,-11){$M_6$}
\put(1,-11){$D_6$}
\put(34,60){$M_5$}
\put(-128,62){$M_3$}
\end{picture}
\caption{If the triangle $V_1V_2V_4$ is right-angled at $V_1$,
then the 9-point circle degenerates into a 5-point (medial)
circle $V_1^3\,M_5^2\,D_6\,M_6\,M_3^2$}
\label{degen}
\end{figure}

Sixty years ago the present writer showed \cite{Gy} that a necessary
and sufficient condition for a triangle to be right-angled is
that the sum of the touch-radii, $r+r_1+r_2+r_3$, should be
equal to the perimeter, $2s$.  Also that another such
condition is $rr_3=r_1r_2$.  This was discovered after noting
that, for a 3,4,5 triangle, the radii are 1,2,3,6, and that
6 is a perfect number.

\clearpage

{\bf The obtuse case}

If the angle at $V_1$ is obtuse, then the orthocentre lies
outside the triangle, there are only two real altitude
points, both of which are exterior to the triangle, on
the altitude through $V_1$.
Only one pair of medial circles intersect in real points.
However, the three altitudes are still radical axes for
the pairs of circles.  This case does not occur in the
context in which we are presently interested. 

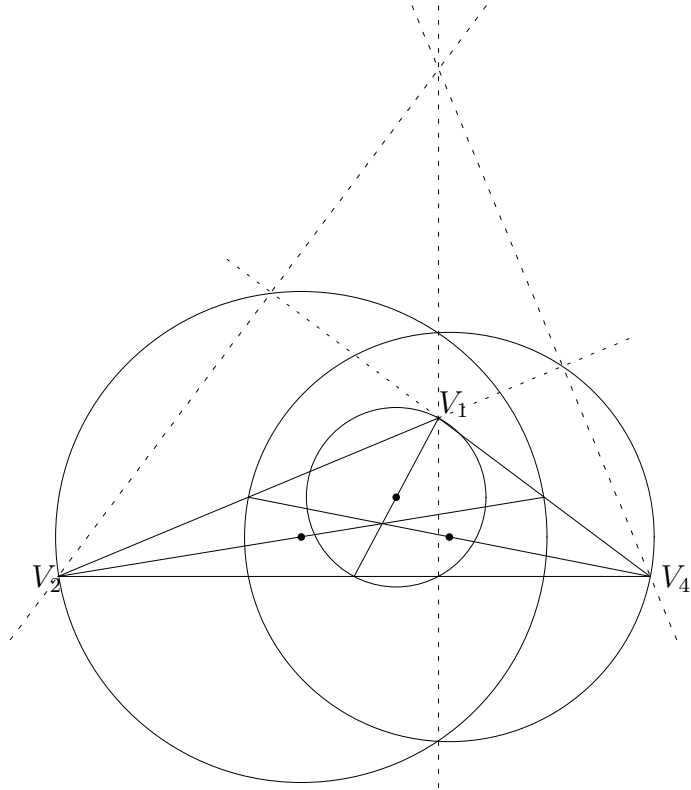
\begin{figure}[h]  
\begin{picture}(400,290)(-250,-80)
\setlength{\unitlength}{2pt}
\drawline(-36,15)(40,0)(-72,0)(0,30)(40,0) 
\drawline(20,15)(-72,0)              
\drawline(-16,0)(0,30)               
\dashline{0.75}(0,30)(36,45)
\dashline{0.75}(0,30)(-40,60)
\put(-8,15){\circle*{1.2}}
\put(2,7.5){\circle*{1.2}}
\put(-26,7.5){\circle*{1.2}}
\put(-8,15){\circle{34}}
\put(2,7.5){\circle{77.466}}
\put(-26,7.5){\circle{92.957}}
\dashline{1}(-81,-12)(9,108)
\dashline{1}(0,-40)(0,108)
\dashline{1}(45,-12)(-5,108)
\put(0,31){$V_1$}
\put(-77,-2){$V_2$}
\put(42,-2){$V_4$}

\end{picture}
\caption{Only one pair of medial circles intersect.  The $V_1$-circle
doesn't intersect the other altitudes in real points, but the
altitudes (dashed) are still radical axes.}
\label{obtuse}
\end{figure}

It is well known (and quadration makes it clear) that

{\bf Theorem} \cite{Gob} There are three times as many
obtuse triangles as there are acute ones.

So that there are only 12 (= 6 + $3\times 2$)
real altitude points in a quadrated triangle, not 24.

\clearpage

\section{The 4-leaf clover theorem.}

If we throw in the 9-point circle to the clover-leaf theorem
we have a four-leaf clover, each of the six pairs of circles
has either an altitude or an edge for radical axis, which
concur in threes at the vertices and orthocentre of the
triangle.  [In the language of Quadration, each pair of
circles has an edge of the (quadrated) triangle as a
radical axis.  See next chapter.]

\begin{figure}[h]  
\begin{picture}(400,340)(-250,-70)
\setlength{\unitlength}{10pt}
\drawline(-18,0)(5,12)
\drawline(0,24)(-4,0)
\drawline(10,0)(-9,12)
\put(-2,12){\circle{24.33}}
\put(-6.5,6){\circle{25.942}}
\put(0.5,6){\circle{22.472}}
\put(-2,7.875){\circle{16.25}}
\put(-2,7.875){\circle{16.3}}
\dashline{0.5}(0,-8)(0,28)
\dashline{0.5}(14,-3)(-18,21)
\dashline{0.5}(-24,-2.5)(12,12.5)
\thicklines
\drawline(-18,0)(0,24)(10,0)(-18,0)
\end{picture}
\caption{The four-leaf clover theorem}
\label{fourleaf}
\end{figure}
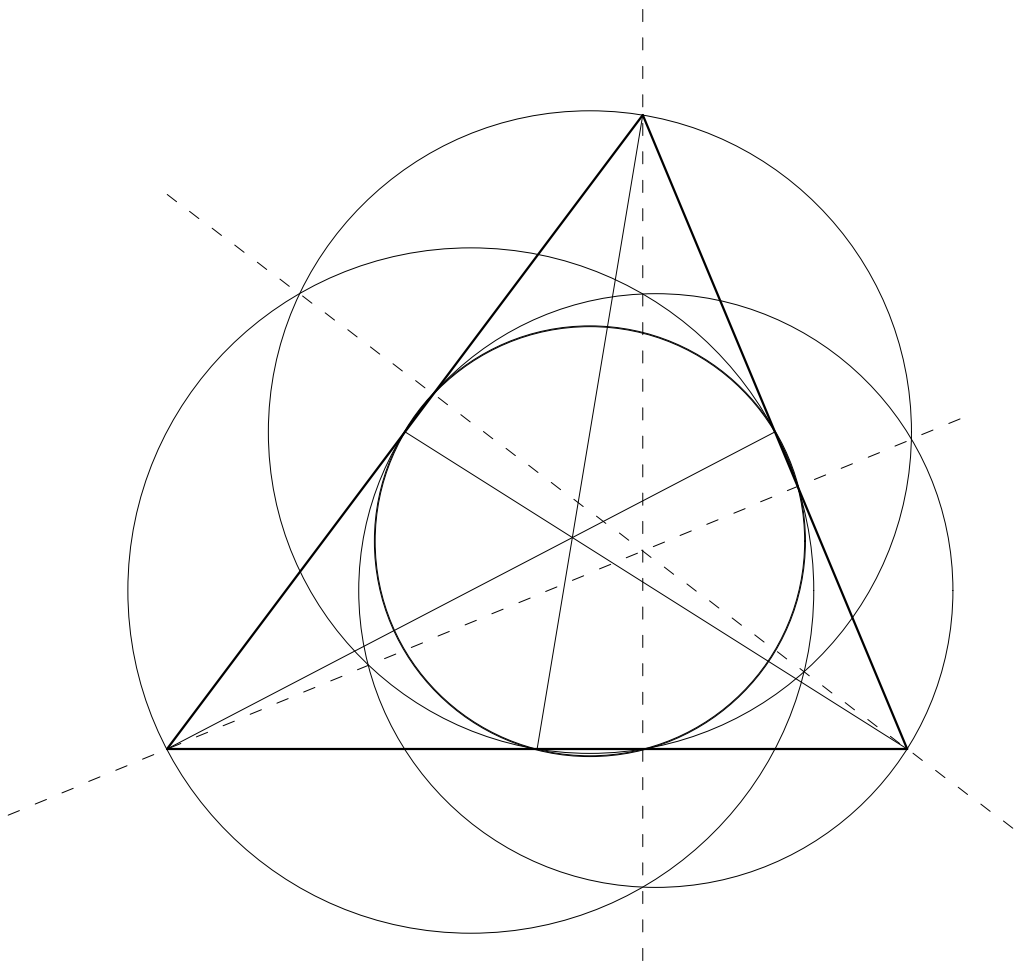

\chapter{A triangle is an orthogonal quadrangle!}

It is reasonable to give equal status to the orthocentre
and the vertices of a triangle, since then each point is the
orthocentre of the triangle formed by the other three.

\section{Quadration}

We then have 4 vertices, 6 edges, 4 triangles, 12 medians,
12 medial circles.  The same 9-point circle serves for
all four triangles.  Note that the edges serve equal time
as edges and altitudes, so that there are twelve altitude
points.

\begin{figure}[h]  
\begin{picture}(400,380)(-250,-105)
\setlength{\unitlength}{2.8pt}
\drawline(-36,48)(40,0)(0,96)(20,15)(-72,0)(40,0)(-36,15)(0,96)(0,30)
(20,48)(-72,0)(0,96)(-16,0)(0,30)(40,0)(0,63)(-72,0)(0,30)(-36,48)(0,96)
\dashline{1}(0,30)(60,55)
\dashline{1}(0,30)(-72,84)
\dashline{1}(0,30)(0,-40)
\put(-8,15){\circle{34}}
\put(2,7.5){\circle{77.466}}
\put(-26,7.5){\circle{92.957}}
\put(-36,31.5){\circle{95.671}}
\put(-18,55.5){\circle{88.64}}
\put(-18,39){\circle{40.249}}
\put(20,31.5){\circle{74.626}}
\put(10,55.5){\circle{83.433}}
\put(10,39){\circle{26.907}}
\put(2,24){\circle{89.889}}
\put(-8,48){\circle{97.324}}
\put(-26,24){\circle{103.769}}
\put(-8,31.5){\circle{65}}
\put(-8,31.5){\circle{64.7}}
\end{picture}
\caption{Twelve medial circles and one 9-point circle. There are
three circles though each of the 4 vertices; three circles
through each of the 6 midpoints, and four through each of the 3 diagonal
points; two circles through each of the 12 altitude points,
and two through each of the 12 midfoot points.}
\label{quadra}
\end{figure}
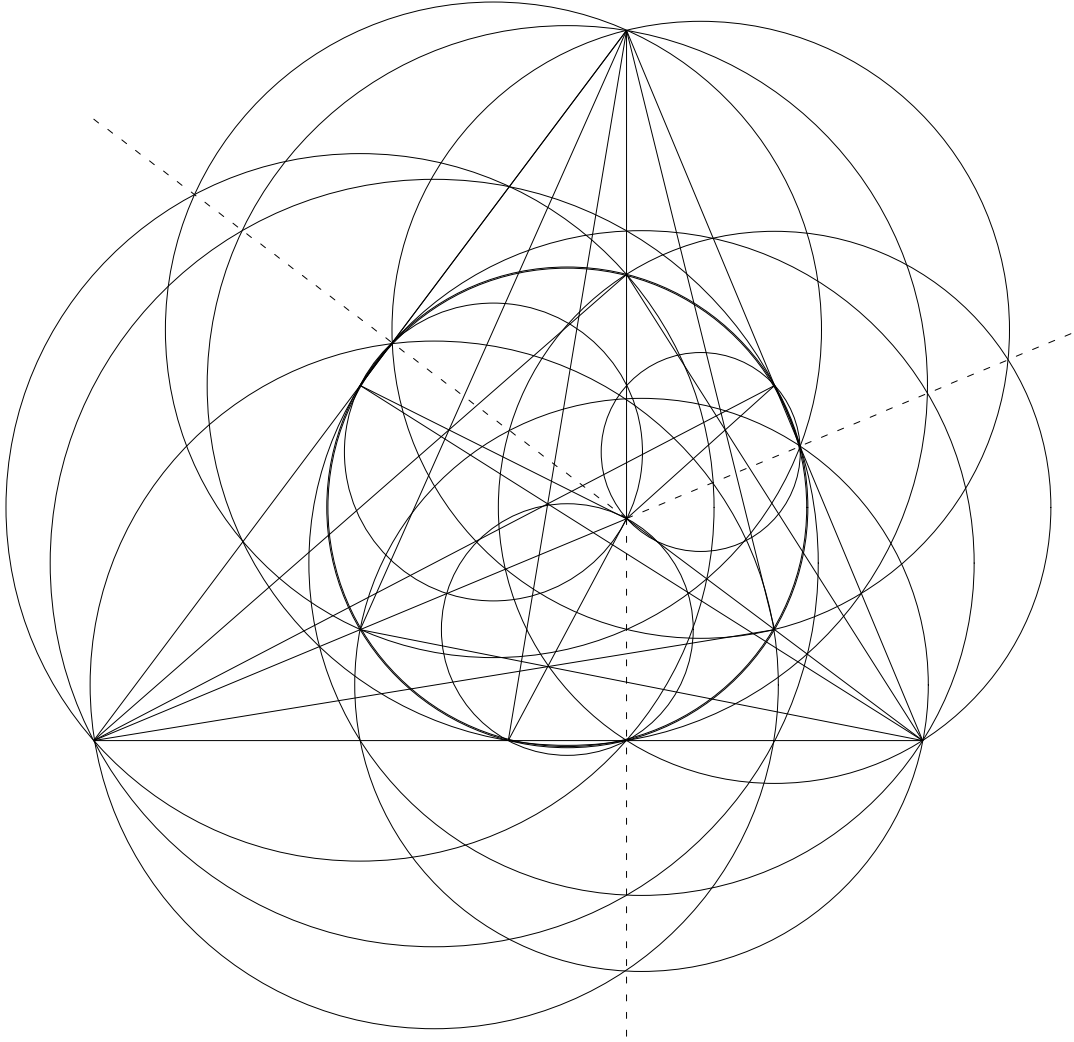


\section{Where did medial circles come from?}

{\bf Ans: The 260?-point sphere.  And where did the 260?-point
sphere come from?  Obviously the number 260 is arbitrary, and
has not been given an updated explanation in what follows.}

\chapter{Ans: The orthocentric tetrahedron}

Although the altitudes of a triangle concur, the
altitudes of a tetrahedron do not usually do so.
In 1827 Jacob Steiner noted that they are
generators of an equilateral hyperboloid.  That is,
a hyperboloid of one sheet formed by rotating a
rectangular hyperbola about its `imaginary' axis.
If two altitudes concur, the hyperboloid degenerates
into a pair of planes and the other two altitudes
also concur, and the tetrahedron is called
semi-orthocentric.  If two altitudes intersect a
third, then all four concur.

Just over two hundred years ago, Gaspard Monge
discovered that the six planes (call them {\bf midplanes}),
each passing through a midpoint of an edge of a
tetrahedron (any tetrahedron) and perpendicular
to the opposite edge, all passed through a single
point, now called the Monge point of the tetrahedron.
The Monge point is the centre of the equilateral
hyperboloid; it specializes to the orthocentre of
the tetrahedron in case the tetrahedron is
orthocentric, which is a very interesting
configuration.  See \cite{HW}, which has
a good bibliography; also \cite{C}, and \cite{B}.

Here we follow \cite{L} and \cite[pp.\,13--14]{R}.

Suppose that a
tetrahedron has an orthocentre and take that as origin, $O$.
Let the vertices $A$, $B$, $C$, $D$ define vectors $OA=$ {\bf a},
$OB=$ {\bf b}, $OC=$ {\bf c}, $OD=$ {\bf d}.  Then {\bf a} is
perpendicular to the plane $BCD$, the scalar products

\smallskip

\centerline{{\bf a}$\cdot(${\bf b} -- {\bf c}) =
{\bf a}$\cdot(${\bf c} -- {\bf d}) =
{\bf a}$\cdot(${\bf d} -- {\bf b}) = 0,}

\smallskip

and by symmetry the six products
{\bf a}$\cdot${\bf b}, {\bf a}$\cdot${\bf c}, {\bf a}$\cdot${\bf d},
{\bf b}$\cdot${\bf c}, {\bf c}$\cdot${\bf d}, {\bf d}$\cdot${\bf b},
are all equal (to $\sigma$, say).

Let the altitudes through $A$, $B$, $C$, $D$ meet the opposite
faces in $O_7$, $O_6$, $O_5$, $O_3$, respectively.
Then the scalar product {\bf a}$\cdot${\bf b}, for example,
is equal to the magnitude of {\bf a} times the projection of
{\bf b} on it.  More generally,
$$\sigma=OA\cdot OO_7=OB\cdot OO_6=OC\cdot OO_5=OD\cdot OO_3$$
Opposite edges are perpendicular.  For example,

\smallskip

\centerline{({\bf b} -- {\bf a})$\cdot$({\bf d} -- {\bf c})=
$\sigma-\sigma-\sigma+\sigma=0$}

\smallskip

Note that this implies that each Monge midplane not only
passes through the midpoint of an edge, but it also contains
that edge.  Moreover, since

\smallskip

\centerline{({\bf b} -- {\bf a})$^2 +$({\bf d} -- {\bf c})$^2=
a^2+b^2+c^2+d^2-4\sigma$}

\smallskip

is symmetrical in $a$, $b$, $c$, $d$, the sum of the squares
of the lengths of opposite edges is constant.

The feet of the altitudes are the orthocentres of the faces.
E.g., $BO_7$ is in the plane $OAB$ and so has shape
$h$\,{\bf a} + $k$\,{\bf b} for some $h$, $k$ and

\smallskip

\centerline{($h$\,{\bf a} + $k$\,{\bf b})$\cdot$({\bf d} -- {\bf c}) =
$h\sigma-h\sigma+k\sigma-k\sigma=0$}

\smallskip

\noindent
so that $BO_7$ is perpendicular to $CD$, and similarly
$CO_7\perp DB$ and $DO_7\perp BC$.

The square of the distance of the point
$\frac{1}{2}$({\bf a} + {\bf b} + {\bf c} + {\bf d}) from $A$ is \\
\centerline{\{{\bf a} -- $\frac{1}{2}$({\bf a} + {\bf b} +
{\bf c} + {\bf d})\}$^2$ = \{$\frac{1}{2}$({\bf b} + {\bf c}
+ {\bf d} -- {\bf a})\}$^2$ =
$\frac{1}{4}(a^2+b^2+c^2+d^2+6\sigma-6\sigma)$.}  By symmetry,
the point $\frac{1}{2}$({\bf a} + {\bf b} + {\bf c} + {\bf d})
is equidistant from the four vertices of the tetrahedron and so is
the circumcentre, $Q$, say.

The centroid, $G=\frac{1}{4}$({\bf a} + {\bf b} + {\bf c} + {\bf d}),
is the midpoint of $OQ$, in a somewhat analogous fashion for
the Euler line of a triangle.  In fact, for the general
tetrahedron, the centroid is the midpoint between the
Monge point and the circumcentre.  The analog of the 9-point
circle for the general tetrahedron is the 12-point sphere,
which passes through the centroids of the faces of the
tetrahedron, and has its centre midway between the Monge
point and the centroid.  But we will soon see that it
is the centroid that is the centre of our
{\bf 260?-point sphere}.  For example, the square of
the distance of the centroid from the midpoint
$\frac{1}{2}$({\bf c} + {\bf d}) of $CD$ is \\
\{$\frac{1}{4}$({\bf a} + {\bf b} + {\bf c} + {\bf d}) --
$\frac{1}{2}$({\bf c} + {\bf d})\}$^2$ =
\{$\frac{1}{4}$({\bf a} + {\bf b} -- {\bf c} -- {\bf d})\}$^2$
$=\frac{1}{16}(a^2+b^2+c^2+d^2+4\sigma-8\sigma)$
and is equal to that from each of the other five midpoints.
So we have six points on a sphere.

Of course, for {\em any} tetrahedron, the joins of the midpoints
of opposite edges concur and bisect each other, but in general
the lengths of these three joins are not equal.

\section{Think inside the box!}

 It helps to visualize properties
of tetrahedra if you put them in boxes.  Through each of the
six edges, draw a plane parallel to the opposite edge.  These
six planes are the faces of a parallelepiped.  In the case of
the regular tetrahedron, the box is a cube, of course.

\newpage

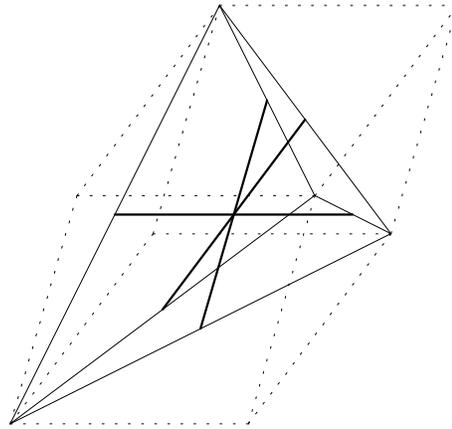
\begin{figure}[h]  
\begin{picture}(400,150)(-150,5)
\setlength{\unitlength}{1.2pt}
\dashline{1}(0,0)(21,72)(96,72)(75,0)(0,0)
(45,60)(66,132)(141,132)(120,60)(45,60)
\dashline{1}(21,72)(66,132)
\dashline{1}(96,72)(141,132)
\dashline{1}(75,0)(120,60)
\drawline(0,0)(66,132)(120,60)(0,0)(96,72)(66,132)
\drawline(120,60)(96,72)
\thicklines
\drawline(81,102)(60,30)
\drawline(108,66)(33,66)
\drawline(48,36)(93,96)
\end{picture}
\caption{Orthocentric tetrahedron in an equilateral box}
\label{box}
\end{figure}

For the more general orthocentric tetrahedron, since
opposite edges are perpendicular, the diagonals of the
parallelogram faces of the box are perpendicular,
so that the faces are rhombuses and the parallelepiped
has all twelve edges equal in length.  The joins of the
midpoints of opposite edges of the tetrahedron have the
same length as the edges of the box, and they each have
their midpoint at the centre of the box, which is the
centroid of the tetrahedron.  The midpoints of the
edges lie on a sphere with centre at this centroid.

This is our 260?-point sphere.  There aren't enough letters
to cope with the rather large number of points we wish to
deal with, so, for the nonce, we introduce the following
notation, which is sufficiently different from that used
elsewhere in the literature (including other parts of
this paper) to avoid confusion.  [Probably better to
change $V_0$ to $V_8$.]  Let the vertices of our
orthocentric tetrahedron be $V_0$, $V_1$, $V_2$, $V_4$,
and give the label $i+j$ ($1\leq i+j\leq6$) to the edge
$V_iV_j$ and the label $i+j+k$ (= 7, 6, 5 or 3) to the
face $V_iV_jV_k$.  For example, the midpoint of the edge
$V_iV_j$ is denoted by $E_{i+j}$ and the common foot of the
perpendiculars from the other two vertices onto that edge
is denoted by $F_{i+j}$.  The orthocentre and centroid of the
tetrahedron are $O$ and $G$ respectively, and the orthocentre
and centroid of the face $V_iV_jV_k$ are $O_{i+j+k}$ and
$G_{i+j+k}$.  The joins of the midpoints of opposite edges
$V_0V_i,\, V_{2i}V_{4i}$ ($i=1$, 2, 4; subscripts mod 7) are
$E_iE_{7-i}$. They concur and bisect
each other at the centroid, $G$.  The common perpendiculars
of opposite edges are $F_iF_{7-i} (i=1, 2, 4)$, and they
concur at the orthocentre, $O$.

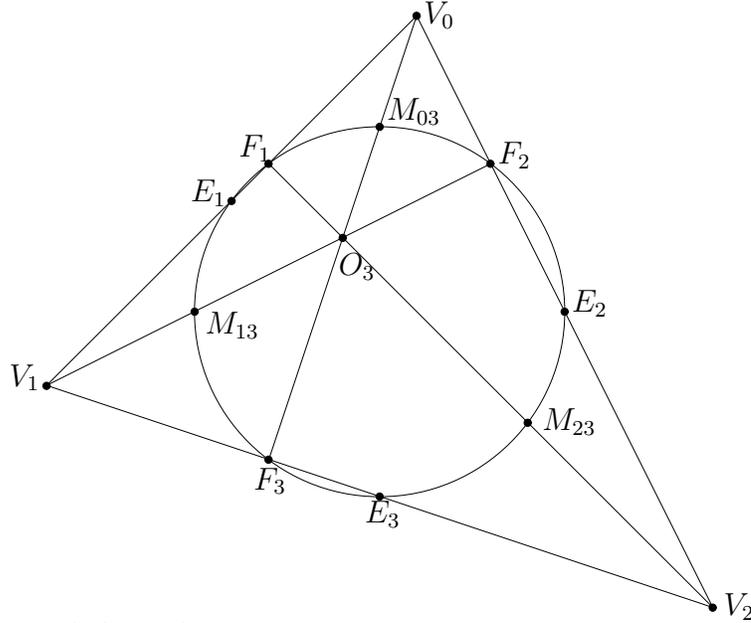
\begin{figure}[h]  
\begin{picture}(400,210)(-225,-158)
\setlength{\unitlength}{14pt}
\drawline(-3,-8)(1,4)(-9,-6)(9,-12)(1,4) 
\drawline(-9,-6)(3,0)   
\drawline(9,-12)(-3,0)   
\put(0,-4){\circle{10}}  
\put(1,4){\circle*{0.2}} 
\put(-9,-6){\circle*{0.2}} 
\put(9,-12){\circle*{0.2}} 
\put(0,-9){\circle*{0.2}} 
\put(5,-4){\circle*{0.2}} 
\put(-4,-1){\circle*{0.2}} 
\put(-3,-8){\circle*{0.2}} 
\put(3,0){\circle*{0.2}} 
\put(-3,0){\circle*{0.2}} 
\put(-1,-2){\circle*{0.2}} 
\put(0,1){\circle*{0.2}} 
\put(-5,-4){\circle*{0.2}} 
\put(4,-7){\circle*{0.2}} 
\put(1.2,3.8){$V_0$} 
\put(-10,-6){$V_1$} 
\put(9.3,-12.2){$V_2$} 
\put(-0.4,-9.7){$E_3$} 
\put(5.2,-4){$E_2$} 
\put(-5.1,-1){$E_1$} 
\put(-3.4,-8.8){$F_3$} 
\put(3.2,0){$F_2$} 
\put(-3.8,0.2){$F_1$} 
\put(-1.1,-3){$O_3$} 
\put(0.2,1.2){$M_{03}$} 
\put(-4.7,-4.6){$M_{13}$} 
\put(4.4,-7.2){$M_{23}$} 
\end{picture}
\caption{A face of an orthocentric tetrahedron and its
9-point circle.}
\label{face}
\end{figure}

Consider the intersection of the sphere with a face of the
tetrahedron (Figure \ref{face}).  It is a circle through
the midpoints of the
edges of the triangle, that is, its 9-point circle.  As
there are four faces we appear to have found 36 points on
the sphere, but two points on each edge have been counted
twice, namely the six midpoints of the edges and the six
feet of the common perpendiculars to pairs of opposite
edges.

\newpage

So far we have only 24 points on our sphere.  They are
$E_i$, $F_i$ ($1\leq i\leq6$) and three on each of the
four 9-point circles.  These turn out to be midfoot
points; we will denote the midpoint of $V_iO_j$ by
$M_{ij}$, ($i\in\{0,1,2,4\}, j\in\{7,6,5,3\}$).
To find more points on the sphere we next consider

\section{Medial planes of an orthocentric tetrahedron.}

The plane through an edge, say $V_1V_4$, and the midpoint,
$E_2$, of the opposite edge $V_0V_2$, (call it the
$E_5$ medial plane --- we economize on notation and
use the same name for the plane and the midpoint of the
edge through which it passes) cuts the tetrahedron
in a triangle, and our sphere in a great circle, since
the centre of the sphere, $G$, is the centroid of the
tetrahedron, i.e., the midpoint of $E_2E_5$, where
$E_5$ is the midpoint of $V_1V_4$.  See the left part
of Figure \ref{section}.  This circle is a 5-point
circle of the triangle, yielding two more midfoot
points (black dots), $S_{15}$, $S_{45}$, on the sphere.
As there are six such medial planes,
we have found 12 points to bring our total back to 36.

The faces and plane sections of an orthocentric
tetrahedron are always acute-angled triangles, so
that the altitude points are real, and we may include
them in our count.  They are indicated by small
circles;\ e.g., in the left of Figure \ref{section}, $T_{15}$
and $T_{65}$, where the medial circle intersects
the altitude through $V_1$; and $T_{45}$, $T_{35}$,
on the altitude through $V_4$.  We now have
$4\times6$ more points on our sphere, giving a total of 60.
In our notation, the first subscript refers to the relevant
vertex and the second to the edge, or plane (here number 5).
When two points are associated with a vertex $V_i$, then
the interior one has the odious subscript $i=0$, 1, 2 or 4,
and the exterior one the evil $7-i=7$, 6, 5 or 3.

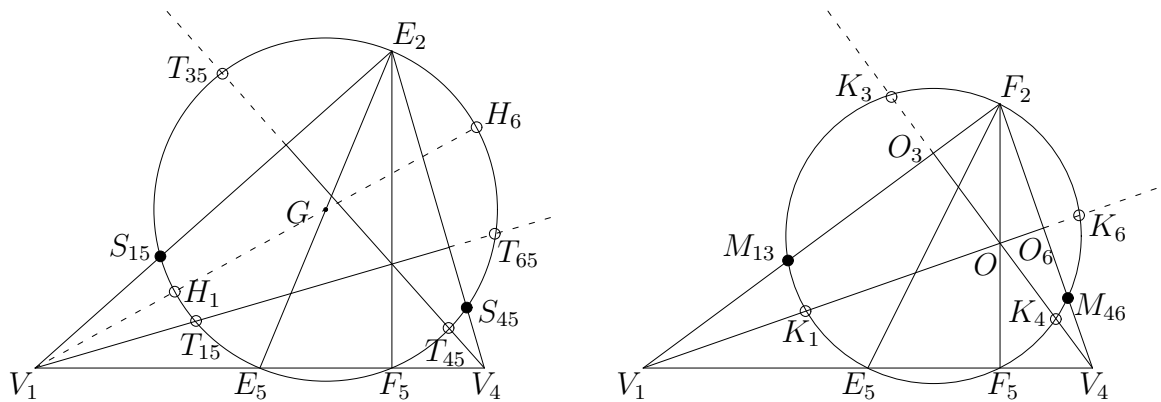
\begin{figure}[h]  
\begin{picture}(400,150)(-160,-10)
\setlength{\unitlength}{5pt}
\drawline(-10,0)(0,24)(-27,0)(7,0)(0,24)(0,0) 
\drawline(-27,0)(4.3344,9.1392)   
\dashline{0.5}(-27,0)(6,18)      
\drawline(7,0)(-8.0069,16.88276)      
\drawline(36,0)(46,20)(19,0)(53,0)(46,20)(46,0)  
\drawline(19,0)(49.2895,10.6013)   
\drawline(53,0)(40.95394,16.2622)      
\put(-5,12){\circle{26}}
\put(41,10){\circle{22.3607}}
\put(-29,-2){$V_1$}
\put(-12,-2){$E_5$}
\put(-1,-2){$F_5$}
\put(6,-2){$V_4$}
\put(0,24.5){$E_2$}
\put(-8,11){$G$}
\put(-21.4,8.5){$S_{15}$}
\put(6.5,3.5){$S_{45}$}
\put(-16,1){$T_{15}$}
\put(8,8){$T_{65}$}
\put(2.5,0.5){$T_{45}$}
\put(-17,22){$T_{35}$}
\put(17,-2){$V_1$}
\put(34,-2){$E_5$}
\put(45,-2){$F_5$}
\put(52,-2){$V_4$}
\put(46,20.5){$F_2$}
\put(-15.8,5){$H_1$}
\put(7,18.5){$H_6$}
\dashline{0.5}(-8.0069,16.88276)(-17,27)  
\dashline{0.5}(4.3344,9.1392)(12,11.375) 
\dashline{0.5}(49.2895,10.6013)(59,14)   
\dashline{0.5}(40.95394,16.2622)(33,27)   
\put(44,7){$O$}
\put(47.4,8.2){$O_6$}
\put(37.5,16){$O_3$}
\put(25,8.5){$M_{13}$}
\put(51.6,4){$M_{46}$}
\put(29.5,2){$K_1$}
\put(52.7,10){$K_6$}
\put(46.6,3.5){$K_4$}
\put(33.5,20.6){$K_3$}
\put(-17.5034,8.4414){\circle*{0.8}}
\put(5.6672,4.5696){\circle*{0.8}}
\put(51.1448,5.3007){\circle*{0.8}}
\put(29.97697,8.1311){\circle*{0.8}}
\put(-5,12){\circle*{0.3}}
\put(7.8,10.15){\circle{0.8}}
\put(-14.8,3.558){\circle{0.8}}
\put(4.3,3.0375){\circle{0.8}}
\put(-12.8,22.275){\circle{0.8}}
\put(50.25,3.7125){\circle{0.8}}
\put(37.8,20.52){\circle{0.8}}
\put(31.3,4.305){\circle{0.8}}
\put(52,11.55){\circle{0.8}}
\put(-16.4126,5.7749){\circle{0.8}}  
\put(6.4126,18.225){\circle{0.8}}    
\end{picture}
\caption{A medial plane and a midplane of an orthocentric
tetrahedron each intersect the tetrahedron in a triangle,
one of whose medial circles is a section of the
260?-point sphere.}
\label{section}
\end{figure}

There are three medial planes through each vertex (1,4,2
through $V_0$; 2,6,3 through $V_2$; 3,1,5 through $V_1$;
and 4,5,6 through $V_4$), e.g., those through $V_1$
intersect in the line $V_1G$ which in turn intersects
the sphere in the {\bf medial points} $H_1$, $H_6$
which are common to all three planes and circles.
The line $V_4H_4GH_3$ is omitted from Figure \ref{section}.
Two such points associated with each of 4 vertices
brings the total to $60+(2\times4)=68$.


\section{Midplanes of an orthocentric tetrahedron}.

We can also consider Monge's midplanes; for example,
the plane through an edge, say $V_1V_4$ again, which
is perpendicular to the opposite edge (see the right
of Figure \ref{section}).  This plane contains $O$,
the orthocentre of the tetrahedron, and $O_6$ and
$O_3$, the feet of the perpendiculars from $V_1$ and
$V_4$ onto their opposite faces, and so contains
the perpendicular, $F_5F_2$, common to the edges
$V_1V_4$ and $V_0V_2$.  The intersection with the
sphere is not a great circle this time, but it is
a circle through the feet of this common perpendicular
and also through the midpoint, $E_5$, of $V_1V_4$,
so that it is a medial circle for the triangle
$V_1V_4F_2$, yielding two midfoot points, $M_{13}$
and $M_{46}$, which we have already seen on the
9-point circles.  There are three midplanes through
each vertex (1,4,2; 2,6,3; etc., as with the medial
planes).  These concur in the four altitudes
$V_i\,O\,O_{7-i}$ which intersect
the sphere in eight altitude points, $K_i$, $K_{7-i}$
($i=0, 1, 2, 4$).  This brings our total of points
on the sphere to $68+8=76$.

The number 76 is somewhat
arbitrary and has varied upwards and downwards several
times during the writing of this paper as we discovered
new points or coincidences between old ones.  For example,
Conway likes to think of the 9-point circle as a 12-point
circle with three inscribed rectangles, the fourth corners
being the reflexions in the 9-point centre of the feet of
the altitudes.  This brings our total up to 88.  Others
would include the points of contact of the three common
tangents of the 9-point circle with the deltoid (envelope of
the Simson-Wallace line), homothetic with all eighteen Morley
triangles (see \cite{M} or \cite[pp.72--79]{Loc2}).
Now we have a nice round hundred-point
sphere.  Yet again, Feuerbach's theorem (due to Brianchon
\& Poncelet \cite{BP}?) tells us that the 9-point circle
also touches the four touch-circles (incircle and excircles)
so this gives us 16 more candidates \ldots, but I digress \ldots .

\section{Seeing the sphere}

It is hard for most of us to visualize three dimensions.
Hopefully, by the time this article appears, there will be a
reference to where you can view the sphere on-line, but
while we are confined to {\sc Monthly} pages, I will
content myself with a stereographic projection from a
``north pole'' of the sphere onto the tangent plane at
the ``south pole'', i.e., we will invert with respect to
a point on the sphere, so that circles invert into circles.
Figure \ref{24point} shows the (inverses of the) four
9-point circles of the faces of the tetrahedron.

\begin{figure}[h]
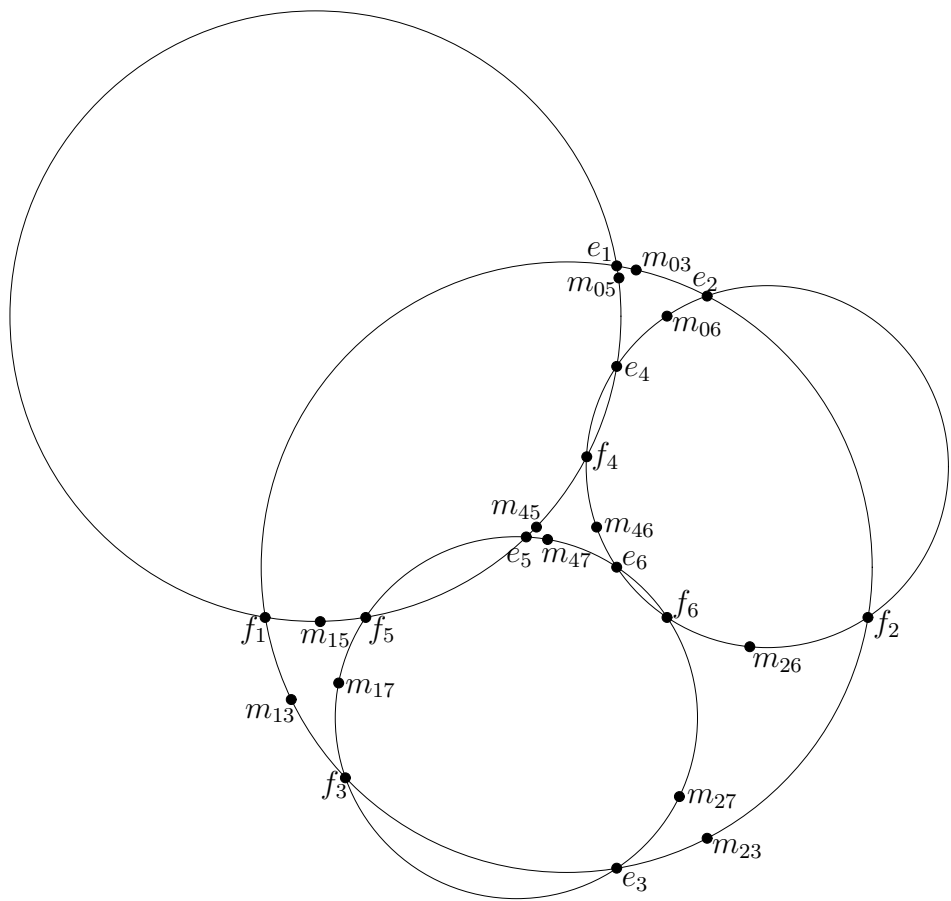
  

\caption{Projection of the four 9-point circles, showing the
first 24 points of the 260?-point sphere.  Lower case labels
denote the projections of the corresponding upper case points
on the sphere.}
\label{24point}
\end{figure}

\clearpage

Figure \ref{bisect} shows the projections of the six medial
circles $E_iE_{7-i}F_i$, $i=1$, 2, \ldots 6, which are the
sections of the 260?-point sphere by the medial planes.

\begin{figure}[h]
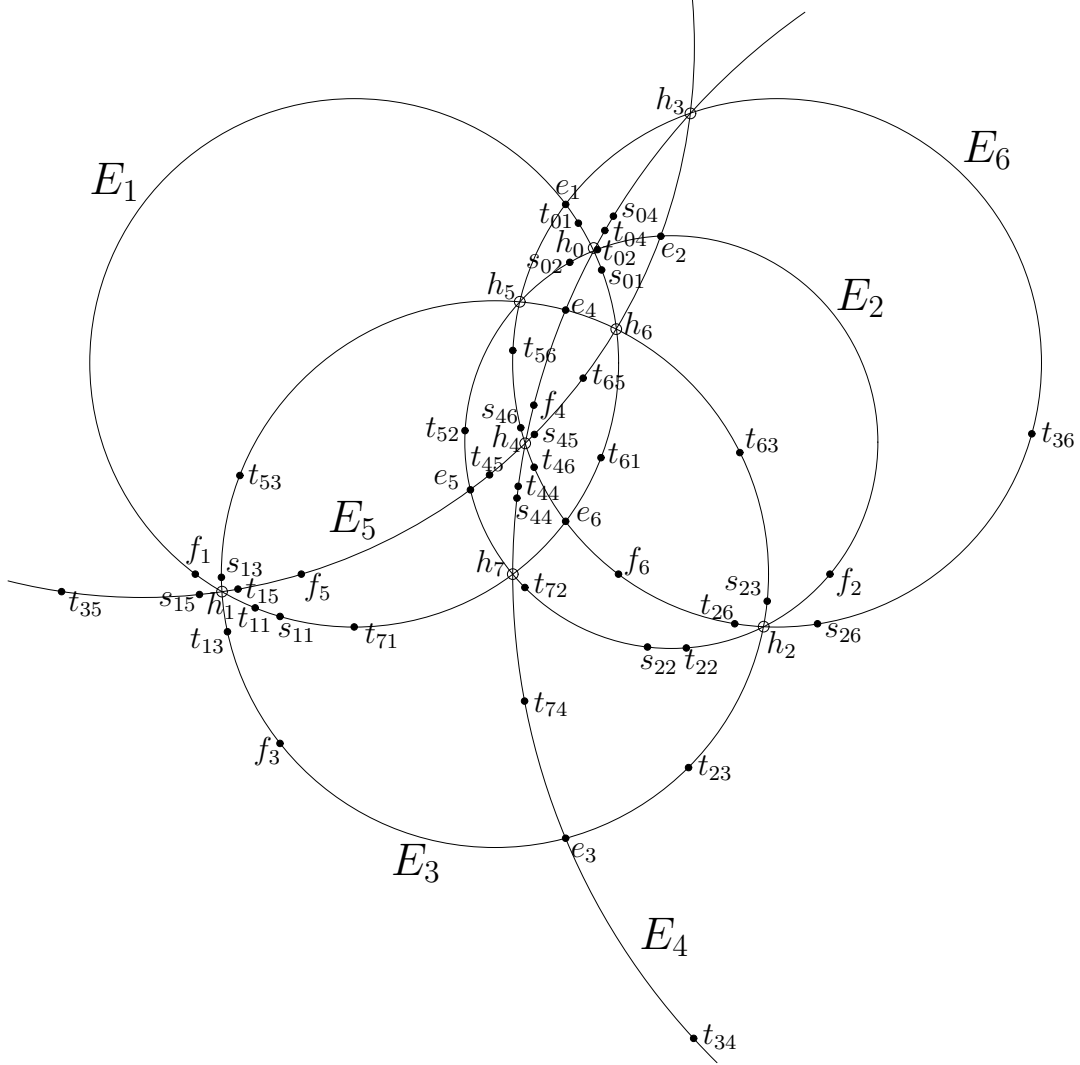
  
%
\caption{Projection of the six medial circles
in the medial planes, $E_i\ (1\leq i\leq6)$ showing the
next 12+24+8 points of the 260?-point sphere.  The eight medial
points, $h_i\,(0\leq i\leq7)$, are represented by small
circles.  $t_{02}$ is a different point from $h_0$.}
\label{bisect}
\end{figure}

\clearpage

Figure \ref{right} shows the projections of the six medial
circles $E_iF_iF_{7-i}$, $i=1$, 2, \ldots 6, which are the
sections of the 260?-point sphere by Monge's midplanes.

\begin{figure}[h]
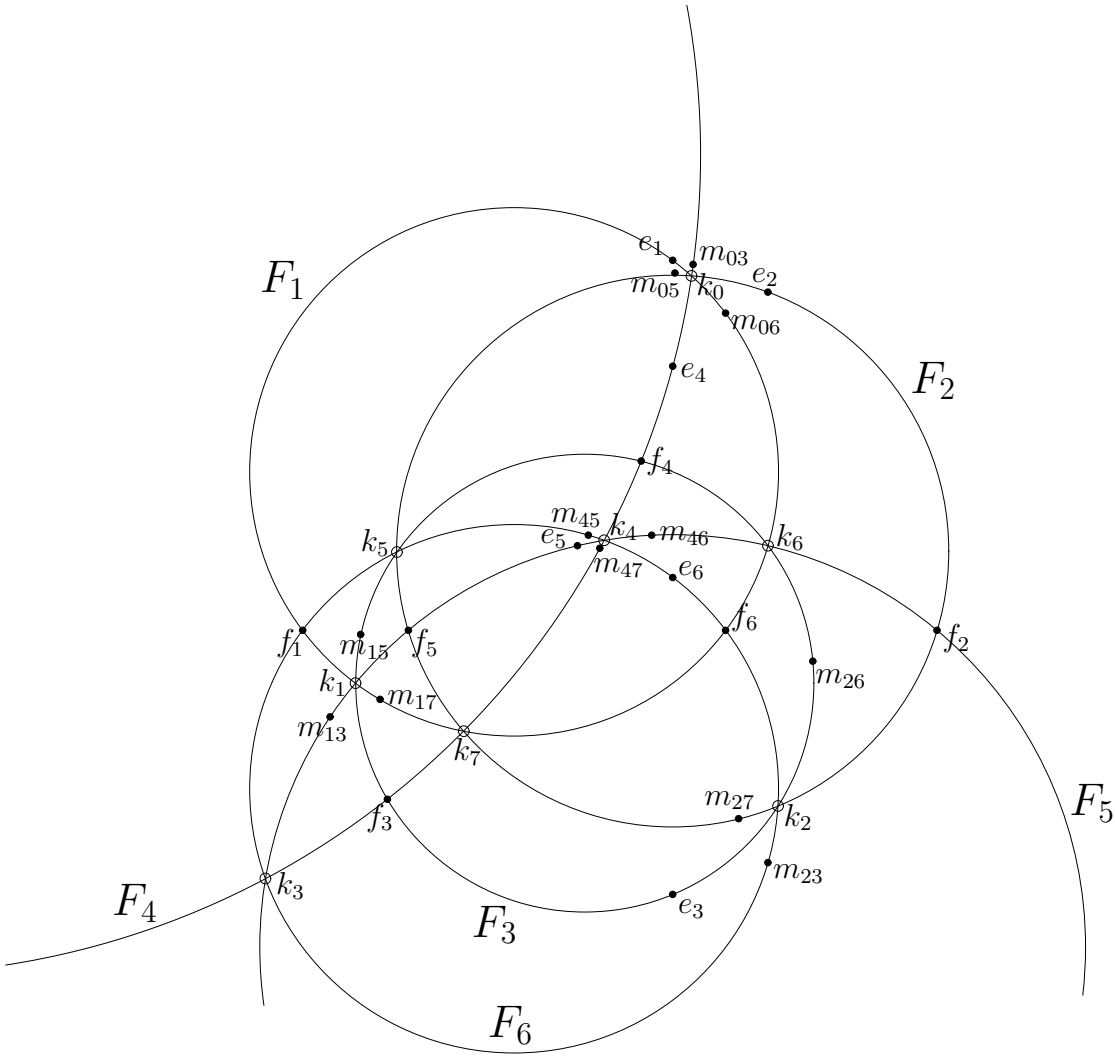
  
%
\caption{Projection of the remaining six medial circles,
the sections of the 260?-point sphere by Monge's midplanes,
$F_i, (1\leq i\leq6)$, showing the final eight points,
$k_i (0\leq i\leq7)$, of the 260?-point sphere, which are
represented by small circles; the black dots are the
original 24 points on the four 9-point circles.}
\label{right}
\end{figure}

{\bf Why were we interested in the orthocentric tetrahedron?}

This is the beginning of our story.

\chapter{The Pavillet tetrahedron}

Axel Pavillet \cite{P} quite
recently discovered that if you erect perpendiculars
$AA_o$, $BB_o$, $CC_o$ (vertical dashed lines in
Figure \ref{pav}) to the plane of the triangle
$ABC$, at $A$, $B$, $C$, of respective lengths $s-a$, $s-b$,
$s-c$, where $a$, $b$, $c$ are the edge-lengths of
the triangle and $s=\frac{1}{2}(a+b+c)$ its semiperimeter,
then the three points $A_o$, $B_o$, $C_o$, together
with the incentre, $I_o$, of $ABC$, form an orthocentric
tetrahedron.  Note that  $s-a$, $s-b$, $s-c$, are the
lengths of the tangents from $A$, $B$, $C$ to the incircle.

\begin{figure}[h]  
\begin{picture}(400,330)(-120,-4)
\setlength{\unitlength}{2.2pt}
\drawline(-13,0)(112,0)(0,84)(-13,0)    
\drawline(14.6283,28.1489)(14.6283,49.093)   
\dashline{1.5}(40,54)(-13,43.75)           
\dashline{1.5}(-7.647,34.588)(112,112.5)             
\dashline{1.5}(22,0)(0,146.5)               
\dashline{3}(-13,0)(-13,43.75)           
\dashline{3}(112,0)(112,112.5)             
\dashline{3}(0,84)(0,146.5)               
\thicklines
\drawline(-13,43.75)(112,112.5)(0,146.5)(-13,43.75)  
\drawline(22,30)(-13,43.75)         
\drawline(22,30)(112,112.5)         
\drawline(22,30)(0,146.5)           
\put(22,30){\circle{60}}
\put(10,50){$O$}
\put(10,25){${\mathbb{G}}_o$}
\put(23,27){$I_o$}
\put(-18,-2){$B$}
\put(-5,82){$A$}
\put(113,-2){$C$}
\put(-2,148){$A_o$}
\put(-19,42){$B_o$}
\put(113,111){$C_o$}
\put(20,-4){$X$}
\put(41,54){$Y$}
\put(-11.6,33){$Z$}
\put(62,-4){$s-c$}
\put(73,30){$s-c$}
\put(114,55){$s-c$}
\put(-26,20){$s-b$}
\put(-2,-4){$s-b$}
\end{picture}
\caption{The Pavillet tetrahedron, $A_oB_oC_oI_o$.  Three of
its altitudes are the three dashed lines through its
orthocentre, $O$.  Check that the three joins $AX$, $BY$,
$CZ$, concur at $G_0$.}
\label{pav}
\end{figure}
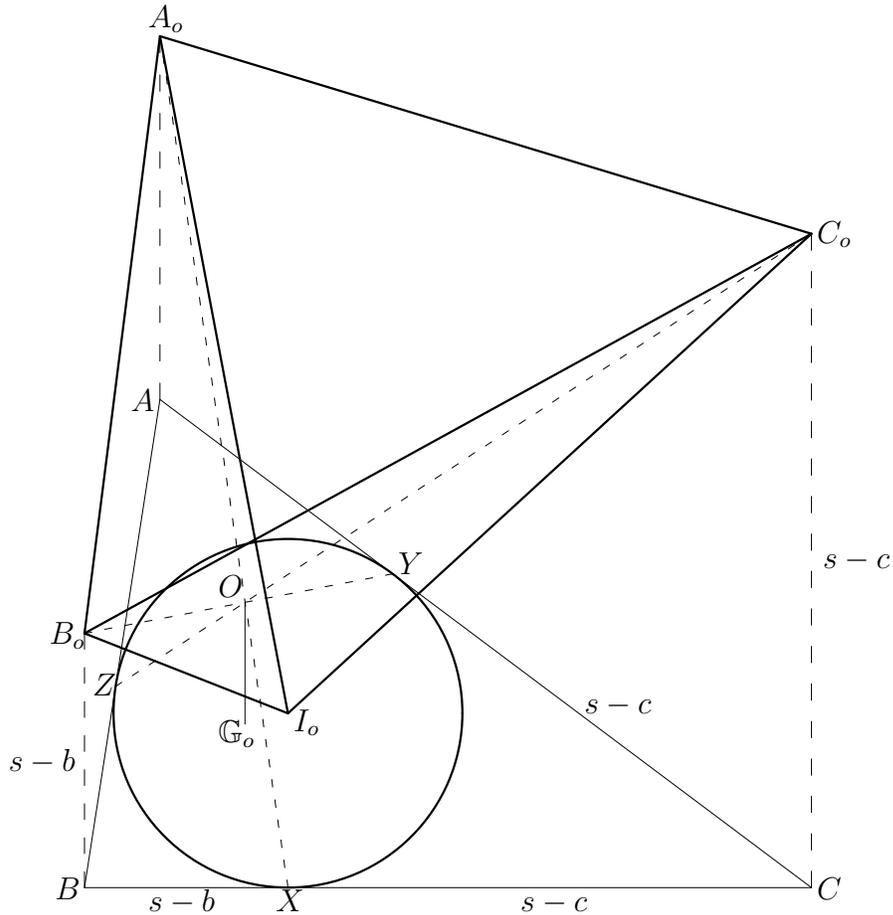

This tetrahedron has many interesting properties
which will hopefully have been revealed \cite{P} by the time this
paper appears.  For example, the projection of its orthocentre,
$O$, onto the plane $ABC$ is the {\bf Gergonne point} of
the triangle $ABC$; that is, the point of concurrence
of the three joins, $AX$, $BY$, $CZ$ of the vertices to
the touchpoints $X$, $Y$, $Z$, of the incircle with the
opposite edges.


In fact the joins $XA_o$, $YB_o$, $ZC_o$, of these points of
contact to the relevant vertices of the Pavillet tetrahedron
are three of its altitudes, concurring in its orthocentre, $O$.
The plane $A_oB_oC_o$ meets the edges $BC$, $CA$, $AB$
respectively in the points $X'$, $Y'$, $Z'$ which are the
harmonic conjugates of $X$, $Y$, $Z$ with respect
to $BC$, $CA$, $AB$.  The line $X'Y'Z'$ is the polar of the
Gergonne point, ${\mathbb{G}}_o$, with respect to the
triangle $ABC$, and
is called the {\bf Gergonne line}.  This is perpendicular to
the projection onto the plane $ABC$ of the fourth altitude
of the tetrahedron, which has been called the {\bf Soddy line}
\cite{O} of the triangle $ABC$, though both Apollonius and Descartes
might wish to claim priority.  This last line contains,
in addition to the the incentre and the Gergonne point,
the interior and exterior Soddy centres, $S$ and $S'$,
and it intersects the Euler line of the triangle $ABC$
in the {\bf deLongchamps point} \cite{Va}, $D$, which is
the reflexion of the orthocentre of the triangle in its
circumcentre, or, alternatively, the orthocentre of the
triangle formed by parallels to $BC$, $CA$, $AB$ passing
respectively through $A$, $B$, $C$.

The so-called Soddy centres are the centres of the circles
which are tangent to each of the three mutually tangent
circles with centres $A$, $B$, and $C$, passing respectively
through the touch-points $Y$ \& $Z$, $Z$ \& $X$, and $X$
\& $Y$ of the incircle with the edges of the triangle $ABC$.
I will call these circles the {\bf tangent circles}
(compare \cite{Ep}), since they have the tangents to the
incircle from the vertices $A$, $B$, $C$ as radii
(lengths $s-a$, $s-b$, $s-c$).

{\bf Theorem}.  {\it The incircle and two Soddy circles of a triangle
form a coaxal system with the Soddy line as line of centres
and the Gergonne line as radical axis.  The circles with
centres $X'$, $Y'$, $Z'$, passing respectively through
$X$, $Y$, $Z$, are members of the orthogonal coaxal system
which has the Gergonne line as line of centres and the
Soddy line as radical axis.} [See Figure \ref{coaxal}.  There
are extraversions of this theorem; compare \S{\bf10} below.
They are depicted in Figure \ref{fourco}.]

\clearpage

\begin{figure}[h]  
\begin{picture}(420,310)(-120,-164)
\setlength{\unitlength}{1.4pt}
\drawline(-60,0)(108,0)   
\drawline(-18,56)(24,0)  
\drawline(-60,0)(16.36,40.73)  
\drawline(-60,-112)(32,95)  
\drawline(-59.9,-112)(32.1,95)  
\drawline(-63,76.05)(227.7,-53.15)   
\drawline(-63,75.95)(227.7,-53.25)   
\put(-4,14){\circle{27.9}}     
\put(-4,14){\circle{28.1}}     
\put(-2.8085,16.68085){\circle{7.14}}     
\put(-2.8085,16.68085){\circle{7.16}}     
\put(-60,-112){\arc{336}{4.68}{6.35}}  
\put(-60,-112){\arc{336.2}{4.68}{6.35}}  
\put(108,0){\circle{224.1}}     
\put(108,0){\circle{224.2}}     
\put(-18,56){\circle{84.1}}     
\put(-18,56){\circle{83.9}}     
\put(16.364,40.727){\circle{61.0}}     
\put(16.364,40.727){\circle{61.2}}     
\put(-36,-58){\arc{215}{4.42}{6.95}}     
\put(-36,-58){\circle*{1.2}}     
\put(-4,14){\circle*{1.2}}     
\put(-1.565,19.478){\circle*{1.2}}     
\put(-2.8085,16.68085){\circle*{1.2}}     
\put(-60,-112){\circle*{1.2}}     
\put(108,0){\circle*{1.2}}     
\put(-18,56){\circle*{1.2}}     
\put(16.36,40.73){\circle*{1.2}}     
\put(-144,0){\circle*{1.2}}     
\put(-35,-62){$D$}
\put(-0.5,19){\scriptsize${\mathbb{G}}_o$}
\put(-2.5,9){\scriptsize$I_o$}
\put(-1.5,13.5){\scriptsize$S_o$}     
\put(-58.5,-116){$S_o'$}     
\put(-66,2){$B$}
\put(-9,31){$A$}
\put(25,2){$C$}
\put(107,2){$X'_o$}
\put(-19,58){$Y'_o$}
\put(13,43.5){$Z'_o$}
\end{picture}
\caption{Two orthogonal coaxal systems.  The Gergonne and
Soddy lines of triangle $ABC$ are the lines of centres of
one system and the radical axes of the other.  The `$o$'
subscripts can be ignored; they are to distinguish from
the extraversions in Figures \ref{extra}, \ref{fourco},
and \ref{finale}.}
\label{coaxal}
\end{figure}
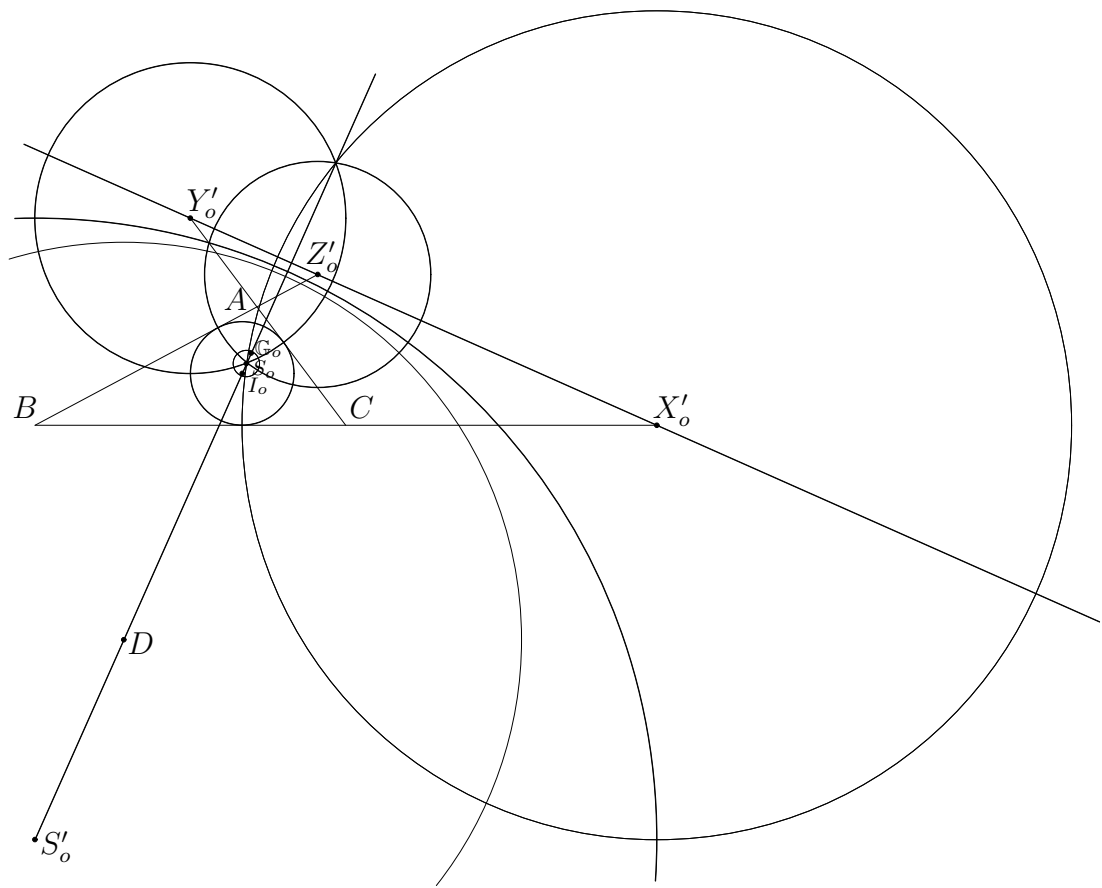

In Figure \ref{coaxal} only part of the outer Soddy circle is shown.
Note that one of the two limiting points on the Soddy line is
close to, but not coincident with, the inner Soddy centre.
Also included is (an arc of) the circle of the system having the
deLongchamps point, $D$, as centre; but no special significance of
this circle is yet known to the present writer.

\clearpage

The original draft of this paper considered three cases of
triangle: a critical case, and more obtuse, and less obtuse ones.
This was influenced by the paper of Veldkamp \cite{V}, published
a quarter of a century ago in a widely read, and oft-cited,
periodical.  It unfortunately contains a mistatement relating
these angle-dependent cases to the ones we now discuss.  The
error was eventually corrected by Hajja \& Yff \cite{HY} in a not so
often read journal.  With hindsight, our original classification
needs to be revised.  The situation is confused by the fact
that the critical case for the triangle $ABC$ can be defined
in a number of equivalent ways, which differ widely in appearance,
and only peripherally involve the size of the largest angle.

The following are equivalent definitions of the critical case.
The edge-lengths of
the triangle are $a$, $b$, $c$; the angles $A$, $B$, $C$;
area $\Delta$; perimeter $2s$; circumradius $R$; inradius $r$;
incentre $I_o$; Gergonne point ${\mathbb{G}}_o$.
\begin{enumerate}
\item Outer Soddy centre at infinity, outer Soddy circle a
straight line;
\item outer Soddy circle coincides with Gergonne line;
\item inner Soddy centre is the midpoint of segment $I_o{\mathbb{G}}_o$;
\item isoperimetric point at infinity (see \S{\bf8});
\item $\tan\frac{A}{2}+\tan\frac{B}{2}+\tan\frac{C}{2}=2$;
\item $(s-b)^2(s-c)^2+(s-c)^2(s-a)^2+(s-a)^2(s-b)^2=2\Delta^2$;
\item $2s=4R+r$;
\item (one of Pavillet's delightful discoveries) the angle
between the fourth altitude (the one through $I_o$) of the
Pavillet tetrahedron, and the Soddy line, is $\pi/4$;
\end{enumerate}
and see \${\bf9}.  It is worth noting that these definitions
are also equivalent to putting radius $\infty$ into Descartes's
theorem (or putting `bend zero' in Soddy's Kiss Precise \cite{S}):
$$\left(\frac{1}{s-a}+\frac{1}{s-b}+\frac{1}{s-c}+
\frac{1}{\infty}\right)^2=2\times\left(\frac{1}{(s-a)^2}+
\frac{1}{(s-b)^2}+\frac{1}{(s-c)^2}+\frac{1}{\infty^2}\right)$$
the other radii, $s-a$, $s-b$, $s-c$, being those of the
tangent circles.

\clearpage

\begin{figure}[h]  
\begin{picture}(420,300)(-180,-150)
\setlength{\unitlength}{0.5pt}
\drawline(-346.15,0)(480,0)    
\drawline(-280,-236.77)(140,276.23)  
\drawline(-207,196.43)(483,-107.57)  
\dashline[+30]{7}(-70,248.8)(152.77,-285.86)  
\put(0,105.23){\circle{72}}        
\put(-86.15,0){\circle{200}}       
\put(238.85,0){\circle{450}}       
\put(9.29,58.3){\circle{23.68}}    
\put(13.85,47.37){\circle{94.74}}  
\dashline[+30]{7}(-346.15,0)(350,290.06)   
\put(4.74,69.23){\circle*{5}}      
\put(9.29,58.3){\circle*{5}}       
\put(13.85,47.37){\circle*{5}}     
\put(152.77,-285.86){\circle*{5}}  
\put(-358,5){$X'$}
\put(-57,132){$Y'$}
\put(32,166){$Z'$}
\put(5,102){$A$}
\put(-90,-22){$B$}
\put(225,-22){$C$}
\put(-9,-22){$X$}
\put(43,89){$Y$}
\put(-50,67){$Z$}
\put(4,73){{\scriptsize ${\mathbb{G}}_o$}}
\put(13,55){{\scriptsize $S_o$}}
\put(20,35){{\scriptsize I}}
\put(159,-285){$D$}
\end{picture}
\caption{The critical case.  The outer Soddy circle coincides with the
Gergonne line, $X'Y'Z'$.  The inner Soddy circle, centre $S_o$, passes
through the incentre, $I$, and the Gergonne point, ${\mathbb{G}}_o$.
The outer Soddy centre is the point at infinity on the Soddy line.}
\label{critical}
\end{figure}
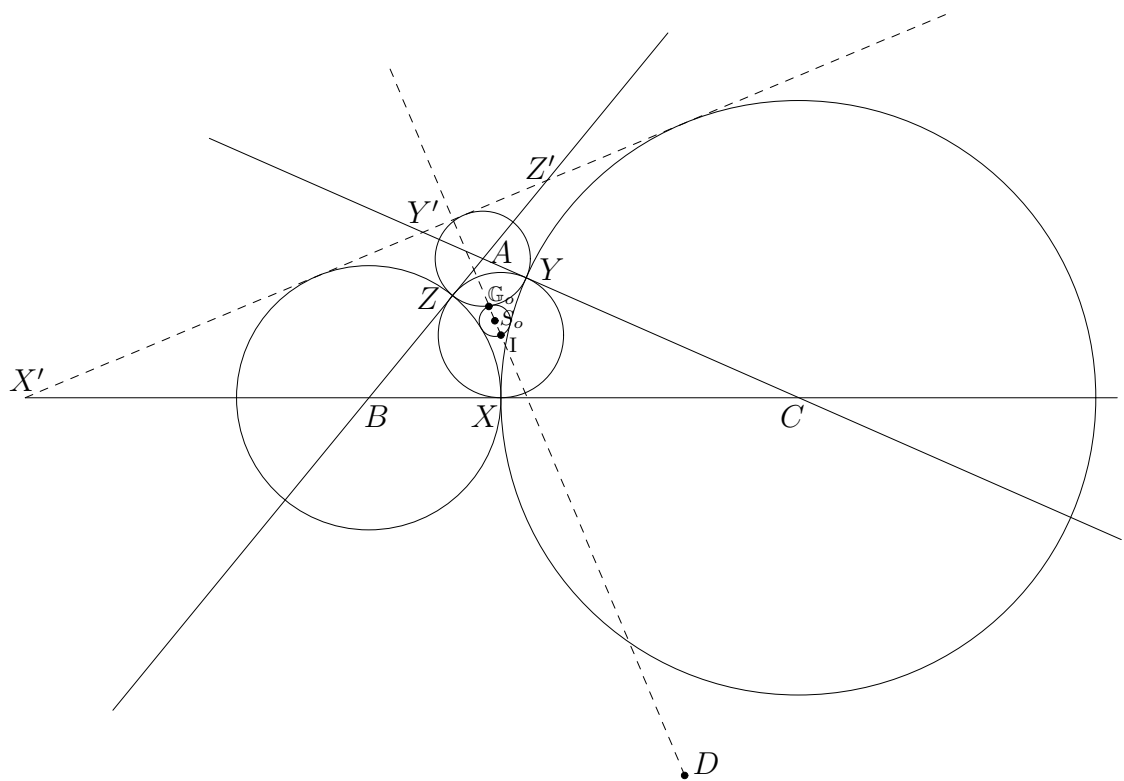

Figure \ref{critical} shows the critical case, in which the
outer Soddy circle is a straight line, coincident with the
Gergonne line.  One Soddy centre is at infinity, and since these
centres form a harmonic range with the Gergonne point,
${\mathbb{G}}_o$, and the incentre, $I_o$, the (inner)
Soddy centre, $S_o$, is the midpoint of their join.
Note that the Soddy line passes through the deLongchamps
point, that is, the reflexion of the orthocentre of $ABC$
in its circumcentre.

We will rename and redefine the other cases as {\bf external}, with
$2s<4R+r$, or {\bf internal}, with $2s>4R+r$.  Respectively,
the outer Soddy circle is touched externally or internally by
the tangent circles.  They correspond roughly to the triangle's
being more obtuse or less obtuse than in the critical case,
but this correspondence cannot be exact, since there are overlaps,
as we shall see in \S9.  This was the main error in \cite{V}.

\clearpage

Figure \ref{more} shows the external case.  The outer Soddy circle,
having $S_e$ as centre, is touched externally by the three tangent
circles.  The insert shows the detail in the neighborhood of the
vertex $A$ and the intersection of the Gergonne line $A'B'C'$ with
the Soddy line, which is perpendicular to it.  This intersection
would be better renamed as the {\it radical origin} of the triangle
$ABC$.  The Soddy line contains the Soddy centres, the
incentre, the Gergonne point and the deLongchamps point.
$S_o$ is not the midpoint of ${\mathbb{G}}_oI_o$ in this case.

\begin{figure}[h]  
\begin{picture}(400,290)(-270,-157)
\setlength{\unitlength}{1.85pt}
\drawline(-120,65)(84,-20)
\drawline(-120,0)(80,0)
\drawline(-120,-75)(80,75)
\thicklines
\drawline(-10,-82)(-10,-39.25)(-52,-39.25)(-52,-82)(-10,-82)  
\drawline(-20,0)(36,0)(0,15)(-20,0)    
\dashline[+30]{3}(-120,-2.72)(80,31.36)    
\dashline[+30]{3}(-10.25,73)(16,-81)       
\put(0,15){\circle{8}}             
\put(-20,0){\circle{42}}           
\put(36,0){\circle{70}}            
\put(0.591,9.396){\circle{3.268}}  
\put(-5.176,43.235){\circle{49.412}} 
\put(0.234,11.496){\circle*{1}}   
\put(1,7){\circle*{1}}            
\put(0.591,9.396){\circle*{1}}    
\put(-5.176,43.235){\circle*{1}}  
\put(16,-81){\circle*{1}}         
\put(-104,0){\circle*{1}}         
\put(-4.645,16.935){\circle*{1}}  
\put(4.706,18.529){\circle*{1}}   
\put(-4,43.5){$S_e$}
\put(17,-82){$D$}
\put(-22,-5){$B$}
\put(34,-5){$C$}
\put(-106,-6){$X'$}
\put(-30.298,-65.512){\circle*{1}}   
\put(-28,-79){\circle*{1}}            
\put(-29.227,-71.812){\circle*{1}}    
\put(-29.227,-71.812){\circle{9.804}}    
\put(-31,-55){\circle{24}}             
\drawline(-52,-70.75)(-10,-39.25)      
\drawline(-52,-46.25)(-10,-63.75)         
\put(-44.935,-49.195){\circle*{1.2}}  
\put(-16.882,-44.413){\circle*{1.2}}   
\put(-29,-66){${\mathbb{G}}_o$}
\put(-26,-81){$I_o$}
\put(-28,-74){$S_o$}
\put(-33,-53){$A$}
\put(-48,-47){$Y'$}
\put(-22,-45){$Z'$}
\put(-91,-100){\arc{126}{5.389}{5.986}}           
\put(77,-100){\arc{210}{3.313}{3.726}}            
\dashline[+30]{3}(-37,-26.2)(-25.75,-92.2)         
\dashline[+30]{3}(-59.2,-51.625)(0.2,-41.5)       
\end{picture}
\caption{The external case.  The insert shows the inner
Soddy circle being touched externally by the $A$-tangent
circle and (arcs of) the $B$- and $C$-tangent circles;
also portions of the Gergonne line, $X'Y'Z'$, and the
(perpendicular) Soddy line, which passes through
the incentre, $I_o$, and the Gergonne point, ${\mathbb{G}}_o$,
of $ABC$ (but not through $A$).}
\label{more}
\end{figure}
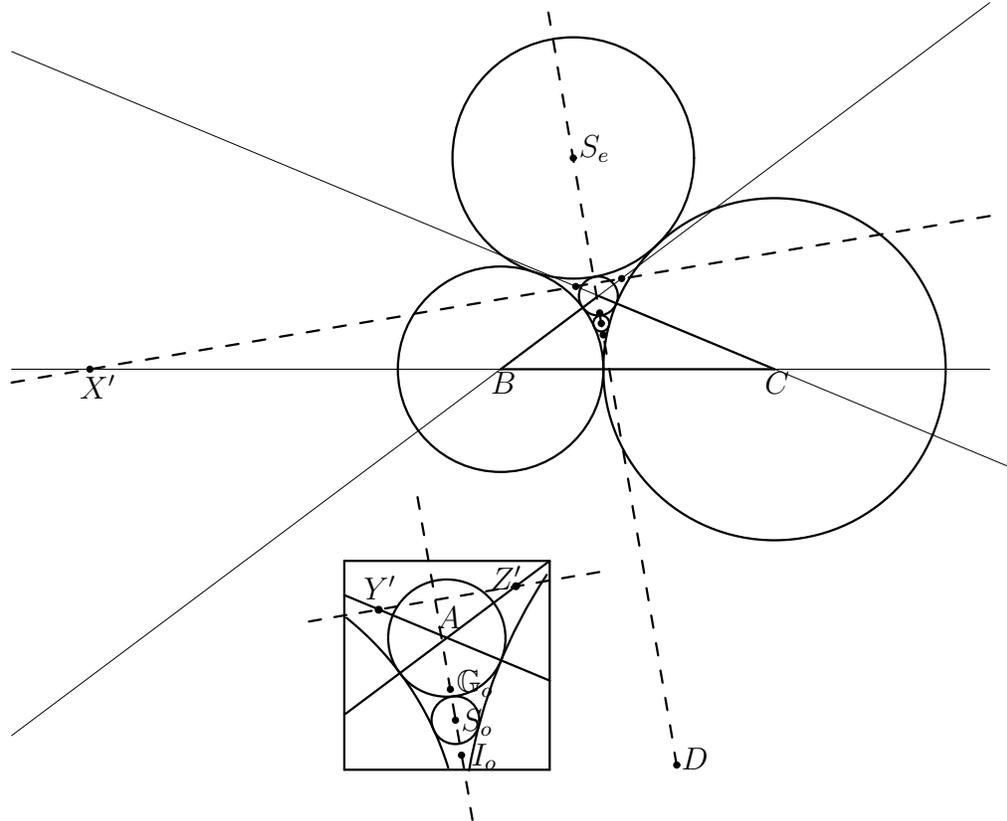

\clearpage

Figure \ref{less} shows the internal case in which the outer
Soddy circle is touched internally by the tangent circles.
Only an arc of this outer circle is shown.  The dashed lines
are the Gergonne and Soddy lines.  The dotted lines are the
lines of centres of the outer Soddy circle with the tangent
circles.  The three unlabelled points on the Soddy line are
the Gergonne point, the inner Soddy centre, and the incentre.
They form a harmonic range with the outer Soddy centre.
Contrary to appearances, the Gergonne point and the incentre
do not lie on the inner Soddy circle.

\begin{figure}[h]  
\begin{picture}(420,270)(-230,-175)
\setlength{\unitlength}{6.5pt}
\drawline(-27,0)(15,0)    
\drawline(-4.091,10.182)(15,0)  
\drawline(-6,0)(4.5,14)           
\dashline[+3]{0.3}(-27,0)(4.5,14)    
\dashline[+3]{0.3}(15,-28)(-2.3196,10.9691)    
\dottedline{0.3}(15,-28)(15,14)
\dottedline{0.3}(15,-28)(-1.1538,10.7692)
\dottedline{0.3}(15,-28)(-10.2,5.6)
\put(0.3913,4.8696){\circle*{0.3}}      
\put(1,3.5){\circle*{0.3}}           
\put(9,-14.5){\circle*{0.3}}       
\put(0,8){\circle{6}}  
\put(-6,0){\circle{14}}  
\put(15,0){\circle{28}}  
\put(0.7021,4.1702){\circle{1.7872}} 
\put(15,-28){\arc{84}{3.2}{5}}  
\put(0.7021,4.1702){\circle*{0.3}} 
\put(15,-28){\circle*{0.3}}  
\put(9,-14.5){\circle*{0.3}}  
\put(-7,-1.5){$B$}
\put(0.4,7.7){$A$}
\put(15.3,-0.6){$C$}
\put(-0.6,-1.3){$X$}
\put(3.4,6.3){$Y$}
\put(-28,0.3){$X'$}
\put(-5.2,10.3){$Y'$}
\put(4,14.3){$Z'$}
\put(7.6,-15.5){$D$}
\put(15.3,-28.5){$S_e$}
\end{picture}
\caption{The internal case.  The three tangent circles touch
the outer Soddy circle internally.}
\label{less}
\end{figure}


\section{Isoperimetric and equal detour points}

In this internal case, with $S$ the centre of the outer
Soddy circle ($S_e$ in the figure), and $\rho$ its radius,
$$SA=\rho-(s-a), \ SB=\rho-(s-b), \ SC=\rho-(s-c),$$
$$SB+SC+BC=\rho-(s-b)+\rho-(s-c)+a=2\rho=SC+SA+CA=SA+SB+AB$$
by symmetry and the perimeters of triangles $SBC$, $SCA$, $SAB$
are equal and $S$ is an {\bf isoperimetric point}.  Such a point
exists only in the internal case.

The inner Soddy circle
is always touched by the tangent circles externally so that,
if $S$ is its centre and $\rho$ its radius,
then
$$SA=\rho+s-a, \ SB=\rho+s-b, \ SC= \rho+s-c,$$
$$SB+SC-BC=\rho+s-b+\rho+s-c-a=2\rho=SC+SA-CA=SA+SB-AB$$
by symmetry, and $S$ is an {\bf equal detour point} in the
sense that if you
traverse the edges of the triangle by detouring through $S$, the
extra distance is the same for each edge.  In the external case
there are two such equal detour points; in the internal
and critical cases, only the inner one.  See X(176) and X(175)
in \cite{K}.

\section{Numbertheoretic interlude}

Andrew Bremner has kindlily supplied the solution in
the critical case for all rational triangles with $2s=4R+r$.
If the area of $ABC$ is $\Delta$, and the perimeter is
$a+b+c=2s$, then $4R=abc/\Delta$ and $r=\Delta/s$, so
the critical case has $abc/\Delta+\Delta/s=2s$, and
\begin{eqnarray*}
\Delta^2+abcs & = & 2\Delta s^2 \\
(s-a)(s-b)(s-c) + abc & = & 2\Delta s \\
-s^2+bc+ca+ab & = & 2\Delta \\
(bc+ca+ab-s^2)^2 & = & 4s(s-a)(s-b)(s-c)
\end{eqnarray*}
which, in terms of $a$, $b$, $c$ only, simplifies (?) to
$$5(a^4+b^4+c^4)-4\left(a^3(b+c)+b^3(c+a)+c^3(a+b)\right)-2(b^2c^2
+c^2a^2+a^2b^2)+4abc(a+b+c)=0$$
He observes that this equation defines a quartic curve
in projective two-dimensional space with arithmetic
genus 3, which is highly singular, and so has geometric
genus 0, and is parametrizable!  For example, 
$$a=8u^2(u^2+25v^2) \qquad b,c = 5(u\pm 5v)^2(u^2\mp2uv+5v^2)$$
or any rational multiple thereof.
This solution satisfies the triangle inequality, except that
the triangle is degenerate when $u=\pm5v$.  Also,
$$\frac{b^2+c^2-a^2}{2bc}+\frac{7}{25}=
\frac{32v^2(6u^2+25v^2)}{25(u^4+6u^2v^2+25v^4}\geq0$$
with equality only if $v=0$, when the triangle is isosceles
with shape (8,5,5).  So, except in this case,
$\cos A < -\frac{7}{25}$ and angle $A$ is always less than
$2\arcsin4/5$.  Examples of triangles with $2s=4R+r$ are
(45,40,13), (160,153,25),
(325,261,126), (425,416,41).

Compare such triangles with those (rational) triangles
having an angle $\arccos-\frac{7}{25}$.
The cosine formula gives us $25(b^2+c^2-a^2)=-14bc$, which
may be written
$$(25a)^2=(24c)^2+(25b+7c)^2$$
with solution  \qquad
$a:b:c=24(p^2+q^2):-7p^2+48pq+7q^2:25(p^2-q^2)$ \\
which satisfies the triangle inequality if $7q>p>q>0$.

The only (shape of) triangle with both
$2s=4R+r$, and angle $\arccos-\frac{7}{25}$, is (8,5,5),
given by $p=2q$.
Other triangles with this angle are (26,25,3), (30,25,11),
(51,38,25), (136,125,29).  They all have $2s<4R+r$.

Armed with numerical examples, and with help from Peter
Moses, we can now clarify the relationship between the
three cases and the size of the greatest angle, say
$A$, of the triangle $ABC$.  The external case comprises
X1, Y1, and Z1, where $A$ may be less than,
greater than, or equal to $2\arcsin\frac{4}{5}$.
X2 is the internal case, where $A$ must exceed
$2\arcsin\frac{4}{5}$ (there are no cases Y2, Z2).  The
critical case comprises X3 and Z3, where $A$
may not exceed $2\arcsin\frac{4}{5}$ (there is no case Y3).
\begin{tabbing}
dun \=  ex \=  a bit more space t \= \kill
X. \> $A<2\arcsin\frac{4}{5}$. \\
   \>             \> X1) \ $2s<4R+r$.  \> Outer Soddy
                         radius negative. Two equal detour points.  \\
 \> \> \> No isoperimetric point. \\
 \> \> \quad (23,22,3). $48=2s<4R+r=\frac{130}{\sqrt{7}}\approx49.135$.
$A=2\arcsin\sqrt{\frac{7}{11}}<2\arcsin\frac{4}{5}$. \\		
      \>          \> X2) \ $2s>4R+r$.  \> Outer Soddy
                         radius positive. Inner equal detour point only. \\
   \> \> \> Isoperimetric point. \\
 \> \> \quad (6,5,5). $16=2s>4R+r=\frac{25}{2}+\frac{3}{2}=14$.
$A=2\arcsin\frac{3}{4}<2\arcsin\frac{4}{5}$. \\
   \>             \> X3) \ $2s=4R+r$.  \> Outer Soddy
	                 radius infinite. Inner equal detour point. \\
   \> \> \> Outer detour \& isoperimetric points at infinity. \\
 \> \> \quad (45,40,13). $98=2s=4R+r=\frac{650}{7}+\frac{36}{7}=98$.
$A=2\arcsin\frac{9}{\sqrt{130}}<2\arcsin\frac{4}{5}$. \\
Y. \> $A>2\arcsin\frac{4}{5}$. \\
   \>             \> Y1) \ $2s<4R+r$.  \> Outer Soddy radius negative.
                     Two equal detour points. \\
   \> \> \> No isoperimetric point. \\
 \> \> \quad (56,39,25). $120=2s<4R+r=130+7=137$.
$A=2\arcsin\frac{7}{\sqrt{65}}>2\arcsin\frac{4}{5}$. \\
Z. \> $A=2\arcsin\frac{4}{5}$. \\
   \>             \> Z1) \ $2s<4R+r$.  \> Outer Soddy radius negative.
                     Two equal detour points. \\
   \> \> \> No isoperimetric point. \\
 \> \> \quad (26,25,3). $54=2s<4R+r=\frac{325}{6}+\frac{4}{3}=
55\frac{1}{2}$.  $A=2\arcsin\frac{4}{5}$. \\
   \>             \> Z3) \ $2s=4R+r$.  \> Outer Soddy
	                 radius infinite.
                     Inner equal detour point. \\
   \> \> \> Outer detour \& isoperimetric point at infinity. \\            		
 \> \> \quad Isosceles triangle (8,5,5). $18=2s=4R+r=
\frac{50}{3}+\frac{4}{3}=18$.  $A=2\arcsin\frac{4}{5}$. \\
\end{tabbing}


\section{There are four (or eight) Pavillet tetrahedra}

[In fact, with extraversion, quadration, twinning and reflexion
in the plane of the triangle, there are 64 Pavillet tetrahedra.
The midpoints of the joins of the orthocentres of a twin pair
are the 32 Gergonne points of the triangle (cf. \S 1.10 and
Figure \ref{16gerg})]

\begin{figure}[h]  
\begin{picture}(420,390)(-150,-112)
\setlength{\unitlength}{1.2pt}
\drawline(0,144)(-17,0)(108,0)(0,144)    
\dashline{2}(9.38,223.45)(0,144)
\dashline{2}(0,144)(-63,228)    
\dashline{2}(108,0)(208,0)    
\dashline{3}(0,64)(-108,228)    
\dashline{3}(0,144)(0,64)         
\dashline{3}(108,0)(108,-100)             
\dashline{3}(-17,0)(-17,225)               
\dashline{1.5}(64,99.69)(64,44.31)   
\dashline{1.5}(208,0)(0,64)            
\dashline{1.5}(9.38,223.45)(108,-100)    
\dashline{1.5}(-17,225)(64,44.31)               
\dashline[30]{1.5}(138.5,-99)(-89.7,228.6)   
\dashline[30]{1.5}(219.7,208.15)(-96.2,-11.9)  
\put(208,200){\arc{400}{1.5}{3.3}}     
\thicklines
\drawline(0,64)(108,-100)(-17,225)(0,64)  
\drawline(208,200)(0,64)        
\drawline(208,200)(108,-100)         
\drawline(208,200)(-17,225)           
\put(48,80){\circle*{2}}     
\put(208,0){\circle*{2}}  
\put(9.38,223.45){\circle*{2}}  
\put(69.54,0){\circle*{2}}    
\put(-4.46,106.23){\circle*{2}}    
\put(64,99.69){\circle*{2}}    
\put(64,44.31){\circle*{2}}    
\put(91,118.5){\circle*{2}}    
\put(93,114){$D$}
\put(66,38){$O_a$}
\put(66,94){${\mathbb{G}}_a$}
\put(209,194){$I_a$}
\put(2,143){$B$}
\put(-25,-4){$A$}
\put(110,-9){$C$}
\put(-29,221){$A_a$}
\put(-8,58){$B_a$}
\put(110,-101){$C_a$}
\put(50,81){$X_a$}
\put(204,-9){$Y_a$}
\put(11,224){$Z_a$}
\put(-85,223){$(X'_a)$}
\put(60,-8){$Y'_a$}
\put(-15,103){$Z'_a$}
\put(140,-100){{\scriptsize Gergonne line}}
\put(-95,-16){{\scriptsize Soddy line}}
\end{picture}
\caption{The Pavillet $a$-tetrahedron, $A_aB_aC_aI_a$}
\label{aflip}
\end{figure}
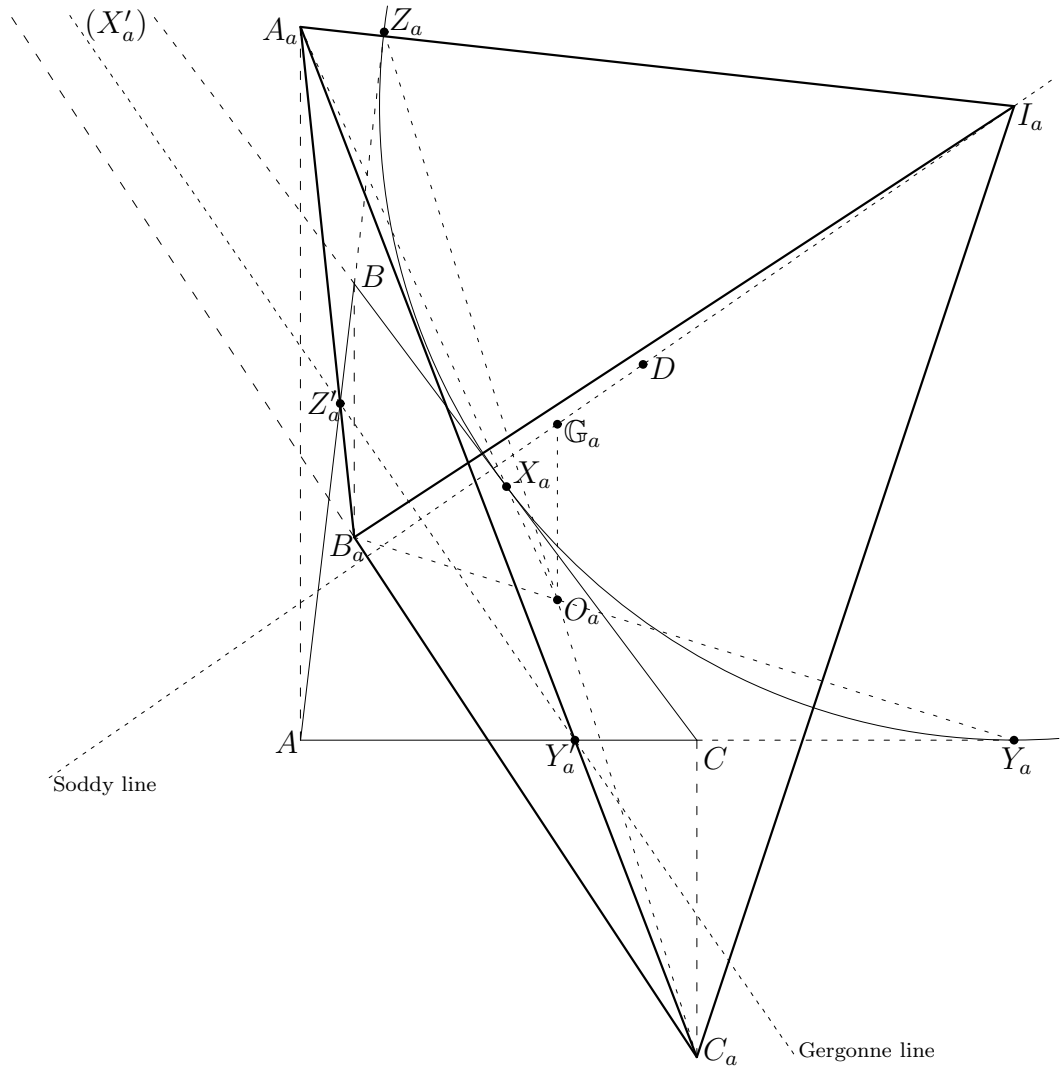

\clearpage

Conway's {\bf extraversion} often gives us several items for
the price of one.  For example, an $a$-flip replaces the angles
$A$, $B$, $C$ of a triangle by $-A$, $\pi-B$, $\pi-C$ respectively,
which changes the sign of $a$, and makes the interchanges

$$s\longleftrightarrow s-a \mbox{\quad and \ }s-b
\longleftrightarrow-(s-c)$$

so if we erect perpendiculars to the plane of $ABC$
of lengths $s$, $c-s$, $b-s$ at $A$, $B$, $C$ respectively
(note that the last two are negative) we arrive at points,
say $A_a$, $B_a$, $C_a$, which form an orthocentric tetrahedron
with $I_a$, the $a$-excentre of $ABC$ (Figure \ref{aflip}).

The projection of the orthocentre of this $a$-flip tetrahedron
onto the plane $ABC$ is the Gergonne point,
${\mathbb{G}}_a$, the point of concurrence of $AX_a$, $BY_a$, $CZ_a$,
where $X_a$, $Y_a$, $Z_a$ are the respective touchpoints of
the $a$-excircle with $BC$, $CA$, $AB$.
The altitudes of the tetrahedron are $A_aX_a$, $B_aY_a$,
$C_aZ_a$, which concur at $O_a$, together with $I_aO_a$,
which projects onto the Soddy line, $I_a{\mathbb{G}}_a$, in the plane
$ABC$.  This last contains the deLongchamps point and the
centres of the two circles which touch the three tangent
circles, centres $A$, $B$, $C$, and radii $s$, $s-c$, $s-b$
respectively.

So, there are four Pavillet tetrahedra and it may be useful to
consider their reflexions in the plane $ABC$ as well, giving
a total of eight.
There are four manifestations of many of the
creatures mentioned above (but only one deLongchamps point).
This was noticed nearly half a century ago by Vandeghen
\cite{Va} who also observed that, if $\odot$ is the circumcentre of
$ABC$, then $\odot I_o$, $\odot I_a$, $\odot I_b$, $\odot I_c$
are the Euler lines for the triangles formed by the touchpoints,
$X_j$, $Y_j$, $Z_j$ ($j=o,a,b,c$) of the four touchcircles.
Figure \ref{extra} shows four Gergonne lines (dotted), four
Soddy lines (dashed), eight Soddy points, four Gergonne
points, four touchcentres, \ldots

\begin{figure}[h]
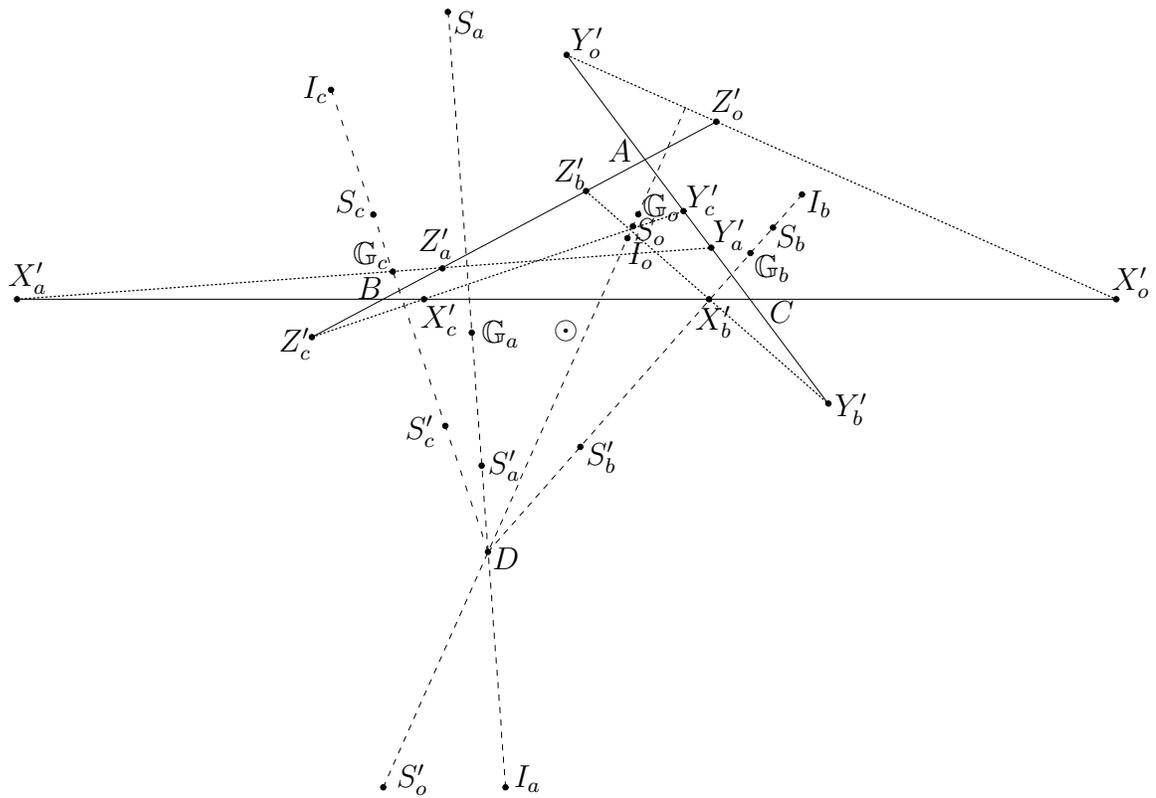
  

\caption{Extraversions.  Four Soddy lines each contain a
touch-centre, $I_j$ ($j=o,a,b,c$); a Gergonne point,
${\mathbb{G}}_j$; two Soddy points, $S_j$ and $S'_j$; and
the deLongchamps point, $D$. Four Gergonne lines,
$X'_jY'_jZ'_j$, each contain three harmonic conjugates
of touchpoints.}
\label{extra}
\end{figure}

\clearpage

Twenty-four of the points in Figure \ref{extra} appear
in Figure \ref{fourco} as black dots (not all labelled).
Each of four Soddy lines contains a touch-centre,
$I_j$, and two Soddy centres, $S_j$, $S'_j$
($j=o$, $a$, $b$, $c$).  The corresponding circles have the
Gergonne line perpendicular to the Soddy line as radical
axis.  Each Gergonne line contains three centres,
$X'_j$, $Y'_j$, $Z'_j$; and the circles intersect in
points on the corresponding Soddy line.  The edges of $ABC$
are shown dashed.

\begin{figure}[h]
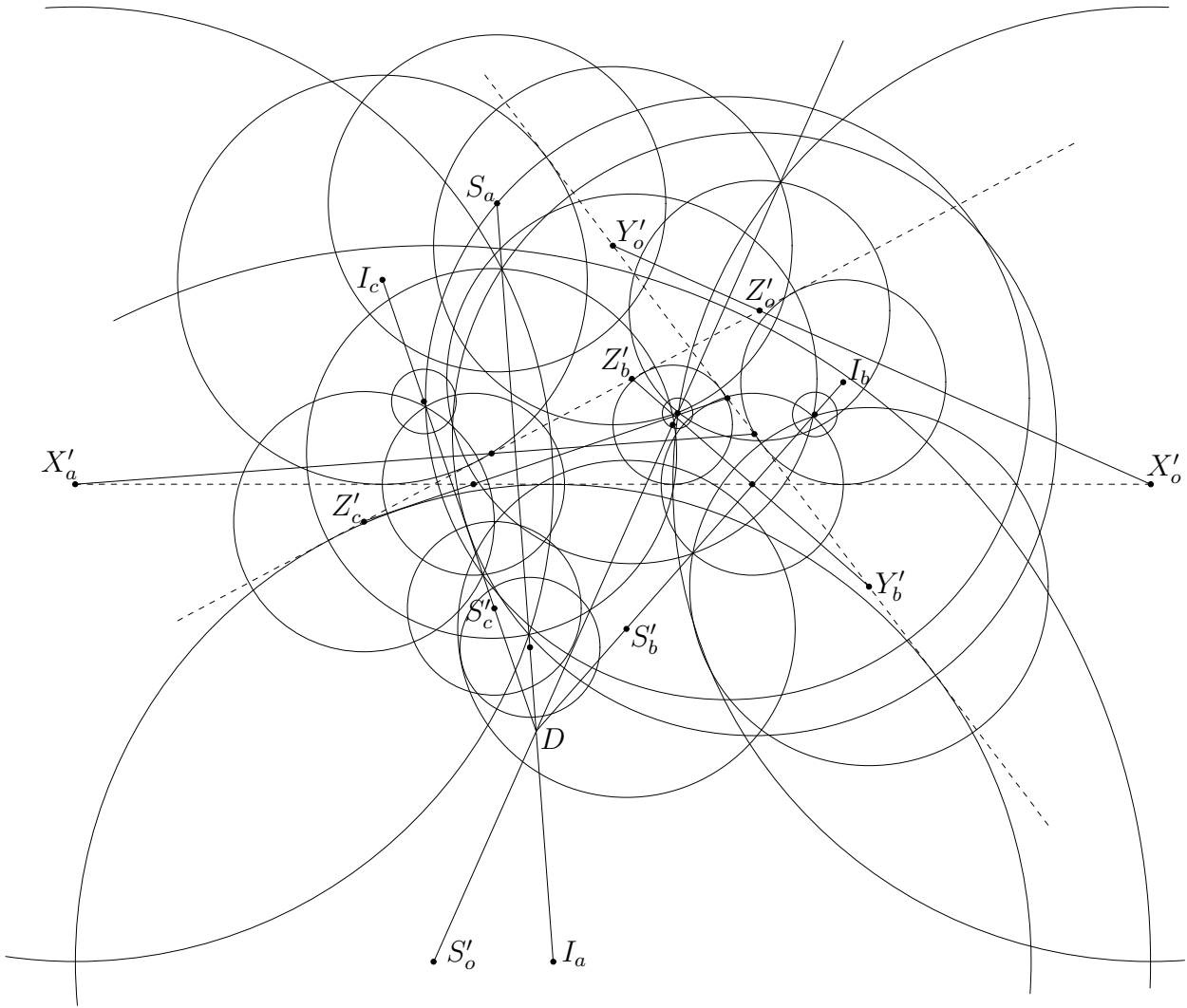
  

\caption{Four pairs of orthogonal coaxal systems.  Three circles
intersect in pairs of points on each of the four Soddy lines (lines
through $D$).  Three circles have centres on each of the four
Gergonne lines (perpendicular to the Soddy lines).  Contrary to
appearances, the $Y'_b$ and $Y'_c$ circles do not intersect on
the edge $BC$.}
\label{fourco}
\end{figure}

\clearpage

\section{L'esprit de l'escalier}

No sooner had I submitted this paper that I realized that
I should have included a picture of all the extraversions
of the Pavillet tetrahedron.  Recall that they are
$$A_oB_oC_oI_o, \quad A_aB_aC_aI_a, \quad A_bB_bC_bI_b,
 \quad A_cB_cC_cI_c$$
where $I_o$, $I_a$, $I_b$, $I_c$ are
the touch-centres of triangle $ABC$ and that
$$\{A_o,B_o,C_o\}, \quad \{A_a,B_a,C_a\}, \quad
\{A_b,B_b,C_b\}, \quad \{A_c,B_c,C_c\}$$
are at respective distances
$$\{s-a,s-b,s-c\}, \{s,-(s\!-\!c),-(s\!-\!b)\},
\{-(s\!-\!c),s,-(s\!-\!a)\}, \{-(s\!-\!b),-(s\!-\!a),s\}$$
vertically above the points $(A,B,C)$.

Had I drawn such a picture
I would have noticed that, since $A_oA_a$, $B_oB_a$,
$C_oC_a$ are equal and parallel (length $a$ and
perpendicular to plane $ABC$) the pairs of corresponding
edges $A_oB_o$ \& $A_aB_a$, and $A_oC_o$ \& $A_aC_a$,
bisect each other, and that their common points,
$E_1$, $E_2$ lie on the respective 260?-point spheres.
Also, $B_oC_o$ is equal and parallel to $B_aC_a$, and the
join of their midpoints is perpendicular to the
plane of $ABC$.
Moreover, the right triangles $A_oAI_o$ and $A_aAI_a$ are
similar, since the ratio of their edge-lengths is $s-a\,:\,a$,
so that $A_oI_o$ is parallel to $A_aI_a$.  All of these
remarks still hold when the letters $ABC$ and $abc$ are
permuted cyclically.

The three points $E_1$, $E_2$, $E_4$ in Figure
\ref{finale} indicate points of intersection of three
260?-point spheres (not the same set of spheres in each case).
The nine dotted lines are diameters of these spheres,
three through each of these points of intersection,
three through the 260?-point centre, $G_o$, and two through
each of the other three 260?-point centres, $G_a$, $G_b$, $G_c$.
They join pairs of midpoints of opposite edges of tetrahedra.

\begin{figure}[h]
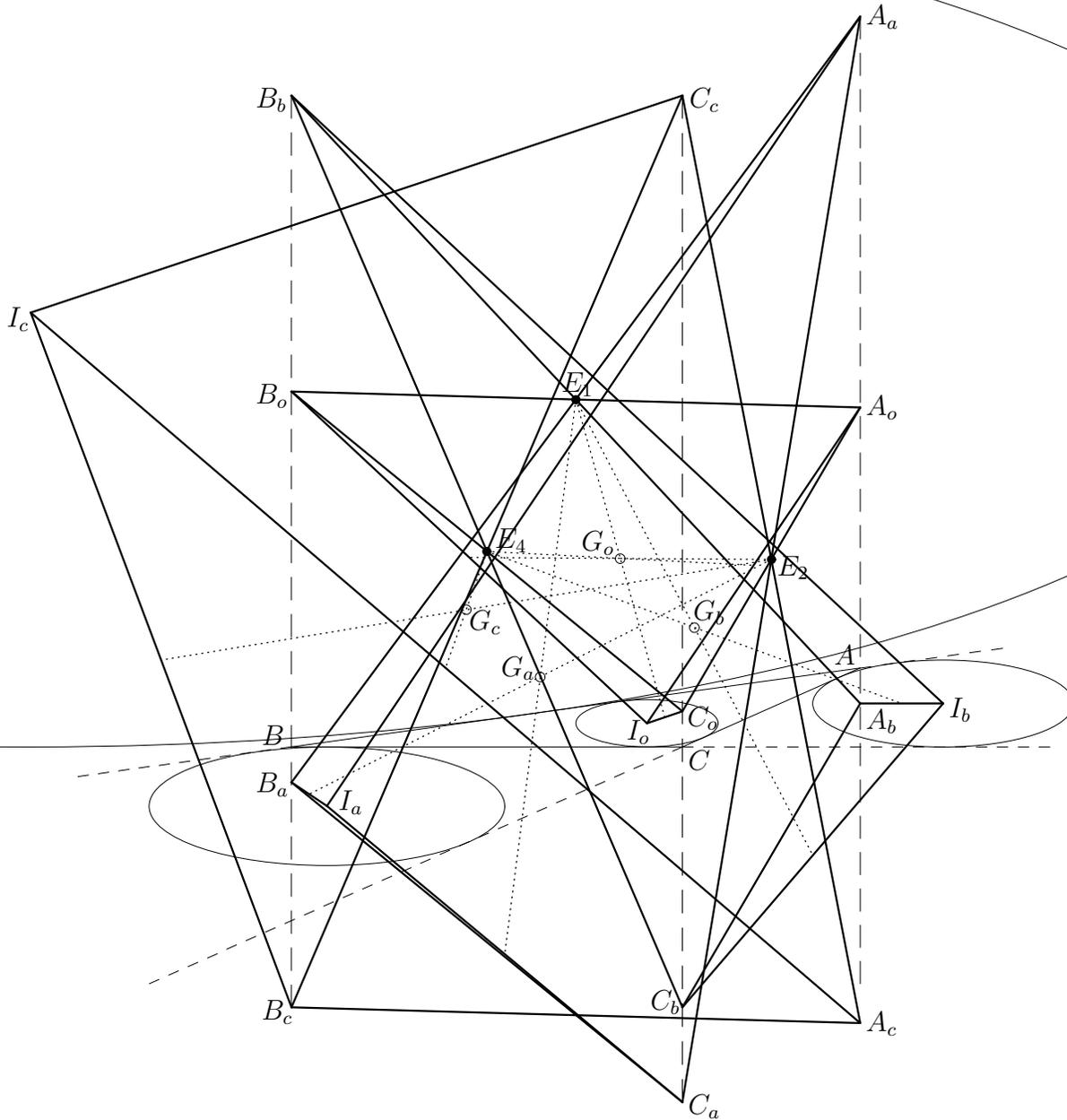
  

\caption{Extraversions of the Pavillet tetrahedron.  The
ellipses are perspective drawings of touch-circles
of the triangle $ABC$.  The (large) $c$-circle is only
partly accommodated.  Small circles indicate the centres,
$G_j$ ($j=o,a,b,c$), of the associated 260?-point spheres.
The points $A_b$, $B_a$ and $C_o$ are not in the plane
$ABC$, but are respectively at distance $s-c$ below,
below, and above it.}
\label{finale}
\end{figure}

\clearpage

\section{So little done; so much to do}

{\bf1.}  What are the various perspectrices of the triangles
of interior or of exterior altitude points relative to the
triangle $ABC$ or to its pedal triangle, $DEF$ ?  All
four triangles have the orthocentre as perspector.

{\bf2.}  Given an orthocentric tetrahedron, are there sixteen planes
(or some subset thereof), four through each vertex, such that
the vertex is an incentre or excentre of a triangle for which
the tetrahedron is a Pavillet tetrahedron?  ANSWER NO!!  Most
orthocentric tetrahedra will not serve as Pavillet tetrahedra.

{\bf3.}  For the potentially up to sixteen triangles in {\bf2.},
are there any incidences among their sixteen Euler lines,
or 64 Soddy lines, or 64 Gergonne lines?  The Soddy lines concur
in fours at sixteen deLongchamps points; do they form any particular
configuration?

{\bf4.}  For these sixteen conjectured triangles, are there
any collinearities among the more important of their triangle
centres?

{\bf5.}  Are there further interesting points on the 260?-point
sphere where it intersects the four opposite pairs of parallel,
congruent, negatively homothetic faces; and the three rectangular
diametral planes; of the octahedron $E_1E_2E_3E_4E_5E_6$ ?

{\bf6.}  Same question, but for the four faces of the
medial tetrahedron, $O_7O_6O_5O_3$.

{\bf7.}  What relations are there between the four (or eight)
260?-point spheres associated with the four (or eight)
Pavillet tetrahedra?

{\bf8.}  Characterize the common points of the three circles,
centres $X'$, $Y'$, $Z'$, passing respectively through the
touchpoints $X$, $Y$, $Z$, of the incircle with the edges
of $ABC$.  These circles form a coaxal system with the Gergonne
line as its line of centres and the Soddy line as radical axis,
so the common points are on the Soddy line.  Do this also for
the extraversions based on the three excircles.

{\bf9.}  The coaxal system which is orthogonal to that mentioned
in {\bf8.}\ contains the incircle and the two Soddy circles and
has the Gergonne line as radical axis.  Investigate the member
of the system having the deLongchamps point as centre [see
Figure \ref{coaxal}].

{\bf10.}  As the deLongchamps point lies on each of the four
Soddy lines (obtained by extraversion; see Figure \ref{extra})
there are in fact four (concentric) circles of the kind
mentioned in {\bf9.}  Investigate their relationship.

{\bf11.}  Examine the critical cases, where the outer Soddy
circle is a straight line, for the three extraversions.
Compare Figures \ref{critical} and \ref{extra}.

{\bf12.}  Investigate tetrahedra whose six edges all touch a
sphere.  It is easy to see that the three sums of the lengths
of pairs of opposite edges are equal.  They have four `tangent
spheres' \cite{Ep} and a `Descartes line' through the centres
of the two spheres which touch all four tangent spheres. A
particular case is the equilateral pyramid, with
$CD=DB=BC\neq AB=AC=AD$, which is also an orthocentric
tetrahedron.

{\bf Acknowledgements}

I am grateful for discussions with several people, including
Andrew Bremner, Fred Lunnon, Peter Moses and Axel Pavillet, and
perhaps (?) for useful suggestions by the referees.

\section{Referees' comments.}

Reviewer \#1: The paper describes some amazing plane geometry
theorems that are new to me, and I thought I knew a fair amount of geometry.
I found this interesting, but somewhat overwhelming. There are few proofs
given, and those that are in the paper are terse almost to the point of
unintelligibility.

I think this paper tries to do too much. The point seems to be to
acquaint readers with some ``new'' geometry theorems, and that is a
worthwhile goal for a paper in the Monthly, but I wonder if discussing
fewer theorems in a little more detail would better serve the readership.

Reviewer \#2:
This paper contains some really interesting results in Euclidean triangle
geometry.  Despite the clear interest of the mathematics, I recommend
against publication on the basis of the exposition.  I don't think many
Monthly readers will be willing to do the work required to read the paper.

The author makes no attempt to explain why his results are interesting or
to put them in any sort of context---he simply launches into technical
proofs.  This appears to be a deliberate strategy and the paper begins
with the words ``I'll tell the story backwards?''.  But I don't think
it works.  There is nothing that would lure the casual reader into the
paper.  The audience is limited to those who already know quite a bit
about the subject.

A Monthly paper should be accessible to a fairly broad audience of
mathematicians, but anyone other than an expert on the kind of Euclidean
geometry discussed here would have to do a tremendous amount of
background work to get very far into the paper.  Of course there's
always the problem of how many definitions to include and which ones
can be assumed, but I don't think the author's solution is a good
one.  He is safe in assuming that his readers will know what the
medians and orthocenter of a triangle are, but the very first theorem
statement on page 2 requires that the reader know what a ``radical
center'' is.  The first line of the proof on the top of page 3 refers
to the ``radical axis''.  There is no hint of definition provided.
Page 4 uses the words ``perspector'' and ``pedal triangle''.  This
means that the readership of the paper is effectively limited to
those who already know what all these words mean, which is far too
restrictive for a Monthly article.  Anyone else would have to dig out
all the definitions from other sources in order to have any clue what
the author is talking about, and I just don't see that happening.

The mathematics in the paper appears to be sound, but some of it
is difficult to follow.  The diagrams are helpful, but some of the
more intricate diagrams are quite obscure.  For example, it is
difficult for me to see what a human being could be expected to
learn from looking at Figure 8.  It looks like a mass of circles
drawn around a triangle and does not provide me with any insight
into what is going on.

Despite the fact that I don't believe the paper is right for the
Monthly in its present form, I hope the author will consider
reworking the exposition and resubmitting the paper.  Euclidean
geometry is a highly suitable topic for the Monthly and the results
in this paper appear to have the right combination of beauty and
surprise.

\chapter{The Lighthouse Theorem}

This is an almost verbatim copy of \cite{G}.  The notation
will be at variance with earlier attempts.  Some errors are
corrected and improvements made in an inchoate paper of
Conway \& the present writer, given as Chapter 7 below.

\section{Introduction}

The story began in 1993,
when Dick Bumby, then Problems Editor of this
{\sc Monthly}, sent me the following problem
to referee:

\begin{quote}
Prove that every prime $p>7$ of the form $3n+1$ can be
written as \\
$p=\sqrt[6]{a^2+4762800b^2}$ for a unique
choice of natural numbers $a$ and $b$.
\end{quote}

\noindent
{\bf Paradox 1.}  The problem was deemed unsuitable for
the {\sc Monthly}, but here it is.

Fortunately the proposer, Joseph Goggins, had included
the motivation for the problem.  He was searching for
integer-edged triangles whose Morley triangle
was also integer-edged.  The Morley triangle
theorem is still not as well-known as it
deserves to be.  Its comparatively short
history is outlined in \S{\bf6.26}.
Figure \ref{morley} shows its simplest form.  The pairs
of angle-trisectors of the triangle $ABC$ meet
at points which form an equilateral triangle
$PQR$, regardless of the shape of the original
triangle.  Goggins found an infinity of triangles
$ABC$ with integer edge-lengths whose Morley
triangle also had an integer edge-length.

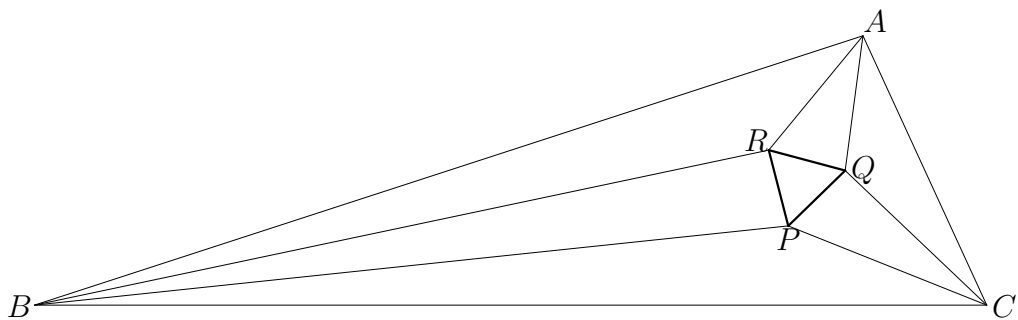
\begin{figure}[h]  
\begin{picture}(400,110)(-40,0)
\setlength{\unitlength}{1.5pt}
\drawline(0,0)(240.1,0)(208.915,67.989)(0,0)
(185.131,39.116)(208.915,67.989)(204.397,33.968)
(240.1,0)(190.05,20.005)(0,0)
\thicklines
\drawline(190.05,20.005)(204.397,33.968)
(185.131,39.116)(190.05,20.005)
\put(209,69){$A$}
\put(-7,-3){$B$}
\put(241.5,-3){$C$}
\put(187,14){$P$}
\put(206,32){$Q$}
\put(179,39){$R$}
\end{picture}
\caption{The simplest form of Morley's Theorem.}
\label{morley}
\end{figure}

The combination of geometry and number theory is dear
to my heart, and I set to work on the general problem.
To trisect angles you need to solve cubic equations,
which have three roots.  Also, just as angle-bisectors
come in pairs, so do (pairs of) trisectors come in
triples.  I soon discovered that there was not one,
but as many as 18 Morley triangles.  I excitedly
rang Coxeter, but he claimed not to know that; nor
did another geometer who was visiting him.  Perhaps
they didn't want to destroy my delight of discovery?

I rang John Conway, and he said, ``Funny thing!  I was
just looking at that last weekend.''  So I dashed off
a paper \cite{CG}, added John's name to it, and sent
it to the {\sc Monthly}.  A while later the editor,
John Ewing, said that he didn't normally intervene
in the refereeing process, but the situation was
unusual.  One referee had said something like,
``O.K., but there are references and other things
missing.''  A second said, ``Isn't this like a paper
I refereed for you a few years ago?''  Indeed it was!
It had been written by the first referee, John Rigby,
who hadn't resubmitted it.

\noindent
{\bf Paradox 2.}  Rigby knew that there are
18 Morley triangles, and he knew that Morley
knew, but few others seemed to know.

Our paper wasn't resubmitted, either.

\noindent
{\bf Paradox 3.}  But here it is!

Goggins's real problem was answered by the following:

\noindent
{\bf Theorem}\cite{BGGG}.  {\it If a rational-edged
triangle has a rational Morley triangle, then
either the original triangle is equilateral
$($and 6 of the 18 Morley triangles
are rational---in fact, congruent to the
original triangle$\,)$, or it is Pythagorean
and belongs to a one-parameter family $($and
just 2 of the 18 Morley triangles
are rational$\,)$, or it belongs to a
two-parameter family of triangles $($and all
18 Morley triangles are rational}\,).

The Pythagorean family has edge-lengths
$$2t(3-t^2)(1-3t^2),\quad(1-t^2)(1-14t^2+t^4),
\quad(1+t^2)^3$$
for rational $t$, while the more general
solution has edge-lengths
$$x_i(x_i^2-1)(x_i^2-9)/(x_i^2+3)^3 \qquad (i=1,2,3)$$
where \quad $3(x_1+x_2+x_3)=x_1x_2x_3$.

We omitted to mention an observation of Goggins that
the lengths of the trisectors are also rational.  In
fact the reader may like to reconstruct the
Conway-Doyle 7-piece jigsaw proof of Morley's
theorem by cutting out an equilateral triangle of
edge 1001, three triangles with edge-length triples
(1001,1716,1859), (1001,9464,9555), (1001,2695,2464)
and three with edge-length triples (9555,2695,12005),
(2464,1716,3740), (1859,9464,10985).  If you want
to be sure that they fit together, use the cosine law
to calculate the angles.

\noindent
{\bf Paradox 4.}  Although \cite{BGGG} is a
numbertheoretic result, published in {\it Acta
Arithmetica}, it was reviewed in {\it MR} under
the heading Geometry.

But our purpose here is to discuss:

\section{The Lighthouse Theorem.}
  {\it Two sets of $n$ lines
at equal angular distances, one set through each
of the points $B$, $C$, intersect in $n^2$ points
that are the vertices of $n$ regular $n$-gons.
The circumcircles of the $n$-gons each pass
through $B$ and $C$.}

That is, they form a {\bf coaxal system} with $BC$
as {\bf radical axis}.  In the exceptional case
that the sets of lines are parallel, one $n$-gon
is at infinity.

\noindent
{\bf Paradox 5.}  Although the theorem is quite
striking and an immediate consequence of the
following well-known propositions of Euclid, it
doesn't seem to be in the literature:

III.20.  The angle at the centre of a circle is
twice that at the circumference.

III.21.  Angles in the same segment are equal.

III.22.  Opposite angles of a cyclic quadrilateral
are supplementary.

The first of these also gives: the angle in a
semicircle is a right angle; and the last leads
to: an exterior angle of a cyclic quadrilateral
is equal to the interior opposite angle. And,
in the limiting case: the angle
between tangent and chord is equal to the
angle in the alternate segment.  For example,
in Figure \ref{liteproof}, the chord $PU$
makes two supplementary angles with the tangent
at $P$ (shown dotted); these are equal
in size to $\angle PVU$ and $\angle PBU$, one in
each of the segments into which $PU$ partitions
the circle.

\medskip

\noindent
{\bf Informal demonstration of the Lighthouse
Theorem for $n=3$.}

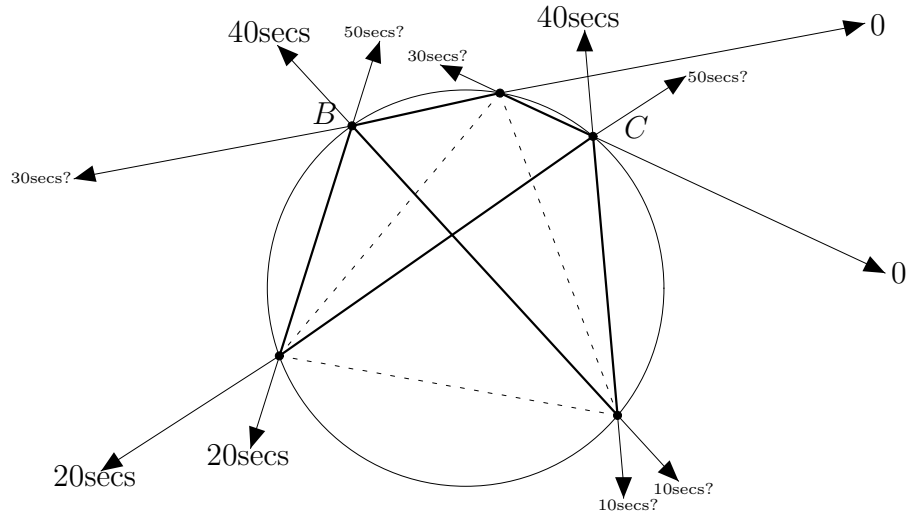
\begin{figure}[h]  
\begin{picture}(400,190)(-225,-80)
\put(0,0){\circle{150}}
\put(-70.48,-25.65){\circle*{3}}
\put(-43.02,61.44){\circle*{3}}
\put(13.02,73.86){\circle*{3}}
\put(48.21,57.45){\circle*{3}}
\put(57.45,-48,21){\circle*{3}}
\put(-58,62){$B$}
\put(60,57){$C$}
\put(145,98.7){\vector(4,1){6}}
\put(-68.14,88.85){\vector(-1,1){3}}
\put(-142.08,42.81){\vector(-4,-1){6}}
\put(77.54,-70.14){\vector(1,-1){3}}
\put(-79.3,-54.68){\vector(-1,-3){2}}
\put(-33.5,90.47){\vector(1,3){1}}
\put(153.78,8.22){\vector(2,-1){5}}
\put(-4.58,82.06){\vector(-2,1){5}}
\put(80.38,78.23){\vector(3,2){3}}
\put(-134.83,-67.20){\vector(-3,-2){3}}
\put(59.76,-74.62){\vector(0,-1){5}}
\put(45.13,92.67){\vector(0,1){5}}
\put(153,96){0}
\put(161,2){0}
\put(-90,93){40secs}
\put(27,99){40secs}
\put(-98,-67){20secs}
\put(-156,-75){20secs}
\put(-22,86){{\tiny 30secs?}}
\put(-172,40){{\tiny 30secs?}}
\put(50,-84){{\tiny 10secs?}}
\put(71,-78){{\tiny 10secs?}}
\put(-46,95){{\tiny 50secs?}}
\put(84,78){{\tiny 50secs?}}
\dashline{2}(13.02,73.86)(57.45,-48.21)(-70.48,-25.65)(13.02,73.86)
\drawline(-43.02,61.44)(-142.08,42.81)
\drawline(-43.02,61.44)(-68.14,88.85)
\drawline(-43.02,61.44)(-33.87,90.47)
\drawline(48.21,57.45)(153.78,8.22)
\drawline(48.21,57.45)(80.38,78.23)
\drawline(48.21,57.45)(45.13,92.67)
\drawline(13.02,73.86)(145.1,98.7)
\drawline(13.02,73.86)(-4.58,82.06)
\drawline(-70.48,-25.65)(-79.63,-54.68)
\drawline(-70.48,-25.65)(-134.83,-67.20)
\drawline(57.45,-48.21)(77.54,-70.14)
\drawline(57.45,-48.21)(59.76,-74.62)
\thicklines
\drawline(-43.02,61.44)(13.02,73.86)
\drawline(-43.02,61.44)(-70.48,-25.65)
\drawline(-43.02,61.44)(57.45,-48.21)
\drawline(48.21,57.45)(13.02,73.86)
\drawline(48.21,57.45)(-70.48,-25.65)
\drawline(48.21,57.45)(57.45,-48.21)
\end{picture}
\caption{Part of the Lighthouse Theorem for $n=3$.}
\label{lite}
\end{figure}

In Figure \ref{lite}, each of two lighthouses at $B$ and $C$
has one doubly-infinite beam, and each rotates
with a uniform angular velocity of one
revolution per minute. It's night-time.  I take
photographs every 20 seconds and superimpose
them.  The locus of the point of intersection of
the beams is a circle (Euclid III.21).  The
point traces out the circle with uniform angular
velocity twice that of the lighthouses (Euclid III.20).
At 20-second intervals the beams will be equally inclined
at angles that are multiples of $\pi/3$, and the
points on the circle will be at angular distances
$2\pi/3$.  In other words, they are the vertices of
an equilateral triangle (the dashed lines in Figure
\ref{lite}).  \qquad $\blacksquare$

\noindent
{\bf Paradox 6.}  I can't tell from the superimposed
photographs whether they were taken every 20 seconds,
or every 10 seconds.  Indeed, one lighthouse might
have been going round twice as fast as the other.

\noindent
{\bf Paradox 7.}  I can't tell the senses of rotation
from my photographs.

This time, ignorance is bliss, because, in Morley's
theorem, it's important to choose the intersections
of the {\bf proximal} trisectors of the angles.  In
this context I will think of the lighthouses as
rotating in opposite senses.  But now the locus
of intersection of their beams is not a circle,
but a rectangular hyperbola.  In either case the
curve passes through the two lighthouses.

\noindent
{\bf Paradox 8.}   There appear to be a circle and a
rectangular hyperbola intersecting in five points.

The Lighthouse Theorem comprises much more than the
bald statement that I made earlier.  I will also
prove:

\section{The Lighthouse Lemma.}
{\it The
$\binom{n}{2}$ edges of an $n$-gon arising in the
Lighthouse Theorem are parallel to those of any
of the other $n\!-\!1$ $n$-gons.}

The Lighthouse Lemma could be expressed by saying
that the $n$-gons are homothetic, but this is so only
if $n$ is odd.  For example, if $n=4$, there are
two parallel squares and two more at an angle
of $\pi/4$ (see Figure \ref{fourlite}).  So, on third
thoughts, I include as edges all $\binom{n}{2}$
joins of vertices to each other (and could include
the joins of vertices to themselves, that is the
tangents to the circumcircles at the vertices) so
that, for $n=4$, the ``diagonals'' of one pair of
squares are parallel to the ``sides'' of the other
pair.

There is not enough room here to deal
with all of the aspects of the Lighthouse
Theorem, but I will at least prove the important

\noindent
{\bf Lighthouse Duplication Theorem.}  {\it The points
$U$, $P$, $V$ are consecutive vertices of one of
the $n$ $n$-gons arising in the Lighthouse Theorem}
(see Figure \ref{lemmaproof}).  {\it $X$ and $Y$
are the intersections
of the next beams with the original ones, making
$XY$ an edge of the neighboring $n$-gon, with
$XY$ parallel to $UV$ by the Lighthouse Lemma.
If $UP$ and $VP$ intersect $XY$ at $Q$ and
$R$, respectively, then $\angle RBP=\angle PBC$
and $\angle QCP=\angle PCB$, where $B$ and $C$ are
the lighthouses.}

It's time to justify some of our statements.

\medskip

\section{Proofs of Lighthouse Results}

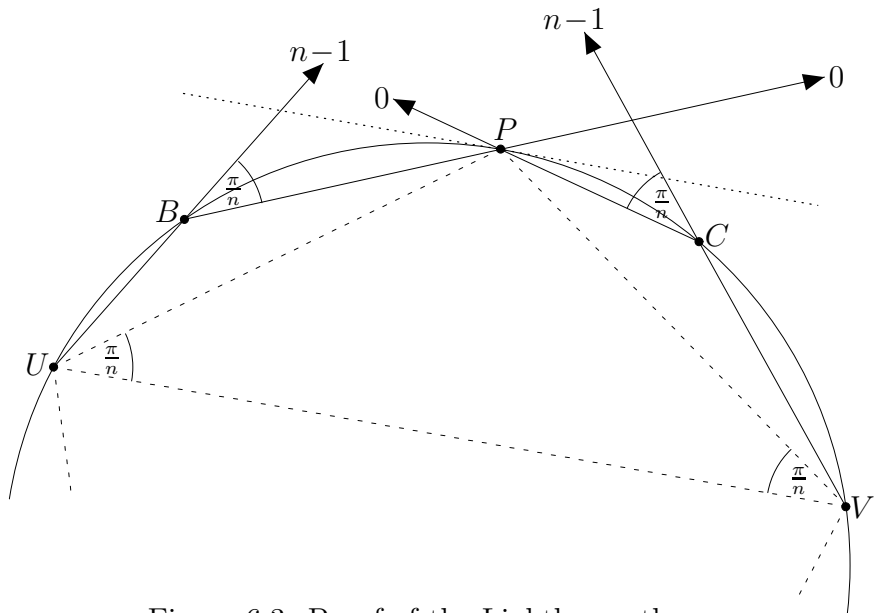
\begin{figure}[h]  
\begin{picture}(400,195)(-218,0)
\put(0,0){\arc{320}{3.3}{6.4}}
\put(-141.27,75.11){\circle*{3}}
\put(-91.77,131.06){\circle*{3}}
\put(27.78,157.56){\circle*{3}}
\put(102.84,122.56){\circle*{3}}
\put(158.44,22.26){\circle*{3}}
\dashline{2}(-134.73,28.55)(-141.27,75.11)(27.78,157.56)(158.44,22.26)(140.78,-10.95)
\dashline{2}(-141.27,75.12)(158.44,22.26)
\drawline(-141.27,75.12)(-91.772,131.06)
\drawline(-91.772,131.06)(-42.27,187.01)
\drawline(-91.772,131.06)(27.78,157.57)
\drawline(27.78,157.57)(-9.75,175.07)
\drawline(102.85,122.57)(27.78,157.57)
\drawline(27.78,157.57)(147.34,184.07)
\drawline(158.44,22.26)(102.85,122.57)
\drawline(102.85,122.57)(61.55,197.8)
\dottedline{3}(-92.22,178.72)(147.78,136.3955)
\put(-42.27,187.01){\vector(1,1){3}}
\put(-9.75,175.07){\vector(-2,1){3}}
\put(147.34,184.07){\vector(4,1){3}}
\put(61.55,197.8){\vector(-1,2){2}}
\put(-20,173){0}
\put(-52,191){$n\!-\!1$}
\put(152,181){0}
\put(44,203){$n\!-\!1$}
\put(-152,72){$U$}
\put(-103,130){$B$}
\put(25,161){$P$}
\put(105,121){$C$}
\put(160,18){$V$}
\put(-141.27,75.11){\arc{60}{5.85}{6.46}}
\put(-91.77,131.06){\arc{60}{5.464}{6.072}}
\put(102.84,122.56){\arc{60}{3.58}{4.2}}
\put(158.44,22.26){\arc{60}{3.32}{3.95}}
\put(-123,76){{\footnotesize$\frac{\pi}{n}$}}
\put(-77,140){{\footnotesize$\frac{\pi}{n}$}}
\put(85,135){{\footnotesize$\frac{\pi}{n}$}}
\put(137,30.5){{\footnotesize$\frac{\pi}{n}$}}
\end{picture}
\caption{Proof of the Lighthouse theorem.}
\label{liteproof}
\end{figure}

{\it Proof of the $($naive$)$ Lighthouse Theorem.}
In Figure \ref{liteproof}, $BP$ and $CP$ are
beams from lighthouses at
$B$ and $C$, and $U$, $P$, $V$ are consecutive
vertices of a regular $n$-gon inscribed in the
circle $BPC$.  Then $\angle PUV=\angle PVU=\pi/n$.
Since $UBPV$ is cyclic, the exterior angle between
the directions $UB$ and $BP$ is also $\pi/n$.
Similarly for the angle between the directions
$VC$ and $CP$.  The regular $n$-gon is generated
by beams through $B$ and $C$ which are equally
spaced by angles that are multiples of $\pi/n$.
Generally, the $n^2$ intersections form $n$ regular
$n$-gons whose circumcircles pass through $B$ and
$C$ and form a coaxal system. \quad $\blacksquare$

{\it Proof of the Lighthouse Lemma.}
I write $\beta$ and $\gamma$ for the magnitudes of
angles $CBP$ and $BCP$, respectively, and refer to
them as the {\bf phases} of the two lighthouses.
Number the beams from each lighthouse from 0 to
$n\!-\!1$ cyclically towards the baseline $BC$,
so that they are counted in opposite senses.  In
Figure \ref{lemmaproof}, $BPY$ and $UBX$ are
beams 0 and $n\!-\!1$
from lighthouse $B$, while $CPX$ and $VCY$ are
beams 0 and $n-1$ from lighthouse $C$.

Now $\angle XBY=\pi/n=\angle XCY$, so that
$BCYX$ is cyclic.  Then, by Euclid III.21,
$\angle YXC=\angle YBC=\beta$, implying
that $XY$ makes an angle $\gamma-\beta$ with
$BC$.  Now the difference between angles $BUV$
and $XBC$ is $(\pi/n+\gamma)-(\pi/n+\beta)=\gamma-\beta$,
making $XY$ parallel to $UV$ (and also parallel
to the tangent at $P$ to the circle $UBPCV$).
Any edge of any $n$-gon thus makes an angle with
any other edge of any $n$-gon that is a multiple
of $\pi/n$. \qquad $\blacksquare$

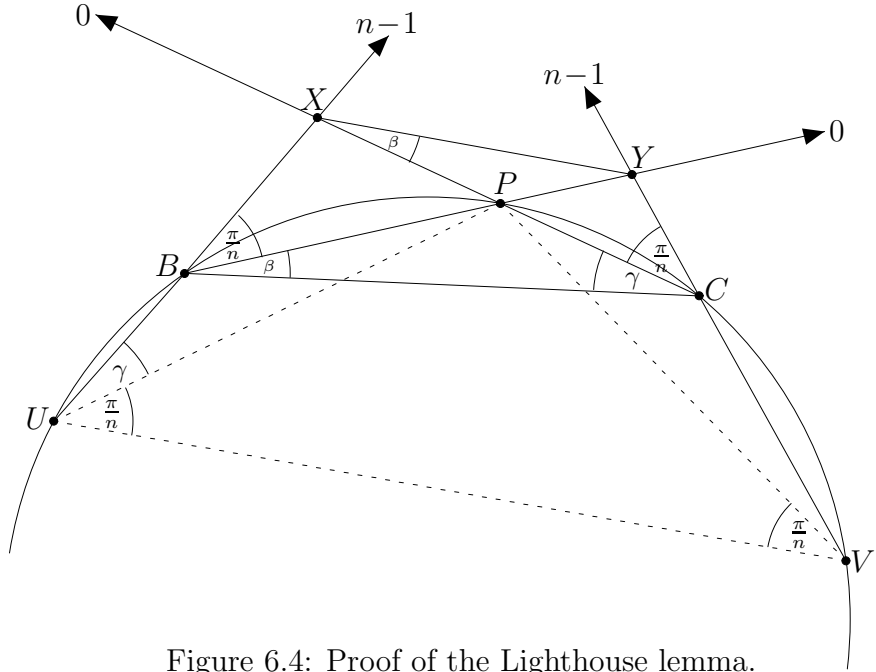
\begin{figure}[h]  
\begin{picture}(400,225)(-215,5)
\put(0,0){\arc{320}{3.3}{6.4}}
\put(-141.27,75.11){\circle*{3}}
\put(-91.77,131.06){\circle*{3}}
\put(27.78,157.56){\circle*{3}}
\put(102.84,122.56){\circle*{3}}
\put(158.44,22.26){\circle*{3}}
\put(-41.5,190){\circle*{3}}
\put(77.5,168.5){\circle*{3}}
\dashline{2}(-141.27,75.11)(27.78,157.56)(158.44,22.26)
\dashline{2}(-141.27,75.12)(158.44,22.26)
\drawline(-141.27,75.12)(-91.772,131.06)
\drawline(-91.772,131.06)(-17.52,218)
\drawline(-91.772,131.06)(27.78,157.57)
\drawline(27.78,157.57)(-122.34,227.57)
\drawline(102.85,122.57)(27.78,157.57)
\drawline(27.78,157.57)(147.34,184.07)
\drawline(158.44,22.26)(102.85,122.57)
\drawline(102.85,122.57)(61.55,197.8)
\drawline(-91.772,131.06)(102.85,122.57)
\drawline(-41.5,190)(77.5,168.5)
\put(-17.52,218){\vector(1,1){3}}
\put(-122.34,227.57){\vector(-2,1){3}}
\put(147.34,184.07){\vector(4,1){3}}
\put(61.55,197.8){\vector(-1,2){2}}
\put(-133,225){0}
\put(-27,222){$n\!-\!1$}
\put(152,181){0}
\put(44,203){$n\!-\!1$}
\put(-152,72){$U$}
\put(-103,130){$B$}
\put(25,161){$P$}
\put(105,121){$C$}
\put(160,18){$V$}
\put(-48,193){$X$}
\put(77,171){$Y$}
\put(-141.27,75.11){\arc{60}{5.85}{6.46}}
\put(-91.77,131.06){\arc{60}{5.464}{6.072}}
\put(102.84,122.56){\arc{60}{3.58}{4.2}}
\put(158.44,22.26){\arc{60}{3.32}{3.95}}
\put(-141.27,75.11){\arc{80}{5.44}{5.82}}
\put(-91.77,131.06){\arc{80}{6.07}{6.32}}
\put(102.84,122.56){\arc{80}{3.2}{3.58}}
\put(-41.5,193){\arc{80}{6.54}{6.8}}
\put(-123,76){{\footnotesize$\frac{\pi}{n}$}}
\put(-77,140){{\footnotesize$\frac{\pi}{n}$}}
\put(85,135){{\footnotesize$\frac{\pi}{n}$}}
\put(137,30.5){{\footnotesize$\frac{\pi}{n}$}}
\put(-62,132){\tiny$\beta$}
\put(-119,92){\footnotesize$\gamma$}
\put(-15,179){\tiny$\beta$}
\put(75,127){\footnotesize$\gamma$}
\end{picture}
\caption{Proof of the Lighthouse lemma.}
\label{lemmaproof}
\end{figure}

Observe that all edges of all $n$-gons make angles
with $BC$ that differ from $\gamma-\beta$ only
by multiples of $\pi/n$.  Our proofs still
hold if $\beta$ or $\gamma$ (or both) are changed
by a multiple of $\pi/n$, although the pictures
look different.

\newpage

{\it Proof of the Lighthouse Duplication Theorem.}
The notation of Figure \ref{dupproof} is as
in Figure \ref{lemmaproof}, and
the edges $UP$ and $VP$ intersect the edge $XY$ at
$Q$ and $R$, respectively.  From the Lighthouse Lemma,
$\angle QRP=\pi/n$ and is therefore equal to
$\angle PBX$.  Hence $RXBP$ is cyclic and
$\angle RBP=\angle RXP=\angle YXC=\angle YBC=\beta$.
Accordingly, $R$ lies on a beam through $B$ with
double the phase, $2\beta$.  Similarly, $Q$ lies on a beam
through $C$ with double the phase, $2\gamma$.
\qquad $\blacksquare$

\begin{figure}[h]
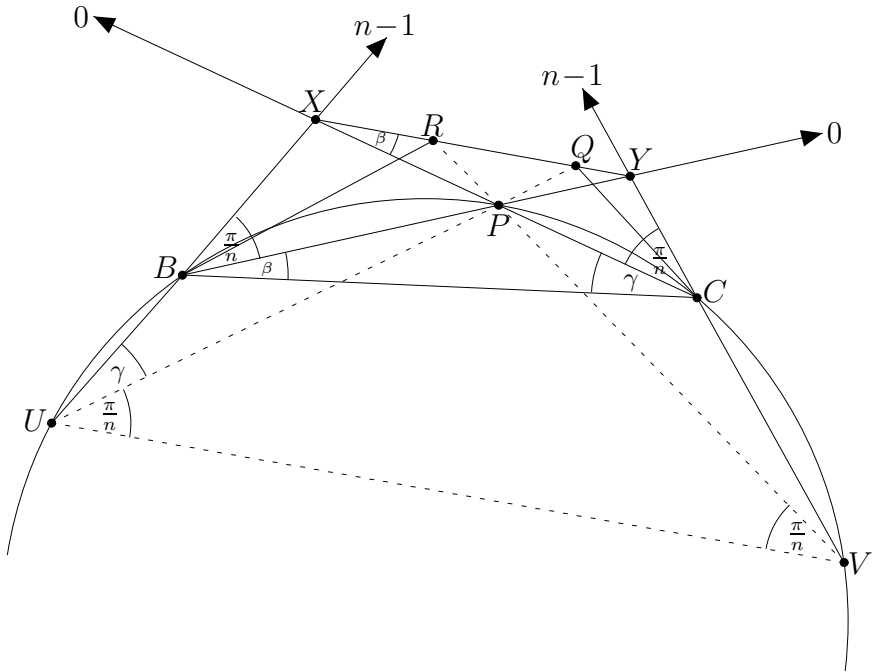
  

\caption{Proof of the Lighthouse Duplication Theorem.}
\label{dupproof}
\end{figure}

\noindent
{\bf Note.}  To avoid repetition when I prove
Morley's theorem, it is convenient to observe that
$QYCP$ and $BCYX$ are also cyclic; that
$\angle PRQ=\angle PQR=\pi/n$; that
$\angle BPX=\angle CPY=\beta+\gamma$; that
$\angle BPR=\angle CPQ=\pi/n+\gamma$; and that
$\angle BRP=\angle CQP=(n\!-\!1)\pi/n-\beta-\gamma$.
This last angle may be written $(n\!-\!2)\pi/n+\alpha$,
where $\alpha+\beta+\gamma=\pi/n$.

The Lighthouse Duplication Theorem has many
ramifications.  The edges of the $n$-gons form
$n$ families of $\binom{n}{2}$ parallel lines\footnote
{The regular $n$-gons are thought of as having
every pair of vertices joined, yielding $\binom{n}{2}$ edges.}.
Two families intersect in $\binom{n}{2}^2$
points, so that the complete configuration
contains $\binom{n}{2}^3$ points, though here
some points (for example, the vertices of the $n$-gons)
are counted by multiplicity.  We now know that
several of these points lie on additional beams
through the lighthouses, which in turn generate
new $n$-gons (or at least $(n/2)$-gons if $n$
is even) whose edges intersect again, and so on
indefinitely.  Also, as we shall see in the proof
of Morley's theorem, some sets of points lie on new
lines, beams from additional lighthouses.

\medskip

\section{Morley's Theorem}
For $n=3$ the Law of
Small Numbers \cite{Gsm} intervenes, because
$\binom{n}{2}$ is no larger than $n$, and the
process is in some sense closed, and we have

\smallskip

\noindent
{\bf The Morley Miracle.}  {\it The nine edges of the
equilateral triangles of the Lighthouse Theorem
for $n=3$ are the Morley lines of a triangle.}

\noindent
{\bf Paradox 9.}  A proof of the
complete Morley theorem, without even having a
triangle to start with.

In fact the 9 edges of the 3 triangles from
one pair of lighthouses form 27 equilateral
triangles of which 18 are genuine Morley
triangles.  The other 9 comprise 3 sets of
what Conway has called the {\bf Guy Faux
triangles}, one set from each of 3 pairs
of lighthouses situated at $B$ \& $C$, $C$
\& $A$, and $A$ \& $B$, where $A$ is a
mystery point, yet to be determined.
A good way to count the Morley triangles is to
notice that they are each formed from two sides
of a GF-triangle together with a third side
taken from a different GF-triangle:
$3\times3\times2=18$.

\smallskip

\noindent
{\bf Yet another proof of Morley's theorem.}
I will in fact prove more: what was probably
known to Morley, and certainly to Rigby \cite{Rj}
(see also \cite{R2}).

\noindent
{\bf Theorem.}  {\it The circumcircles of
the nine GF-triangles meet in threes, not only
at the vertices of the triangle $ABC$, but
also at nine other points.}

More precisely, any choice of two GF-circles
(that is, circumcircles of GF-triangles), one
from each of two of the three families through $B\&C$,
$C\&A$, $A\&B$, pass through two pairs of the
six Morley points on one of the 9 Morley lines.
These, together with the circle from the third
family that passes through the remaining two
points on the Morley line, concur in an
additional point: a configuration of
27+9+3 points, 9 circles and 9 lines, each circle
passing through 3+3+2 points and each line through
6+0+0 points (see Figure \ref{9lines}).

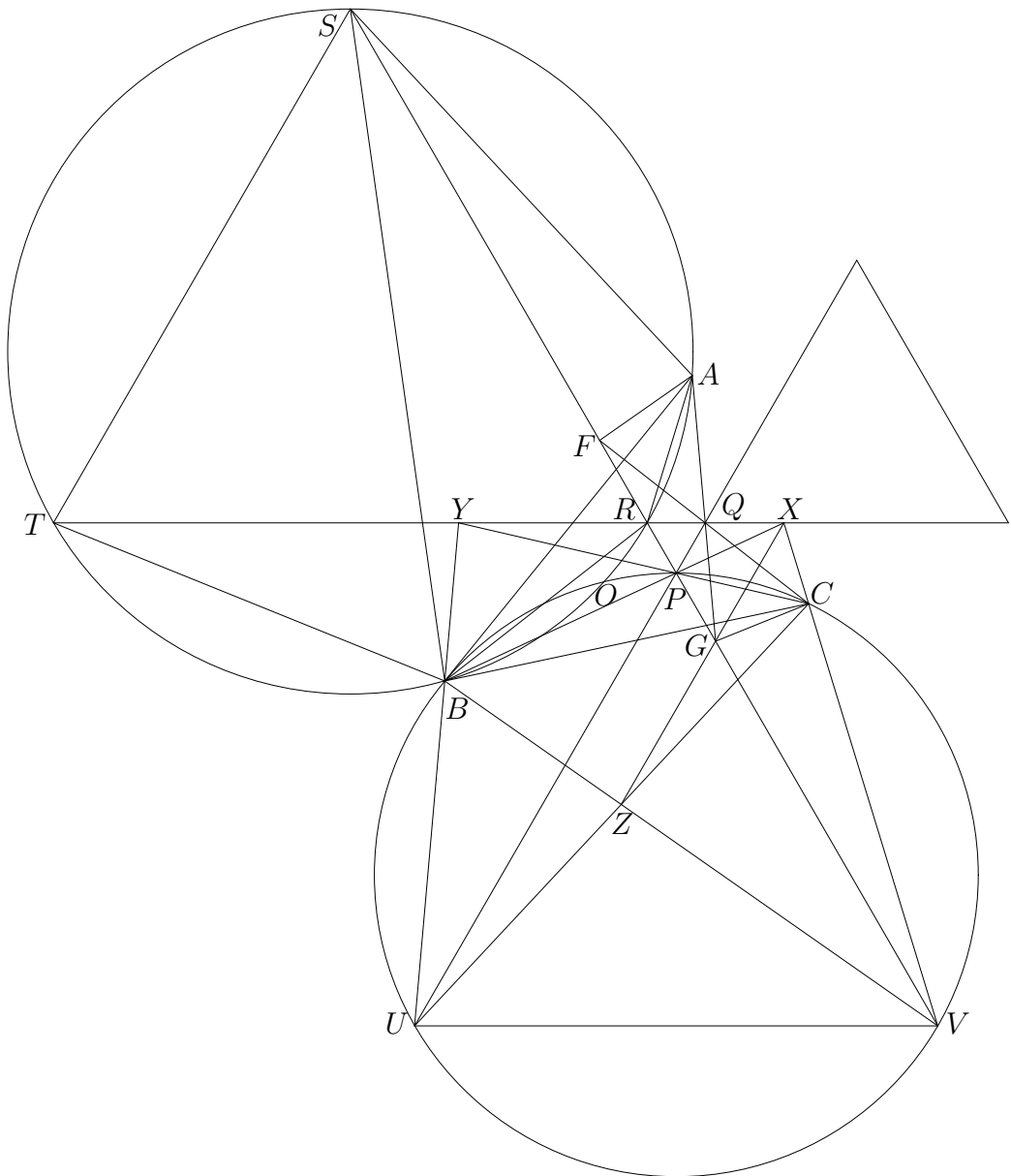
\begin{figure}[h]  
\setlength{\unitlength}{6pt}
\begin{picture}(65,72)(11,-7.5)
\put(59.45016,9.83048){\circle{39.32193}}
\put(38.23853,43.9175){\circle{44.62714}}
\drawline(18.9144,32.76071)(81.0856,32.76071)
\drawline(38.23853,66.23107)(76.47706,0)
\drawline(42.42326,0)(71.21163,49.86292)
\drawline(44.38898,22.46831)(60.49775,42.36098)
\drawline(44.38898,22.46831)(68.06896,27.50164)
\drawline(44.38898,22.46831)(57.56265,32.76071)
\drawline(44.38898,22.46831)(66.4611,32.76071)
\drawline(18.9144,32.76071)(38.23853,66.23107)
\drawline(81.0856,32.76071)(71.21163,49.86292)
\drawline(42.42326,0)(76.47706,0)
\drawline(42.42326,0)(45.28945,32.76071)
\drawline(42.42326,0)(68.06896,27.50164)
\drawline(44.38898,22.46831)(76.47706,0)
\drawline(76.47706,0)(66.4611,32.76071)
\drawline(54.46056,38.13369)(60.49775,42.36098)
\drawline(54.46056,38.13369)(68.06896,27.50164)
\drawline(60.49775,42.36098)(57.56265,32.76071)
\drawline(44.38898,22.46831)(38.23853,66.23107)
\drawline(60.49775,42.36098)(38.23853,66.23107)
\drawline(44.38898,22.46831)(18.9144,32.76071)
\drawline(45.28945,32.76071)(68.06896,27.50164)
\drawline(55.875275,14.42552)(66.4611,32.76071)
\drawline(62.01188,25.05442)(68.06896,27.50164)
\drawline(62.01188,25.05442)(60.49775,42.36098)
\put(60.8,41.8){$A$}
\put(44.4,20){$B$}
\put(68.2,27.5){$C$}
\put(52.7,37.1){$F$}
\put(60,24.1){$G$}
\put(54.1,27.4){$O$}
\put(58.6,27.1){$P$}
\put(62.4,33.3){$Q$}
\put(55.3,33){$R$}
\put(36.1,64.5){$S$}
\put(17,32){$T$}
\put(40.5,-0.5){$U$}
\put(77,-0.5){$V$}
\put(66,33){$X$}
\put(44.8,33){$Y$}
\put(55.2,12.5){$Z$}
\end{picture}
\caption{Proof of Morley's theorem.}
\label{morlproof}
\end{figure}

\clearpage

\noindent
{\bf Nothing that Euclid couldn't have done.}
In Figure \ref{morlproof} the labels are
as in earlier figures,
but in the particular case $n=3$, so that $PQR$ is
a typical potential Morley triangle.  The point
$A$ is defined by $\angle RBA=\beta$ and
$\angle RAB=\pi/3-(\beta+\gamma)=\alpha$, say.
[Warning: it's best not to connect the edge $CA$,
lest one assume things that, while
true, have not yet been established.]

Draw the circle $BRA$ and let $PR$ and $QR$ cut it
again at $S$ and $T$ respectively.  From triangle
$ABR$, $\angle ARB=\pi-\alpha-\beta$.  From the
Note following the proof of the Duplication Theorem,
$\angle BRP=\pi/3+\alpha$, so that
$$\angle QRA=2\pi-\pi/3-(\pi/3+\alpha)-
(\pi-\alpha-\beta)=\pi/3+\beta.$$  Hence $RST$ is a GF-triangle
generated by lighthouses at $A$ and $B$ with
phases $\alpha$ and $\beta$.  By the Duplication
Theorem, $\angle RAQ=\alpha$ and
$\angle RQA\!=\!\pi\!-\!\alpha\!-\!(\pi/3\!+\!\beta)\!=\!\pi/3\!+\!\gamma$.

Let $CQ$ meet $PR$ in $F$.  The Note also asserts
that $\angle CQP=\pi/3+\alpha$, whence $\angle RQF=
\pi-\pi/3-(\pi/3+\alpha)=\pi/3-\alpha$ and
$\angle RFQ=\alpha=\angle RAQ$, which shows that
the quadrilateral $RQAF$ is cyclic.  It follows that
$\angle CFA=\angle QFA=\angle QRA=
\pi/3+\beta$ and $\angle FAQ=\pi/3$.

Let $XZ$, an edge of the GF-triangle $XYZ$, meet
$PR$ in $G$.  Since $XZ$ is parallel to $PU$, the
quadrilateral $GXQP$ has angles $\pi/3$ and $2\pi/3$
and so is cyclic.  Moreover, $C$ lies on the same
circle, because the Note gives $\angle CQP=\pi/3+\alpha$
which is also the measure of $\angle CXP=
\angle CXB=\angle CYB$.  So $\angle GQP=\angle GXP=
\pi/3-\gamma$ and then $\angle RQA=\pi/3+\gamma$
implies that $A$, $Q$ and $G$ are collinear and
$AF$ and $AG$ are beams through $A$ inclined at $\pi/3$.
Also $\angle GCF=\angle GCQ=\angle GXQ=\pi/3$
so that $CF$ and $CG$ are beams through $C$ inclined
at $\pi/3$.  Hence the quadrilateral $AFGC$ is cyclic,
$\angle CAG= \angle CFG=\alpha$, and $\angle ACF=\angle AGF=
\angle QGP=\angle QXP=\gamma$, implying that
$FG$ is an edge of a GF-triangle, $EFG$ say,
generated by lighthouses at $C$ and $A$ with
phases $\gamma$ and $\alpha$, respectively.

Finally, let the GF-circles $ARB$ and $BPC$ meet
again at $O$.  Then
$\angle AOB=\angle ARB=\pi-\alpha-\beta$ and
$\angle BOC=\angle BPC=\pi-\beta-\gamma$.  This
means that
$$\angle COA=2\pi-(\pi-\alpha-\beta)-(\pi-\beta-\gamma)=
\pi/3+\beta=\angle CFA=\angle CGA$$
and that the GF-circles $BPC$, $C{\bf GF}A$,
$ARB$ concur at $O$. \qquad $\blacksquare$

The point $O$ is associated with the Morley
line $VGPRFS$ in the sense that the three circles
through $O$ intersect the line in the pairs of points
$\{P,V\};\ \{G,F\};\ \{R,S\}$.

\newpage

\begin{figure}[h]
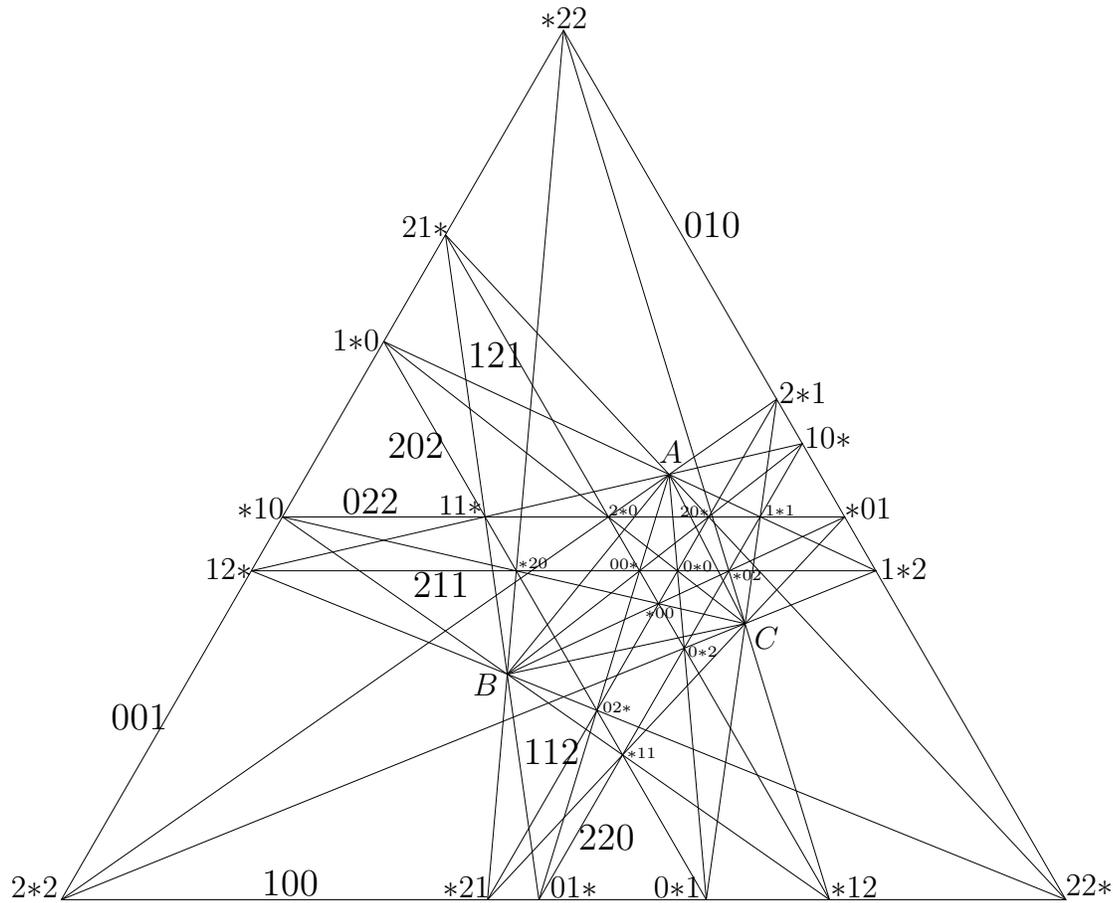
  
\setlength{\unitlength}{3.8pt}

\caption{The eighteen Morley triangles.}
\label{18tri}
\end{figure}

Figure \ref{18tri} shows the complete configuration of
18 trisectors, 27 Morley points, 9 Morley
lines, 18 Morley triangles and 9 GF-triangles.
The figure is labelled according to Conway's scheme,
in which the 27 Morley points each receive a three-digit
label comprising two numbers 0, 1 or 2 and a star.
The (1st, 2nd, 3rd) position of the star indicates
the vertex ($A$, $B$, $C$) of the triangle that
is {\em not} responsible for the point.  For
example, 1$\ast$2 is the intersection of beam 1 from
a lighthouse at $A$ with beam 2 from a lighthouse
at $C$.  Remember that the beams from the
$(A,C)$-lighthouse pair are numbered starting
from 0 for the internal trisector proximal to
the edge $CA$ and continuing across $CA$ with
the beams counted in opposite senses.

\newpage

Figure \ref{scheme} is a schematic diagram of
the 9 Morley lines and 18 Morley points.
The six points on any Morley line come two
from each of the three pairs of lighthouses,
with two point-labels having a star in the first
position, two having it in the second, and two
with it in the third.  On any Morley line
only two of the digits 0,1,2 occur in any one
of the three positions.  The line-label
(large in Figure \ref{18tri}) consists of the three
digits which do {\em not} occur in the three
positions.  Each line-label has a repeated digit
and a different digit, and the sum of the three
digits is congruent to 1 modulo 3.  For example,
the points labelled $\ast00,\ \ast12,\ 0\!\ast\!2,\
2\!\ast\!0,\ 00\ast,\ 21\ast$, lie on the Morley
line with label 121.

\begin{figure}[h]
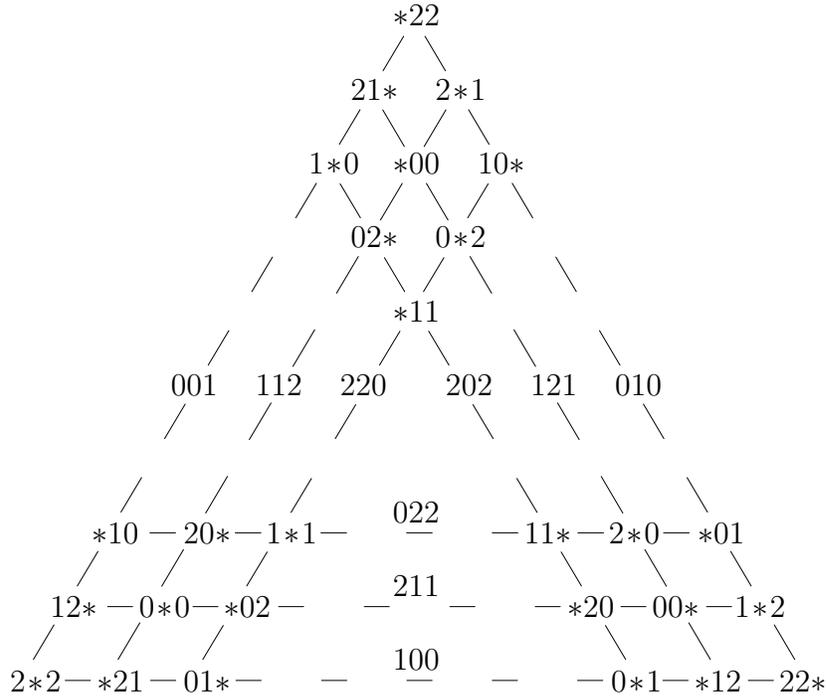
  
\setlength{\unitlength}{4pt}

\caption{Schematic diagram of Morley points, lines and triangles,
and GF-triangles.}
\label{scheme}
\end{figure}

The GF-triangles are found by choosing any
vertex, say 0$\ast$1 in the bottom row of
Figure \ref{scheme}, and then taking the only other two
points on the Morley lines through that vertex
that have stars in the same position (here
the $B$ position): 2$\ast$2 and 1$\ast$0.  The
other 18 triangles are genuine Morley
triangles whose vertices have their stars in
the three different positions.

\begin{figure}[h]  
\setlength{\unitlength}{4pt}
%
\caption{Nine Morley lines and nine GF-circles.}
\label{9lines}
\end{figure}

\clearpage

Figure \ref{9lines} shows the 9 Morley lines
and 9 GF-circles,
with more detail shown in Figure \ref{large9}.
There are 27 Morley points, 9 points of
concurrence of triads of GF-circles, and 3
vertices of the original triangle.  Each
circle passes through 3+3+2 of these points.
Each of the 9 Morley lines passes through 6
Morley points.  The label for each of the 9
points of concurrence of triads of GF-circles
is a three-digit number that differs from the
label of the corresponding Morley line in
only one digit: the digit that differs from the
other two is changed to the third possibility,
so that the sum of the digits of the label of
an associated point is congruent to 2 modulo 3,
as follows:

\begin{center}

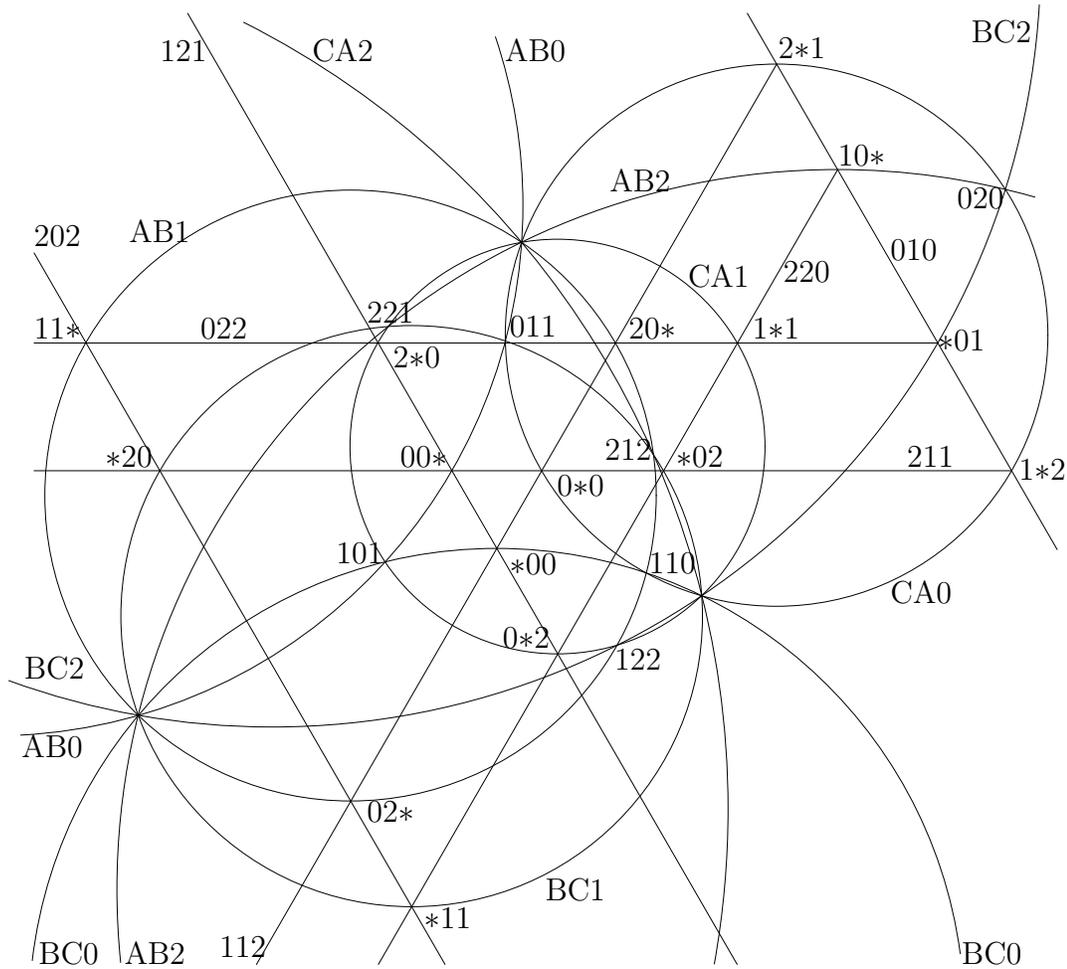
\caption{Enlargement of part of Figure 6.9.}
\label{large9}
\end{figure}

\clearpage

Table \ref{mpoints} facilitates the location of
all the points, lines and circles.  The
circumcircle $ABC$ belongs to each of the
coaxal systems of GF-circles.  The three lines
of centres of the systems are the perpendicular
bisectors of the edges of $ABC$ and they concur
at its circumcentre.  If you omit the Morley lines
from Figure \ref{9lines} you are left with a pleasing
configuration of 9 circles and 12 points
with 5 points on each circle, 6 circles
through each of 3 points, and 3 circles through
each of the remaining 9 points.


\begin{table}[ht]   
\begin{center}
\begin{tabular}{cccc}
GF-circle & Morley points & Morley lines
 & associated points \\[3pt]
BC0 & $\ast$00 $\ast$21 $\ast$12 & 100 121 112
 & 200 101 110 \\
BC1 & $\ast$11 $\ast$02 $\ast$20 & 211 202 220
 & 011 212 221 \\
BC2 & $\ast$22 $\ast$10 $\ast$01 & 022 010 001
 & 122 020 002 \\[3pt]
CA0 & 2$\ast$1 0$\ast$0 1$\ast$2 & 211 010 112
 & 011 020 110 \\
CA1 & 0$\ast$2 1$\ast$1 2$\ast$0 & 022 121 220
 & 122 101 221 \\
CA2 & 1$\ast$0 2$\ast$2 0$\ast$1 & 100 202 001
 & 200 212 002 \\[3pt]
AB0 & 21$\ast$ 12$\ast$ 00$\ast$ & 211 121 001
 & 011 101 002 \\
AB1 & 02$\ast$ 20$\ast$ 11$\ast$ & 022 202 112
 & 122 212 110 \\
AB2 & 10$\ast$ 01$\ast$ 22$\ast$ & 100 010 220
 & 200 020 221 \\
\end{tabular}
\end{center}

\begin{center}
\begin{tabular}{cccc}
Morley & GF-circles & Morley points & associated \\
 line  &            &               &   point \\
211 & BC1 CA0 AB0 &  $\ast$20 $\ast$02 \ 0$\ast$0
 1$\ast$2 \ 12$\ast$ 00$\ast$ & 011 \\
121 & BC0 CA1 AB0 &  $\ast$12 $\ast$00 \ 2$\ast$0
 0$\ast$2 \ 00$\ast$ 21$\ast$ & 101 \\
112 & BC0 CA0 AB1 &  $\ast$00 $\ast$21 \ 2$\ast$1
 0$\ast$0 \ 02$\ast$ 20$\ast$ & 110 \\[3pt]
022 & BC2 CA1 AB1 &  $\ast$10 $\ast$01 \ 1$\ast$1
 2$\ast$0 \ 20$\ast$ 11$\ast$ & 122 \\
202 & BC1 CA2 AB1 &  $\ast$20 $\ast$11 \ 1$\ast$0
 0$\ast$1 \ 11$\ast$ 02$\ast$ & 212 \\
220 & BC1 CA1 AB2 &  $\ast$11 $\ast$02 \ 0$\ast$2
 1$\ast$1 \ 01$\ast$ 10$\ast$ & 221 \\[3pt]
100 & BC0 CA2 AB2 &  $\ast$21 $\ast$12 \ 2$\ast$2
 0$\ast$1 \ 01$\ast$ 22$\ast$ & 200 \\
010 & BC2 CA0 AB2 &  $\ast$01 $\ast$22 \ 2$\ast$1
 1$\ast$2 \ 22$\ast$ 10$\ast$ & 020 \\
001 & BC2 CA2 AB0 &  $\ast$22 $\ast$10 \ 1$\ast$0
 2$\ast$2 \ 12$\ast$ 21$\ast$ & 002
\end{tabular}
\vspace{-20pt}
\end{center}
\caption{Morley points \& lines; GF-circles \& associated points.}
\label{mpoints}
\end{table}

\section{Conway's extraversion.}
The best insight into the
Morley configuration is provided by Conway's
{\bf extraversion}.  An ``A-flip'' of a triangle, for
example, replaces angles $A$, $B$, $C$
with $-A$, $\pi-B$, $\pi-C$, respectively.  It takes
us from THE Morley triangle, 000, with vertices $\ast$00,
0$\ast$0, 00$\ast$, to triangle 011 with vertices
$\ast$11, 0$\ast$1, 01$\ast$.  Repeated use of
A-, B- and C-flips generates the toroidal map of
Figure \ref{18torus}, where similarly
labelled Morley triangles
are to be identified.  There are 9 hexagonal regions,
bounded by $9\times6/2=27$ edges, meeting 3 at each
of the 18 Morley triangles.

\begin{figure}[h]  
\begin{picture}(400,360)(-200,-175)
\setlength{\unitlength}{1pt}
\put(12,0){\bf 000}     
\put(12,50){\bf 011}    
\put(-33,-40){\bf 101}  
\put(57,-40){\bf 110}   
\put(-78,-130){\bf 221} 
\put(57,140){\bf 221}   
\put(102,-130){\bf212}  
\put(-33,140){\bf 212}  
\put(12,-130){\bf 200}  
\put(12,-180){\bf 222}  
\put(57,-90){\bf 102}   
\put(-33,-90){\bf 120}  
\put(147,-40){\bf 122}  
\put(-123,-40){\bf 122} 
\put(102,50){\bf 020}   
\put(102,0){\bf 012}    
\put(147,90){\bf 222}   
\put(57,90){\bf 210}    
\put(-78,50){\bf 002}   
\put(-78,0){\bf 021}    
\put(-33,90){\bf 201}   
\put(-123,90){\bf 222}  
\put(-123,-90){\bf 111} 
\put(147,-90){\bf 111}  
\put(12,180){\bf 111}   

\multiput(-22,102)(90,0){2}{\line(0,1){30}}
\multiput(-67,12)(90,0){3}{\line(0,1){30}}
\multiput(-112,-78)(90,0){4}{\line(0,1){30}}
\put(23,-168){\line(0,1){30}}
\put(-18,151){\line(3,2){36}}
\multiput(-67,61)(90,0){3}{\line(3,2){36}}
\multiput(-108,-29)(90,0){3}{\line(3,2){36}}
\multiput(-63,-119)(90,0){3}{\line(3,2){36}}
\put(64,151){\line(-3,2){36}}
\multiput(-71,61)(90,0){3}{\line(-3,2){36}}
\multiput(-26,-29)(90,0){3}{\line(-3,2){36}}
\multiput(-71,-119)(90,0){3}{\line(-3,2){36}}
\multiput(-20,114)(90,0){2}{$A$}
\multiput(-65,24)(90,0){3}{$A$}
\multiput(-110,-66)(90,0){4}{$A$}
\put(25,-156){$A$}
\put(-10,161){$B$}
\multiput(-55,71)(90,0){3}{$B$}
\multiput(-100,-19)(90,0){3}{$B$}
\multiput(-55,-109)(90,0){3}{$B$}
\put(35,156){$C$}
\multiput(-100,66)(90,0){3}{$C$}
\multiput(-55,-24)(90,0){3}{$C$}
\multiput(-100,-114)(90,0){3}{$C$}
\end{picture}
\caption{The 18 Morley triangles related by extraversion.}
\label{18torus}
\end{figure}
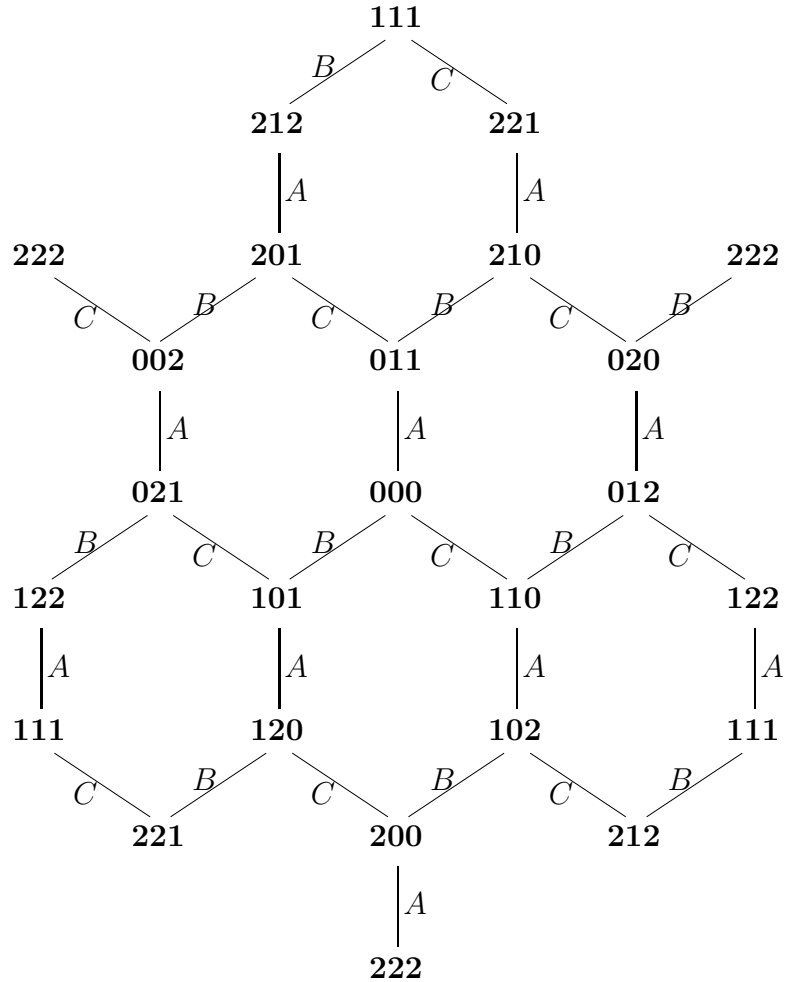

\section{The Morley Group}

The least obvious of the
five abstract groups of order 18 is the semidirect
product of $C_3\times C_3$ with $C_2$.  One rarely
meets such things in the heat of battle, but here
it is: a set of relations is $x^2=y^3=z^3=1$,
$yz=zy$, $yxy=x=zxz$.  One way to realize it is
with paths in Figure \ref{18torus}: $y=AB$,
$z=AC$, $x=ABA$.  It deserves a proper name.

\section{Where are the Morley centres?}
There is some debate as to what constitutes a
``centre'' of a triangle \cite{K}, but for
the equilateral
triangle, one point, and hence at least one
Morley centre (namely 000) is beyond doubt; but
222 and 111 should also be considered.  Label
the centre of a Morley triangle with the three-digit
number that gives the vertices on replacing
the digits in turn with a star.  For example,
120 is the centre of the Morley triangle
with vertices $\ast20,\ 1\!\ast\!0,\ 12\ast$.
In Figure \ref{crate}
the Morley lines are dotted and the 18 Morley
centres form a ``crate'', which, when extended
by the 9 GF-centres (smaller unlabelled dots)
appears to consist of 27 cuboids, 18 of which
have 5 Morley vertices and 3 GF-vertices, and
the other 9 have 6 Morley vertices and 2
GF-vertices.
\begin{figure}[h]
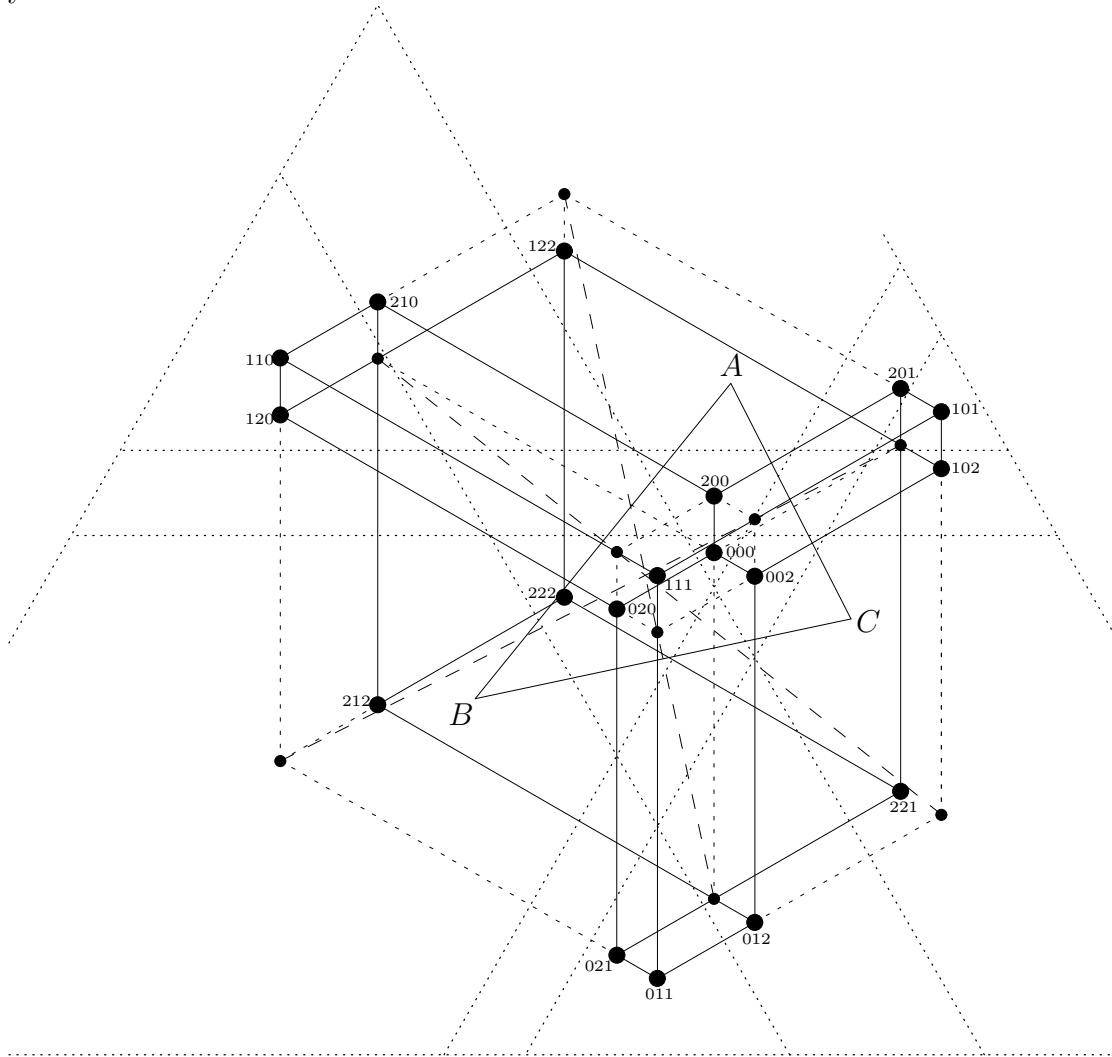
  
\setlength{\unitlength}{6pt}

\caption{Eighteen Morley centres and nine GF-centres.}
\label{crate}
\end{figure}

\clearpage

 There appear to be 108
body-diagonals that bisect each other in fours
at 27 points that form a half-sized ``crate''.
The 9 GF-centres are collinear in threes---the
dashed lines in Figure \ref{crate} add to the
visual confusion.

\noindent
{\bf Paradox 10}.  Figure \ref{crate} is a
drawing of an impossible object.

\section{The Lighthouse Theorem for $n=2$}

For $n=2$, $\binom{n}{2}<n$, and the
ramifications that occur for larger $n$ do not
appear, but even so, the Lighthouse Theorem
still has a lot to tell us.  Paradox 6 no longer
arises, while Paradoxes 7 and 8 are not so
apparent, because we expect a rectangular hyperbola and a circle to
intersect in four points.  In
Figure \ref{altitudes} two pairs of perpendicular
lines through lighthouses at $E$ and $F$ intersect
at the vertices of two 2-gons, say $AH$ and
$BC$, segments that are perpendicular to
each other. In other words:

\smallskip

\noindent
{\bf Theorem.}  {\it The altitudes of a triangle
$(ABC)$ concur $($at the orthocentre} $H$).

\smallskip

\noindent
{\bf Paradox 11.}  Another triangle theorem
without having a triangle.

More symmetrically, each of the 4 points $A$, $B$, $C$,
$H$ is the orthocentre of the triangle formed by
the other three.  This is a special case of the theorem
that any conic through the intersections of two
rectangular hyperbolas is a rectangular hyperbola.

But more: if the two 2-gons meet at $D$, then the Lighthouse
Duplication Theorem tells us that $\angle DEH=\angle HEF$
and $\angle DFH=\angle HFE$ (i.e., $EH$ \& $AC$ and $FH$
\& $AB$ are the pairs of angle-bisectors at vertices
$E$ and $F$ of triangle $DEF$).  The
angle-bisectors of triangle $DEF$ concur in four
points, the incentre $H$ and the excentres $A$, $B$, $C$,
{\em provided} that $DH$ \& $BC$ are the
angle-bisectors of angle $D$.  But
$\angle BDF=\angle BHF$ ($BDHF$ cyclic) = $\angle EHC$
(vertically opposite) = $\angle EDC$ ($EHDC$ cyclic).
\qquad $\blacksquare$

\smallskip

\noindent
{\bf Paradox 12.}  A third triangle theorem without
starting from a triangle.

Note that the theorems are quite general.  Any
triangle $DEF$ is obtained from lighthouses at
$E$ and $F$ whose beams are phased to pass
simultaneously through $D$.

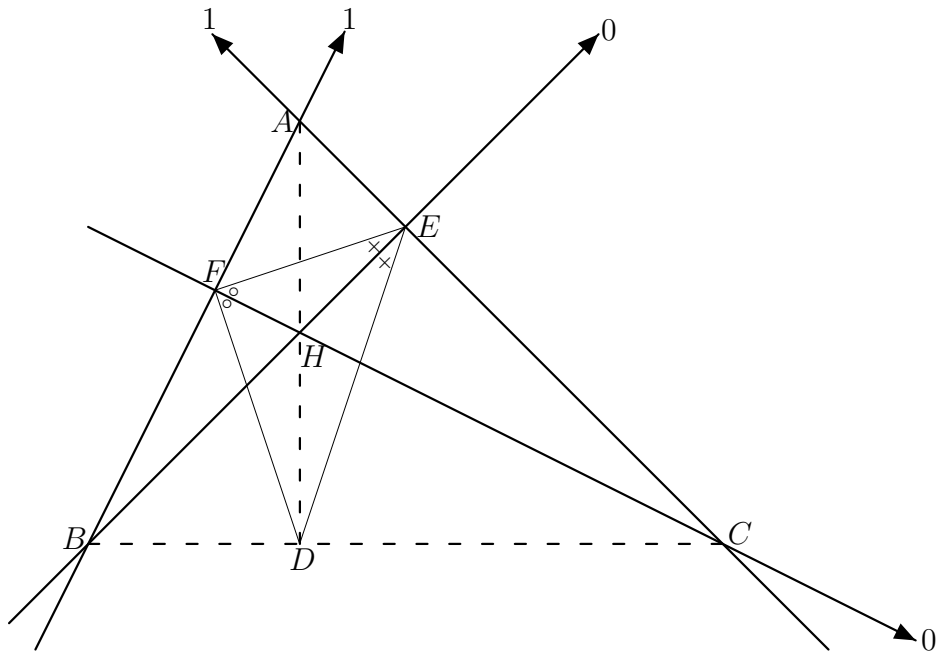
\begin{figure}[h]  
\begin{picture}(320,230)(-80,4)
\drawline(80,40)(120,160)
\drawline(80,40)(48,136)
\drawline(120,160)(48,136)
\thicklines
\drawline(-30,10)(190,230)
\put(190,230){\vector(1,1){3}}
\put(194,231){0}
\drawline(280,0)(50,230)
\put(50,230){\vector(-1,1){3}}
\put(43,235){1}
\drawline(0,160)(310,5)
\put(310,5){\vector(2,-1){3}}
\put(315,0){0}
\drawline(-20,0)(95,230)
\put(95,230){\vector(1,2){2}}
\put(96,235){1}
\dashline{5}(80,40)(80,200)
\dashline{5}(0,40)(240,40)
\put(69,196){$A$}
\put(-10,38){$B$}
\put(242,40){$C$}
\put(76,30){$D$}
\put(124,156){$E$}
\put(43,139){$F$}
\put(80,107){$H$}
\put(49,129){$~^{\circ}$}
\put(101,146){$~^{\times}$}
\put(46.5,124.5){$~^{\circ}$}
\put(105,140){$~^{\times}$}

\end{picture}
\caption{The altitudes of a triangle concur: also the angle-bisectors.}
\label{altitudes}
\end{figure}

Figure \ref{altitudes} contains six cyclic
quadrangles and three rectangular hyperbolas
and the circumcircles of $ABC$, $BHC$, $CHA$, $AHB$
are all congruent, the last three being reflexions
of the first in $BC$, $CA$, $AB$ respectively.

\section{The Thrice Sixteen Theorem}

Figure \ref{thrice16} shows a cyclic quadrangle 0123.  The
incentres of triangles 123, 023, 013, 012 are
the points that are respectively labelled  00,
11, 22, 33.  The excentres are labelled with
the other twelve two-digit base-4 numbers.
The midpoints of the segments joining these
in-(\& ex-)centres are labelled, somewhat
arbitrarily, as follows:

\small

\begin{center}
\begin{tabular}{lc@{\hspace{5pt}}c@{\hspace{5pt}}c
  @{\hspace{5pt}}c@{\hspace{5pt}}c@{\hspace{5pt}}c
  @{\hspace{5pt}}c@{\hspace{5pt}}c@{\hspace{5pt}}c
  @{\hspace{5pt}}c@{\hspace{5pt}}c@{\hspace{5pt}}c}
the point &   04  &   05  &   06  &   14  &   15  &   16
     &   24  &   25  &   26  &   34  &   35  &   36 \\
is mid of & 10-13 & 10-12 & 22-23 & 00-03 & 00-02 & 21-20
     & 01-02 & 11-13 & 00-01 & 11-12 & 01-13 & 02-03 \\
\& mid of & 20-23 & 31-33 & 32-33 & 30-33 & 21-23 & 31-30
     & 31-32 & 30-32 & 10-11 & 21-22 & 20-22 & 12-13
\end{tabular}
\end{center}

\normalsize

The Thrice Sixteen Theorem states that the $4\times4$
incentres of the four triangles 123, 023,
013, 012 lie in fours in a rectangular array on two
perpendicular sets of 4 parallel lines.  Second, the
$4\times6$ midpoints of segments joining these centres
coincide in twelve pairs, one pair (the dots in
Figure \ref{thrice16}) at each end of six diameters of the
{\bf 16-point circle} 0123.  Third, this circle
is the {\bf nine-point circle} for each of the
$4\times4=16$ triangles with vertices
$\{ab,ac,ad\}$, where $a$, $b$, $c$,
$d\in\{0,1,2,3\}$ and $b$, $c$, $d$ are
distinct.
Each of the 16 centres
is the orthocentre of one of these triangles,
implying that the 16 circumcentres are the
reflexions of the 16 incentres in the centre of
the 16-point circle.  Note that the four centres
in any of the four rows or four columns come one each
from the four concyclic triangles, so that the
first digits of the two-digit labels form a latin
square.

\begin{figure}[h]
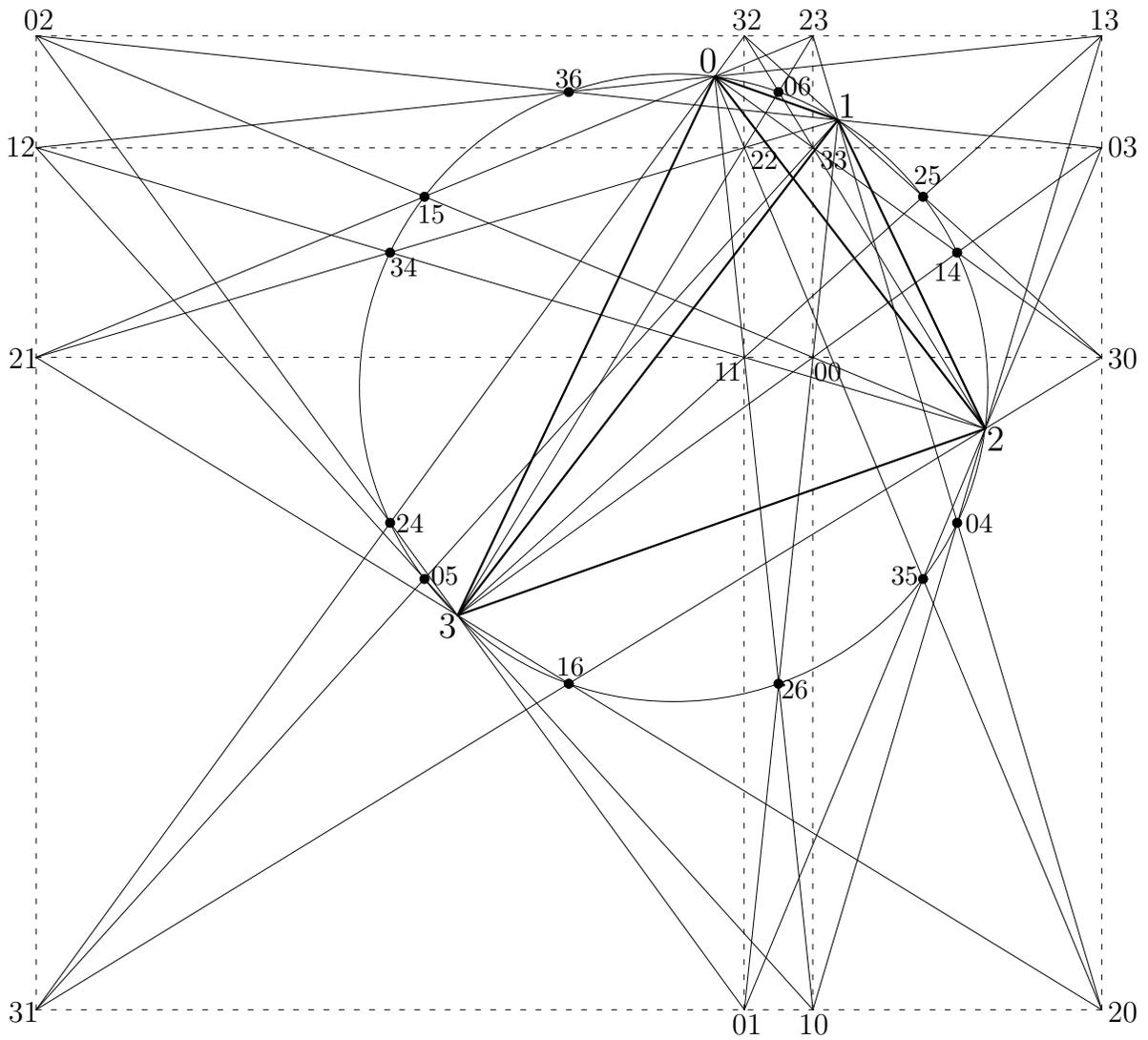
  
\setlength{\unitlength}{1.2pt}

\caption{The Thrice Sixteen Theorem.}
\label{thrice16}
\end{figure}

\clearpage

\noindent
{\bf Exercise for the reader.}  Prove the Thrice Sixteen
Theorem.  Hint: use Euclid III.21 and the Lighthouse
Theorem for $n=2$ with pairs of lighthouses at each of
the six pairs of points (0,1), (2,3), (0,2), (3,1),
(0,3), (1,2) to show that the join 00-11 is
perpendicular to 22-32 and many other perpendicularities.
Use the Lighthouse Lemma to prove that 00-11 is parallel
to 23-32 and many other parallelisms.

\noindent
Prove also that there are twelve sets of six
concyclic points:

\begin{center}
\begin{tabular}{c@{\hspace{2pt}}c@{\quad}|@{\quad}c
 @{\hspace{2pt}}c@{\quad}|@{\quad}c@{\hspace{2pt}}c}
Sets of points & Centre & Sets of points & Centre & Sets of points & Centre \\[3pt]
 \{0,1,33,23,32,22\} & 06 & \{0,2,31,12,33,10\} & 05
 & \{0,3,20,10,23,13\} & 04 \\
 \{1,0,30,20,31,21\} & 16 & \{1,2,03,33,00,30\} & 14
 & \{1,3,23,02,21,00\} & 15 \\
 \{2,0,32,13,30,11\} & 25 & \{2,1,01,31,02,32\} & 24
 & \{2,3,11,00,10,01\} & 26 \\
 \{3,0,21,11,22,12\} & 34 & \{3,1,22,01,20,03\} & 35
 & \{3,2,12,03,13,02\} & 36
\end{tabular}
\end{center}

There is a fourth manifestation of 16:
combine these 12 circles with the 16-point circle
counted with multiplicity four as the circumcircle of
the triangles 123, 230, 301, 012.  In Figure
\ref{thrice16} each of the sextuples
(1,2,3,04,05,06), (2,3,0,14,15,16),
(3,0,1,24,25,26), (0,1,2,34,35,36) lies on the 16-point
circle. In Figure \ref{12circles} the
incentre grid is dashed; the
lighthouse beams are dotted; and the 16-point circle is
drawn with a thick line.

\begin{figure}[h]
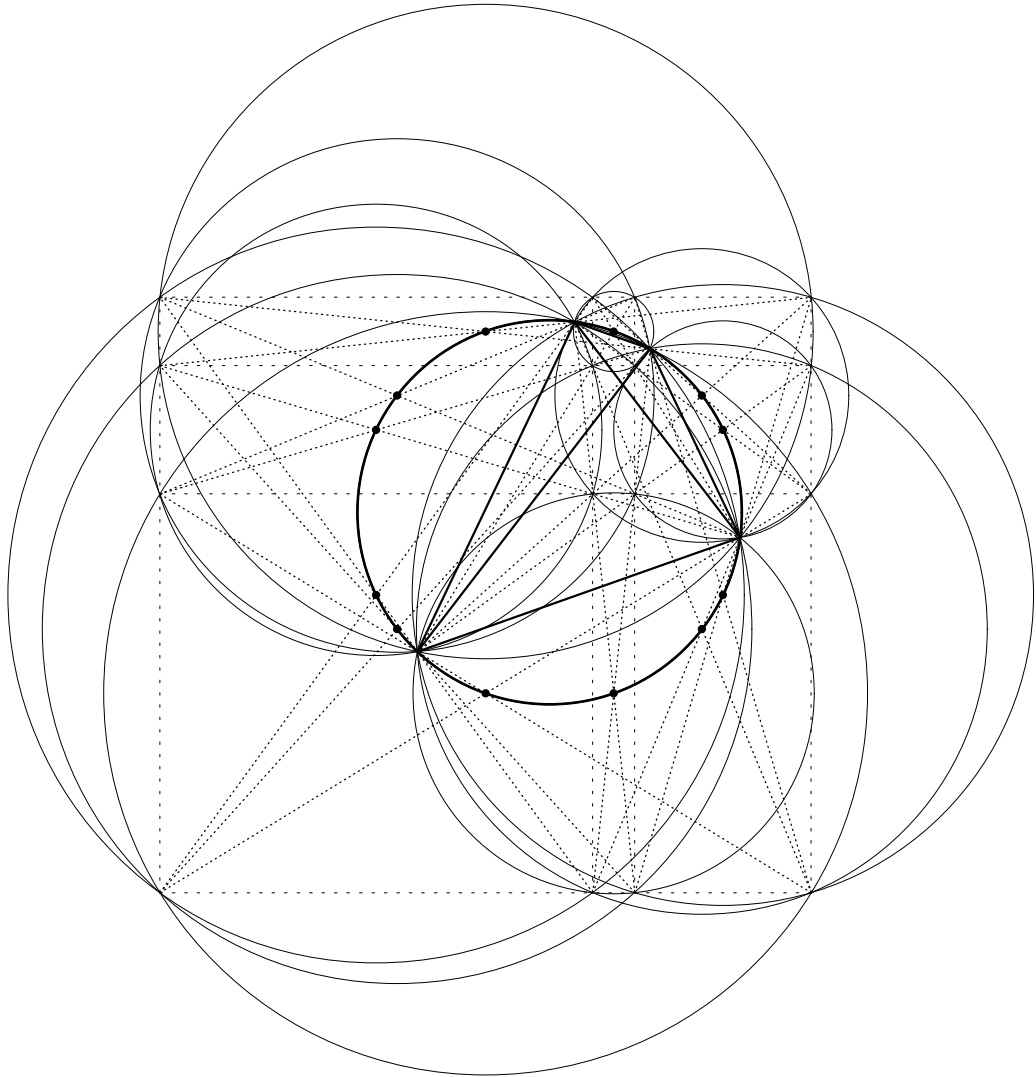
  
\setlength{\unitlength}{0.7pt}
%
\caption{Twelve circles and four times one other.}
\label{12circles}
\end{figure}

\clearpage

Further manifestations of 16 come from the 4 sets of
4 orthocentric points $a0,a1,a2,a3$ ($a=0$,1,2,3).
Each of the 16 incentres $ab$ is the orthocentre of
the triangle $a\bar{b}$, where $\bar{b}$ denotes the
3-element complement of $b$ in the set $\{0,1,2,3\}$.
These 16 triangles have a common nine-point circle:
the 16-point circle, with centre $N$, say.
If $O_{ab}$ is the circumcentre of triangle $a\bar{b}$,
then, since $N$ is the midpoint of the segment
of the Euler line that joins the orthocentre to
the circumcentre, the 16 circumcentres form a
$4\times4$ grid that is congruent to the grid
of orthocentres, being its reflexion in $N$.
Moreover the 16 circumcircles are congruent,
for each of their radii is twice that of the
nine-point radius.
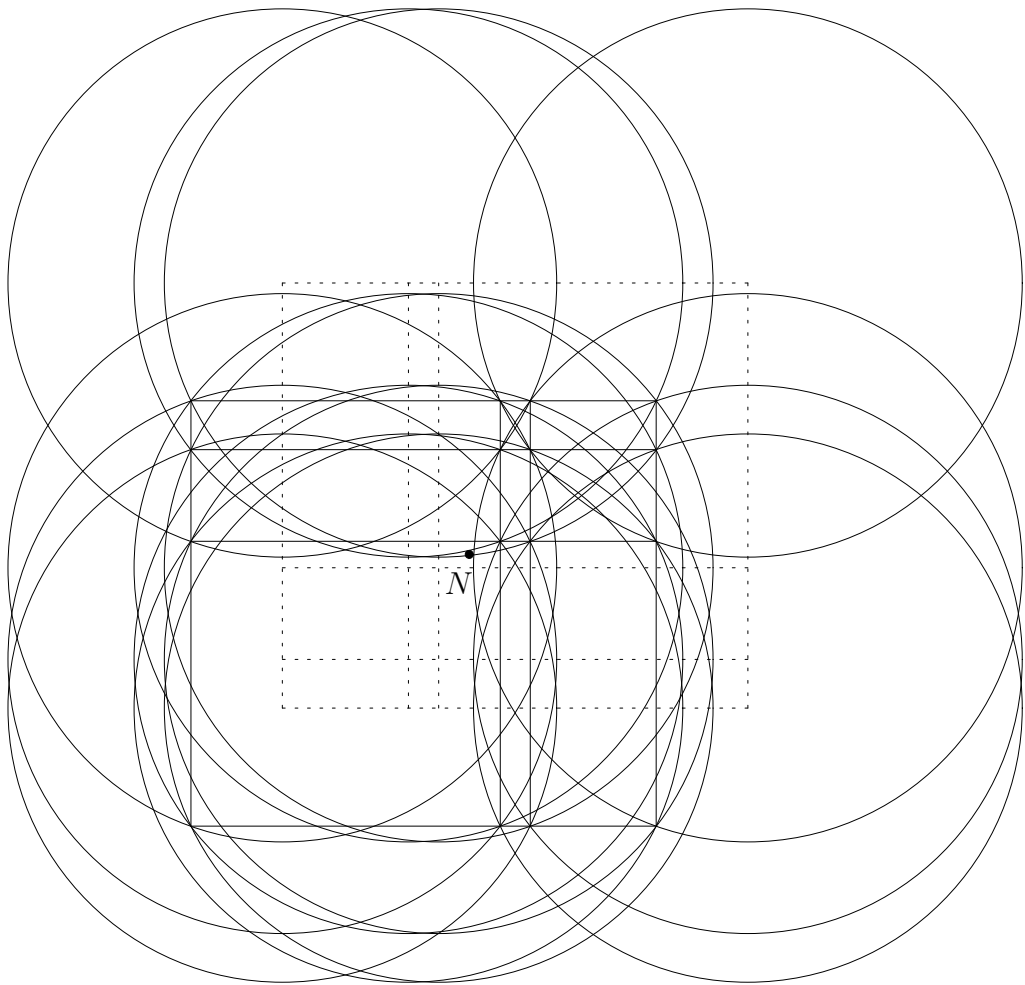
\begin{figure}[h]  
\setlength{\unitlength}{0.5pt}
\begin{picture}(500,780)(-190,-120)
\drawline(0,0)(352,0)
\drawline(0,215.7059)(352,215.7059)
\drawline(0,285.0033)(352,285.0033)
\drawline(0,322)(352,322)
\drawline(0,0)(0,322)
\drawline(233.9467,0)(233.9467,322)
\drawline(256.6193,0)(256.6193,322)
\drawline(352,0)(352,322)
\put(210.6415,205.6773){\circle*{6}}
\put(192,176){$N$}
\dashline{3}(69.283,89.355)(69.283,411.355)  
\dashline{3}(164.664,89.355)(164.664,411.355)
\dashline{3}(187.336,89.355)(187.336,411.355)
\dashline{3}(421.283,89.355)(421.283,411.355)
\dashline{3}(69.283,89.355)(421.283,89.355)  
\dashline{3}(69.283,126.351)(421.283,126.351)
\dashline{3}(69.283,195.649)(421.283,195.649)
\dashline{3}(69.283,411.355)(421.283,411.355)
\put(69.283,89.355){\circle{415.1075}}
\put(69.283,126.351){\circle{415.1075}}
\put(69.283,195.649){\circle{415.1075}}
\put(69.283,411.355){\circle{415.1075}}
\put(164.664,89.355){\circle{415.1075}}
\put(164.664,126.351){\circle{415.1075}}
\put(164.664,195.649){\circle{415.1075}}
\put(164.664,411.355){\circle{415.1075}}
\put(187.336,89.355){\circle{415.1075}}
\put(187.336,126.351){\circle{415.1075}}
\put(187.336,195.649){\circle{415.1075}}
\put(187.336,411.355){\circle{415.1075}}
\put(421.283,89.355){\circle{415.1075}}
\put(421.283,126.351){\circle{415.1075}}
\put(421.283,195.649){\circle{415.1075}}
\put(421.283,411.355){\circle{415.1075}}
\end{picture}
\caption{Sixteen orthocentres, circumcentres and circumcircles.}
\label{16ortho}
\end{figure}

\clearpage

  The 16 centroids form a similar
grid, with one-third the linear dimensions.  This
is not included in Figure \ref{16ortho}, which shows
the orthocentre grid dashed and the circumcentre
grid solid, together with the nine-point
centre $N$ and the 16 congruent circumcircles.

\noindent
{\bf Take five!}  (Further exercise for the reader).
Draw a figure with five concyclic points, forming
$\binom{5}{3}=10$ triangles with five $4\times4$
grids of incentres coinciding as 40 pairs.

\section{There are 72 Morley triangles!}

[In fact, with {\bf twinning} there are 144.] Does
THE Morley triangle, together with its twin, form a
Star of David inscribed in the double-deltoid\,?  Not
quite!?  Needs investigation!  Preliminary experiment:
The axes of symmetry of the deltoid pass through
the 50-point (9-point) centre.  The axes of symmetry of
THE Morley triangle are parallel to these, but not
coincident with them.  The angles they make with the edges
of the triangle are $|B-C|/3$, $|C-A|/3$, $|A-B|/3$.

Repeat the proof of Morley's theorem with the
lighthouses at phases $\frac{\pi}{6}-\gamma$ and
$\frac{\pi}{6}-\beta$ in place of $\beta$ and $\gamma$.
This yields 18 Morley triangles for a triangle with
base angles $3(\frac{\pi}{6}-\gamma)=\frac{\pi}{2}-C$
and $\frac{\pi}{2}-B$.  This is the triangle $BHC$,
where $H$ is the orthocentre of $ABC$.  The edges
of these Morley triangles make angles with $BC$
that differ only by multiples of $\frac{\pi}{3}$
from $(\frac{\pi}{6}-\gamma)-(\frac{\pi}{6}-\beta)=
\gamma-\beta$ and so are parallel to the edges of
the Morley triangles of $ABC$.  With 18 more from
each of the triangles $CHA$ and $AHB$ this makes
a grand total of 72. In the rational
case \cite{BGGG} all 72 have rational edges!  Note
that while $BC$, $CA$, $AB$ are rational, $AH$,
$BH$, $CH$ are not: these, as well as the
altitudes and the area, are rational multiples
of $\sqrt{3}$.

\noindent
{\bf Buy three: get two free!}

\smallskip

\noindent
{\bf Theorem.}  {\it If a triangle is inscribed in a
rectangular hyperbola, then its orthocentre also
lies on the hyperbola, at the opposite end of
the diameter of the hyperbola through the fourth point
of intersection of the hyperbola with the
circumcircle of the triangle.}

I prefer synthetic proofs, but this analytic one
is so attractive that I can't resist:

{\it Proof.}  Choose the axes and scale so
that the rectangular hyperbola has
equation $xy=1$ and the triangle has vertices
$(t,1/t)$, $(u,1/u)$, $(v,1/v)$.  The line
through the first two has slope
$-1/tu$.  The line through the third,
with perpendicular slope, has equation
$y-1/v=tu(x-v)$, which intersects the
hyperbola again in $(-1/tuv,-tuv)$.  By
symmetry, the altitudes of the
triangle concur at this point of the hyperbola.
If the equation to the circumcircle of the
triangle is $x^2+y^2+2gx+2fy+c=0$, then the
$x$-coordinates of its points of intersection
with the hyperbola are given by
$x^4+2gx^3+cx^2+2fx+1=0$.  The
product of the roots is 1, so the fourth
point has $x$-coordinate $1/tuv$, at the opposite
end of the diameter of the hyperbola through
the orthocentre of the triangle. \qquad $\blacksquare$


In case there was ever any doubt about Paradox 8, here
is a list of the 27 Morley points.  They lie in threes
on 9 GF-circles that pass through pairs of the
vertices $ABC$.  They also lie in
threes on 9 rectangular hyperbolas that pass
through pairs of $ABC$.  The following table
shows what lies on what.  Rows of three points are
on a GF-circle through the pair of vertices listed at
their head.  Columns of three points are on a
GF-hyperbola through these vertices.

\begin{center}
\begin{tabular}{c@{\hspace{5pt}}c@{\hspace{5pt}}c@{\hspace{20pt}}
c@{\hspace{5pt}}c@{\hspace{5pt}}c@{\hspace{20pt}}
c@{\hspace{5pt}}c@{\hspace{5pt}}c}
\multicolumn{3}{c}{$B$ \quad $C$ \quad~} &
\multicolumn{3}{c}{$C$ \quad $A$ \quad~} &
\multicolumn{3}{c}{$A$ \quad $B$ \quad~} \\
$\ast$00 & $\ast$12 & $\ast$21 & 0$\ast$0 & 1$\ast$2 & 2$\ast$1 & 00$\ast$ & 12$\ast$ & 21$\ast$ \\
$\ast$11 & $\ast$20 & $\ast$02 & 1$\ast$1 & 2$\ast$0 & 0$\ast$2 & 11$\ast$ & 20$\ast$ & 02$\ast$ \\
$\ast$22 & $\ast$01 & $\ast$10 & 2$\ast$2 & 0$\ast$1 & 1$\ast$0 & 22$\ast$ & 01$\ast$ & 10$\ast$
\end{tabular}
\end{center}

In fact each 3-by-3 array may be thought of as an
affine geometry with 9 points.  Two sets of three
parallel ``lines'' are the rows and columns, i.e.,
the GF-circles and the GF-hyperbolas.  The other
two sets of three parallel lines are the broken
diagonals and correspond to sets of points on
beams through $B$ \& $C$, or $C$ \& $A$, or $A$ \& $B$.

\section{So little done --- so much to do!}

Nine rectangular hyperbolas, each through 5 points,
$\binom{5}{3}$ triangles inscribed in each.  Where
are the ninety orthocentres?  Where do the circumcircles
meet the hyperbolas again?  Where are the circumcentres?
The nine-point circles?  The Euler lines?  The ``buy 3,
get 2 free'' theorem gives us a good start.

For simplicity, look only at the three hyperbolas
associated with the lighthouses at $B$ and $C$.
Their centres are all at the midpoint, $M$, of $BC$.
Simplify the labels of the Morley points by
omitting the star and reading the two digits
as a ternary number.  For example, $\ast$21 is 7.
Figure \ref{3hyp} (see also Figure \ref{large3hyp})
depicts the nine points whose star is in the first ($A$)
position as  0, 1, 2, $\ldots$, 8.

Consider just the nine triangles
$BCa$ where $0\leq a\leq8$.  Each is inscribed
in hyperbola 0 or 1 or 2 according as
$a$ belongs to $\{0,8,4\}$ or $\{5,1,6\}$ or $\{7,3,2\}$.
Their nine orthocentres $a^{\prime}$ lie
respectively on these hyperbolas, and the nine
sets $\{B,C,a,a^{\prime}\}$ are each orthocentric:

$$aa^{\prime}\perp BC, \quad aB\perp a^{\prime}C,
  \quad aC\perp a^{\prime}B.$$

\noindent
The circumcircles of the triangles $Caa^{\prime}$
$Baa^{\prime}$, $BCa^{\prime}$ are the reflexions
of the circumcircle of $BCa$ in its respective
edges $Ca$, $aB$, $BC$, and so each has the same
circumradius, twice that of the radius of the
common nine-point circle of the four triangles.
Each of the nine points $a^{\prime\prime}$ is at
the opposite end of the diameter through $a$ of
the appropriate GF-circle $BC057$, $BC813$, $BC462$.
They form regular hexagons with the original
points $a$ inscribed in the GF-circles and are
part of a manifestation of the Lighthouse Theorem
for $n=6$.  Each is also at the opposite end of
the diameter through $a^{\prime}$ of the
appropriate hyperbola $BC084$, $BC516$, $BC732$.
These nine diameters concur at $M$.

\begin{figure}[h]
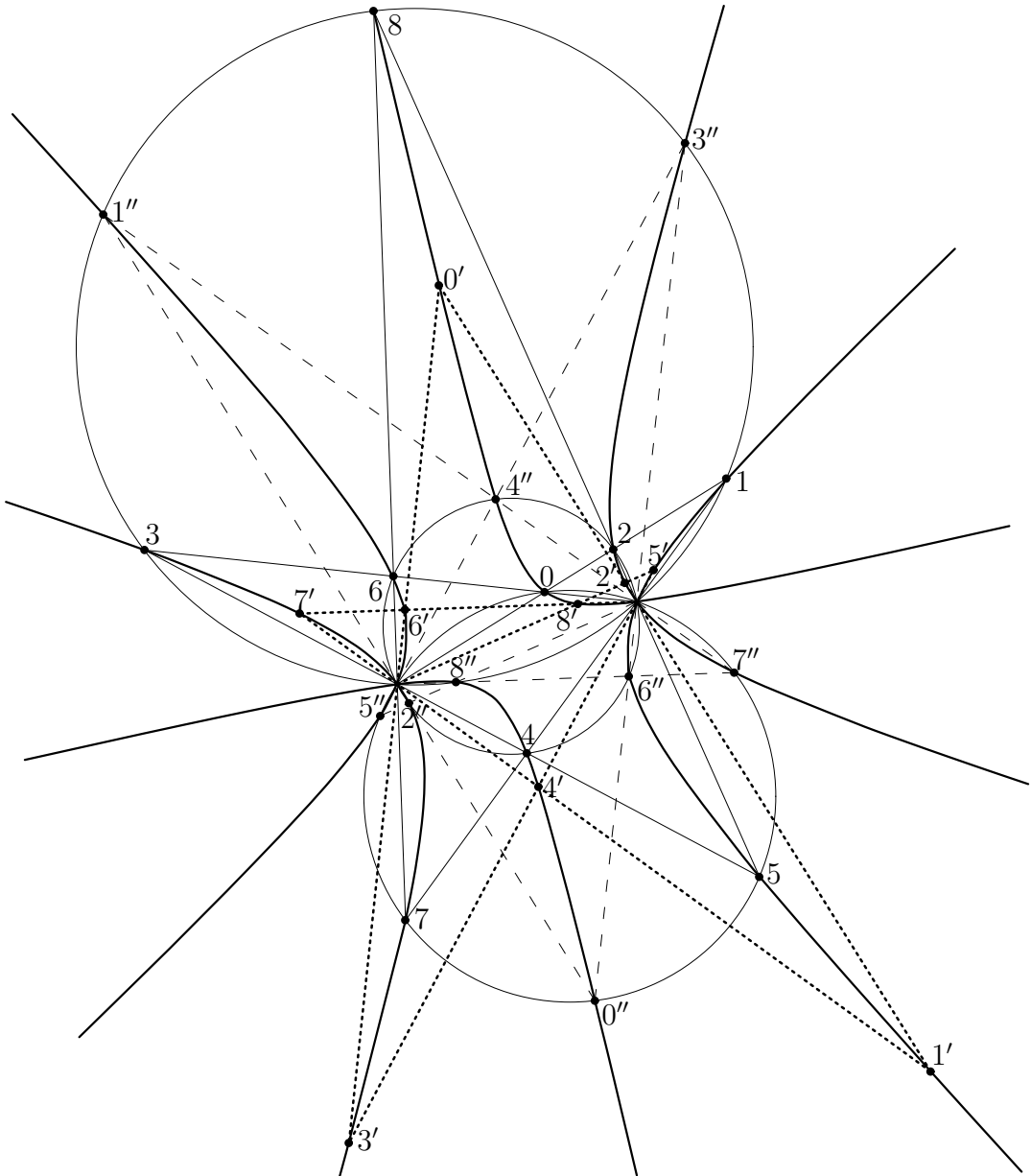
  

\caption{The three GF-hyperbolas through $B$ and $C$.}
\label{3hyp}
\end{figure}

\clearpage

The nine points $a^{\prime\prime\prime}$ are not
shown in Figures \ref{3hyp} or \ref{large3hyp},
though Figure \ref{further}
contains three specimens.  They may be variously
described as the second intersection of
$aa^{\prime}$ with the appropriate GF-circle;
or as the reflexion of $a^{\prime}$ in $BC$
(i.e., as being generated by lighthouses whose
beams are the reflexions in $BC$ of the original
beams through $B$ and $C$); or as the reflexion
of $a^{\prime\prime}$ in the perpendicular
bisector of $BC$.

\begin{figure}[h]
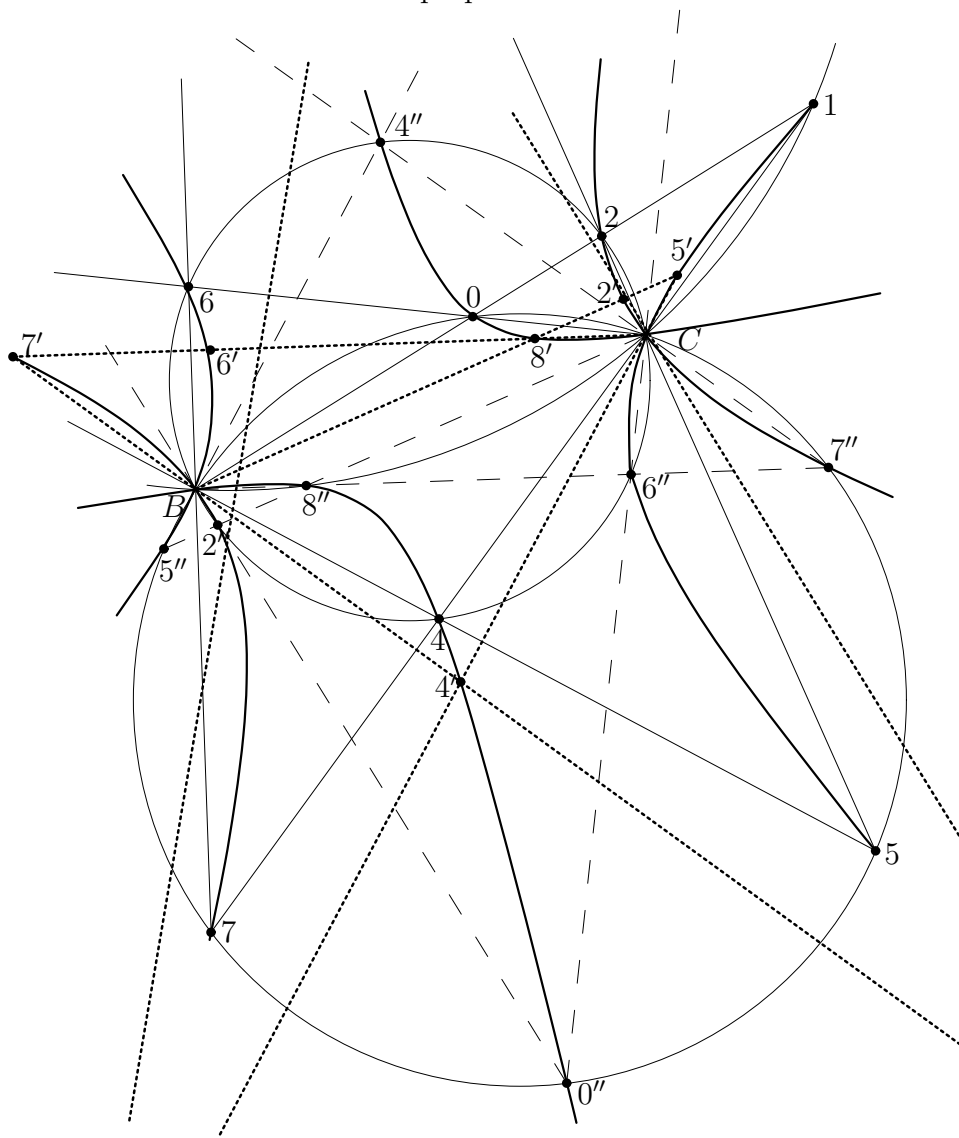
  

\caption{Enlargement of part of Figure \ref{3hyp}.}
\label{large3hyp}
\end{figure}

\clearpage

\section{The Four Nines Theorem.}
{\it The four sets
of nine points, $\{a\}$, $\{a^{\prime}\}$,
$\{a^{\prime\prime}\}$, $\{a^{\prime\prime\prime}\}$
are the congruent configurations}

\begin{center}
\begin{tabular}{ccc@{\qquad}ccc@{\qquad}ccc@{\qquad}ccc}
0 & 1 & 2 & $0^{\prime}$ & $3^{\prime}$ & $6^{\prime}$ &
$0^{\prime\prime}$ & $1^{\prime\prime}$ & $2^{\prime\prime}$ &
$0^{\prime\prime\prime}$ & $3^{\prime\prime\prime}$ &
 $6^{\prime\prime\prime}$ \\
3 & 4 & 5 & $1^{\prime}$ & $4^{\prime}$ & $7^{\prime}$ &
$3^{\prime\prime}$ & $4^{\prime\prime}$ & $5^{\prime\prime}$ &
$1^{\prime\prime\prime}$ & $4^{\prime\prime\prime}$ &
 $7^{\prime\prime\prime}$ \\
6 & 7 & 8 & $2^{\prime}$ & $5^{\prime}$ & $8^{\prime}$ &
$6^{\prime\prime}$ & $7^{\prime\prime}$ & $8^{\prime\prime}$ &
$2^{\prime\prime\prime}$ & $5^{\prime\prime\prime}$ &
 $8^{\prime\prime\prime}$
\end{tabular}
\end{center}

\noindent
{\it whose rows are points on Beams through $B$ and
whose Columns are points on beams through} $C$.

\begin{figure}[h]
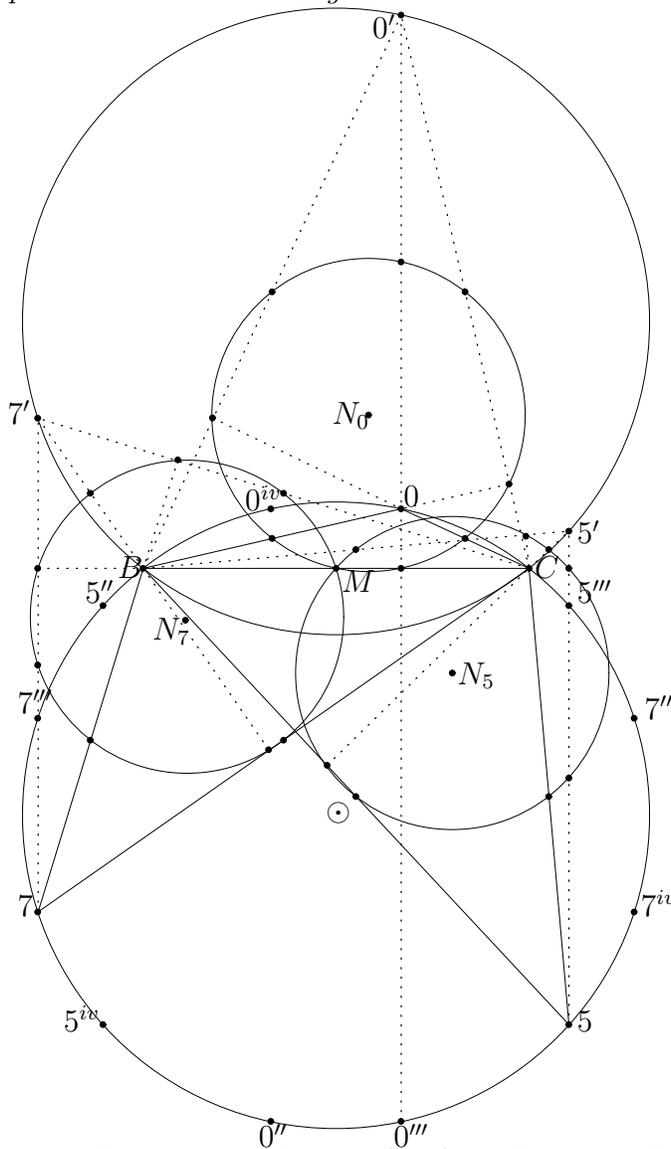
  
%
\caption{Further fun from Figure \ref{3hyp}.}
\label{further}
\end{figure}

\clearpage

In Figure \ref{further} some of the relations
between these configurations can be seen as
\begin{enumerate}
\item The triples of points
$00^{\prime}0^{\prime\prime\prime}$,
$55^{\prime}5^{\prime\prime\prime}$,
$77^{\prime}7^{\prime\prime\prime}$ are
collinear and perpendicular to $BC$, with
$00^{\prime}$, $55^{\prime}$, $77^{\prime}$
equal in length.
\item The joins $00^{\prime\prime}$, $55^{\prime\prime}$,
$77^{\prime\prime}$ are diameters of the
GF-circle 057.
\item Each of $0^{\prime}0^{\prime\prime}$,
$5^{\prime}5^{\prime\prime}$,
$7^{\prime}7^{\prime\prime}$ has $M$
as midpoint.
\item Each of $0^{\prime\prime}0^{\prime\prime\prime}$,
$5^{\prime\prime}5^{\prime\prime\prime}$,
$7^{\prime\prime}7^{\prime\prime\prime}$ is
parallel to $BC$.
\end{enumerate}

Again there are ramifications and more sets
of nine points.  The astute reader would demand
a set $\{a^{iv}\}$ of reflexions of
$\{a^{\prime\prime}\}$ in $BC$ or of
$\{a^{\prime\prime\prime}\}$ in $M$, forming
a four-group of configurations with $\{a\}$,
$\{a^{\prime\prime}\}$, $\{a^{\prime\prime\prime}\}$.
The points $0^{iv}$, $5^{iv}$ and $7^{iv}$
have crept into Figure \ref{further} which confines
itself to illustrating triangles $BC0$, $BC5$
and $BC7$, with common circumcentre $\odot$. The
Euler lines $\odot N_0 0^{\prime}$,
$\odot N_5 5^{\prime}$, $\odot N_7 7^{\prime}$
are not drawn.  The nine-point centres $N_0$,
$N_5$, $N_7$ form an equilateral triangle.
The nine-point circles are congruent, each
passes through $M$, and they intersect at
angles $\pi/3$, forming a pleasing cloverleaf.

\smallskip

\noindent
{\bf Only a tithe.}  I have considered one-third
of the GF-hyperbolas and three-tenths of the
triangles inscribed in each.  The ``buy three,
get two free'' theorem gives us 180 bonus points from
the 90 triangles. What further coincidences,
collinearities, concyclicities are there
among these points?

\section{Bifaux and Skewfaux.}
When you state the
Morley theorem, you must be careful to
specify the intersections of the {\em proximal}
trisectors.  What if you make a mistake and
use the {\em distal} (duplicated) beams?  You
get a new set of GF-triangles whose edges are
the Morley lines of a triangle, $A^{\prime}BC$
say, with base angles $6\beta$ and $6\gamma$,
and the original $A$ is an incentre of
$A^{\prime}BC$.  Or you might get mixed up
and choose one proximal and one distal trisector,
and arrive at the Morley lines for a triangle
$A^{\prime\prime}BC$ or $A^{\prime\prime\prime}BC$,
with just one base angle doubled, and
$A^{\prime\prime}$ and $A^{\prime\prime\prime}$
are the intersections of $BA$ with $CA^{\prime}$
and of $CA$ with $BA^{\prime}$.

\begin{figure}[h]  
\begin{picture}(400,295)(-40,-55)
\setlength{\unitlength}{1pt}
\put(90,-30){\circle*{3}}   
\put(60,60){\circle*{3}}   
\put(0,0){\circle*{3}}   
\put(180,0){\circle*{3}}   
\put(0,240){\circle*{3}}   
\put(108,-36){\circle*{3}}   
\put(45,-45){\circle*{3}}   
\put(360,0){\circle*{3}}   
\put(25.714,77.143){\circle*{3}}   
\put(67.5,67.5){\circle*{3}}   
\put(60,-60){\circle*{3}}   
\put(72,144){\circle*{3}}   
\put(0,180){\circle*{3}}   
\put(60,120){\circle*{3}}   

\drawline(0,0)(360,0)
\drawline(180,0)(25.714,77.143)
\drawline(0,0)(67.5,67.5)
\drawline(90,-30)(0,240)
\drawline(180,0)(45,-45)
\drawline(0,0)(108,-36)
\drawline(0,240)(108,-36)
\drawline(0,240)(45,-45)
\drawline(72,144)(45,-45)
\drawline(108,-36)(0,180)
\drawline(60,-60)(180,0)
\drawline(0,0)(0,240)
\drawline(0,240)(180,0)
\drawline(60,-60)(0,0)
\drawline(360,0)(45,-45)
\drawline(360,0)(0,180)
\drawline(60,-60)(72,144)
\drawline(60,-60)(0,180)
\drawline(360,0)(25.714,77.143)
\drawline(180,0)(0,180)
\drawline(60,120)(60,-60)
\drawline(0,0)(72,144)
\thicklines
\drawline(0,0)(180,0)(60,60)(0,0)
\put(88,-26){$\odot$}
\put(46,55){$A$}
\put(-10,-4){$B$}
\put(179,3){$C$}
\put(2,240){$A^{\prime}$}
\put(104,-47){$B^{\prime}$}
\put(33,-55){$C^{\prime}$}
\put(358,3){$\alpha$}
\put(17,73){$\beta$}
\put(70,70){$\gamma$}
\put(56,-69){$\alpha^{\prime}$}
\put(74,146){$\beta^{\prime}$}
\put(-10,176){$\gamma^{\prime}$}
\put(56,124){$H$}
\put(60,60){\arc{9}{0.75}{2.05}}
\put(60,60){\arc{9}{2.65}{3.95}}
\put(60,60){\arc{9}{5.15}{6.45}}
\put(60,60){\arc{13}{0.75}{2.05}}
\put(60,60){\arc{13}{2.65}{3.95}}
\put(60,60){\arc{13}{5.15}{6.45}}
\put(6,1){$\bullet$}
\put(5,-6){$\bullet$}
\put(1.5,6){$\bullet$}
\put(166,-0.5){$\circ$}
\put(166,-5.5){$\circ$}
\put(168,5){$\circ$}
\end{picture}
\caption{Inside-out Morley's theorem.}
\label{inout}
\end{figure}
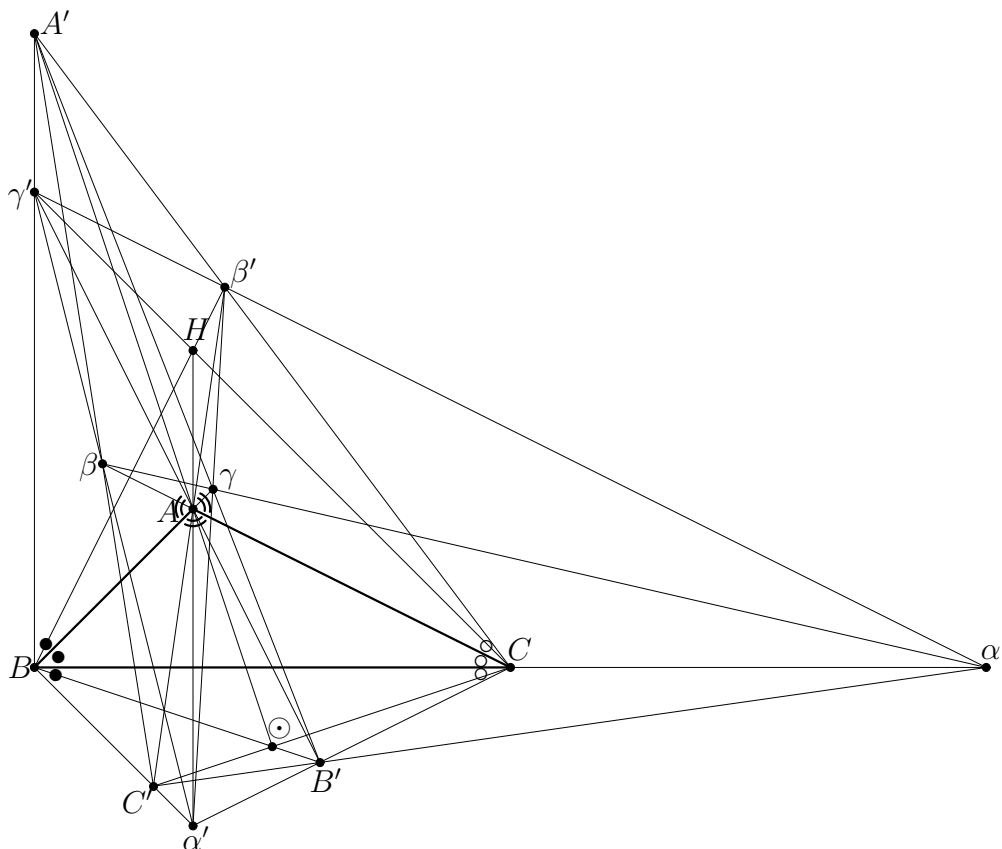

\section{Turning Morley inside-out.}
The construction
of $A^{\prime}$ could be described as the intersection
of the ``distal treblers'' at $B$ and $C$, in
contrast to the {\em proximal trisectors} that are used
in Morley's theorem.  Construct $B^{\prime}$
and $C^{\prime}$ similarly. In Figure \ref{inout}, Arcs
mark angles equal to $A$, Bullets mark angles
equal to $B$, and Circles those equal to $C$.
The (reflex) angle $B'AC'=3A$, angle $C'BA'=3B$
and angle $A'CB'=3C$.  Let the ``proximal
treblers'' meet at $\alpha^{\prime}$,
$\beta^{\prime}$ and $\gamma^{\prime}$, and let
$BC$, $B^{\prime}C^{\prime}$ meet at $\alpha$,
and so forth, as in Figure \ref{inout}.

\noindent
{\bf Theorem.} 
(a) $AA^{\prime}$, $BB^{\prime}$, $CC^{\prime}$
{\it concur at the circumcentre of} $ABC$.

\noindent
(b) $A\alpha^{\prime}$, $B\beta^{\prime}$,
$C\gamma^{\prime}$ {\it concur at the orthocentre of} $ABC$.

\noindent
(c) {\it The triads of points} $(\alpha,\beta,\gamma)$,
$(\alpha,\beta^{\prime},\gamma^{\prime})$,
$(\alpha^{\prime},\beta,\gamma^{\prime})$,
$(\alpha^{\prime},\beta^{\prime},\gamma)$ {\it are
each collinear}.

{\it Sketch of proof.}  (a)  $AB$, $AC$ are angle-bisectors
of $\angle A^{\prime}BC$, $\angle A^{\prime}CB$,
so that $A$ is an incentre of $A^{\prime}BC$ and
$AA^{\prime}$ is a bisector of $\angle BA^{\prime}C$
and makes angles $\frac{\pi}{2}-C$, $\frac{\pi}{2}-B$
with $AB$, $AC$ and hence passes through $\odot$,
the circumcentre of $ABC$.  Similarly for $BB^{\prime}$,
$CC^{\prime}$.

(b)  $\alpha^{\prime}$ is the reflexion of $A$ in $BC$,
so that $A\alpha^{\prime}$ is perpendicular to $BC$.

(c)  $AA'$, $BB'$, $CC'$ concur at $\odot$, making
it the perspector of triangles $ABC$,
$A^{\prime}B^{\prime}C^{\prime}$, and Desargues
tells us that $(\alpha,\beta,\gamma)$ are collinear.
Similarly $\odot$ is the perspector of triangles
$A^{\prime}BC$, $AB^{\prime}C^{\prime}$, and of
two other pairs. \quad $\blacksquare$

\section{Desargues distended.}
The last part of
this proof reminds us that we may swap a pair of
vertices of two triangles in Desargues's theorem
and produce a new perspectrix.  Since
any of the 10 points may serve as perspector,
there are 40 perspectrices and a configuration
of 25 points and 55 lines.  Nine
lines through each of the original 10 points
and six through each of 15 new points.
There are four points on each of 15 new lines
and three points on each of the 10 old lines
and on each of the 30 new perspectrices.

Relabel Figure \ref{inout} with Alex Fink's beautiful
labelling, as in Figure \ref{distend}:

\begin{center}
\begin{tabular}{lccccccccccccc}
Fig.\ref{inout} labels & $\odot$ & $A$ & $B$ & $C$
  & $A^{\prime}$ & $B^{\prime}$ & $C^{\prime}$
  & $\alpha$ & $\beta$ & $\gamma$
  & $\alpha^{\prime}$ & $\beta^{\prime}$ & $\gamma^{\prime}$ \\
Alex's labels &    0    &  8  &  7  &  1
  &       9      &      3       &      2
  &    4     &     5   &    6
  &        04         &       05         &        60 \\
\end{tabular}
\end{center}

\noindent
and see the duality between the Desargues
configuration and the Petersen graph
(Figure \ref{petersen}), each of which has
automorphism group $S_5$.

\begin{figure}[h] 
\begin{picture}(125,125)(-70,5)
\setlength{\unitlength}{4pt}
\put(0,22.27){\circle*{0.4}}   
\put(-10,0){\circle*{0.4}}   
\put(8.09,16.39){\circle*{0.4}}   
\put(-16.18,19.02){\circle*{0.4}}   
\put(5,6.88){\circle*{0.4}}   
\put(0,30.777){\circle*{0.4}}   
\put(-5,6.88){\circle*{0.4}}   
\put(16.18,19.02){\circle*{0.4}}   
\put(-8.09,16.39){\circle*{0.4}}   
\put(10,0){\circle*{0.4}}   
\thicklines
\drawline(0,22.27)(5,6.88)(-8.09,16.39)(8.09,16.39)(-5,6.88)(0,22.27)
\drawline(-10,0)(-16.18,19.02)(0,30.777)(16.18,19.02)(10,0)(-10,0)
\drawline(0,22.27)(0,30.777)
\drawline(-10,0)(-5,6.88)
\drawline(8.09,16.39)(16.18,19.02)
\drawline(-16.18,19.02)(-8.09,16.39)
\drawline(5,6.88)(10,0)
\put(-1.5,22){0}
\put(-12,-1){1}
\put(7,17){2}
\put(-18,18){3}
\put(5.5,6.5){4}
\put(-0.8,31.2){5}
\put(-6.7,6.5){6}
\put(17,18){7}
\put(-8,17){8}
\put(11,-1){9}
\put(22,15){\parbox{300pt}{
\small
\begin{tabular}{|l@{\hspace{5pt}}c
  @{\hspace{9pt}}c@{\hspace{9pt}}c
  @{\hspace{9pt}}c@{\hspace{9pt}}c|}\hline
 (a)  &   6  &   8  &   0  &   2  &   4  \\
 (b)  &  012 &  234 &  456 &  678 &  890 \\
 (c)  &   1  &   3  &   5  &   7  &   9  \\
 (d)  &  369 &  581 &  703 &  925 &  147 \\
 (e)  & 9023 & 1245 & 3467 & 5689 & 7801 \\
 (f) & 61 57 48 & 83 79 60 & 05 91 82 & 27 13 04 & 49 35 26 \\ \hline
\end{tabular}
\normalsize}}
\end{picture}
\caption{The Petersen graph and Desargues-Petersen labels.}
\label{petersen}
\end{figure}
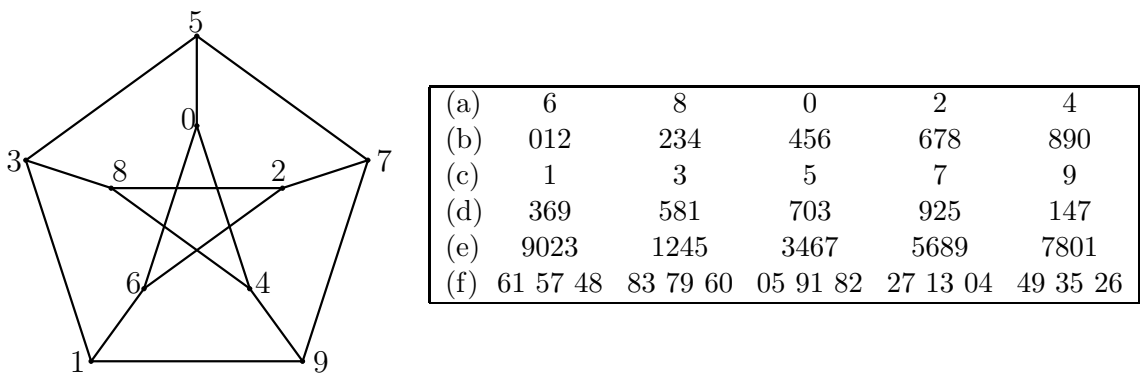

In the table, rows (a) \& (c) are perspectors
in the Desargues configuration, and (b) \& (d)
are the corresponding perspectrices.  Vertices (b)
\& (d) in the Petersen graph are adjacent to the
corresponding vertices in rows (a) \& (c).
Row (e) shows the maximal independent sets of the
Petersen graph; they comprise the $5\cdot\binom{4}{2}$
pairs of points which have a common Desargues line;
the edges of the five complete quadrangles in the
Desargues configuration.  The other
$\binom{10}{2}-30=15$ pairs
of points, which are independent in the Desargues
configuration, correspond to the 15 edges of the
Petersen graph.  These are given as five triples
of pairs in row (f) of the table.  The first pair
of each triple is a copy of the entries in (a) \& (c)
and, in the Petersen graph, is the edge perpendicular
to the 2nd \& 3rd pairs.  These pairs can also be
read as two-digit labels of the 15 points that
amplify the Desargues configuration.  For example,
read the middle triple as ``the joins of 9 \& 1 and
of 8 \& 2 meet in the point whose two-digit label is
05'' etc.

\begin{figure}[h]
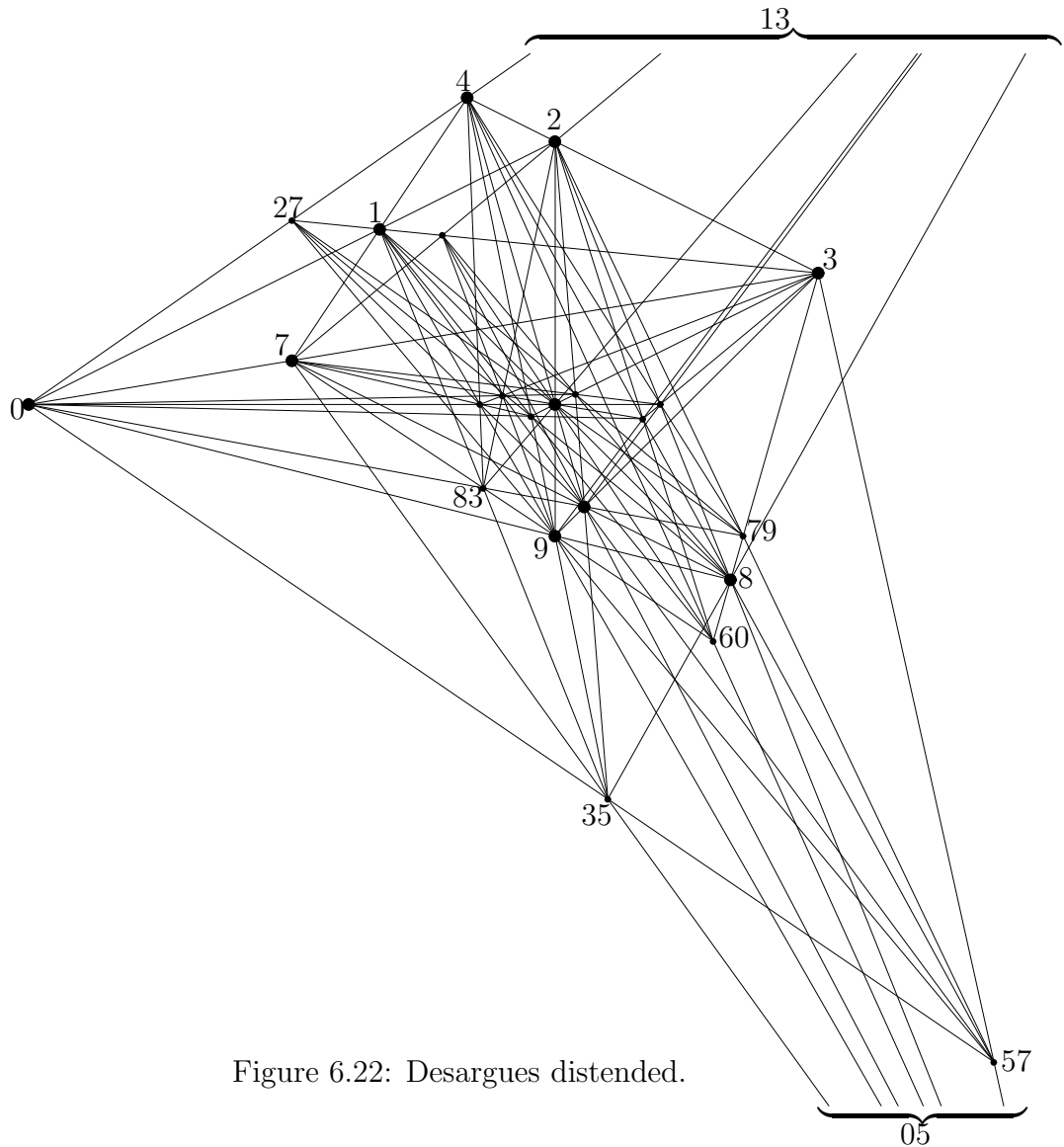
  

\caption{Desargues distended. \qquad\qquad\qquad}
\label{distend}
\end{figure}

\vspace*{-8pt}

\noindent
With perspector 0, triangles 871 \& 932 have perspectrix   4\, \ 5\, \ 6 \\
\rule{97pt}{0pt} triangles 971 \& 832 have perspectrix   4 05 60 \\
\rule{97pt}{0pt} triangles 831 \& 972 have perspectrix   04 5 60 \\
\rule{97pt}{0pt} triangles 872 \& 931 have perspectrix   04 05 6 \\
which is a description of the relabelled Figure \ref{inout},
shown now as Figure \ref{distend}

\clearpage

\section{The Lighthouse Theorem when $n=4$.}

When $n>3$, think of the $n$-gons as
complete graphs on $n$ vertices, two-dimensional
representations of regular $(n\!-\!1)$-dimensional
simplexes.  There are $n$ sets of $\binom{n}{2}$
parallel lines that intersect in $\binom{n}{2}^3$
points, though I now count points with
multiplicity.  The $n^2$ vertices are each
counted $\binom{n\!-\!1}{2}$ times.  If $n$
is even the ``diameters'' concur, and for certain
$n$ there are other concurrencies, originally
enumerated by G.~Bol \cite{B} and often
rediscovered (see \cite{R1} for a good account).

\begin{figure}[h]
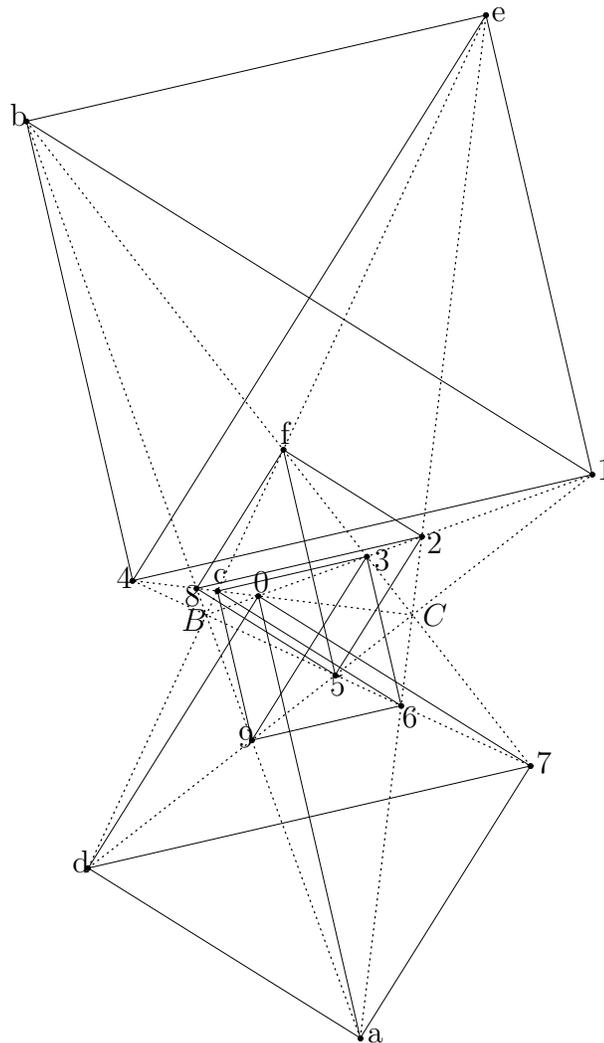
  

\caption{The four 4-gons formed from lighthouses with $n=4$.}
\label{fourlite}
\end{figure}

\clearpage

Figures \ref{fourlite} and \ref{dup4}
are for $n=4$.  In Figure \ref{fourlite}
the beams are dotted and the quadrangles
solid.  In Figure \ref{dup4} the original beams are
omitted, the edges of the quadrangles are
dotted, and the ``duplicated'' beams are solid.
I number the beams as before and will label
the intersections of the lighthouse beams
with the hexadecimal numbers
0, 1, $\ldots$, 9, a, b, $\ldots$, f=15
which should be thought of as two-digit
quaternary numbers $\beth\gimel$, where
$\beth$ and $\gimel$ are the numbers of
the beams from $B$ and $C$.  Quadrangle $q$
$(0\leq q\leq3)$ then has vertices whose
quaternary digits add to $q$ modulo 4:

\smallskip

\centerline{$q=0$ is 00=0 13=7 22=a 31=d;\qquad
$q=1$ is 01=1 10=4 23=b 32=e;}
\centerline{$q=2$ is 02=2 11=5 20=8 33=f;\qquad
$q=3$ is 03=3 12=6 21=9 30=c.}

\begin{figure}[h]
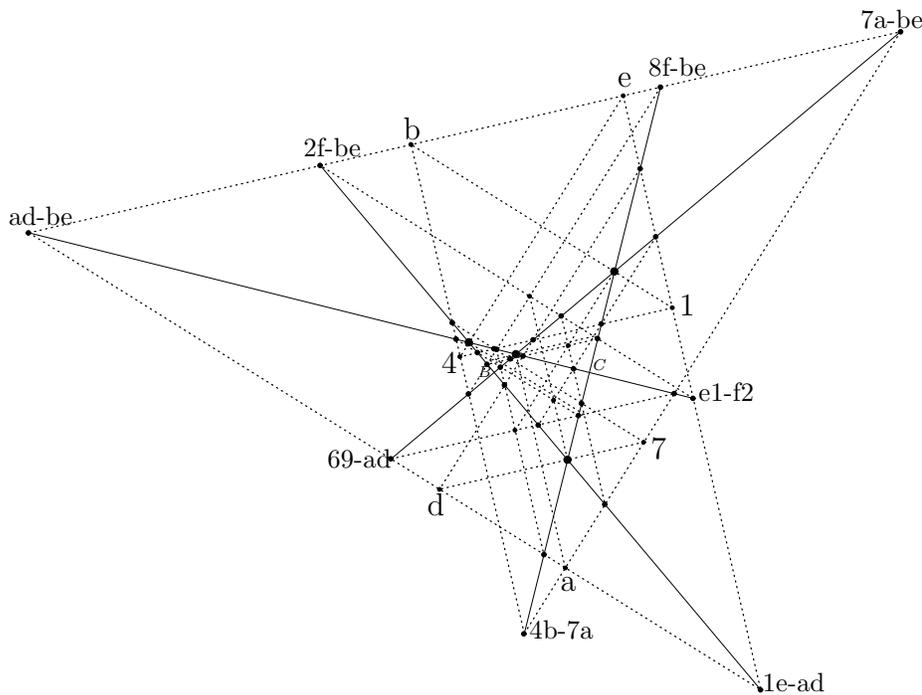
  

\caption{The duplicate beams through edge-intersections
when $n=4$.}
\label{dup4}
\end{figure}

Edges of a quadrangle are denoted by concatenating
the point-labels in these cyclic orders,
starting at the ``earliest'' end.  For instance,
if $q=2$ the edges of ``length'' 1 are 25, 58, 8f
and f2, and the two of ``length'' 2 (the
``diagonals'') are 28 and 5f.  The Lighthouse
Duplication Theorem tells us that the following
columns of ten intersections lie on beams
through $B$ with phases $2\beta$ and
$2\beta+\frac{\pi}{2}$ and through $C$ with
phases $2\gamma$ and $2\gamma+\frac{\pi}{2}$
(hexadecimal numbers are converted
back to base-4 numbers, better to show the pattern).

\begin{table}[h] 
\begin{center}
\begin{tabular}{cccc}
$2\beta$ & $2\beta+\frac{\pi}{2}$
   & $2\gamma$ & $2\gamma+\frac{\pi}{2}$ \\[3pt]
00.13--30.03 & 00.13--10.23 & 00.13--01.10 & 00.13--03.12 \\
01.10--31.00 & 01.10--11.20 & 02.11--03.12 & 01.10--02.11 \\
02.11--32.01 & 02.11--12.21 & 10.23--11.20 & 10.23--13.22 \\
03.12--33.02 & 03.12--13.22 & 12.21--13.22 & 11.20--12.21 \\[2pt]
10.23--20.33 & 20.33--30.03 & 20.33--21.30 & 20.33--23.32 \\
11.20--21.30 & 21.30--31.00 & 22.31--23.32 & 21.30--22.31 \\
12.21--22.31 & 22.31--32.01 & 30.03--31.00 & 30.03--33.02 \\
13.22--23.32 & 23.32--33.02 & 32.01--33.02 & 31.00--32.01 \\[3pt]
00.22--02.20 & 30.12--32.10 & 00.22--02.20 & 03.21--01.23 \\
03.21--01.23 & 11.33--13.31 & 30.12--32.10 & 11.33--13.31 \\[-19pt]
\end{tabular}
\end{center}
\caption{10 points on each of 4 duplicated beams.}
\label{10on4}
\end{table}

The first eight rows indicate intersections of
sides of quadrangles, each  of which meets a pair
of adjacent sides of the two adjacent
quadrangles.  Here ``side'' means the join of two
adjacent vertices, and ``adjacent'' means
``neighboring'' in the cyclic order 0123.
The last two rows contain (duplicates of) the
four intersections of diagonals with the
opposite diagonal of the opposite quadrangle:
the set of four orthocentric points 0a-28,
39-1b, c6-e4, 5f-7d, (the larger dots in Figure
\ref{dup4}).  These constitute an imbedding
of the Lighthouse Theorem for $n=2$. Compare
the labelling of
Figures \ref{altitudes} and \ref{dup4}:

\begin{center}
\begin{tabular}{lccccccc}
Fig.\ref{altitudes} labels & $E$ & $F$
 &  $H$  &  $A$  &  $B$  &  $C$  & $D$ \\
Fig.\ref{dup4}   labels    & $B$ & $C$
 & 0a-28 & 5f-7d & 39-1b & c6-e4 & --- \\
Relabel as                 & $B$ & $C$
 &  $I$  & $I_A$ & $I_B$ & $I_C$ & $A$
\end{tabular}
\end{center}

Figure \ref{largedup4} does not show the new
label $A$, which is for the point of
intersection of the line joining
5f-7d to 0a-28 with the line joining
39-1b to c6-e4.  By the Lighthouse Duplication
Theorem these two lines are perpendicular
and are the angle-bisectors of
$\angle BAC$ and the beams $BA$ and $CA$ have phases
$4\beta$ and $4\gamma$.

The original beams are angle-bisectors of the four
triangles $BCI$, where $I$ is written collectively
for the four incentres of $ABC$.  The original 16 points
are the incentres of these triangles, given as
follows in the order: incentre, excentre
opposite $I$, excentre opposite $B$, excentre
opposite $C$:

\begin{center}
\begin{tabular}{llcccc}
Incentres of triangle  $BCI$  & ($I$ \ = 0a-28) & 0 & a & 2 & 8 \\
Incentres of triangle $BCI_A$ & ($I_A$ = 5f-7d) & 5 & f & 7 & d \\
Incentres of triangle $BCI_B$ & ($I_B$ = 39-1b) & 3 & 9 & 1 & b \\
Incentres of triangle $BCI_C$ & ($I_C$ = c6-e4) & c & 6 & e & 4
\end{tabular}
\end{center}

\noindent
These are, of course, the 16 points of Figure
\ref{fourlite}, similarly labelled.

\begin{figure}[h]
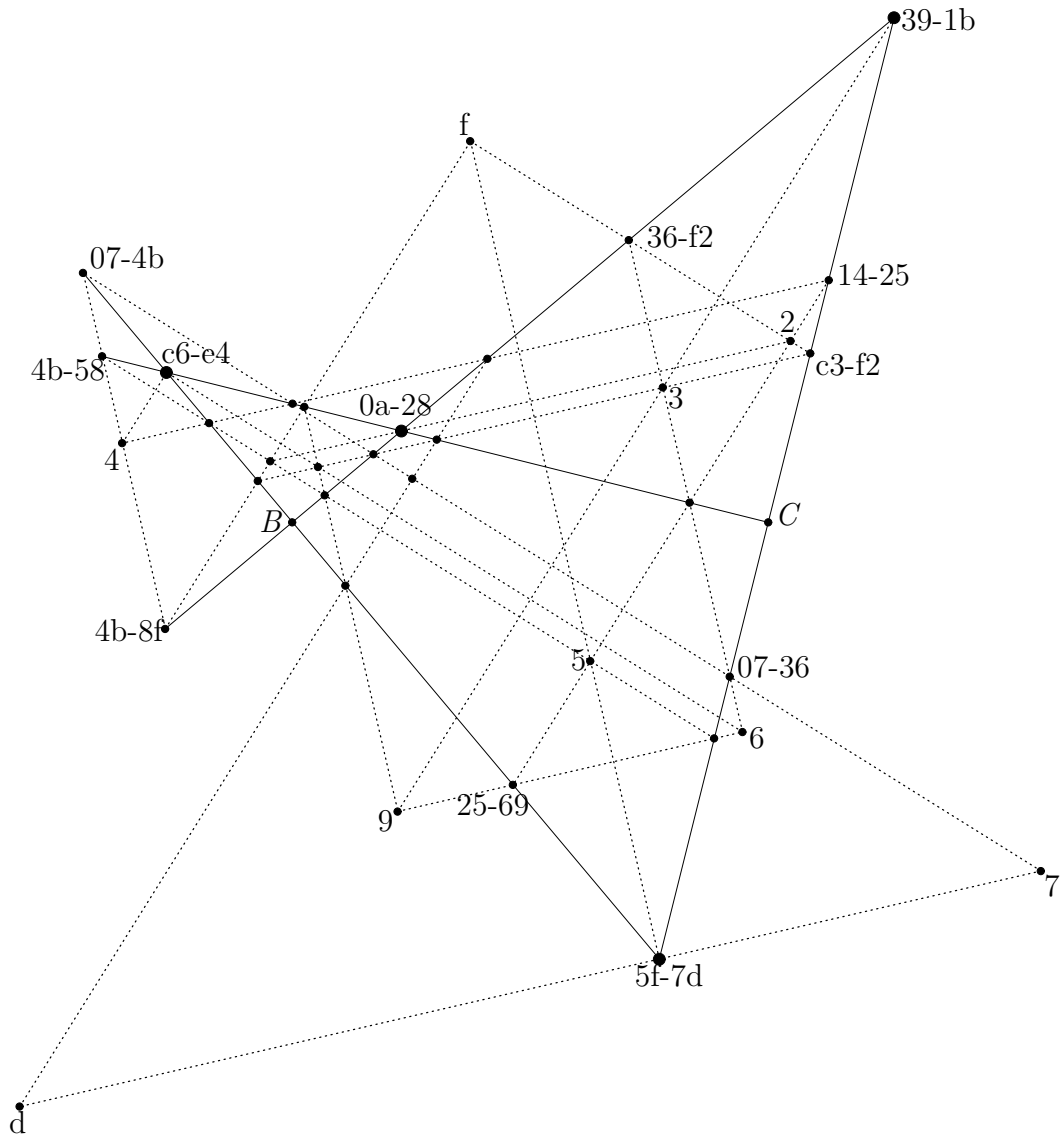
  
%
\caption{Enlargement of part of Figure 24.}
\label{largedup4}
\end{figure}

\clearpage

\section{There's a Malfatti miracle, too!}

In an earlier draft of this paper I asked:

\begin{quote}
``Is there any significance in the fact that
W.~E.~Philip's proof of Morley's theorem
\cite{T+M}, \cite{J} and Steiner's construction
for the Malfatti circles both start with
the incircle whose centre is lit by beams
of phases $B/n(=\beta)$ and $C/n(=\gamma)$ ?''
\end{quote}

Indeed there is!  Malfatti \cite{Mal} asked
for three circular cylinders of maximum
total volume to be cut from a triangular prism.
The complete solution to his problem, begun
by Lob \& Richmond \cite{LR}, was obtained
only comparatively recently; see
\cite[Vol.2, pp.\,245--247]{Ev}, \cite{Gol1},
\cite{Gol2}, \cite{GL}, \cite[pp.\,145--147]{Og},
\cite{Rog} and Guggenheimer's
review of Zalgaller \& Los$'$ \cite{ZL}.

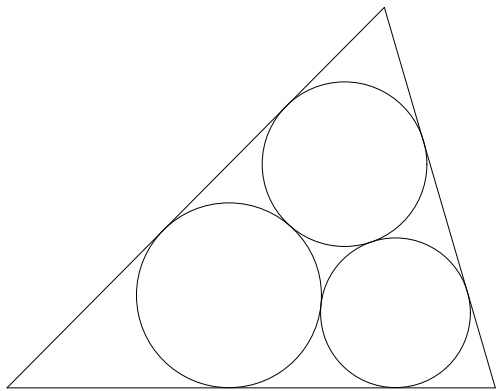
\begin{figure}[h]  
\begin{picture}(160,130)(-265,4)
\setlength{\unitlength}{1.2pt}
\put(-12.5714,70.5){\circle{51.8571}}
\put(-49,29.1667){\circle{58.3333}}
\put(3.5,23.625){\circle{47.25}}
\drawline(0,120)(35,0)(-119,0)(0,120)
\end{picture}
\caption{The popular version of the Malfatti problem.}
\label{popmalf}
\end{figure}

\noindent
{\bf Paradox 13.}  Malfatti misconstrued his
own problem.

This is lucky, because his interpretation
is a much more interesting problem. He
thought that the solution was given by
three cylinders, each touching the other
two, and each touching two faces of the prism,
with a cross-section like Figure \ref{popmalf}.

\noindent
{\bf Paradox 14.}  This problem is not
due to Malfatti.

It occurs in Japanese temple geometry
\cite[pp.\,28, 103--106]{FP}, \cite{Hir},
\cite{N2}, \cite{Rot}, where it is often
ascribed to Ajima.  The isosceles triangle case
was considered by Jacques Bernoulli \cite{Ber};
see Dieudonn\'e's review of Fiocca \cite{Fio}.
But it was posed 417 years before Malfatti, and in
his own country; see the review of Simi
\& Rigatelli \cite{STR}.

\noindent
{\bf Paradox 15.}  Almost all writers state
that Steiner gave no proof for his solution
to the Malfatti problem, but he did give one!

Steiner's proof is scattered over several
sections of two separate articles
\cite{Ste1}, \cite{Ste2}.  It is indirect, via the
more general problem in which the edges of
the triangle are replaced by three circles.
It refers to figures which are not in the
text, but are tucked away at the end of the
first volume of {\it Crelle} on separate
plates from different fascicules.

The literature on Malfatti's problem is
extensive, widely scattered, and not always
aware of itself.  Many of the published proofs
follow Hart's justification \cite{Har} of
Steiner's construction.  Others use
trigonometry, or are concerned with
generalizations to circles, spheres and
sections of quadrics.  Some references,
chronologically since Malfatti, are
\cite{A}, \cite{Cay1}, \cite{Q}, \cite{Sch},
\cite{Cay2}, \cite{Cay3}, \cite{Cat}, \cite{Sim},
\cite{Cay4}, \cite[pp.\,154--155]{C}, \cite{Der},
\cite[{\bf I} book III pp.\,311--314]{Rou}, \cite{Pam1},
\cite{Pam2}, \cite[{\bf IV} ch.II Ex.8 pp.\,65--69]{Bak},
\cite{Rog1}, \cite{Pro}, \cite{N1}, \cite{Rog}.

Here is a pot-pourri of facts.  I'll hide
my proofs because they use what
Swinnerton-Dyer calls ``the sort
of mathematics that no gentleman does in
public.''

\begin{figure}[h]
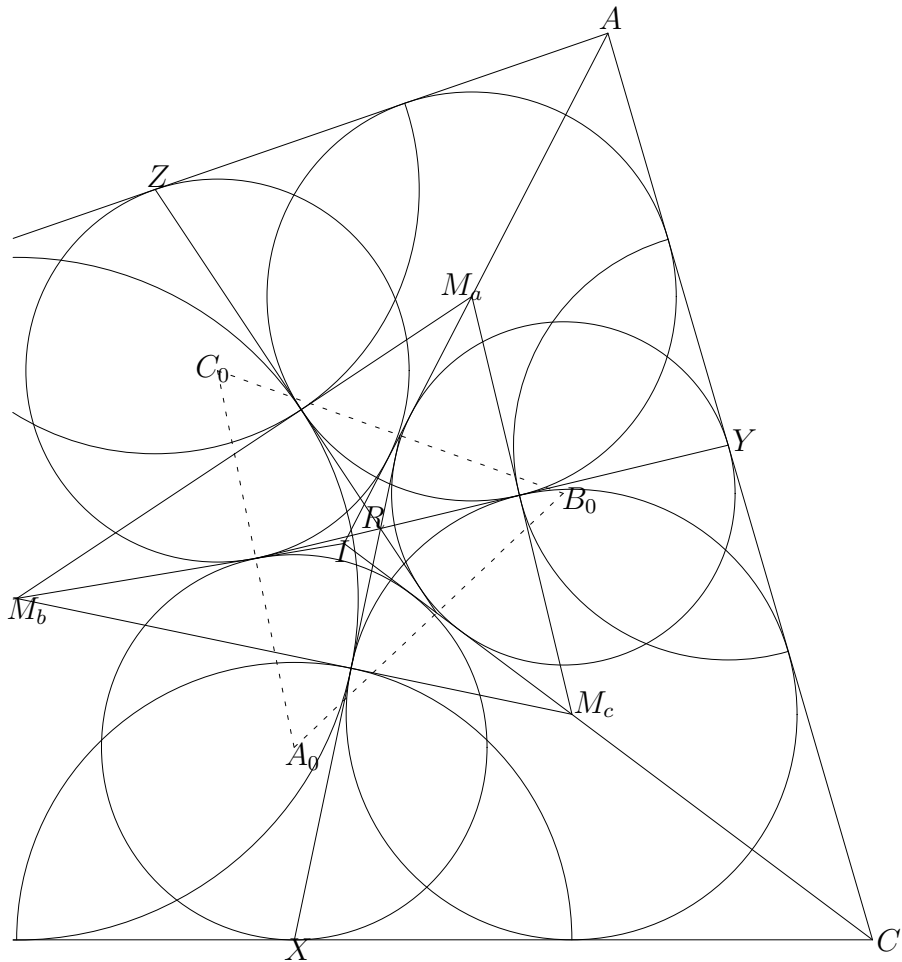
  

\caption{Construction of the Malfatti circles.}
\label{constructmalf}
\end{figure}

\clearpage

Most constructions start with incircles
of triangles $BIC$, $CIA$, $AIB$, where
$I$ is the incentre of $ABC$.  Even when
$ABC$ is exceptionally scalene, the radii
of these incircles remain surprisingly
equal.  It is hard to draw a picture which
clearly distinguishes the points.  In
Figure \ref{constructmalf} vertex $B$ is
off left. The three incircles touch
$BC$, $CA$, $AB$ at $X$, $Y$, $Z$ respectively
and their centres are $A_0$, $B_0$, $C_0$.

\begin{figure}[h]  
%
\caption{Detail of the Malfatti construction.}
\label{detailmalf}
\end{figure}

\clearpage

More detail is shown in Figure \ref{detailmalf}
where the perpendiculars (dotted lines) from
$A_0$, $B_0$, $C_0$ onto $AI$, $BI$, $CI$
indicate their points of contact.  Letters
in parentheses at the ends of lines are
labels of points on those lines, but too
distant to appear in the figure.
Lines $B_0C_0$, $C_0A_0$ and $A_0B_0$
are dashed and they intersect $AI$, $BI$
and $CI$ at $X'$, $Y'$ and $Z'$ respectively.
Then, surprisingly, $XX'$, $YY'$ and $ZZ'$
are the reflexions of $AI$, $BI$ and $CI$ in
$B_0C_0$, $C_0A_0$ and $A_0B_0$ respectively,
and they are not only transverse common
tangents to pairs of the incircles, but also
transverse common tangents to pairs of the
sought-after Malfatti circles.  They concur
in $R$, the radical centre of these circles.
The quadrilaterals $YRZA$, $ZRXB$, $XRYC$
are inscribable in the sense that they each
have an incircle which touches all four
sides.  These incircles are the
Malfatti circles, which touch $XX'$, $YY'$,
$ZZ'$ in pairs at $X_1$, $Y_1$, $Z_1$
respectively.

A variant construction uses the circle centre
$X$ and radius $r(1+u)/2$, where $r$ is the
inradius of the original triangle $ABC$,
and $u=\tan A/4$.  This circle, just half of
which is shown in Figure \ref{constructmalf},
passes through the points of contact of the
Malfatti circles with $BC$ and with each other.
Similarly for circles with centres $Y$, $Z$
and radii $r(1+v)/2$, $r(1+w)/2$, where
$v=\tan B/4$ and $w=\tan C/4$.

Note that $u$, $v$, $w$ are related, because $A+B+C=\pi$,
$\tan\frac{A+B}{4}=\tan\frac{\pi-C}{4}$ and
$(u+v)/(1-uv)=(1-w)/(1+w)$ or $1+uvw=u+v+w+vw+wu+uv$

\medskip

\noindent
{\bf Paradox 16.}  This identity may be
written in more than fifty different ways,
but never in the one that's needed at any
particular moment.

Andrew Bremner \cite{BG} has shown that the
radii of all 96 Malfatti circles are rational
just if (two of) $u$, $v$ and $w$ are rational.

\section{The light dawns.}
The Malfatti miracle comes
in four parts.  The first part is the relevance of
The Lighthouse Theorem.  The trick is to concentrate,
not on the Malfatti circles, but on the three
incircles with which the construction begins.
The 16 incentres of the four triangles BIC are
waiting for us in the previous section!  They
are generated by two 4-beam lighthouses at
$B$ and $C$ with phases $\beta=B/4$ and
$\gamma=C/4$ and have hexadecimal labels 0
to f in Figures \ref{fourlite}, \ref{dup4}
and \ref{largedup4}, where the missing vertex
$A$ is on the verge of making its ghostly appearance.

\noindent
{\bf Paradox 17}.  With no triangle to start
from, the Lighthouse Theorem for $n=4$ generates
the full set of 32 Malfatti solutions.

It remains to set up two pairs of 4-beam
lighthouses, one pair at $C$ \& $A$ with
phases $\gamma$ \& $\alpha=A/4$  and
another at $A$ \& $B$ with phases
$\alpha$ \& $\beta$.  These have
angle-quadrisecting beams of angles
$BCA$ \& $CAB$ and of angles $CAB$ \&
$ABC$ that intersect at the 16 incentres
of the four triangles $CIA$ and the 16
incentres of the four triangles $AIB$.
One application of the Lighthouse
Duplication Theorem yields the pairs of
angle-bisectors of triangle $ABC$, and
a second application produces the edges!

\section{Extraversion again.}

Recall that an
$A$-flip replaces angles $A$, $B$, $C$ by
$-A$, $\pi-B$, $\pi-C$ respectively.
That is, it replaces $(u,v,w)$  by
$\left(-u,\frac{1-v}{1+v},\frac{1-w}{1+w}\right)$.

\begin{figure}[h]
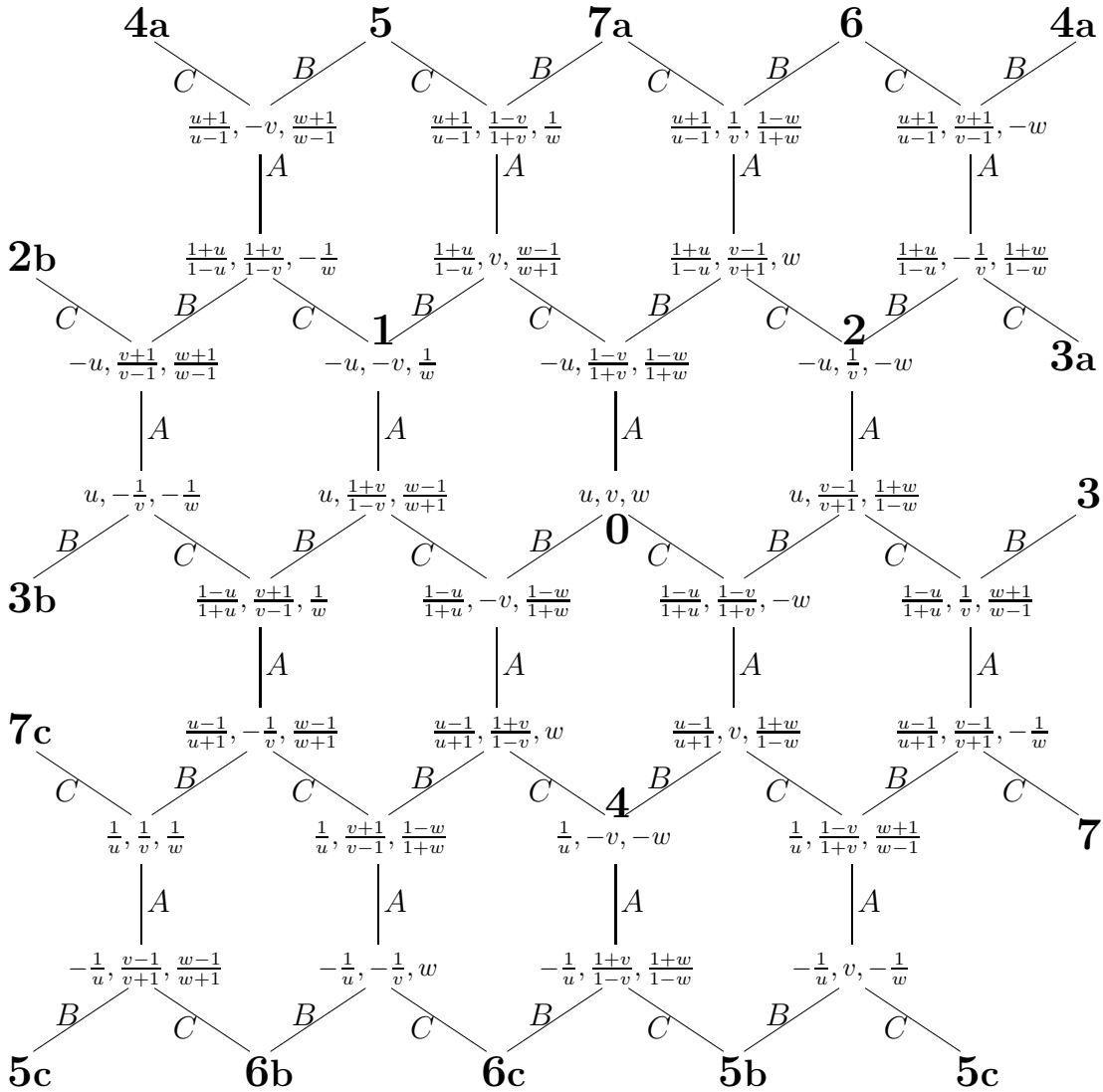
  

\caption{There are 32 Malfatti solutions.}
\label{32torus}
\end{figure}

\clearpage

So start at $(u,v,w)$ in the middle of
Figure \ref{32torus} and keep flipping.
You'll find that you're making a toroidal map
like Figure \ref{18torus}, but this one has
16 hexagonal regions, bounded by 48 edges,
16 each of $A$-flips, $B$-flips and $C$-flips,
which meet three at each of 32 vertices which
represent the 32 Malfatti solutions.

Eight {\bf ordinary} solutions (Table \ref{quarter})
are labelled with a {\bf bold} digit:

\begin{table}[h]   
\begin{tabbing}
\qquad \qquad \= {\bf0}$(u,v,w)$ \qquad \=
 {\bf3}$(u,-\frac{1}{v},-\frac{1}{w})$ \qquad \=
 {\bf5}$(-\frac{1}{u},v,-\frac{1}{w})$ \qquad \=
 {\bf6}$(-\frac{1}{u},-\frac{1}{v},w)$ \\
 \> {\bf7}$(\frac{1}{u},\frac{1}{v},\frac{1}{w})$
 \> {\bf4}$(\frac{1}{u},-v,-w)$
 \> {\bf2}$(-u,\frac{1}{v},-w)$
 \> {\bf1}$(-u,-v,\frac{1}{w})$ \\[-18pt]
\end{tabbing}
\vspace{-10pt}
\caption{Tangents of quarter-angles of ordinary solutions.}
\label{quarter}
\end{table}

The other 24 {\bf flipped} solutions are obtained
via an $A$-flip, $B$-flip or $C$-flip and are
labelled by respectively appending {\bf a, b} or
{\bf c} to the ordinary solution digit.  For
example, {\bf2b} is a $B$-flip away from
{\bf2}$(-u,\frac{1}{v},-w)$, and the tangents
of its quarter-angles are
$$\left(\frac{1-(-u)}{1+(-u)},
-\!\left(\frac{1}{v}\right),
\frac{1-(-w)}{1+(-w)}\right)=
\left(\frac{1+u}{1-u},-\frac{1}{v},
\frac{1+w}{1-w}\right)$$

\section{The Malfatti group.}

What is the group of
the 32 solutions?  It is generated by
$A$, $B$, $C$: notice the following
relations:
$$A^2=B^2=C^2=(ABC)^2=(BC)^4=(CB)^4=(CA)^4=
(AB)^4=I,\ ABC=CBA.$$
I thank Alex Fink for locating the group.  It
has no name (nor do most of the other 50
groups of order 32), so let's give it one.
It has Hall-Senior number and Magma number
34.  An official set of generator relations
is $g_2^2=g_4$, $g_3^2=g_5$,
$g_1^{-1}g_2g_1=g_2g_4$, $g_1^{-1}g_3g_1=g_3g_5$
which are realized here by $g_1=A$, $g_2=AB$,
$g_3=AC$, $g_4=(AB)^2=(BA)^2$, $g_5=(AC)^2=(CA)^2$.
The centre is \{{\bf0,3,5,6}\} (the {\bf evil}
numbers).  The order 4 elements are
{\bf 1a,1b,1c,2a,2b,2c,4a,4b,4c,7a,7b,7c}
(the $A$-, $B$- and $C$-flips of the
{\bf odious} numbers).  There are 19 elements
of order 2.  The numbers of subgroups of orders
\{2,4,8,16\} are \{19,31,30,7\}, the 7 maximal
subgroups being C4$\times$C4 and 6 copies
of D8$\times$C2.  It is a semidirect product
of C4$\times$C4 with C2.

\section{Which incircles do I draw?}
Table \ref{incent}
displays the incentres,  $A_i$, $B_j$, $C_k$ of the
appropriate respective  triangles $BIC$, $CIA$, $AIB$,
as triples $ijk$, where $i$, $j$, $k$ are the
hexadecimal numbers which appear in Figures
\ref{fourlite} to \ref{largedup4} and the accompanying
list of incentres at the end of {\bf \S7}.

\begin{table}[h]   
\begin{center}
\begin{tabular}{ccccccccc}
       & {\bf0}  & {\bf1}  & {\bf2}  & {\bf3}
    & {\bf4}  & {\bf5}  & {\bf6}  & {\bf7} \\
ordinary  &   000   &   280   &   802   &   a82
    &    028  &   2a8   &   82a   &   aaa  \\
A-flip &   541   &   d43   &   7c1   &   fc3
    &    feb  &   7e9   &   d6b   &   569  \\
B-flip &   154   &   17c   &   bfe   &   bd6
    &    3d4  &   3fc   &   97e   &   956  \\
C-flip &   415   &   ebf   &   43d   &   e97
    &    c17  &   6bd   &   c3f   &   695 \\[-15pt]
\end{tabular}
\end{center}
\caption{Incentres appropriate to the 32 solutions.}
\label{incent}
\end{table}

The pattern becomes memorable if you write
the solution number in binary and replace the
{\bf1}-digits by {\bf2}s.  Then recall the
base-4, beam-number, forms of the hexadecimal
numbers.  For example, think of solution
{\bf6} = {\bf110} as 220.  Its incentres are
\{8,2,a\} or, in base 4, \{20,02,22\} which
are beam numbers from lighthouses at $BC,CA,AB$,
or 2, 2, 0 from $A$, $B$, $C$ respectively.

To find the beam numbers for a flipped solution,
for example {\bf6b}, fix the $B$ beam number and
ad{\bf v}ance or {\bf d}ecrease the other two beam
numbers according as the solution number is
e{\bf v}il or o{\bf d}ious. As {\bf6} is evil,
2,2,0 become 2+1,2,0+1 = 3,2,1 and the incentres
are illuminated by beams 2 \& 1, 1 \& 3, 3 \& 2
from $B$ \& $C$, $C$ \& $A$, $A$ \& $B$ respectively,
so that they are the points 21=9, 13=7, 32=e.
An odious example is {\bf4c}: solution number
{\bf4} is {\bf100} in binary, giving beam numbers
2,0,0 which $C$-flip to $2\!-\!1,0\!-\!1,0$ =
1,3,0 and illuminate incentres 30,01,13 = c,1,7.

\smallskip

\noindent
{\bf Query}.  Triangle $ABC$ has 4 incentres, $I$.
Each of triangles $IBC$, $ICA$, $IAB$ has 4
incircles. There are $4^4$ choices of a triple of incircles.
Where are the 256 radical centres of these triples?
Are there coincidences? collinearities? concyclicities?

\section{In or out, near or far?}

What do the solutions
look like?  The eight ordinary solutions have all
their three circles in the interior angles.  The
flipped solutions have two exterior circles, but
keep the appropriate (A-, B- or C-)circle interior,
except that it becomes ``vertically opposite''
(outside the triangle, and denoted by lower case
letters in Table \ref{nearfar}) if the solution number
``contains'' 4, 2 or 1 respectively [4, 5, 6, 7
contain 4; while 2, 3, 6, 7 contain 2; and 1, 3, 5, 7
contain 1].  Call a circle ``Near'' if it
touches the edges at points nearer to its vertex
than the other two circles do.  Otherwise it's
``Far''---the proximity of the three circles is as
follows---in Table \ref{nearfar} I've written the
solution numbers in binary so that you can see
the pattern:

\begin{table}[h]    
\begin{center}
\begin{tabular}{ccccccccc}
solution & {\bf000}  & {\bf001}  & {\bf010}  & {\bf011}
    & {\bf100}  & {\bf101}  & {\bf110}  & {\bf111} \\
ordinary &   NNN   &   FFN   &   FNF   &   FNN
    &    NFF  &   NFN   &   NNF   &   FFF  \\
A-flip &   FNN   &   NNF   &   NFN   &   NNN
    &    fNN  &   nNF   &   nFN   &   nNN  \\
B-flip &   NFN   &   NNF   &   NfN   &   NnF
    &    FNN  &   NNN   &   FnN   &   NnN  \\
C-flip &   NNF   &   NNf   &   NFN   &   NFn
    &    FNN  &   FNn   &   NNN   &   NNn
\vspace{-15pt}
\end{tabular}
\end{center}
\caption{Positions of the Malfatti circles.}
\label{nearfar}
\end{table}

\newpage

There's not room to draw 32 separate
pictures, but turning Malfatti inside out
yields eight solutions in one handy package.
It's well known that the simplest construction for
an incircle or circumcircle of a triangle is to
draw the circle first, and the same goes for Malfatti!
Start with three mutually touching circles and draw the
direct common tangents to pairs of them.  You may
choose one from each pair of tangents to get
8 triangles for which the circles are Malfatti circles.

\begin{figure}[h]  
\begin{picture}(320,230)(-300,0)
\setlength{\unitlength}{2pt}
\drawline(-119,0)(35,0)(0,120)(-119,0) 
\put(-49,29.17){\circle{58.33}}        
\put(3.5,23.625){\circle{47.25}}       
\put(-12.5714,70.5){\circle{51.86}}    
\put(1,118,){$A_0$}
\put(-124,-3){$B_0$}
\put(36,-3){$C_0$}
\drawline(22.2297,43.7838)(-58.4915,61.0169)   
\put(-65.5,61){$B_1$}
\put(23,43){$C_2$}
\drawline(-12.5588,0)(-42.9180,76.7213)   
\put(-13,-5){$C_4$}
\put(-50,76){$A_1$}
\drawline(16.0417,65)(-34.78125,0)   
\put(17,65){$A_2$}
\put(-34,-5){$B_4$}
\put(-18,17.5){$A_4$}
\put(2.5,42.5){$B_2$}
\put(-34,56.5){$C_1$}
\end{picture}
\caption{Eight triangles from three Malfatti circles.}
\label{insideout}
\end{figure}
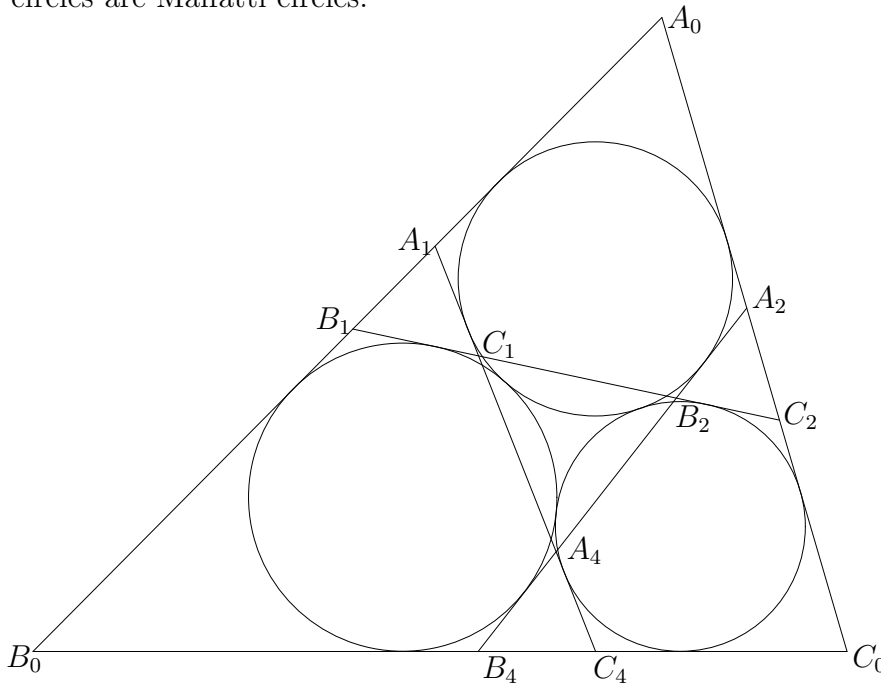

In Figure \ref{insideout} the eight triangles $ABC$
are manifestations of the solutions listed in the
first row of Table \ref{malsol}.  Note the
relation of the subscripts to the solution
numbers.  Other solutions
are found by varying the relative sizes of the
circles, in a kind of extraversion. For example,
if the A-circle is made small, then the solutions
are as in the third row of Table \ref{malsol}: they
are ordinary or A-flip solutions whose numbers
contain 4.  If the A-circle is enlarged, then the
solutions are as in the last row:
all solution numbers appear; solutions are ordinary
if they contain 1 \& 2, A-flips if they do not
contain 1 or 2, B-flips if 2 is the only number
contained or missing, and C-flips if 1 is the
only number present or missing.  Other solutions
can be obtained by permuting the sizes of the
circles.

\begin{table}[h]    
\begin{center}
\begin{tabular}{lcccccccc}
$ABC_{{\mbox{\tiny sub}}}$
   & 000 & 111  & 222  & 012 & 444 & 104 & 240 & 421 \\
        & {\bf0} & {\bf1c} & {\bf2b} & {\bf3a}
       & {\bf4a} & {\bf5b} & {\bf6c} & {\bf7} \\
A small & {\bf4} & {\bf5a} & {\bf6a} & {\bf7a}
       & {\bf4a} & {\bf5}  & {\bf6}  & {\bf7} \\
A large & {\bf3} & {\bf5b} & {\bf6c} & {\bf0a}
       & {\bf7}  & {\bf1c} & {\bf2b} &{\bf4a} \\[-18pt]
\end{tabular}
\end{center}
\caption{Malfatti solutions from given circles.}
\label{malsol}
\end{table}

\noindent
{\bf There are more than 32 Malfatti solutions.}
There are four degenerate solutions, consisting
of an incircle taken with multiplicity three,
and 24 semi-degenerate solutions consisting of a
repeated incircle (4 possibilities) and a third
circle touching two sides (3 choices) and
touching the incircle either proximally or
distally.  These solutions are more relevant to
Malfatti's original problem of maximizing the volume
of three cylinders cut from a triangular prism
\cite{ZL}.

But Steiner mentions 48 others! Twenty-four of these he
attributes to Gergonne \cite{Ger}.  They have two of the
circles touching an edge of the triangle ``in einem und
demselben Punkt'' ($X$ in Figure \ref{oddball}).  The
radii of the circles are no longer rational functions
of $u$, $v$ and $w$.

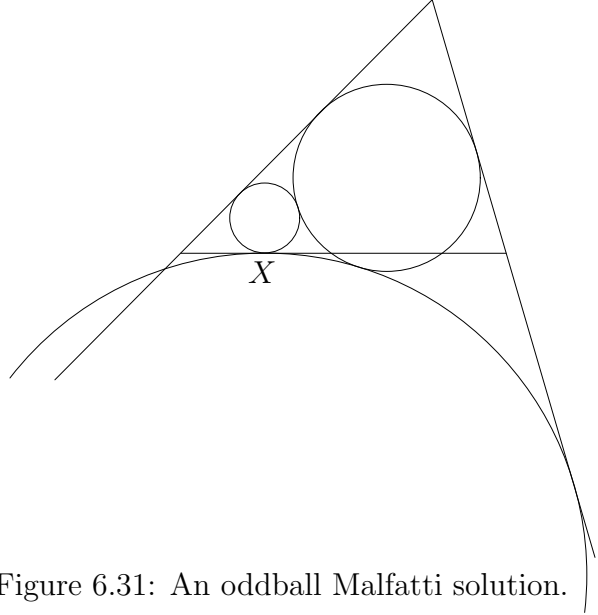
\begin{figure}[h]  
\begin{picture}(180,190)(-285,-105)
\setlength{\unitlength}{0.8pt}
\put(-79.2203,-152.2937){\arc{304.5875}{3.8}{6.4}}
\put(-79.2203,16.5749){\circle{33.1497}}
\put(-21.4472,35.5514){\circle{88.4799}}
\drawline(-119,0)(35,0)
\drawline(-178.5,-60)(0,120,)(77,-144)
\put(-87,-14){$X$}
\end{picture}
\caption{An oddball Malfatti solution.}
\label{oddball}
\end{figure}

\newpage

\section{Radpoints, guylines, oddpoints \& tielines}

[But see Chapter 7 below.]

These are Conway's names for members of a remarkable
configuration.  ``Radpoints'' are radical centres
of triads of Malfatti circles, so that there are 32
of them.   Their areal (barycentric) coordinates are
$$\left(\frac{u(1-u^2)}{1+u^2},\frac{v(1-v^2)}{1+v^2},
\frac{w(1-w^2)}{1+w^2}\right)$$
where $u$, $v$, $w$ are still the tangents of the
quarter-angles of triangle $ABC$ and they
range over the functions shown in Figure \ref{32torus}.
Denote the shapes of the coordinates by $I$ (identity),
$R$ (reciprocal), $S$ (switch), and $T$ (twist),
with lower case letters for their negatives.
``Switch'' and ``twist'' are closely related to ``flip''.

\begin{table}[h]    
\begin{center}
\begin{tabular}{l@{\hspace{4pt}}c@{\hspace{4pt}}c
  @{\hspace{4pt}}c@{\hspace{4pt}}c@{\hspace{4pt}}c
  @{\hspace{4pt}}c@{\hspace{4pt}}c@{\hspace{4pt}}c}
function & $x$ & $-x$ & $\frac{1}{x}$ & $-\frac{1}{x}$
  & $\frac{1-x}{1+x}$ & $\frac{x-1}{x+1}$
     & $\frac{1+x}{1-x}$ & $\frac{x+1}{x-1}$ \\[3pt]
radcoord & $\frac{x(1-x^2)}{1+x^2}$ & $\frac{x(x^2-1)}{x^2+1}$
   & $\frac{x^2-1}{x(x^2+1)}$ & $\frac{1-x^2}{x(1+x^2)}$
   & $\frac{2x(1-x)}{(1+x)(1+x^2)}$ & $\frac{2x(x-1)}{(x+1)(x^2+1)}$
   & $\frac{2x(x+1)}{(x-1)(x^2+1)}$ & $\frac{2x(1+x)}{(1-x)(1+x^2)}$ \\[3pt]
symbol & $I_x$ & $i_x$ & $R_x$ & $r_x$ & $S_x$ & $s_x$ & $T_x$ & $t_x$ \\[-18pt]
\end{tabular}
\end{center}
\caption{Shapes of coordinates of radpoints.}
\label{shapes}
\end{table}

\noindent
Table \ref{coords} shows the areal coordinates $\{x,y,z\}$
of all 32 radpoints, with braces and commas omitted.
As the coordinates are homogeneous, I may change all
three signs when necessary to produce a majority of
upper case letters.

\begin{table}[h]   
\begin{center}
\begin{tabular}{ccccccccc}
 & {\bf0} & {\bf1} & {\bf2} & {\bf3} & {\bf4} & {\bf5} & {\bf6} & {\bf7} \\
 {\bf o} & $I_uI_vI_w$ & $I_uI_vr_w$ & $I_ur_vI_w$
     & $i_uR_vR_w$ & $r_uI_vI_w$ & $R_ui_vR_w$
     & $R_uR_vi_w$ & $R_uR_vR_w$ \\
 {\bf a} & $i_uS_vS_w$ & $I_uT_vs_w$ & $I_us_vT_w$
     & $I_uT_vT_w$ & $r_uT_vT_w$ & $R_uS_vt_w$
     & $R_ut_vS_w$ & $R_uS_vS_w$ \\
 {\bf b} & $S_ui_vS_w$ & $T_uI_vs_w$ & $T_ur_vT_w$
     & $S_uR_vt_w$ & $s_uI_vT_w$ & $T_uI_uT_w$
     & $t_uR_uS_w$ & $S_uR_uS_w$ \\
 {\bf c} & $S_uS_vi_w$ & $T_uT_vr_w$ & $T_us_vI_w$
     & $S_ut_vR_w$ & $s_uT_vI_w$ & $t_uS_vR_w$
     & $T_uT_vI_w$ & $S_uS_vR_w$
\vspace{-15pt}
 \end{tabular}
 \end{center}
\caption{Coordinates of the 32 radpoints.}
\label{coords}
\end{table}

\noindent
Now it is easy to pick out collinearities.
I'll use the sets of coordinates ``III'',
``iSS'', ``SRt'', etc., with the subscripts
omitted, as symbols for ``the radical centre of
the solutions {\bf0}, {\bf0a}, {\bf3b}'' etc.,
while A, B, C still denote the vertices
of the triangle.  For example, the triples of
points \{A, iSS, RSS\}, \{B, III, IrI\}
and \{C, SRt, SRS\} are each
collinear, as you can see from the vanishing of
the determinants
$$%
$$

\noindent
The respective equations to the lines are
$y/S_v=z/S_w$, $x/I_u=z/I_w$ and $x/S_u=y/R_v$.
The triples always comprise a vertex, a
radpoint from an evil solution, and a
radpoint from an odious solution.

\noindent
{\bf There are 64 guylines.}  This is the second
part of the Malfatti miracle:

{\it There are 64 such collinearities!  Sixteen
``vertical'' ones through each of $A$, $B$ and $C$
and four ``nails'' through each of the four
Nagel points} $N_o$, $N_a$, $N_b$, $N_c$.

\noindent
{\bf What is a Nagel point?}  Choose three of the four
incircles of $ABC$ [if you choose {\bf onne}, as
Conway says for ``one and only one'', then you get
a Gerg{\bf onne} point].  Let them [it] touch $BC$,
$CA$, $AB$ at $A'$, $B'$, $C'$ respectively.
Then $AA'$, $BB'$, $CC'$ concur in a Nagel
[Gergonne] point.  The Nagel and Gergonne points are
listed in Table \ref{nagerg} with four representations
of their areal coordinates.

\begin{table}[h]    
\begin{center}
\begin{tabular}{c@{\hspace{2pt}}c@{\hspace{-2pt}}c@{\hspace{2pt}}c@{\hspace{2pt}}c}
$N_o$ & $\{s\!-\!a,s\!-\!b,s\!-\!c\}$
  & $\left\{\cot\frac{A}{2},\cot\frac{B}{2},\cot\frac{C}{2}\right\}$
  & $\{I\!-\!R,I\!-\!R,I\!-\!R\}$
  & $\left\{\frac{1\!-\!u^2}{u},\frac{1\!-\!v^2}{v},\frac{1\!-\!w^2}{w}\right\}$ \\[4pt]
$N_a$ & $\{s,c\!-\!s,b\!-\!s\}$
  & $\left\{-\cot\frac{A}{2},4\tan\frac{B}{2},4\tan\frac{C}{2}\right\}$
  & $\{I\!-\!R,T\!-\!S,T\!-\!S\}$
  & $\left\{\frac{1\!-\!u^2}{2u},\frac{2v}{v^2\!-\!1},\frac{2w}{w^2\!-\!1}\right\}$ \\[4pt]
$N_b$ & $\{c\!-\!s,s,a\!-\!s\}$
  & $\left\{4\tan\frac{A}{2},-\cot\frac{B}{2},4\tan\frac{C}{2}\right\}$
  & $\{T\!-\!S,I\!-\!R,T\!-\!S\}$
  & $\left\{\frac{2u}{u^2\!-\!1},\frac{1\!-\!v^2}{2v},\frac{2w}{w^2\!-\!1}\right\}$ \\[4pt]
$N_c$ & $\{b\!-\!s,a\!-\!s,s\}$
  & $\left\{4\tan\frac{A}{2},4\tan\frac{B}{2},-\cot\frac{C}{2}\right\}$
  & $\{T\!-\!S,T\!-\!S,I\!-\!R\}$
  & $\left\{\frac{2u}{u^2\!-\!1},\frac{2v}{v^2\!-\!1},\frac{1\!-\!w^2}{2w}\right\}$ \\[7pt]
$G_o$ & $\left\{\frac{1}{s\!-\!a},\frac{1}{s\!-\!b},\frac{1}{s\!-\!c}\right\}$
  & $\left\{\tan\frac{A}{2},\tan\frac{B}{2},\tan\frac{C}{2}\right\}$
  & $\{T\!-\!S,T\!-\!S,T\!-\!S\}$
  & $\left\{\frac{u}{1\!-\!u^2},\frac{v}{1\!-\!v^2},\frac{w}{1\!-\!w^2}\right\}$ \\[4pt]
$G_a$ & $\left\{\frac{1}{s},\frac{1}{c\!-\!s},\frac{1}{b\!-\!s}\right\}$
  & $\left\{4\tan\frac{A}{2},\cot\frac{B}{2},\cot\frac{C}{2}\right\}$
  & $\{T\!-\!S,I\!-\!R,I\!-\!R\}$
  & $\left\{\frac{2u}{1\!-\!u^2},\frac{v^2\!-\!1}{2v},\frac{w^2\!-\!1}{2w}\right\}$ \\[4pt]
$G_b$ & $\left\{\frac{1}{c\!-\!s},\frac{1}{s},\frac{1}{a\!-\!s}\right\}$
  & $\left\{\cot\frac{A}{2},4\tan\frac{B}{2},\cot\frac{C}{2}\right\}$
  & $\{I\!-\!R,T\!-\!S,I\!-\!R\}$
  & $\left\{\frac{u^2\!-\!1}{2u},\frac{2v}{1\!-\!v^2},\frac{w^2\!-\!1}{2w}\right\}$ \\
\end{tabular}
\end{center}
\caption{Coordinates of Nagel and Gergonne points.}
\label{nagerg}
\end{table}

The Nagel points are the incentres of the dilated
or anticomplementary triangle [with centre $G$,
the centroid, amplify the triangle $ABC$ by a
factor $-2$; shown dashed in Figure \ref{32nag}
with its angle-bisectors dotted].  They are a
set of orthocentric points generated by a pair
of two-beam lighthouses at two of the vertices
of $A'B'C'$.  Their nine-point centre is the
orthocentre of $ABC$. The Gergonne points are
not, in general, orthocentric.

Table \ref{radonguy} lists  the 64 guylines
({\bf0}$\leq${\bf X}$\leq${\bf7}).  Each radpoint
lies on 4 guylines, one through each of $A$, $B$
and $C$, and one through a Nagel point.  This
table is easy to memorize, since,
in each box, the four pairs of numbers have the
same nim-sum [binary addition without carry,
i.e., vector addition over GF(2), or XOR].  If
this nim-sum is nim-added to the row and column
numbers, the total is always 7.  Each guyline
may be given a three-digit label using the evil
digits 0, 3, 5, 6 taken from the row, column,
and entry in the table.  For example, the line
\{C\,2b\,6b\} has label 656.  Given a line label,
its points may be read as ``The first digit
gives the vertex (A=3, B=5, C=6) or Nagel point
(=0), the second specifies an ordinary (0) or
flipped (3,5,6) solution, the third is the evil
solution number, the odious one being that which
makes the nim-total 7.''  For example, 536
passes through B (=5) and through two A-flip
points (A=3), namely 6a and 7a (since the
nim-sum of 5, 3, 6 and 7 is 7).

\begin{table}[h]    
\begin{center}
\begin{tabular}{c@{\hspace{-1pt}}c@{\hspace{1pt}}|@{\hspace{2pt}}
c@{\hspace{3pt}}c@{\hspace{3pt}}c@{\hspace{2pt}}|@{\hspace{2pt}}c
@{\hspace{3pt}}c@{\hspace{3pt}}c@{\hspace{2pt}}|@{\hspace{2pt}}
c@{\hspace{3pt}}c@{\hspace{3pt}}c@{\hspace{2pt}}|@{\hspace{2pt}}c
@{\hspace{3pt}}c@{\hspace{3pt}}c@{\hspace{2pt}}|}
\cline{3-14}
  &        & \multicolumn{3}{c}{{\bf Xo} (=0)}
           & \multicolumn{3}{c}{{\bf Xa} (=3)}
           & \multicolumn{3}{c}{{\bf Xb} (=5)}
           & \multicolumn{3}{c}{{\bf Xc} (=6)} \\ \cline{3-14}
  &      & No & 0o III & 7o RRR & Na & 0a iSS & 4a rTT
   & Nb & 0b SiS & 2b TrT & Nc & 0c SSi & 1c TTr \\
N & (=0) &    & 3o iRR & 4o rII &    & 3a ITT & 7a RSS
   &    & 3b SRt & 1b TIs &    & 3c StR & 2c TsI \\
  &      &    & 5o RiR & 2o IrI &    & 5a RSt & 1a ITs
   &    & 5b TIT & 7b SRS &    & 5c tSR & 4c sTI \\
  &      &    & 6o RRi & 1o IIr &    & 6a RtS & 2a IsT
   &    & 6b tRS & 4b sIT &    & 6c TTI & 7c SSR \\ \cline{3-14}
  &      &  A & 0o III & 4o rII &  A & 0a iSS & 7a RSS
   &  A & 0b SiS & 1b TIs &  A & 0c SSi & 2c TsI \\
A & (=3) &    & 3o iRR & 7o RRR &    & 3a ITT & 4a rTT
   &    & 3b SRt & 2b TrT &    & 3c StR & 1c TTr \\
  &      &    & 5o RiR & 1o IIr &    & 5a RSt & 2a IsT
   &    & 5b TIT & 4b sIT &    & 5c tSR & 7c SSR \\
\end{tabular}
\end{center}
\caption{32 radpoints in pairs on 64 guylines.}
\label{radonguy}
\end{table}

\begin{figure}[h]
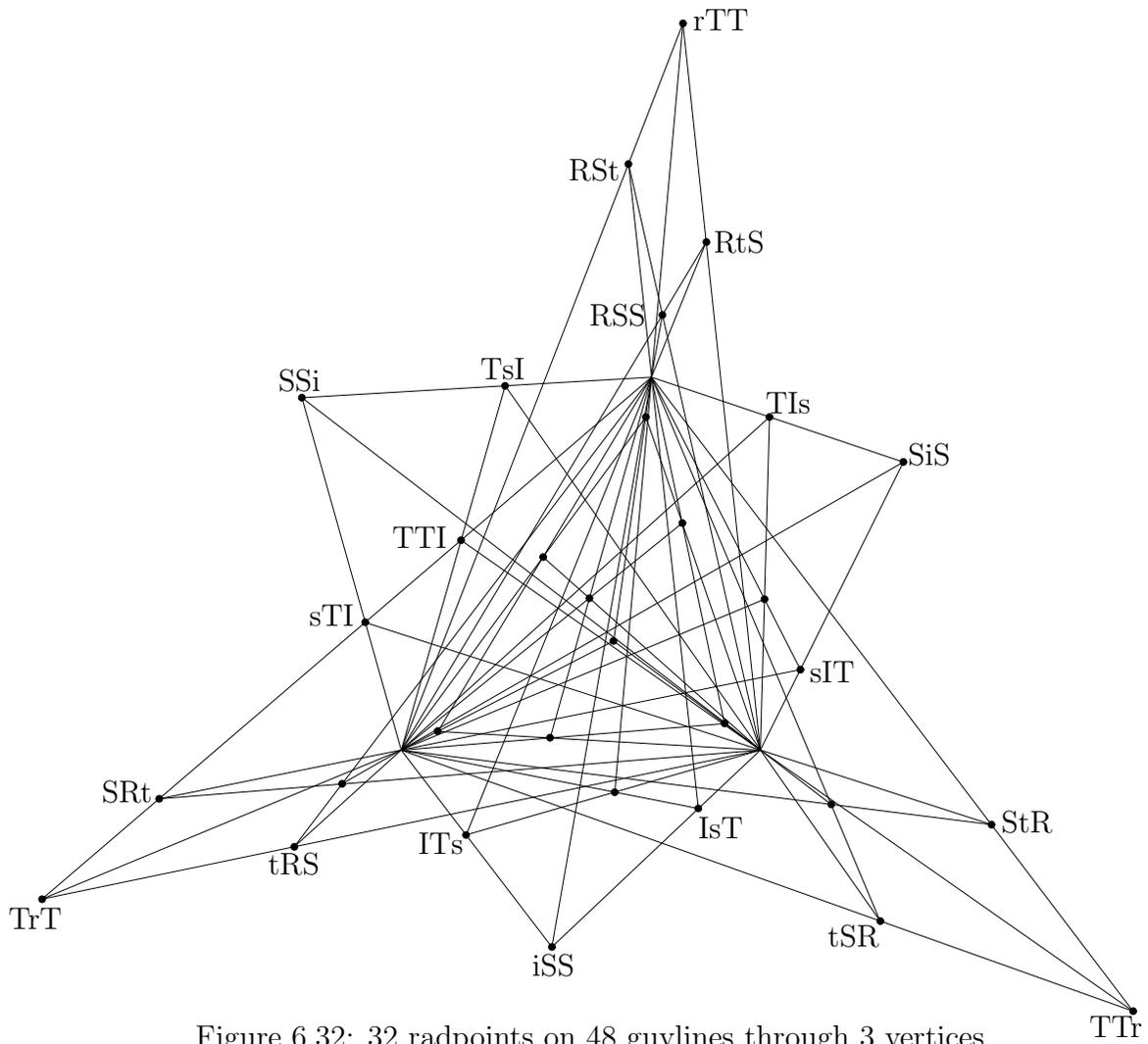
  
%
\caption{32 radpoints on 48 guylines through 3 vertices.}
\label{32rad}
\end{figure}

\clearpage

Figure \ref{32rad} shows the 48 ``vertical'' guylines.
Fortuitously \cite{G}, A 3b 2b and A 6c 4c with
equations $y/R_v=z/T_w$ and $y/T_v=z/I_w$ coincide
in Figure \ref{32rad}, which is drawn with
$v=\tfrac{1}{4}$, $w=\tfrac{1}{3}$, giving
$17y+50z=0$ in each case.

\begin{figure}[h]
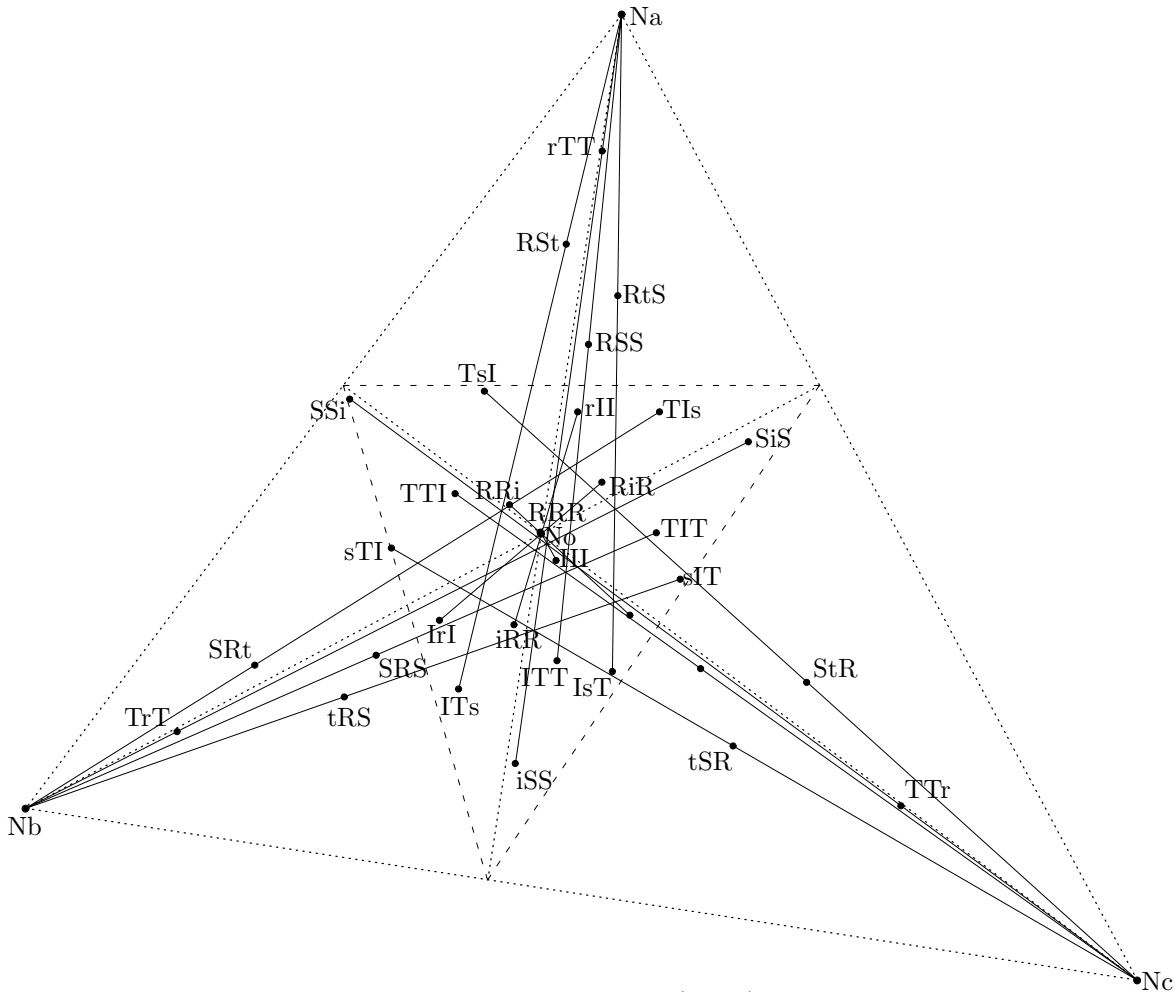
  
%
\caption{32 radpoints on 16 guylines (nails) through 4 Nagel points.}
\label{32nag}
\end{figure}

Figure \ref{32nag} shows the remaining 16 guylines:
four ``nails'' through each of the four Nagel points.
The triangle $ABC$ is not shown, but the dilated
triangle, twice its size, is shown dashed, with
its angle-bisectors, shown dotted, concurring in
the four Nagel points.

\clearpage

\noindent
{\bf There are 32 oddpoints}.  The third part
of the Malfatti miracle is:

{\it The 48 ``vertical'' guylines concur in threes
at 32 oddpoints}.

It is easy to see, for example, that the point
TII (with coordinates $\{T_u,I_v,I_w\}$) lies
on each of the three lines \{A III rII\},
\{B TTI TsI\}, \{C TIT TIs\}.  Table \ref{48in32}
gives the coordinates of the oddpoints and the
labels of the guylines therethrough.

\begin{table}[h]    
\begin{center}
\begin{tabular}{|@{\hspace{1.8pt}}c@{\hspace{1.8pt}}|@{\hspace{1.8pt}}c
@{\hspace{1.8pt}}|@{\hspace{1.8pt}}c@{\hspace{1.8pt}}|@{\hspace{1.8pt}}c
@{\hspace{1.8pt}}|@{\hspace{1.8pt}}c@{\hspace{1.8pt}}|@{\hspace{1.8pt}}c
@{\hspace{1.8pt}}|@{\hspace{1.8pt}}c@{\hspace{1.8pt}}|@{\hspace{1.8pt}}c
@{\hspace{1.8pt}}|} \hline
TTT & TTs & TsT & tSS & sTT & StS & SSt & SSS \\
 {\footnotesize 333 555 666} & {\footnotesize 336 556 666}
  & {\footnotesize 335 555 665} & {\footnotesize 330 556 665}
   & {\footnotesize 333 553 663} & {\footnotesize 336 550 663}
     & {\footnotesize 335 553 660} & {\footnotesize 330 550 660} \\[3pt] \hline
tRR & TIr & TrI & TII & sII & SRi & SiR & SRR \\
 {\footnotesize 303 565 656} & {\footnotesize 305 565 655}
  & {\footnotesize 306 566 656} & {\footnotesize 300 566 655}
   & {\footnotesize 300 560 650} & {\footnotesize 306 560 653}
     & {\footnotesize 305 563 650} & {\footnotesize 303 563 653} \\[3pt] \hline
RtR & ITr & IsI & RSi & rTI & ITI & iSR & RSR \\
 {\footnotesize 363 505 636} & {\footnotesize 363 503 633}
  & {\footnotesize 360 500 630} & {\footnotesize 360 506 635}
   & {\footnotesize 366 506 636} & {\footnotesize 366 500 633}
     & {\footnotesize 365 503 630} & {\footnotesize 365 505 635} \\[3pt] \hline
RRt & IIs & IrT & RiS & rIT & iRS & IIT & RRS \\
 {\footnotesize 353 535 606} & {\footnotesize 350 530 600}
  & {\footnotesize 353 533 603} & {\footnotesize 350 536 605}
   & {\footnotesize 355 535 605} & {\footnotesize 356 530 603}
     & {\footnotesize 355 533 600} & {\footnotesize 356 536 606} \\ \hline
\end{tabular}
\end{center}
\caption{48 guylines concur in threes at 32 oddpoints.}
\label{48in32}
\end{table}

\section{The SHOESTRING duality.}
There are very simple relations
between the coordinates of the radpoints and of the oddpoints.
To get from one radpoint to the other on the same guyline, or
to get from one oddpoint to the other on the same tieline
(soon to be defined in the fourth part of the miracle) just swap

\centerline{$I\longleftrightarrow R$ and $S\longleftrightarrow T$}

\noindent
in all three coordinates if the line is through a Nagel or
Gergonne point, but only in the $x$-, $y$- or $z$-coordinate
if the line is through $A$, $B$ or $C$ respectively.  To get from
a radpoint to an oddpoint on the same line, swap

\centerline{$I\longleftrightarrow S$ and $R\longleftrightarrow T$}

\noindent
in all three coordinates if the line is through a Nagel or Gergonne
point, but only in the one appropriate coordinate if the line
is through a vertex.

There are many other aspects
to this duality, which may be mnemonically listed:

\begin{center}
\begin{tabular}{|c@{\hspace{-1.3pt}}c@{\hspace{-1.3pt}}c@{\hspace{-1.3pt}}
  c@{\hspace{-1.3pt}}c@{\hspace{-1.3pt}}c@{\hspace{-1.3pt}}c|
  c@{\hspace{2.7pt}}c|c@{\hspace{2.7pt}}c|c@{\hspace{2.7pt}}c|
  c@{\hspace{2.7pt}}c@{\hspace{2.7pt}}c@{\hspace{2.7pt}}
  c@{\hspace{2.7pt}}c|}\hline
\rule{0pt}{12pt} & & & & & & & \multicolumn{2}{c|}{coordinates}
 & \multicolumn{2}{c|}{points} &  \multicolumn{2}{c|}{lines}
  & \multicolumn{5}{c|}{numbers} \\[3pt]\hline
\rule{0pt}{12pt} & H & & E & & RIN & & {\footnotesize REcIp}
 & {\footnotesize IdEN} & {\footnotesize Rad} & {\footnotesize NagEl}
  & {\footnotesize RIcH} & {\footnotesize NaIl} & {\footnotesize EvIl}
   & {\footnotesize zERo} & {\footnotesize tHREE} & {\footnotesize pENtE}
    & {\footnotesize HEx} \\
S & & O & & ST & & G & {\footnotesize SwiTch}
 & {\footnotesize TwiST} & {\footnotesize Odd} & {\footnotesize GerG.}
  & {\footnotesize Tie} & {\footnotesize Tail} & {\footnotesize OdiOuS}
   & {\footnotesize Seven} & {\footnotesize TeSSara} & {\footnotesize TwO}
    & {\footnotesize One}  \\[3pt] \hline
\end{tabular}
\end{center}

\noindent
In particular, numerical labels may be given to
radpoints and oddpoints by assigning 0 to a
coordinate I; 7 to a coordinate S; 3, 5, 6
respectively to $x$-, $y$-, $z$-coordinates R;
and 4, 2, 1 to coordinates T.  Radpoints will
have evil labels, oddpoints odious ones.  The
reader will notice many simple relationships
between the point labels and the labels of the
lines on which they lie.

\noindent
{\bf There are 64 tielines.}  The fourth part of the
Malfatti miracle is:

{\it The oddpoints lie in pairs on 64 tielines,
sixteen ``vertical'' ones through each of $A$,
$B$, $C$ and four ``tails'' through each of
the four Gergonne points.}

\begin{table}[h]   
\begin{center}
\begin{tabular}{c@{\hspace{-1pt}}c@{\hspace{1pt}}|@{\hspace{1pt}}
c@{\hspace{3pt}}c@{\hspace{3pt}}c@{\hspace{1pt}}|@{\hspace{1pt}}c
@{\hspace{3pt}}c@{\hspace{3pt}}c@{\hspace{1pt}}|@{\hspace{1pt}}
c@{\hspace{3pt}}c@{\hspace{3pt}}c@{\hspace{1pt}}|@{\hspace{1pt}}c
@{\hspace{3pt}}c@{\hspace{3pt}}c@{\hspace{1pt}}|}
\cline{3-14}
  &        & \multicolumn{3}{c}{{\bf Xo} (=7)}
           & \multicolumn{3}{c}{{\bf Xa} (=4)}
           & \multicolumn{3}{c}{{\bf Xb} (=2)}
           & \multicolumn{3}{c|}{{\bf Xc} (=1)} \\ \cline{3-14}
  &      & Go & 7o SSS & 0o TTT & Ga & 7a SRR & 3a TII
   & Gb & 7b RSR & 5b ITI & Gc & 7c RRS & 6c IIT \\
G & (=7) &    & 4o sTT & 3o tSS &    & 4a sII & 0a tRR
   &    & 4b rTI & 6b iSR &    & 4c rIT & 5c iRS \\
  &      &    & 2o TsT & 5o StS &    & 2a TrI & 6a SiR
   &    & 2b IsI & 0b RtR &    & 2c IrT & 3c RiS \\
  &      &    & 1o TTs & 6o SSt &    & 1a TIr & 5a SRi
   &    & 1b ITr & 3b RSi &    & 1c IIs & 0c RRt \\ \cline{3-14}
  &      &  A & 7o SSS & 3o tSS &  A & 7a SRR & 0a tRR
   &  A & 7b RSR & 6b iSR &  A & 7c RRS & 5c iRS \\
A & (=4) &    & 4o sTT & 0o TTT &    & 4a sII & 3a TII
   &    & 4b rTI & 5b ITI &    & 4c rIT & 6c IIT \\
  &      &    & 2o TsT & 6o SSt &    & 2a TrI & 5a IsT
   &    & 2b IsI & 3b RSi &    & 2c IrT & 0c RRt \\
\end{tabular}
\end{center}
\caption{32 oddpoints in pairs on 64 tielines.}
\label{oddontie}
\end{table}

\noindent
{\bf Help wanted.}    In Table \ref{oddontie}
I have retained the solution numbers, but what,
if any, is their relevance?  The duality between
the radpoints and the oddpoints is clear from
Tables \ref{coords} and \ref{48in32} and
Tables \ref{radonguy} and \ref{oddontie}.
The radpoints are radical centres, but what
relation do the oddpoints have to the solutions
I've listed them against?  From \cite[X(1488) and
X(2089)]{K} we discover the happy alliterations
that SSS is the Second Stevanovic point and
that TTT is the Third mid-arc point. At the latter
reference you can confirm that these two points
are collinear with the Gergonne point, $G_o$.
Perhaps the construction given in \cite{K}
generalizes to give all 32 oddpoints?

\newpage

\begin{quote}

Let $A'$, $B'$, $C'$ be the first points of
intersection of the angle bisectors of triangle
$ABC$ with its incircle. Let $A''B''C''$ be the
triangle formed by the lines tangent to the
incircle at $A'$, $B'$, $C'$.  Then $A''B''C''$
is perspective to the intouch triangle of $ABC$,
and the perspector is X(2089).

\end{quote}

\begin{flushright}

Darij Grinberg, Hyacinthos \#8072, 2003-10-01

\end{flushright}

Another fortuity \cite{G} is the concurrence of
three tielines \{A TTT sTT\}, \{B ITI IsI\},
\{C RRS RRt\} in Figure \ref{32odd} in a point
which is not an oddpoint.  In general these
lines do not concur.

\begin{figure}[h]
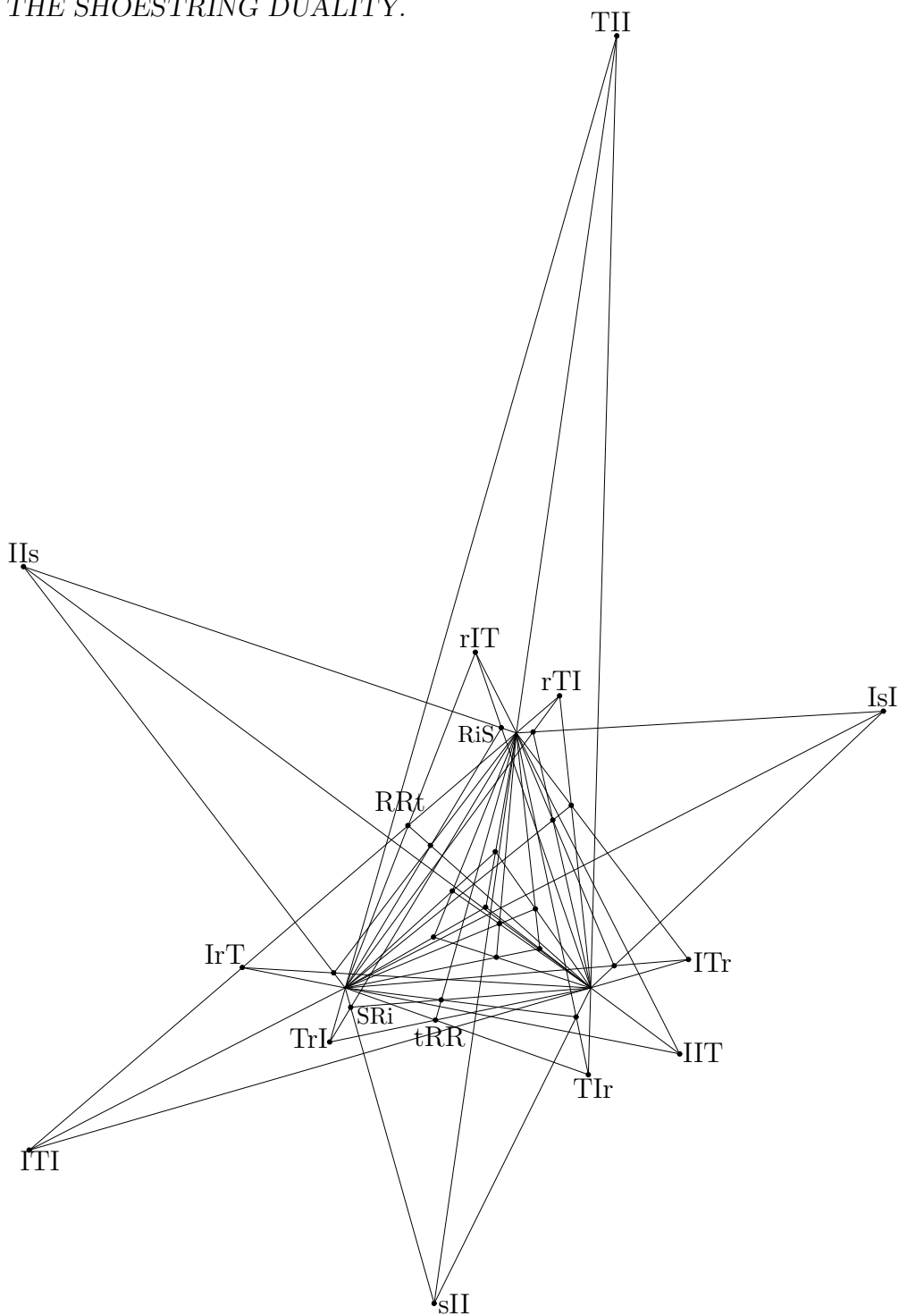
  

\caption{32 oddpoints on 48 tielines through 3 vertices.}
\label{32odd}
\end{figure}

\begin{figure}[h]  
%
\caption{32 oddpoints on 16 tielines (tails) through 4 Gergonne points.}
\label{32ger}
\end{figure}

\clearpage

\noindent
{\bf The truth is out there.} I leave $n=5$
(which, with its $\binom{5}{2}^3=1000$ points,
is particularly beautiful) and larger values
of $n$ to the reader.  There's much, much more
to be discovered.  A glimpse of the richness
for $n=6$ is seen in Figures \ref{3hyp} to
\ref{further}.  Since 2 and 3 divide 6, there
are three imbeddings of Lighthouse-2 and two
imbeddings of Lighthouse-3.  So there are two
complete Morley configurations together with
their interaction via $n=2$ with its sets of
orthocentric points. Are there further
incidences among the 18 GF-circles?

\section{History and Literature.}
The Morley triangle theorem has an interesting history.
It's a very Euclidean theorem, but it was discovered
comparatively recently---by Frank Morley around the
turn of the previous century, though he didn't
publish anything about it until many years later
\cite{Mor},\cite{Mor2}. The first published proofs appeared
in 1909: a trigonometric one by Satyanarayana
\cite{Sat} and a synthetic one by Naraniengar
\cite{Nar}.  Naraniengar's proof, which hinges on
a lemma, can be found in Coxeter \& Greitzer
\cite[pp.\,47--50]{C+G}. It was rediscovered by
Child \cite{Chi} and by others.  A direct proof
seems to be elusive, though the one given here in
Section {\bf3} may qualify.  There is an ``inside-out'' proof
due to Bricard \cite{Bri}, expounded in
\cite[pp.\,23--25]{Cox}; it uses a lemma that is
a corollary to the angle-bisector theorem.
A similar proof is due to Bottema \cite[p.\,34]{Bot1}.
For other proofs, see \cite{C+D}, \cite{Gro},
\cite{Ban}, \cite{T+M}, \cite{Gam}, \cite{Ros},
\cite{T+M}.  There are few hints that there is more
than one Morley triangle, but Honsberger
\cite[p.\,98]{Hon} asks the reader to show that
Morley's theorem holds also in the case of the
trisection of the exterior angles of a triangle.
He gives a proof, for just one more Morley triangle,
on pp.~163--164.  This is quoted in \cite{Kut},
an article drawn to my attention by Coxeter,
that is deserving of wider circulation.

There {\it is} a statement, rarely remembered or
repeated, of the existence of 27 points lying six
by six on nine (Morley) lines, at \cite[p.\,254, \S421]{J}.
This is attributed to \cite{T+M} where there are
two proofs.  The first is quite neat, attributed
to W.~E.~Philip, and given in \cite{J}.  Figure 3
on \cite[p.\,124]{T+M} is very like our Figure \ref{18tri}.
Glanville Taylor \& Marr also give some attention
to nonequilateral choices of three of the 27
Morley points.  Johnson's next section, his {\bf\S422},
gives the theorem that the incentres of the four
triangles of a cyclic quadrangle form the vertices
of a rectangle, and its
generalization to the first sixteen of the
Thrice Sixteen Theorem.  This is attributed to
Fuhrmann \cite[p.\,50]{F}.

I hope that Rigby's paper \cite{R} will eventually see the
light of day.  As he says there:

\begin{quote}
$\ldots$ I then consulted Morley's paper on the subject;
this contains results that rarely seem to be quoted,
so I have included an account of aspects of the theorem
that are apparently not generally known, especially the
connection with cardioids inscribed in the triangle.
\end{quote}

Honsberger \cite[p.\,93]{Hon} quotes Morley in this regard:

\begin{quote}
If a variable cardioid touch the sides of a triangle
the locus of its center, that is, the center of the
circle on which the equal circles roll, is a set of
9 lines which are 3 by 3 parallel, the directions being
those of the sides of an equilateral triangle.  The meets
of these lines correspond to double tangents; they
are also the meets of certain trisectors of angles,
internal and external, of the first triangle.
\end{quote}

These are, of course, the 9 Morley lines of the triangle.
Rigby generalizes this:

\noindent
{\bf Theorem}\cite{R}.  {\it The locus of the centres of
epicycloids $C(m,n)$ touching a given triangle consists
of $(n+m)^2$ lines $($with an exception if the triangle
is equilateral}).

Here the {\bf epicycloid} $C(m,n)$ is the envelope of
the line joining the points with polar coordinates
$(1,m\theta)$ and $(1,n\theta)$, where $0<m<n$ and
$m\perp n$.  The origin, $O$, is the {\bf centre}
of the epicycloid.  E.g., $(m,n)=(1,2)$ gives the
{\bf cardioid} and (1,3) the {\bf nephroid}.
If $m$ is allowed to be negative, $C(m,n)$ is a
{\bf hypocycloid}; for example $(m,n)=(-1,2)$ gives a
{\bf deltoid} and $(-1,3)$ an {\bf astroid}.
Better known are the connexions with the Steiner
deltoid (the envelope of the Simson-Wallace lines of the
triangle) and the nine-point circle, their common
tangents being parallel to the edges of the
Morley triangles.  See \cite[pp.\,345--349]{Bak2}
\cite[Ex.7, p.\,115]{Cox},
\cite[pp.\,226--231]{Dor}, \cite{Gu},
\cite[pp.\,72--79]{Loc2}, \cite{Mac}.
A referee has also pointed to \cite{Gam} and
\cite[ch.\,IV pp.\,173--194]{Leb}.
An interesting generalization, wherein the angles
of the triangle are partitioned in the ratios
1\,:\,$n$\,:\,1, and hence implicating
the Lighthouse Theorem for $n+2$, is given by
Fox \& Goggins \cite{F+G}.

\noindent
{\bf Query.}  If the lighthouses are rotating
with different angular velocities, is there a
connexion between the locus of the intersection
of the beams with the epicycloids considered
by Morley and Rigby?

The story began at about the time that I discovered
a ``lost notebook'' of Conway, and we embarked on
writing {\it The Triangle Book} \cite{CS}.
In the interim, Steve Sigur has taken over, and,
by the time you're reading this, the book will
probably be available.  Also, since the story began,
news has been spreading, and many people, including
D.~J.~Newman \cite{New} \cite[163--166]{Gale}
and Alain Connes, have taken an interest in
Morley's theorem (see
\verb+http://www.cut-the-knot.com/triangle/Morley+).

\chapter{More on Morley \& Malfatti}  

This is an unfinished paper by John Conway and myself.

In \cite{G} the second author discussed the geometry
of the points and lines associated with the Morley and
Malfatti configurations.  Here we shall introduce a
simple notation which has enabled us to establish the
incidences between the points and lines discussed in
\cite{G}.  In particular,
it shows immediately that the ``tielines'' and ``guylines''
are identical; we also correct several errors in sign.

It is simplest to define various points and lines by their
coordinates, and only later establish their relation to
the actual geometry.  We start by defining 64 points
$\langle\,ijk\,\rangle$ $\langle$pointed brackets for
points$\rangle$, where $i$, $j$, $k$ are integers
mod 4, by their barycentric (areal) coordinates
$$\langle\,ijk\,\rangle:\langle\,f(A+i\pi),f(B+j\pi),
f(C+k\pi)\,\rangle$$
where  $f(\theta)=\tan\frac{\theta}{4}\cos\frac{\theta}{2}$
and $A$, $B$, $C$ represent both the vertices and the
angles of our triangle.

It is immediately apparent that the points
$$A\quad\langle\,0jk\,\rangle\quad\langle\,1jk\,\rangle
\quad\langle\,2jk\,\rangle\quad\langle\,3jk\,\rangle$$
are collinear; we call this the {\bf guyline} $[\,*jk\,]$
[\,linear brackets for lines\,]
and define further guylines $[\,i*k\,]$ and $[\,i\,j\,*\,]$ through
$B$ and $C$ respectively.  There are 48 such guylines,
16 through each vertex.

It is easy to see that if $t=\tan\frac{\theta}{4}$, then
the tangents $t_n$ of $\frac{\theta+n\pi}{4}$ are given by
$$t_0=t\qquad t_1=\frac{1+t}{1-t}\qquad t_2=-\frac{1}{t}
\qquad t_3=\frac{t-1}{t+1}$$
If we define $u_n$, $v_n$, $w_n$ to be the values of these
when $\theta=A$, $B$, $C$ respectively, then it is also
easy to see that
$$\langle\,ijk\,\rangle=\left\langle\,
u_i\frac{1-u_i^2}{1+u_i^2},\,v_j\frac{1-v_j^2}{1+v_j^2},\,
w_k\frac{1-w_k^2}{1+w_k^2}\,\right\rangle$$

We check that $$t_0\frac{1-t_0^2}{1+t_0^2}\,+\,
t_2\frac{1-t_0^2}{1+t_0^2}=\left(t+\frac{1}{t}\right)
\frac{1-t^2}{1+t^2}=\frac{1-t^2}{t}
=\frac{2}{\sin\tfrac{\theta}{2}}$$
from which it follows that the line joining 000 to 222
passes through
$$\left\langle\csc\frac{A}{2},\,\csc\frac{B}{2},\,
\csc\frac{C}{2}\right\rangle$$
which are the coordinates of the original Nagel point, $N_0$.
We call this line a {\bf Nail} (German: Nagel = nail).  Similarly
$$t_1\frac{1-t_1^2}{1+t_1^2}+t_3\frac{1-t_3^2}{1+t_3^2}
=\frac{2}{\cos\tfrac{\theta}{2}}$$
showing that the line joining $\langle111\rangle$ to
$\langle333\rangle$ passes through
$$\left\langle\sec\frac{A}{2},\,\sec\frac{B}{2},\,
\sec\frac{C}{2}\right\rangle,$$
the original Gergonne point, $G_0$.  We call this line
a {\bf peG}.

\section{Extraversion}

All the formulas of triangle geometry remain true if the angles
$A$, $B$, $C$ are replaced by
$$~\quad i\pi+A,j\pi+B,k\pi+C\mbox{\qquad if
$i\!+\!j\!+\!k\equiv0$}$$
$$\mbox{or\quad }i\pi-A,j\pi-B,k\pi-C\mbox{\qquad if
$i\!+\!j\!+\!k\equiv2$}$$
The first author has called this process {\bf extraversion}
since it involves turning inside out and produces
``extra versions'' of triangle constructs.

\begin{figure}[h]
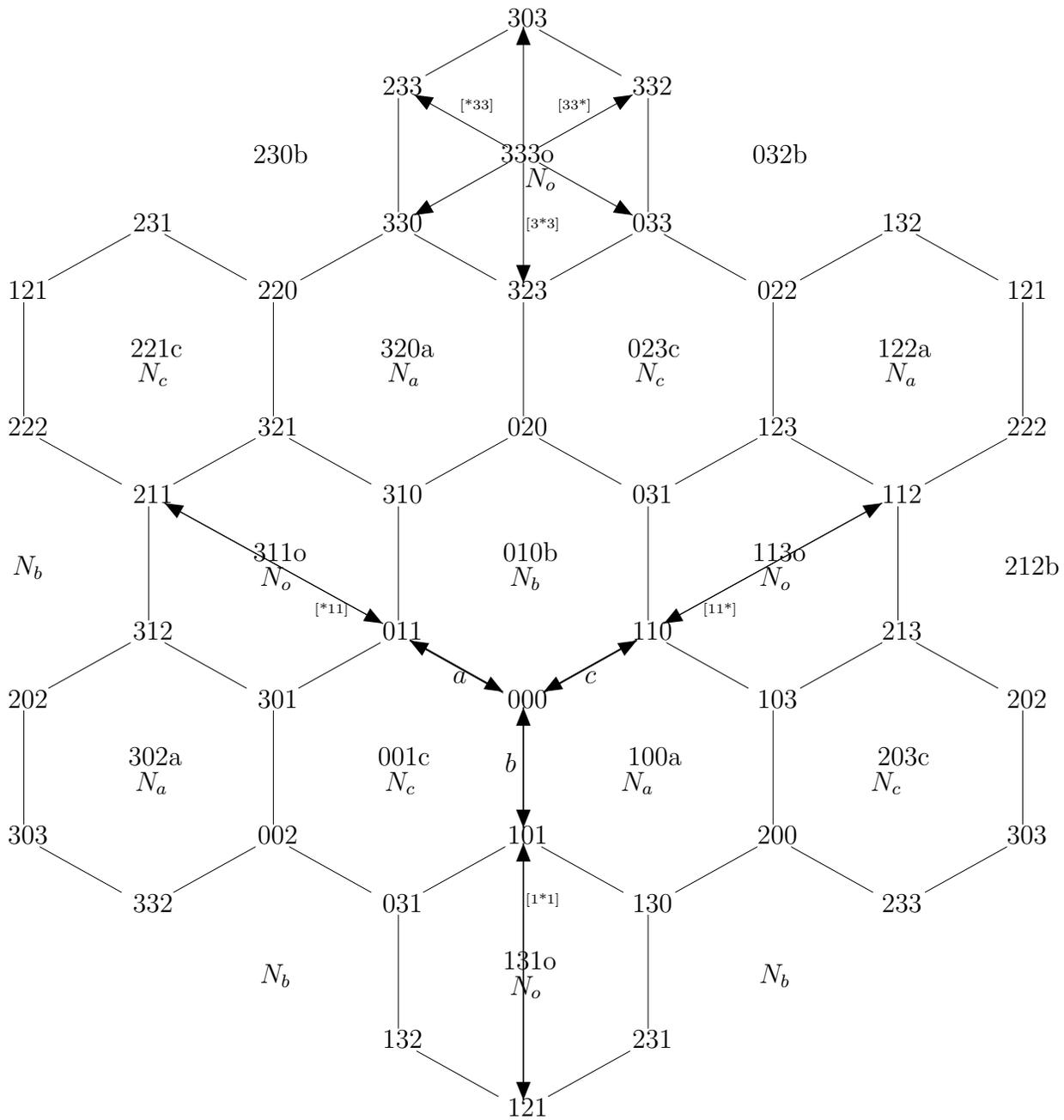
  

\caption{There are 32 Malfatti radpoints}
\label{malf}
\end{figure}

\clearpage

The above operators form the extraversion group, which
is generated by the particular operators called the
$a$-flip, $b$-flip, and $c$-flip, obtained by respectively
taking $ijk=011$, 101, 110 (see the double-arrows labelled
$a$, $b$, $c$ in Figure \ref{malf}).  For the
Malfatti case, the values of $i$, $j$, $k$ are
only relevant mod 4, and the extraversions of 000
are exhibited in Figure \ref{malf}.  Note that the
opposite edges of the figure are identified, forming
a torus, partitioned into 16 hexagonal regions.

\section{Extraverting the Malfatti configuration}

In \cite{G} it is shown that 000 is the radical
centre of the original three Malfatti circles.
Figure \ref{malf} shows that this is one of 32
triples of extra Malfatti circles, indexed by the
triples $ijk$ with $i\!+\!j\!+\!k$ even, that correspond
to the vertices of Figure \ref{malf}.  We call
these (Malfatti or Ajima) {\bf radpoints}, and
we have shown that the points $A$, $0jk$, $1jk$,
$2jk$, $3jk$ are collinear, and since two of these
are radpoints, we see that

\smallskip

\centerline{{\bf Each guyline is the join of two radpoints.}}

\smallskip

The two radpoints on the guyline [\,3*3\,] are the ends
(303 and 323) of the vertical diagonal of the
topmost hexagon in Figure \ref{malf}.  By symmetry
it follows that

\smallskip

\centerline{{\bf The 48 guylines correspond to the
48 diagonals of the hexagons.}}

\smallskip

The remaining two points on the typical guyline have
$i\!+\!j\!+\!k$ odd, so are not radpoints; we call them
{\bf oddpoints}.  However, they are easily
constructed from the radpoints.  For instance,
$\langle\,333\,\rangle$ is the intersection of
the guylines [\,*33\,], [\,3*3\,] and [\,33*\,], which we have
already constructed from the radpoints.  We see
that these guylines are the three diagonals of the
topmost hexagon in Figure \ref{malf}.  Extraverting
these, we obtain:

\smallskip

\centerline{{\bf The 16 points $\langle\,ijk\rangle\,$ with
$i\!+\!j\!+\!k\equiv1$ correspond to the hexagons.}}

\smallskip

We call these the {\bf 1-points}.  Each oddpoint is
the intersection of three guylines.  For the
particular 1-point $\langle\,333\,\rangle$ the guylines
are the diagonals of just one hexagon (the topmost
one), while for the particular 1-point 111 they are
the diagonals of three hexagons.  Extraverting this
we conclude that the oddpoints fall into two
distinct orbits of size 16, namely

\smallskip

\centerline{{\bf1-points}, with $i\!+\!j\!+\!k\equiv1$, whose
guylines are the diagonals of one hexagon;}

\smallskip

\centerline{{\bf3-points}, with $i\!+\!j\!+\!k\equiv3$, whose
guylines are the diagonals of three hexagons.}

\newpage

In Figure \ref{malf} $ijkx$ labels the hexagon whose
1-point is $ijk$, and whose Nail and peG pass through
$N_x$ and $G_x$ respectively.  The four concepts
corresponding to the hexagon $ijkx$ are:

The 1-point $\langle \,ijk\,\rangle$ on guylines $[\,*jk\,]$, $[\,i*k\,]$, $[\,ij*\,]$; \\
The 3-point $\langle \,i+2,j+2,k+2\,\rangle$ ``antipodal'' to this,
through the points \\
\mbox{ \hspace{80pt}} $\langle\,*,j+2,k+2\,\rangle,
\langle \,i+2,*,k+2\,\rangle, \langle \,i+2,j+2,*\,\rangle$; \\
The Nail $[\,ijkx\,]$ through $N_x$ and the points
$\langle \,ijk\,\rangle\pm\langle\,111\,\rangle$; \\
The peG $[\,i+2,j+2,k+2,x\,]$ through $G_x$ and
$\langle \,ijk\,\rangle$ and $\langle \,i+2,j+2,k+2\,\rangle.$

These four are exemplified in Figure \ref{malf} by

The 1-point $\langle\,333\,\rangle$ on diagonals
$[\,*33\,]$, $[\,3*3\,]$, $[\,33*\,]$ of the topmost hexagon; \\
the 3-point $\langle\,111\,\rangle$, visualized as the
intersection of the diagonals $[\,*11\,]$, $[\,1*1\,]$, $[\,11*\,]$
of hexagons 311, 131, 113; \\
the Nail $[\,333o\,]$ through $N_o$ and the points
$\langle\,000\,\rangle$, $\langle\,222\,\rangle$ furthest from
the (topmost) hexagon $333o$; and \\
the invisible peG $[\,111o\,]$ through $G_o$ and the points
$\langle\,333\,\rangle$ $\langle\,111\,\rangle$.  [This last one
is wrong, I think.  We remark that the last two points
are the vertices of Figure \ref{malf} furthest from
the hexagon.

Figure \ref{summ} summarizes the situation.

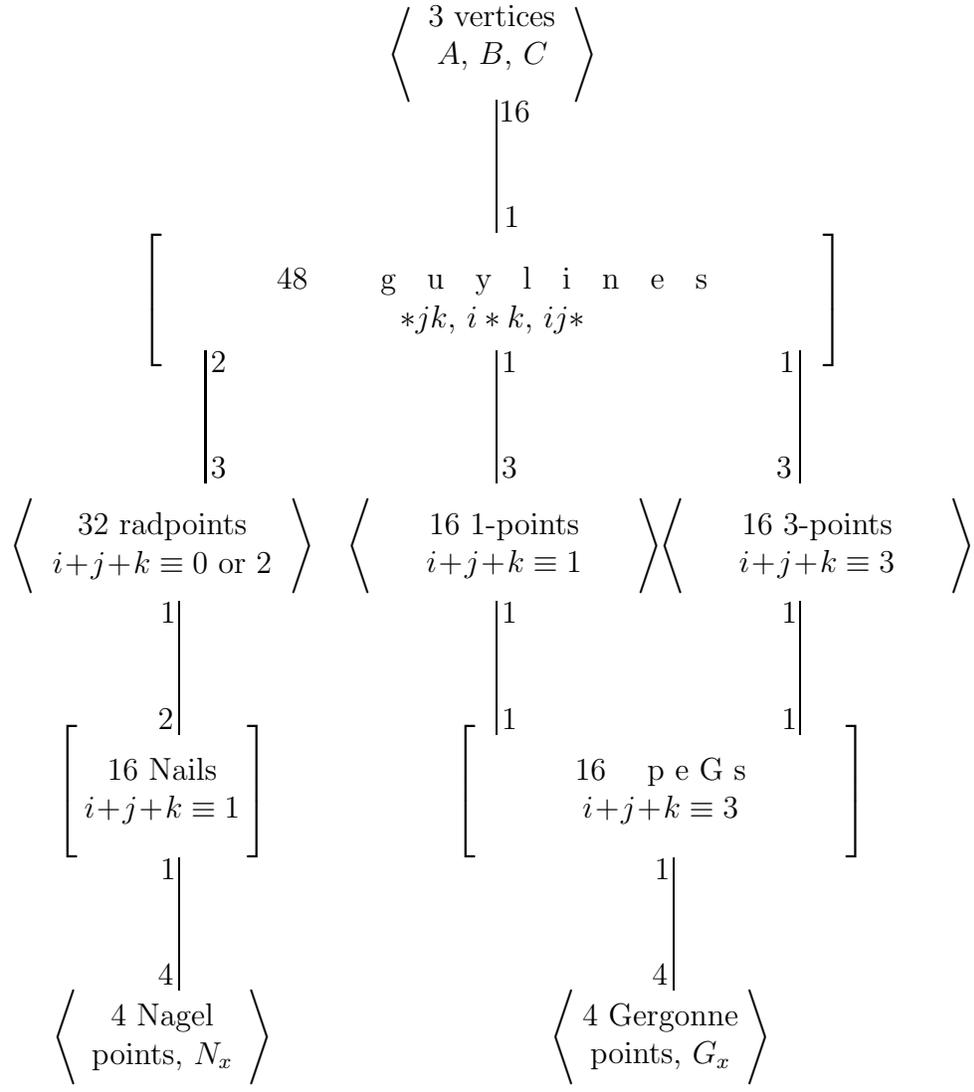
\begin{figure}[h]  
\begin{picture}(440,400)(-50,-200)
\setlength{\unitlength}{1pt}
\put(0,0){\begin{tabular}{ccc@{\hspace{-1pt}}c}
\multicolumn{4}{c}{$\left\langle\parbox{60pt}
{\begin{center}
3 vertices \\ $A$, $B$, $C$
\end{center}}\right\rangle$} \\[60pt]
\multicolumn{4}{c}{$\left[\parbox{250pt}
{\begin{center}
48 \qquad g\quad u\quad y\quad l\quad i\quad n\quad e\quad s \\
 $*jk$, $i*k$, $ij*$
 \end{center}}\right]$} \\[60pt]
\multicolumn{2}{c}{$\left\langle\parbox{96pt}
{\begin{center}
32 radpoints \\
$i\!+\!j\!+\!k\equiv0$ or 2
\end{center}}\right\rangle$} &
$\left\langle\parbox{100pt}{\begin{center}
16 1-points \\ $i\!+\!j\!+\!k\equiv1$
\end{center}}\right\rangle$ &
$\left\langle\parbox{100pt}{\begin{center}
16 3-points \\ $i\!+\!j\!+\!k\equiv3$
\end{center}}\right\rangle$ \\[60pt]
\multicolumn{2}{c}{$\left[\parbox{64pt}
{\begin{center}
16 Nails \\
$i\!+\!j\!+\!k\equiv1$ \end{center}}\right]$} &
\multicolumn{2}{c}{$\left[\parbox{140pt}
{\begin{center}
16 \quad p\ e\ G\ s \\ $i\!+\!j\!+\!k\equiv3$
\end{center}}\right]$} \\[60pt]
\multicolumn{2}{c}{$\left\langle\parbox{64pt}
{\begin{center}
4 Nagel \\ points, $N_x$
\end{center}}\right\rangle$} &
\multicolumn{2}{c}{$\left\langle\parbox{64pt}
{\begin{center}
4 Gergonne \\ points, $G_x$
\end{center}}\right\rangle$} \\
\end{tabular}}
\put(70,-162){\line(0,1){50}}
\put(257,-162){\line(0,1){50}}
\put(70,-65){\line(0,1){50}}
\put(190,-65){\line(0,1){50}}
\put(305,-65){\line(0,1){50}}
\put(80,30){\line(0,1){50}}
\put(190,30){\line(0,1){50}}
\put(305,30){\line(0,1){50}}
\put(190,125){\line(0,1){50}}
\put(191,167){16}
\put(193,127){1}
\put(82,72){2}
\put(82,32){3}
\put(192,72){1}
\put(192,32){3}
\put(297,72){1}
\put(296,32){3}
\put(63,-23){1}
\put(62,-63){2}
\put(192,-23){1}
\put(192,-63){1}
\put(298,-23){1}
\put(298,-63){1}
\put(63,-120){1}
\put(62,-160){4}
\put(250,-120){1}
\put(249,-160){4}
\end{picture}
\caption{Summary of Figure \ref{malf}.  The congruences are mod 4.}
\label{summ}
\end{figure}

\clearpage

\section{Oddpoints, mid-arc points, Stevanovic points}

The next page \& a bit are ramblings.  But I can now
describe the 32 ``3rd mid-arc points'' --- see Section 5,
between the Kimberling descriptions of ``3rd mid-arc
point'' and ``2nd Stevanovic point''.

There are other geometrical constructions for our oddpoints.
Our original 1-point, $\langle\,333\,\rangle$, is what
Kimberling calls the 2nd Stevanovic point, $X_{1488}$ in
\cite{K}, so the 16 1-points are extraversions of that.
$$\left\langle\frac{\sin A}{1+\sin A/2},\frac{\sin B}{1+\sin B/2},
\frac{\sin C}{1+\sin C/2}\right\rangle$$
$\big[$I still need to get Kimberling's $X_{1488}$ and $X_{2089}$
constructed, generalized (extraverted), and identified with
the 1-points and 3-points, respectively.  Presumably, as
there are only 16 1-points and 16 3-points, extraversion
doesn't yield 32 different cases.  If you separate the
operations `negative' and `reciprocal' in
$t\longrightarrow-\frac{1}{t}$, does this swap 1-points
with 3-points??  Here are the (clearly related)
barycentrics for $X_{2089}$\,:
$$\left\langle\frac{\sin A}{1-\sin A/2},\frac{\sin B}{1-\sin B/2},
\frac{\sin C}{1-\sin C/2}\right\rangle$$
So it looks as though you can get $X_{2089}$ from $X_{1488}$
by negating $A$, $B$ and $C$, i.e., by making an $a$-flip,
a $b$-flip, and a $c$-flip, n'est-ce pas?

Let's answer some questions!

From $t_3=\frac{t-1}{t+1}$ we have
$$u_3\frac{1-u_3^2}{1+u_3^2}=\frac{u-1}{u+1}\,\cdot\,
\frac{4u}{2u^2+2}=\frac{2\tan\tfrac{A}{4}(\tan\tfrac{A}{4}-1)}
{\sec^2\tfrac{A}{4}(\tan\tfrac{A}{4}+1)}=\frac{2\sin\tfrac{A}{4}
\cos\tfrac{A}{4}(\sin\tfrac{A}{4}-\cos\tfrac{A}{4})}
{\sin\tfrac{A}{4}+\cos\tfrac{A}{4}}$$
$$=\frac{\sin\tfrac{A}{2}(\sin^2\tfrac{A}{4}-\cos^2\tfrac{A}{4})}
{(\sin\tfrac{A}{4}+\cos\tfrac{A}{4})^2}=
-\frac{\sin\tfrac{A}{2}\cos\tfrac{A}{2}}
{1+2\sin\tfrac{A}{4}\cos\tfrac{A}{4}}
=-\frac{1}{2}\,\cdot\,\frac{\sin A}{1+\sin\tfrac{A}{2}}$$
so the point $\langle\,333\,\rangle$ is indeed the
2nd Stevanovic point, $X_{1488}$, since the factor
$-\frac{1}{2}$ can be omitted by homogeneity.  But
does this factor sabotage things when we come to
extraversions??  Let's do $\langle\,111\,\rangle$
first:

From $t_1=\frac{1+t}{1-t}$ we have
$$u_1\frac{1-u_1^2}{1+u_1^2}=\frac{1+u}{1-u}\,\cdot\,\frac{4u}{2u^2+2}
=\frac{2\tan\tfrac{A}{4}(1+\tan\tfrac{A}{4})}
{\sec^2\tfrac{A}{4}(1-\tan\tfrac{A}{4})}=\frac{2\sin\tfrac{A}{4}
\cos\tfrac{A}{4}(\cos\tfrac{A}{4}+\sin\tfrac{A}{4})}
{\cos\tfrac{A}{4}-\sin\tfrac{A}{4}}$$
$$=\frac{\sin\tfrac{A}{2}(\cos^2\tfrac{A}{4}-\sin^2\tfrac{A}{4})}
{(\cos\tfrac{A}{4}-\sin\tfrac{A}{4})^2}=
\frac{\sin\tfrac{A}{2}\cos\tfrac{A}{2}}
{1-2\sin\tfrac{A}{4}\cos\tfrac{A}{4}}
=\frac{1}{2}\,\cdot\,\frac{\sin A}{1-\sin\tfrac{A}{2}}$$
and we can omit the homogeneous $\frac{1}{2}$ to confirm
that $\langle\,111\,\rangle$ is $X(2089)$ that Kimberling
calls the 3rd mid-arc point.

Here are Kimberling's descriptions, slightly amplified
with formulas in our notation.  Note that there are no
obvious connexions between the constructions and the
Malfatti circles.$\big]$

\medskip

X(2089) = 3rd MID-ARC POINT [\,my TTT,
 i.e. $\langle\,111\,\rangle$\,]

Trilinears: $(\cos B/2 + \cos C/2 - \cos A/2) \sec A/2 :
(\cos C/2 + \cos A/2 - \cos B/2) \sec B/2 :
(\cos A/2 + \cos B/2 - \cos C/2) \sec B/2
= 1/(1 - \sin A/2) : 1/(1 - \sin B/2) : 1/(1 - \sin C/2)$

Barycentrics: $(\sin A)f(A,B,C):(\sin B)f(B,C,A):(\sin C)f(C,A,B)$,

$$\bigg[\mbox{\,i.e.,\ }\bigg\langle\,\ldots,\,\frac{v(1+v)}{(1-v)(1+v^2)},\,\ldots\,\bigg\rangle$$
We seem to have 16 examples of Desargues's theorem. $\bigg]$

Let $A'$, $B'$, $C'$ be the first points of intersection of the
angle bisectors of triangle $ABC$ with its incircle.
Let $A''B''C''$ be the triangle formed by the lines tangent
to the incircle at $A'$, $B'$, $C'$. Then $A''B''C''$ is perspective
to the intouch triangle of $ABC$, and the perspector is X(2089).
(Darij Grinberg, Hyacinthos \#8072, 10/01/03)

X(2089) lies on these lines: 1,167 \  2,178 \  7,1488

X(2089) = X(7)-Ceva conjugate of X(174)

X(2089) = X(173)-cross conjugate of X(174)

\section{Third mid-arc points}

Each of the four touch-circles (incircle \& three
excircles) is cut by three angle-bisectors in two
diametrically opposite points, so there are four
sets of three pairs of parallel tangents there.
These form four sets of eight triangles, each
triangle being homothetic to the respective
touch-triangle (whose vertices are the points of
contact of the edges of the original triangle
with the respective touch-circle).  So each member
of the four sets of eight triangles is in perspctive
with the respective touch-triangle, and we have 32
perspectors, each of which has a claim to be a
``third mid-arc point''.  Until my co-author gives
me a better notation, I will denote these points
by symbols $X_{abc}$ where $X=O, A, B$, or $C$
according to which touch-circle it belongs and
$a, b, c$ are $i$ or $o$, for `inner' and `outer',
according as the tangents to the touch-circle is
nearer to or further from the respective vertices
$A, B, C$.  For example, $O_{iii}$ is the third
mid-arc point as originally defined by Grinberg,
and is our point 111.

In fact 16 of the points are our ``1-points'':

\begin{center}
\begin{tabular}{cccc}
$O_{iii}=111$ & $O_{ioo}=133$ & $O_{oio}=313$ & $O_{ooi}=331$ \\
$A_{ooo}=300$ & $A_{oii}=322$ & $A_{ioi}=102$ & $A_{iio}=120$ \\
$B_{ooo}=030$ & $B_{oii}=012$ & $B_{ioi}=232$ & $B_{iio}=210$ \\
$C_{ooo}=003$ & $C_{oii}=021$ & $C_{ioi}=201$ & $C_{iio}=223$ \\
\end{tabular}
\end{center}

The other 16 have barycentric coordinates

$$\Big\langle \frac{u}{1+u^2}, \quad \frac{v}{1+v^2},
 \quad \frac{w}{1+w^2} \Big\rangle $$

where $u$ (``0'') extraverts to $\frac{1+u}{1-u}$ (``1''),
$-\frac{1}{u}$ (``2''), and $\frac{u-1}{u+1}$ (``3'').

LATER: I discover that this is Kimberling's X(174), the
Yff centre of congruence.   Unfortunately, I'm unable
to understand the description in ETC:  {description of
X(173), congruent isoscelizers point, is appended
below, in case it helps]

``In notes dated 1987, Yff raises a question concerning
certain triangles lying within ABC: can three
isoscelizers (as defined in connection with X(173),
P(B)Q(C), P(C)Q(A), P(A)Q(B) be constructed so that
the four triangles P(A)Q(A)A, P(B)Q(B)B, P(C)Q(C)C,
ABC are congruent? After proving that the answer is
yes, Yff moves the three isoscelizers in such a way
that the three outer triangles, P(A)Q(A)A, P(B)Q(B)B,
P(C)Q(C)C stay congruent and the inner triangle, ABC,
shrinks to X(174).

``Let D be the point on side BC such that (angle BID)
= (angle DIC), and likewise for point E on side CA
and point F on side AB. The lines AD, BE, CF concur
in X(174). [Seiichi Kirikami, Jan. 29, 2010]

``Generalization: if I is replaced by an arbitrary
point P = p : q : r (trilinears), then the lines
AD, BE, CF concur in the point
K(P) = f(p,q,r,A) : f(q,r,p,B) : f(r,p,q,C),
where f(p,q,r,A) = (q2 + r2 + 2qr cos A)-1/2.
Moreover, if P* is the inverse of P in the
circumcircle, then K(P*) = K(P). [Peter Moses,
Feb. 1, 2010, based on Seiichi Kirikami's
construction of X(174)] ''

``X(173) = CONGRUENT ISOSCELIZERS POINT

``Let P(B)Q(C) be an isoscelizer: let P(B) on
sideline AC and Q(C) on AB be equidistant from A,
so that AP(B)Q(C) is an isosceles triangle.
Line P(B)-to-Q(C), P(C)-to-Q(A), P(A)-to-Q(B)
concur in X(173). (P. Yff, unpublished notes, 1989)

``The intouch triangle of the intouch triangle of
triangle ABC is perspective to triangle ABC, and
X(173) is the perspector. (Eric Danneels,
Hyacinthos 7892, 9/13/03) ''

Alternatively, the coordinates are

$$\Big\langle \sin\frac{A}{2}, \quad \sin\frac{B}{2},
 \quad \sin\frac{C}{2} \Big\rangle $$

\noindent
which extravert to

$$\Big\langle-\sin\frac{A}{2}, \quad \cos\frac{B}{2},
 \quad \cos\frac{C}{2} \Big\rangle \mbox{,\ etc.} $$

Note that interchange of ``0'' and ``2'', or of ``1''
and ``3'' only changes the sign of the coordinate, so
that, for example, ``301'' is the same as ``123''.
The 16 points are, until I have a better notation, the
``0-points'' (= the ``2-points'')

\begin{center}
\begin{tabular}{cccc}
$O_{ooo}=``000$'' & $O_{oii}=``022$''
 & $O_{ioi}=``202$'' & $O_{iio}=``220$'' \\
$A_{iii}=``233$'' & $A_{ioo}=``211$''
 & $A_{oio}=``031$'' & $A_{ooi}=``013$'' \\
$B_{iii}=``323$'' & $B_{ioo}=``301$''
 & $B_{oio}=``121$'' & $B_{ooi}=``103$'' \\
$C_{iii}=``332$'' & $C_{ioo}=``310$''
 & $C_{oio}=``130$'' & $C_{ooi}=``112$'' \\
\end{tabular}
\end{center}

One way to describe these 16 points, is as constituting
four quadrangles, one associated with each touch-circle,
whose diagonal point triangles each coincide with the
original triangle, $ABC$.  We have the following
collinearities:  [Done in haste: needs checking]

\begin{center}
\begin{tabular}{ccc}
$A O_{ooo} O_{oii}$ & $B O_{ooo} O_{ioi}$ & $C O_{ooo} O_{iio}$  \\
$A O_{iio} O_{ioi}$ & $B O_{oii} O_{iio}$ & $C O_{oii} O_{ioi}$  \\
$A A_{iii} A_{ioo}$ & $B A_{iii} A_{oio}$ & $C A_{iii} A_{ooi}$  \\
$A A_{oio} A_{ooi}$ & $B A_{ioo} A_{ooi}$ & $C A_{oio} A_{ioo}$  \\
$A B_{iii} B_{ioo}$ & $B B_{iii} B_{oio}$ & $C B_{iii} B_{ooi}$  \\
$A B_{oio} B_{ooi}$ & $B B_{ioo} B_{ooi}$ & $C B_{oio} B_{ioo}$  \\
$A C_{iii} C_{ioo}$ & $B C_{iii} C_{oio}$ & $C C_{iii} C_{ooi}$  \\
$A C_{oio} C_{ooi}$ & $B C_{ioo} C_{ooi}$ & $C C_{oio} C_{ioo}$  \\
\end{tabular}
\end{center}

\section{Second Stevanovic point(s)?}

I first quote from ETC, and interpolate a $\Bigg[\qquad \Bigg]$

X(1488) = 2nd STEVANOVIC POINT  [\,my SSS, sh'd've been sss,
i.e. $\langle\,333\,\rangle$\,]

Trilinears: $f(A,B,C):f(B,C,A):f(C,A,B)$, where $f(A,B,C)=
1/[1 + \sin(A/2)]$

Barycentrics  $(\sin A)f(A,B,C):(\sin B)f(B,C,A):(\sin C)f(C,A,B)$,
$$\bigg[\mbox{\,i.e.,\ }\bigg\langle\,\ldots,\,
\frac{v(1-v)}{(1+v)(1+v^2)},\,\ldots\,\bigg\rangle$$
The points $A''$, $B''$, and $C''$ in the following
construction, have coordinates
$$\bigg\langle\,0,\frac{v(1-v)}{(1+v)(1+v^2)},
\frac{w(1-w)}{(1+w)(1+w^2)}\,\bigg\rangle,\
\bigg\langle\,\frac{u(1-u)}{(1+u)(1+u^2)},\,0,
\frac{w(1-w)}{(1+w)(1+w^2)}\,\bigg\rangle$$
$$\mbox{and\quad}\bigg\langle\,\frac{u(1-u)}{(1+u)(1+u^2)},\,
\frac{v(1-v)}{(1+v)(1+v^2)},\,0\,\bigg\rangle$$
so that the first of them, for example, lies on the join
of $I_o$ to the incentre of $BI_AC$.  $\bigg]$

Let $U$ be the $A$-excenter of triangle $ABC$; let $A'$
be the incenter of triangle $UBC$, and define $B'$, $C'$
cyclically.  Let $A'' = IA'\cap BC$, and define $B''$, $C''$
cyclically. The lines $AA''$, $BB''$, $CC''$ concur
in X(1488). (Milorad R. Stevanovic, Hyacinthos \#7185, 5/21/03.
See also X(1130) and X(1489).)

X(1488) lies on these lines: 1,166 \  7,2089 \  57,173 \
145,188 \  557,1274 \  558,1143

X(1488) = X(1)-cross conjugate of X(174)

In my desperation to generalize the 2nd Stevanovic point
I tried ``inverting'' his construction in the following
way.  If $A'$ is the incentre of triangle $BIC$, let
$A'I_a$ cut $BC$ in $A''$ and define $B''$, $C''$
cyclically.  Then $AA''$, $BB''$, $CC''$ concur.
Where?  In ETC X(483), our 000 !  Are the radpoints
some sort of inverses of 32 ``2nd Stevanovich points'' ?
One can certainly see the connexion with Steiner's
Malfatti construction, which starts with the incircles
of $BIC, CIA, AIB$.  I must get in touch with my coauthor.

\newpage

\section{Notation at last!}

2016-06-06. This section §7.6, was originally headed
``More on Morley'' but its original content now appears
elsewhere, and has been replaced by recent thoughts on
notation.

To avoid coincidences which do not occur in general, we
assume that our triangles are scalene, that is, neither
right-angled nor isosceles.

The labels of the vertices will be three chosen from the
set \{{\bf 1, 2, 4, 7}\}. The fourth member, the nim-sum
of those chosen, is the label of the orthocentre.  See
``quadration'' at the beginning of the next chapter (Chap. 8).

On those occasions where it is desirable to distinguish
between acute and obtuse triangles, the acute triangle
will be {\bf 124} and the related obtuse triangles will be
{\bf 247, 417, 127}, the vertex with the obtuse angle being
{\bf7} in each case.

After quadration our triangle will have four vertices and
six edges. The midpoints of the six edges {\bf24, 41, 12,
 71, 72, 74} are respectively labelled ${\bf6, 5, 3,
 \bar6, \bar5, \bar3}$ where these last three,
`bar\,6', `bar\,5' and `bar\,3'
are the negatives of 6, 5 and 3. They are the midpoints
of pairs of edges which contain an obtuse angle.

These six midpoints are also the midpoints of the edges
${\bf \bar7\bar1, \bar7\bar2, \bar7\bar4,
\bar2\bar4, \bar4\bar1, \bar1\bar2,}$
of the twin triangles -- see §8.3 below.

The diagonal points of the quadrangle {\bf1247} are the
intersections of the pairs of edges
${\bf17 \& 24, 27 \& 41, 47 \& 12}$
and are denoted by ${\bf D_6, D_5 , D_3}$ respectively.

Similarly, the diagonal points of the twin quadrangle
 ${\bf\bar1\bar2\bar4\bar7}$ are the
intersections of the pairs of edges
${\bf\bar1\bar7 \& \bar2\bar4,
\bar2\bar7 \& \bar4\bar1,
\bar4\bar7 \& \bar1\bar2}$
and are denoted by ${\bf D_{\bar6}, D_{\bar5}, D_{\bar3}}$ respectively.

It will be noted that the midpoint of the join of any point
with its negative is always the same point. We will call this
the {\bf Centre} of the triangle, and denote it by {\bf0}.
In particular, the joins of the six points
${\bf6, 5, 3, D_6, D_5, D_3}$ with their negatives are
all diameters of the same circle, often called the
nine-point circle, but which we will call the
{\bf Central circle} of the triangle.

It has also been called the fifty-point centre, but it is not
advisable to attach a number, since this will be a function
of time and taste. In addition to the six midpoints and the
six diagonal points, there are 32 points of contact with
touch-circles, six points of contact with the double deltoid
and with the Steiner Star of David, and eight vertices of
the four conics which are enveloped by Droz-Farny lines.

\begin{figure}[h]
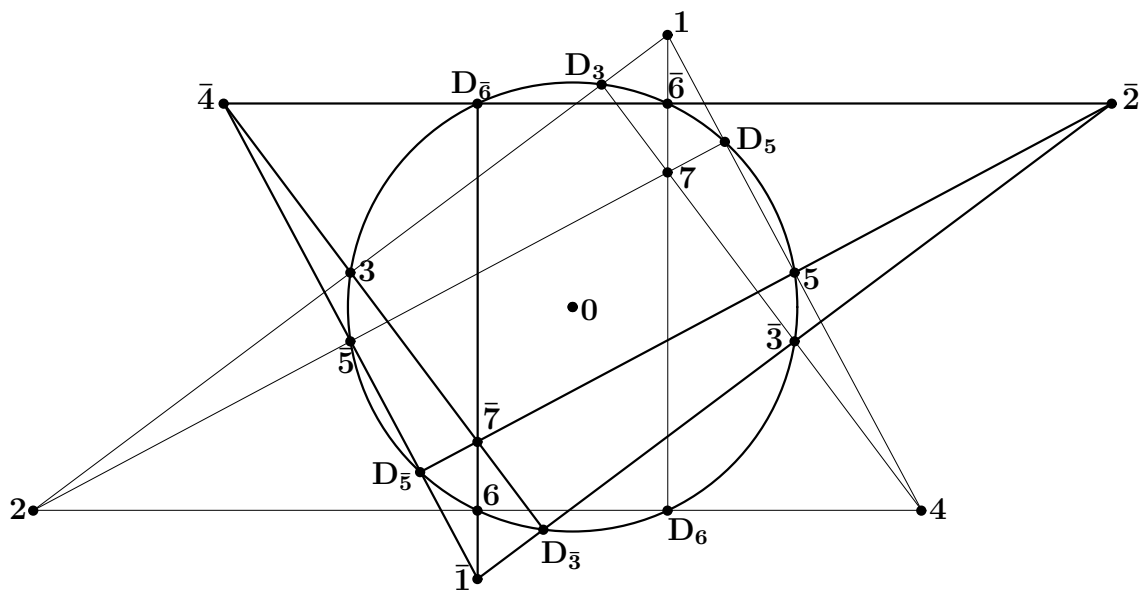
 

\caption{Revised notation}
\label{notation}
\end{figure}

\chapter{Eight Vertices, but only One Centre} 

\section{Quadration}

\centerline{We'll start with something we all know: a triangle}

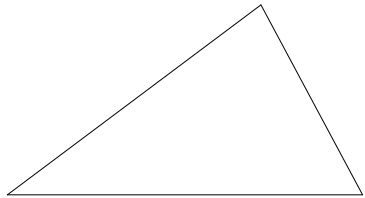
\vspace{4mm}

\begin{figure}[h] 
%
\label{ortho}
\end{figure}

\centerline{has an orthocentre.}
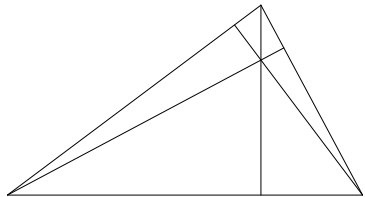
\begin{figure}[h] 
%
\label{perps}
\end{figure}
\begin{center}
That is: the altitudes (perpendiculars from the vertices to the opposite edges) concur.
\end{center}

\centerline{Here are nine proofs of this.}

\bigskip

\begin{center}

1. Euclid first showed that the perpendicular bisectors
of the edges of a triangle concur; then drew a twice-sized
triangle whose perpendicular bisectors were the altitudes
of the original triangle.

\end{center}

\newpage

\begin{figure}[h]
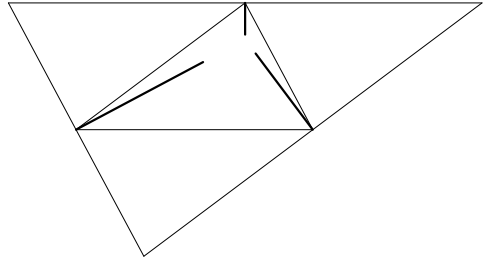
 
%
\caption{Perpendicular bisectors of edges of big
triangle are altitudes of smaller one.}
\label{proof1}
\end{figure}


2. A conic through the four intersections of
two rectangular hyperbolas is also a rectangular hyperbola:
The pair of thin lines and the pair of thick lines are each
perpendicular and form rectangular hyperbolas.  Hence the
pair of dashed lines form a rectangular hyperbola, and are
perpendicular.


\begin{figure}[h]
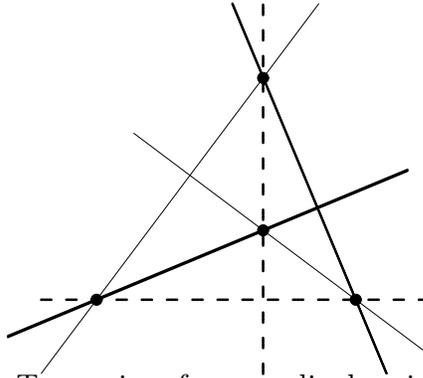
 
%
\caption{Two pairs of perpendiculars induce a third.}
\label{proof2}
\end{figure}


3.  Draw a rectangular hyperbola through the vertices.
Choose axes and scale so that its equation is $xy=1$. The
vertices are $(t,1/t)$, for $t=t_1$, $t=t_2$, $t=t_3$.
The slope of the join of $t_2$ and $t_3$ is $-1/t_2t_3$.
The perpendicular through $t_1$ meets the
hyperbola again at $t=-1/t_1t_2t_3$, whose symmetry
shows that the three perpendiculars concur.

\begin{figure}[h]
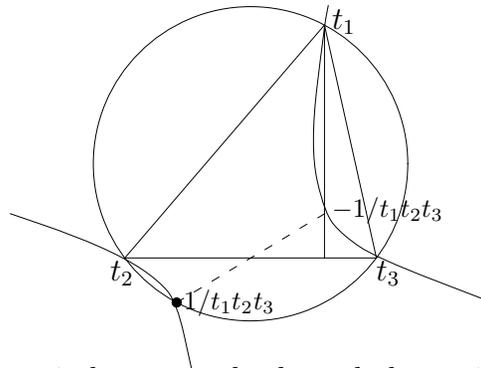
 
%
\caption{Bonus!  The circumcircle meets the hyperbola
again at $t=1/t_1t_2t_3$, the opposite end of
the diameter of the hyperbola through the orthocentre.}
\label{proof3}
\end{figure}

\newpage


\vspace*{40pt}

\begin{figure}[h]  
%
\caption{The four 4-gons formed from lighthouses with $n=4$.}
\label{proof4}
\end{figure}

5. Proof number 5 uses vectors.  Take as origin, {\bf o},
the intersection of the perpendiculars from the vertices
{\bf b} and {\bf c} onto the opposite edges.  That is
{\bf b} is perpendicular to {\bf a-c} and {\bf c} is
perpendicular to {\bf b-a}.

\begin{figure}[h]  
%
\caption{}
\label{proof5}
\end{figure}

\begin{center}
{\bf b$\cdot$(a-c)} = 0 = {\bf c$\cdot$(b-a)} \\

{\bf b$\cdot$a} = {\bf b$\cdot$c} = {\bf c$\cdot$b} = {\bf c$\cdot$a} \\

{\bf (b-c)$\cdot$a} = 0 \\

and {\bf a} is perpendicular to {\bf b-c}.

\end{center}

\newpage

The other four proofs are ``clover-leaf theorems''.
Given three circles, the {\bf radical axes} of pairs
concur in the {\bf radical centre}. What is the radical axis
of a pair of circles?  Not very good definition:
locus of points whence the tangents to the two circles
are equal in length.  Better to define the {\bf power} of
a point w.r.t.\ a circle, that is {\bf the square of the
distance from the point to the centre of the circle minus
the square of the radius of the circle}.  The power is
negative for points inside the circle, positive for points
outside the circle, and zero for points on the circle.
Then the radical axis of two circles is {\bf the locus of
points whose powers w.r.t.\ the two circles are equal.}
If the circles intersect, it is their common chord.  If they
are concentric, it is the line at infinity.

\medskip

6. Draw the three ``edge-circles'', the circles having
the edges of the triangle as diameters.

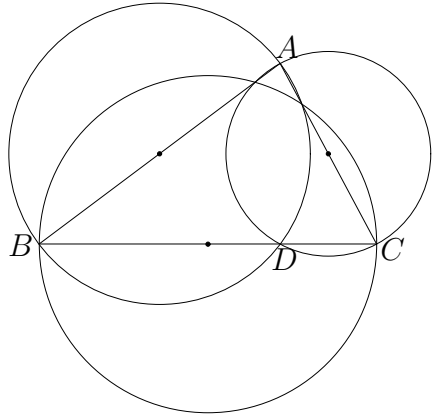
\begin{figure}[h] 
\begin{picture}(440,140)(-190,-82)
\setlength{\unitlength}{0.38pt}
\drawline(-204,-77)(132,-77)(36,103)(-204,-77) 
\put(-84,13){\circle{300}}
\put(84,13){\circle{204}}
\put(-36,-77){\circle{336}}
\put(-84,13){\circle*{4}}
\put(84,13){\circle*{4}}
\put(-36,-77){\circle*{4}}
\put(31,109){$A$}
\put(-235,-89){$B$}
\put(136,-93){$C$}
\put(27,-103){$D$}
\end{picture}
\caption{Edge-circles}
\label{proof6}
\end{figure}

The circles on $AB$, $AC$ as diameters both intersect $BC$
in the point $D$ where $\angle ADB$ = $\angle ADC$ = $\pi/2$
since the angle in a semicircle is a right angle.  $AD$ is
the radical axis of the circles on $AB$, $AC$ as diameters.
The radical axes are the altitudes of the triangle and hence
concur.

\bigskip

7. In Fig.\ref{proof7} the altitudes from $B$ and $C$
meet at $P$.  The circles with $BP$ and $CP$ as diameters
cut the edges $AB$, $AC$ respectively at the feet of the
altitudes from $C$ and $B$, since the angle in a semicircle
is a right angle.  It is clear that the radical axis of the
$A-$ and $B-$circles and that of the $A-$ and $C-$circles
are altitudes of the triangle.  It remains to show that
$A$, $P$, $D$ are collinear

\medskip

Have I wandered up a cul-de-sac?  Does this not lead to a proof?

\bigskip

We shall see, after we've introduced quadration, that these
three diameters may also be considered as edges of the
triangle.

\begin{figure}[h]
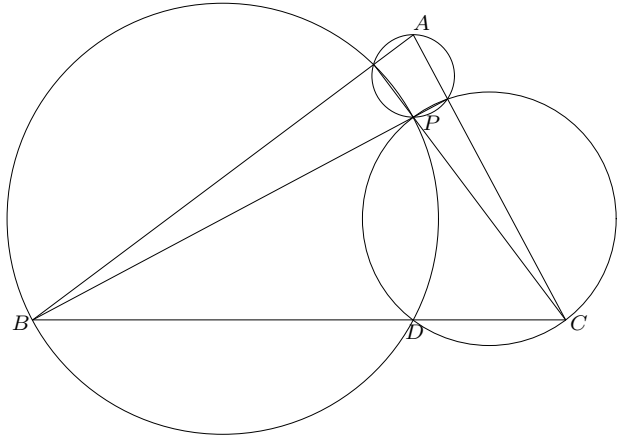
 

\caption{More edge-circles}
\label{proof7}
\end{figure}

\newpage

8.  Draw the reflexions of the circumcircle in the three
edges.  From the dashed lines in Fig.\ref{proof8} it will
be seen that the line of centres of a pair of such circles
is parallel to an edge of the triangle.  The radical axis
will therefore be perpendicular to that edge, and, since
it passes through the opposite vertex, is an altitude of
the triangle.  Since the radical axes concur, so do the
altitudes.

\begin{figure}[h]
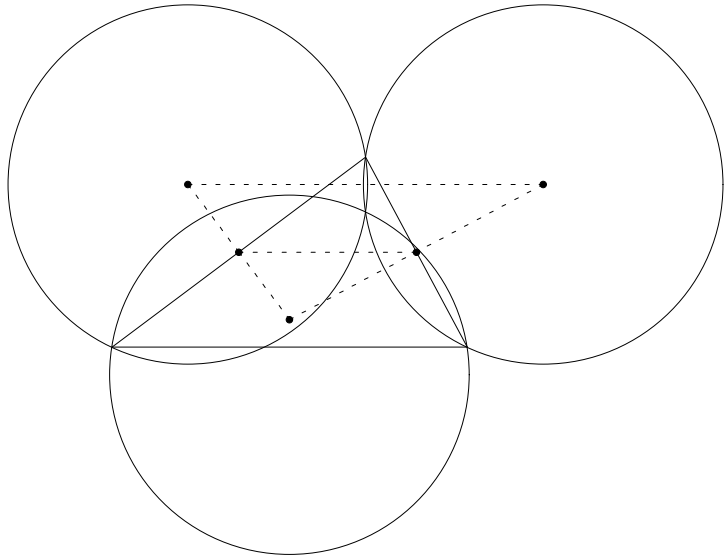
 
%
\caption{Reflected circumcircles}
\label{proof8}
\end{figure}

\newpage

9.  This what I call is THE Clover Leaf Theorem.  Draw the
{\bf medial circles}; that is, the circles having the
medians as diameters.  Here is one:

\begin{figure}[h]
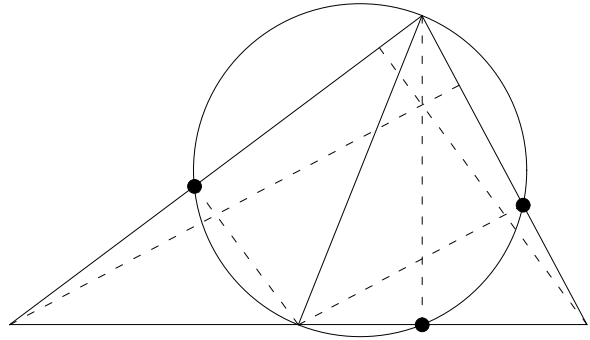
 
%
\caption{A medial circle.}
\label{medial}

\end{figure}

The medial circles pass through some interesting points:

The foot of the altitude, of course; but also what I call
the {\bf midfoot points}, half way between a vertex and
the foot of an altitude.

Now let's find the radical axis of two medial circles.
It must be perpendicular to the line of centres,
which is parallel to an edge.

\begin{figure}[h]
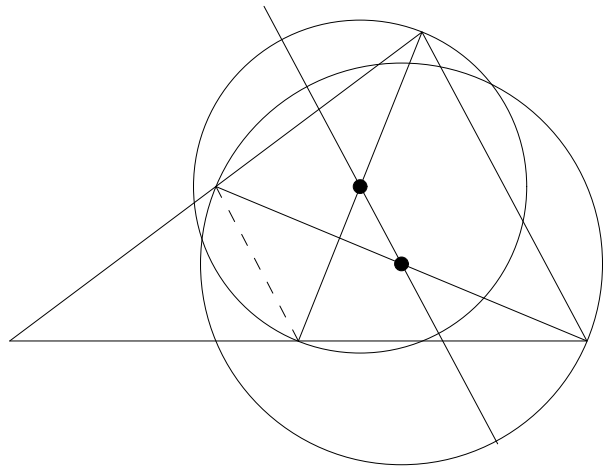
 
%
\caption{Line of centres parallel to an edge.}
\label{perp}
\end{figure}

\clearpage

We'll show that the opposite vertex, $B$ in Fig.\ref{radical},
lies on the radical axis.

\begin{figure}[h]
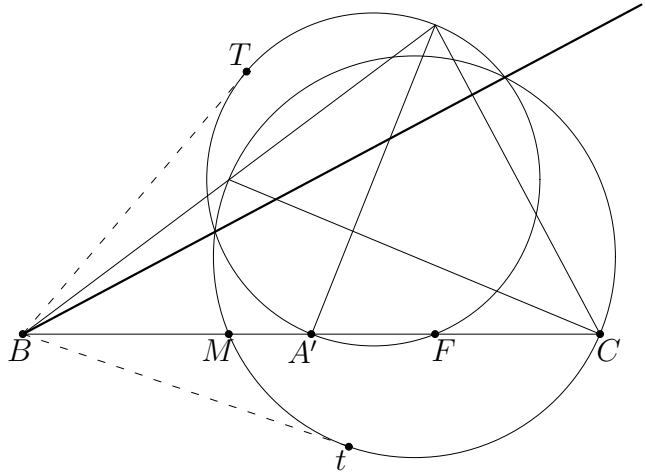
 
%
\caption{$BT^2=BA'\cdot BF=BC\cdot BM=Bt^2$}
\label{radical}
\end{figure}


The radical axes are the {\bf altitudes}, not the medians!
If you like to throw in the nine-point circle, you have the
Four Leaf Clover Theorem.  The $\binom{4}{2}=6$ radical axes
are the six edges of the orthocentric quadrangle.

\begin{figure}[h]
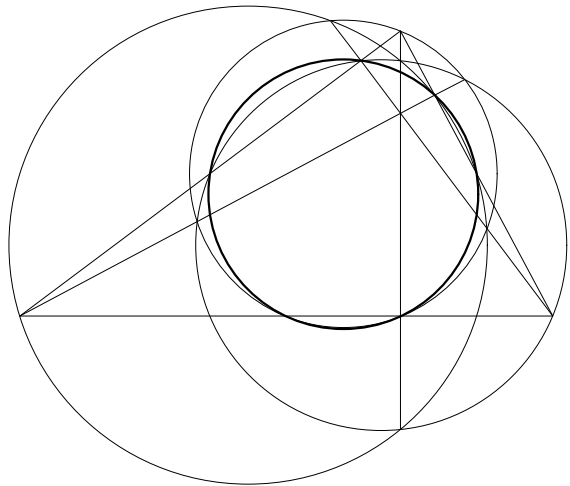
 
%
\caption{The Four Leaf Clover Theorem}
\label{4leaf}
\end{figure}

\newpage

\section{Quadration}


\begin{center}

Here's something else we all know; \\

\medskip
(though some of you may not know that you know!) \\

\bigskip
{\bf The vertices and the orthocentre should have equal status.} \\

\medskip

Each of the four points is the orthocentre of the triangle formed by the
other three! \\

\bigskip
{\Large A triangle is an orthocentric quadrangle!} \\

\bigskip
{\bf A triangle has (at least!) four vertices and six edges.} \\

\bigskip
{\bf Theorem:} There are three times as many obtuse triangles as acute ones. \\

\bigskip

(There are four more proofs in {\it Math.\ Mag.}, {\bf66}(1993) 175--179.) \\

\bigskip
Reminder:  The joins of the midpoints of opposite edges of any quadrangle \\
(not just orthocentric ones) concur and bisect each other.
\end{center}

\section{Circumcentres.  Twinning}

\begin{center}
Let's look for circumcentres. \\

\medskip
Draw the perpendicular bisectors of the edges (there are now six!) \\
\end{center}

\begin{figure}[h] 
\begin{picture}(440,90)(-170,-47)
\setlength{\unitlength}{0.65pt}
\drawline(-204,-77)(132,-77)(36,103)(-204,-77) 
\drawline(-204,-77)(57.59,62.51)               
\drawline(132,-77)(11.04,84.28)                
\drawline(36,-77)(36,103)                      
\thicklines 
\drawline(204,77)(-132,77)(-36,-103)(204,77)   
\drawline(204,77)(-57.59,-62.51)               
\drawline(-132,77)(-11.04,-84.28)              
\drawline(-36,77)(-36,-103)                    
\put(0,0){\circle*{3}}
\put(3,-5){$O$}
\end{picture}
\caption{Hello, twins!}
\label{twins}
\end{figure}

\vspace*{-5pt}

\begin{center}
We get a configuration congruent to the original orthocentric
quadrangle, \\
forming an involution with it, having fixed point $O$. \\

Our triangle now has {\bf eight} vertices. \\

\bigskip
They are {\bf all} orthocentres and {\bf all} circumcentres\,!! \\
\bigskip

There are {\bf twelve edges} which share {\bf six mid-points}.\\
\bigskip
These, with the {\bf six diagonal points}, form three rectangles; \\
whose six diagonals are diameters of a circle, centre $O$. \\

\bigskip
You probably call it the nine-point circle.  \\
\bigskip

I call it the {\bf fifty-point circle.} [NOT ANY MORE! See \S7.6]

\bigskip

Note that $O$ is THE unique centre of the triangle.  \\
\bigskip
Any other candidate would have a twin with equal right.

\end{center}

\begin{figure}[h] 
\begin{picture}(440,120)(-240,-65)
\setlength{\unitlength}{0.7pt}
\drawline(-204,-77)(132,-77)(36,103)(-204,-77)  
\drawline(-204,-77)(57.59,62.51)                
\drawline(132,-77)(11.04,84.28)                 
\drawline(36,-77)(36,103)                       
\drawline(204,77)(-132,77)(-36,-103)(204,77)    
\drawline(204,77)(-57.59,-62.51)                
\drawline(-132,77)(-11.04,-84.28)               
\drawline(-36,77)(-36,-103)                     
\put(0,0){\circle*{5}}
\put(36,-77){\circle*{5}}
\put(-36,77){\circle*{5}}
\put(11.04,84.28){\circle*{5}}
\put(-11.04,-84.28){\circle*{5}}
\put(57.59,62.51){\circle*{5}}
\put(-57.59,-62.51){\circle*{5}}
\put(-36,-77){\circle*{5}}
\put(36,77){\circle*{5}}
\put(84,13){\circle*{5}}
\put(-84,-13){\circle*{5}}
\put(-84,13){\circle*{5}}
\put(84,-13){\circle*{5}}
\put(3,-3){$O$}

\thicklines
\put(0,0){\circle{170}}
\end{picture}
\caption{Twelve of the 50 points}
\label{12point}
\end{figure}
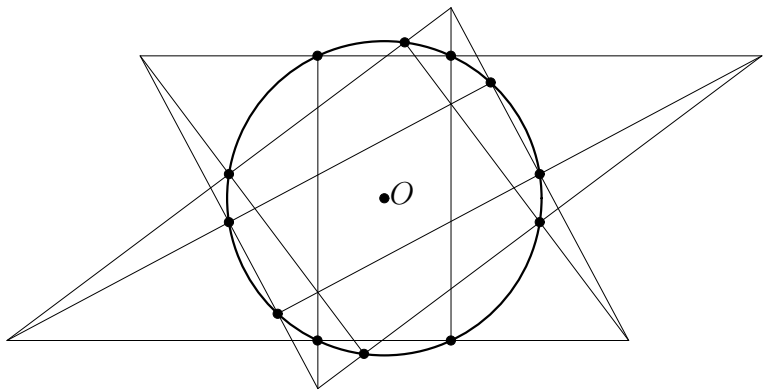

\medskip

\begin{center}

Notice that, although we now have eight circumcircles, they all have
the same radius --- twice that of the fifty-point circle. \\

\medskip

In fact the circumcircles form pairs of reflexions in the twelve edges: \\

\medskip

2 reflexions $\times$ 12 edges \big/ 3 edges in a triangle = 8 circumcircles. \\

\bigskip

For those of you who prefer to do things analytically, \\
here's how quadration appears: \\

\medskip

A triangle with angles $A$, $B$, $C$ and \\
edge-lengths $2R\sin A$, $2R\sin B$, $2R\sin C$, \\
quadrates into three more triangles: \\
\bigskip
\begin{tabular}{cccccc}
$\pi-A$ & $\pi/2-C$ & $\pi/2-B$ & $2R\sin A$ & $2R\cos C$ & $2R\cos B$ \\
$\pi/2-C$ & $\pi-B$ & $\pi/2-A$ & $2R\cos C$ & $2R\sin B$ & $2R\cos A$ \\
$\pi/2-B$ & $\pi/2-A$ & $\pi-C$ & $2R\cos B$ & $2R\cos A$ & $2R\sin C$ \\
\end{tabular}
\end{center}

\section{Touch Circles}

\begin{center}
You will all be familiar with the incircle and three excircles of
a triangle. \\ Their centres, the touch-centres, are where the
angle-bisectors concur:
\end{center}

\begin{figure}[h] 
%
\caption{The incentre and three excentres form an
orthocentric quadrangle whose 50-point centre is the
circumcentre of the original triangle.}
\label{touch}
\end{figure}

\begin{figure}[h] 
%
\caption{32 touch-circles each touch 3 of the 12 edges at
one of 8 points, and each touches the 50-point circle
(which isn't drawn, but you can ``see'' it!)}
\label{32touch}
\end{figure}

\clearpage

\begin{figure}[h]
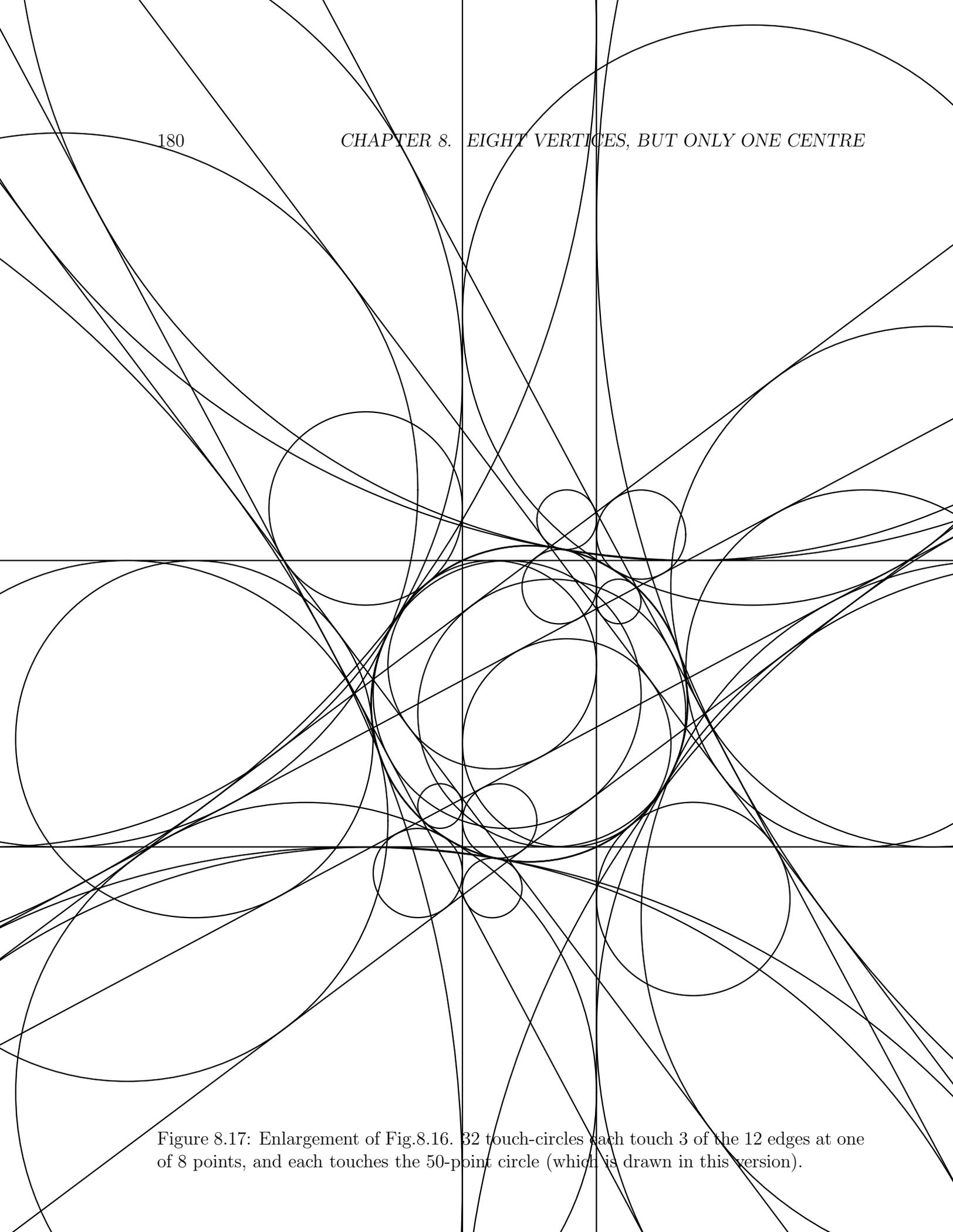
 
%
\caption{Enlargement of Fig.8.16.  32 touch-circles each
touch 3 of the 12 edges at one of 8 points, and each
touches the 50-point circle (which is drawn in this version).}
\label{enlarge}
\end{figure}

\clearpage

\begin{center}

We have 8 triangles, and hence 32 touch-circles.  They all
appear to touch the 50-point circle.  This is Feuerbach's
theorem, though published by \\

\medskip

C.-J.~Brianchon \& J.-V.~Poncelet in Recherches
sur la d\'etermination d'une hyperbole \'equilat\`ere, au
moyen de quatre conditions donn\'ees, {\it Ann.\ des Math.},
{\bf11}(1821) 205--220. \\

\bigskip

As an undergraduate I learned Feuerbach's theorem as the 11-point
conic, the locus of the poles of a line with respect to a
4-point system of conics, but that is another story. \\

\medskip

If the 4 points form an orthocentric quadrangle, the conics
are rectangular hyperbolas; \\
if the line is the line at infinity, the poles are their
centres, which lie on the 9-point circle, \\
the other two points being the circular points at infinity.

\bigskip

I wanted to find a proof of Feuerbach's theorem that was fit
for human consumption.  The best that I've been able to do is
a sort of parody of Conway's {\bf extraversion}. \\

\bigskip

\end{center}

\section{Hexaflexing}

\begin{figure}[h]
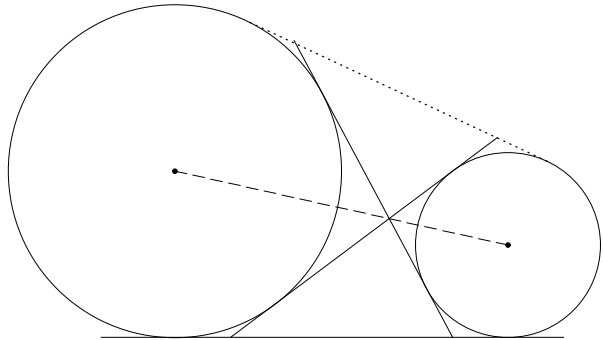
 

\label{common}
\caption{Triangle reflected in angle-bisector}
\end{figure}

\begin{center}

Take one of the $\binom{4}{2}=6$ pairs of touch-circles.
The three triangle edges are common tangents.  Reflect the
triangle in the corresponding angle-bisector.  This
adds the fourth common tangent, dotted in figure.

\end{center}

\newpage

\begin{figure}[h] 
%
\caption{Three pairs of parallel common tangents to the touch-circles}
\label{hexaflex}
\end{figure}

\vspace*{-20pt}

\begin{center}

The 3 points on each of the 4 circles form triangles homothetic
to the original triangle:

\end{center}

\begin{figure}[h]
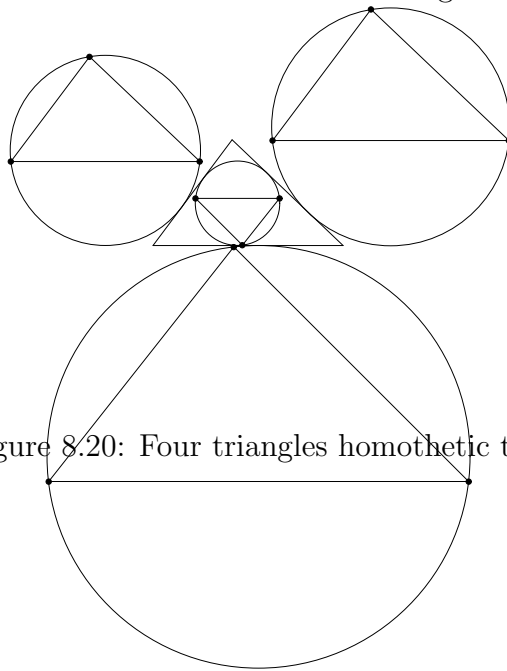
 
%
\caption{Four triangles homothetic to the original triangle}
\label{homothets}
\end{figure}

\clearpage

\ldots and also homothetic to the {\bf medial triangle} (dashed in
the following figure).  Where is (are) the perspector(s)
[centre(s) of perspective]\,?
\begin{figure}[h]
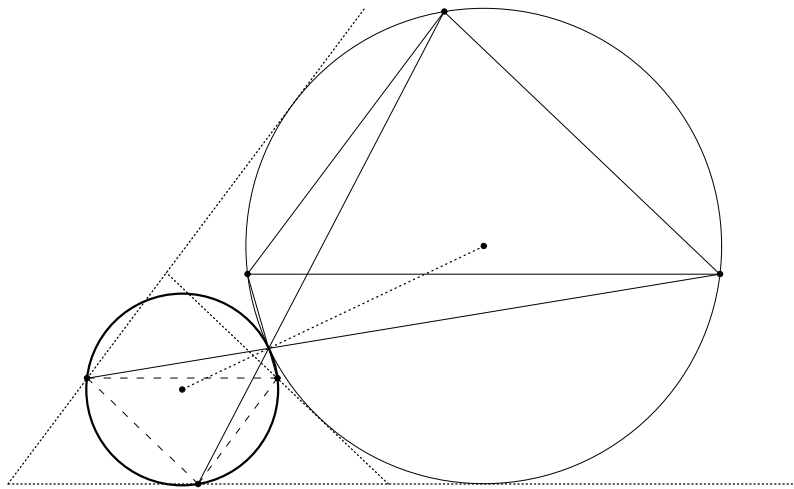
 
%
\caption{A triangle homothetic to the medial triangle}
\label{medial}
\end{figure}

  They lie on the circumcircles of
the homothets, and on the circumcircle of the medial triangle,
which is, of course, the 50-point circle.  That is, all the 32
touch-circles touch the 50-point circle.

\section{Gergonne points}

Join the vertices of a triangle to the touch-points of one
the touch-circles to the opposite edges of the triangle.
By Ceva's theorem, the joins concur.

\begin{figure}[h]
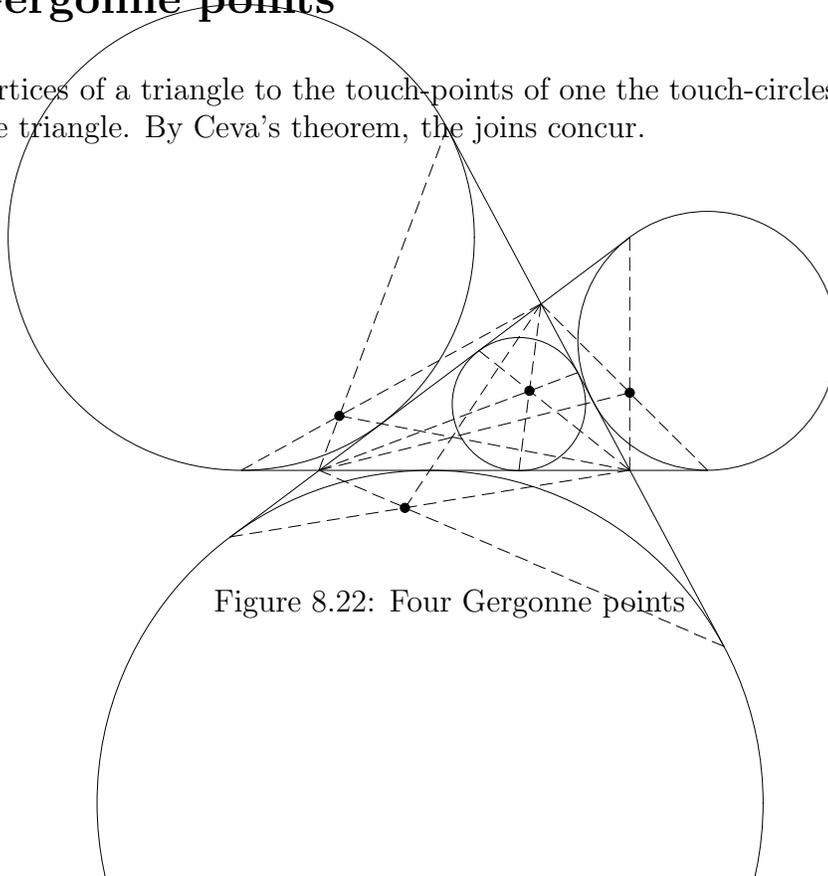
 
%
\label{gergs}
\caption{Four Gergonne points}
\end{figure}

\newpage
\begin{center}
If we join the four Gergonne points to their corresponding
touch-centres, the four segments concur in a \\

\smallskip

de LONGCHAMPS POINT.

\end{center}

\begin{figure}[h]
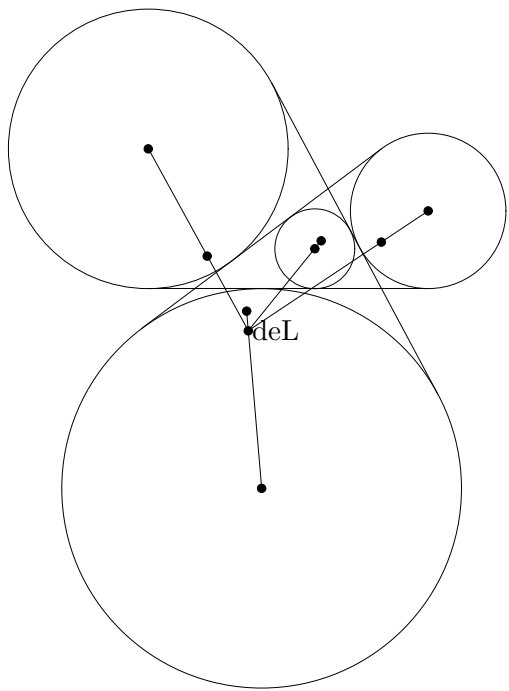
 
%
\caption{Extraversion of Gergonne points in a triangle.
The Gergonne point, touch-centre and deLongchamps point are collinear.}
\label{delong}
\end{figure}


A de Longchamps point lies on the corresponding {\bf Euler line}

\begin{figure}[h]
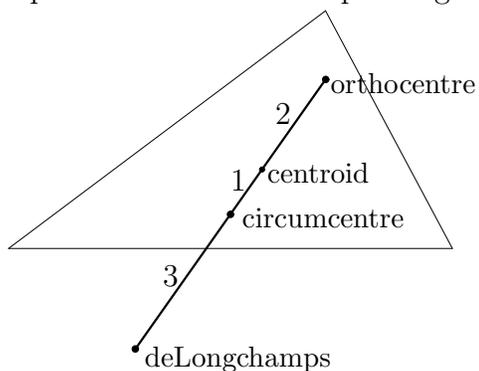
 
%
\caption{A de Longchamps point together with the corresponding
centroid separate the circumcentre and orthocentre harmonically.}
\label{euler}
\end{figure}

\newpage

Here are 32 Gergonne points in one diagram.
Lines connect them to corresponding touch-centres
(unmarked in figure) and concur in fours at 8 deLongchamps
points.

\begin{figure}[h]
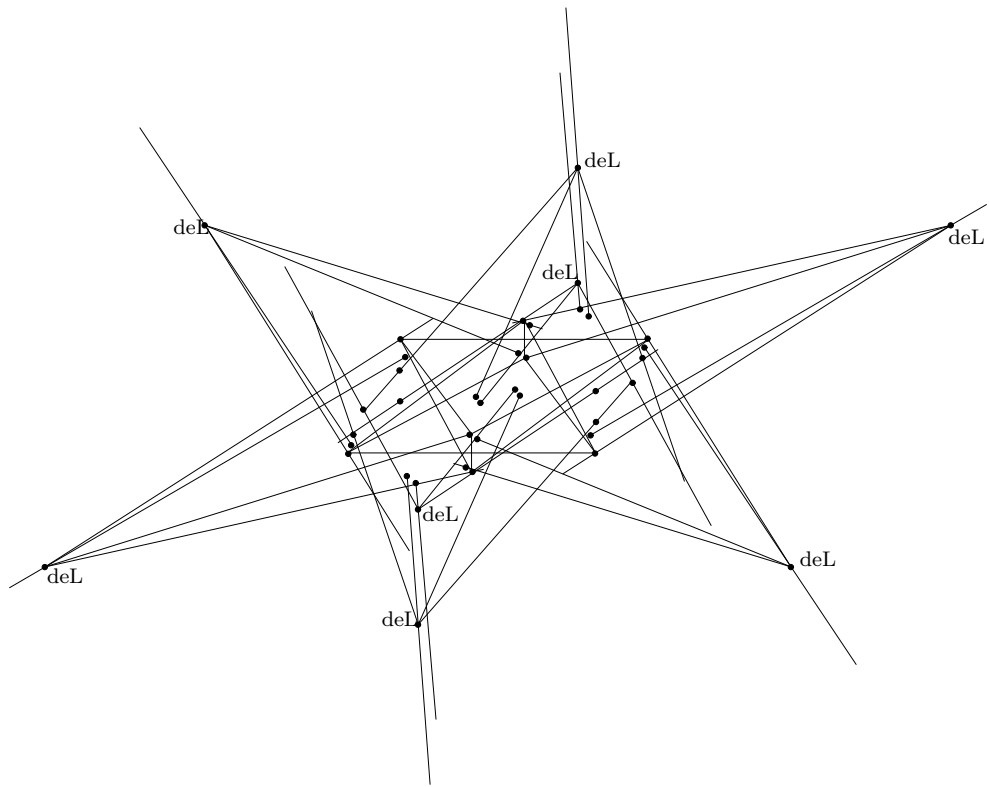
 
%
\caption{Thirty-two Gergonne points.}
\label{32gerg}
\end{figure}

\begin{center}

The 8 de Longchamps points and 8 centroids separate
the 8 circumcentres and \\
8 orthocentres in 8 harmonic
ranges which lie in pairs on

\end{center}

\begin{figure}[h] 
%
\caption{4 EULER LINES}
\label{4euler}
\end{figure}

\newpage

\section{Six more points to make up the fifty.}

\begin{center}

We will find them as the midpoints of the edges of the BEAT twins. \\

The Best Equilateral Approximating Triangles. \\

We start from the Simson Line Theorem. This doesn't appear
in Simson's work, but is due to Wallace, so we'll use his
name instead.

If, from a point on the circumcircle of a triangle, we drop
perpendiculars onto the edges, then their feet are collinear
(the Wallace line).

\end{center}

\begin{figure}[h]
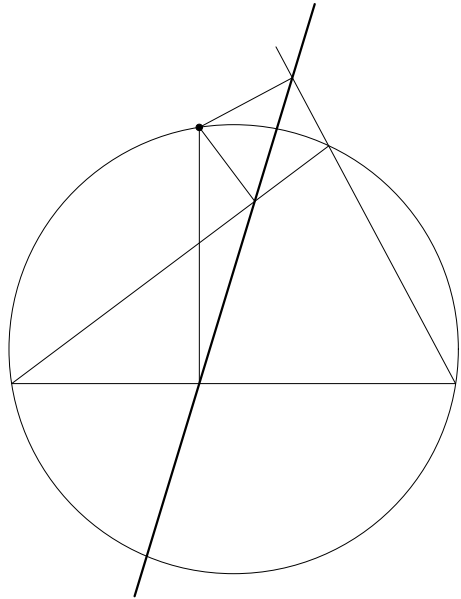
 
%
\caption{The Simson-Wallace line theorem.}
\label{wallace}
\end{figure}

\newpage

\begin{center}

As the point moves round the circumcircle, the Wallace line
rotates in the opposite sense with half the angular
velocity.

\end{center}

\begin{figure}[ht]  
%
\caption{Angle between Wallace lines of $S$ and $S'$ = $\angle UV_1U'$ =
${1\over2}\mbox{arc}UU'={1\over2}\mbox{arc}S'S$.}
\label{angle}
\end{figure}

\newpage

The Wallace line bisects the segment joining the orthocentre
to the point on the circumcircle.

\begin{figure}[ht]
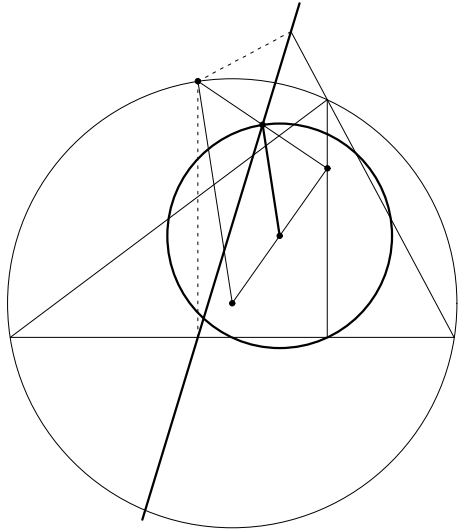
  
%
\caption{The Simson-Wallace line bisects the segment $S_8V_8$}
\label{midpoint}
\end{figure}

This bisection point is on the 50-point circle.

The Wallace line is a diameter of a circle of radius $R$
which encloses the 50-point circle and rolls inside
the concentric {\bf Steiner circle}, radius $\frac{3}{2}R$.

\begin{figure}[ht]
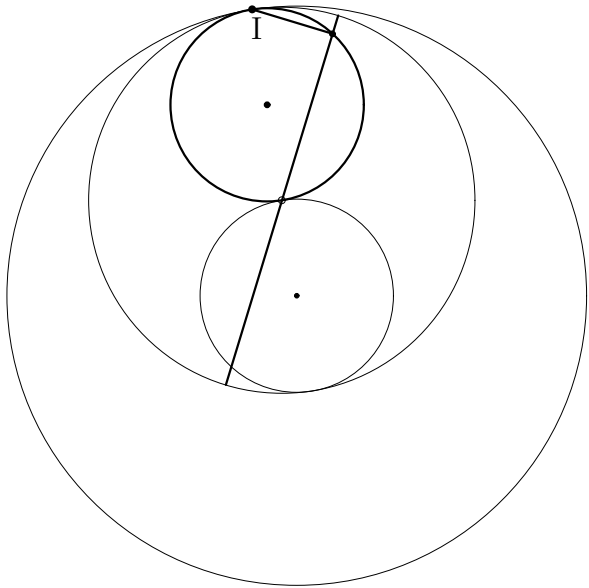
  
%
\caption{The deltoid is a locus and an envelope.}
\label{cycloids}
\end{figure}

\newpage

\ldots and, as the point moves round the circumcircle,
the Wallace line envelops the

\medskip

\centerline{\large STEINER DELTOID}

\begin{figure}[ht]
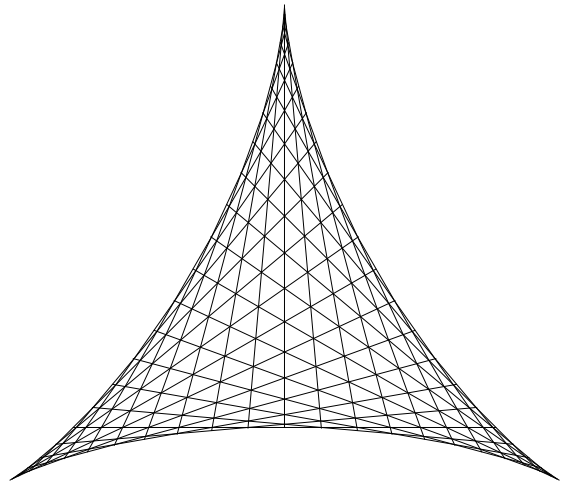
  
%
\caption{The deltoid as envelope.}
\label{deltoid}
\end{figure}

\newpage

\begin{center}

We've found our 3-symmetry! \\

\medskip

Let's try Quadration. \\

\medskip

Three perpendiculars from points on
each of 4 circumcircles coincide
in pairs to give 6 points on the Wallace line. \\

\end{center}

\begin{figure}[h]
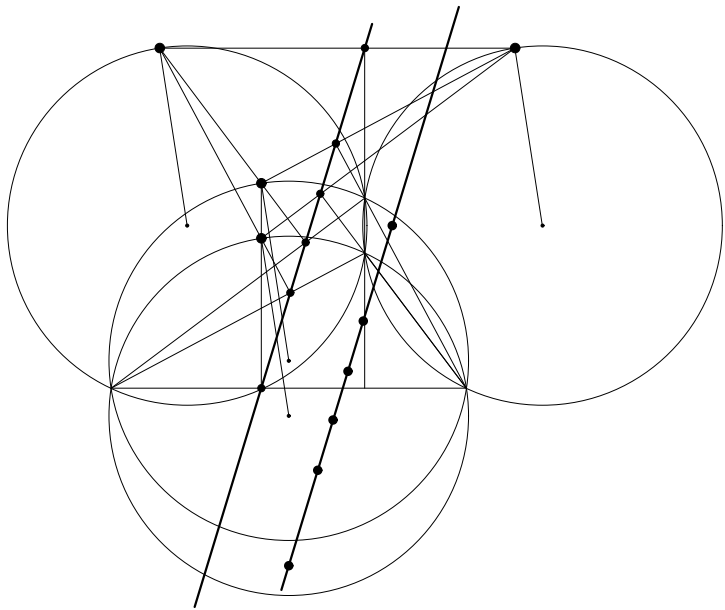
 
%
\caption{Quadration of the Simson-Wallace line theorem.}
\label{quadwall}
\end{figure}

\newpage

\begin{center}

Now let's try Twinning. \\

\medskip

Twinning gives 6 points on a parallel Wallace line,
the reflexion of the first in the 50-point centre. \\

\end{center}

\begin{figure}[h]
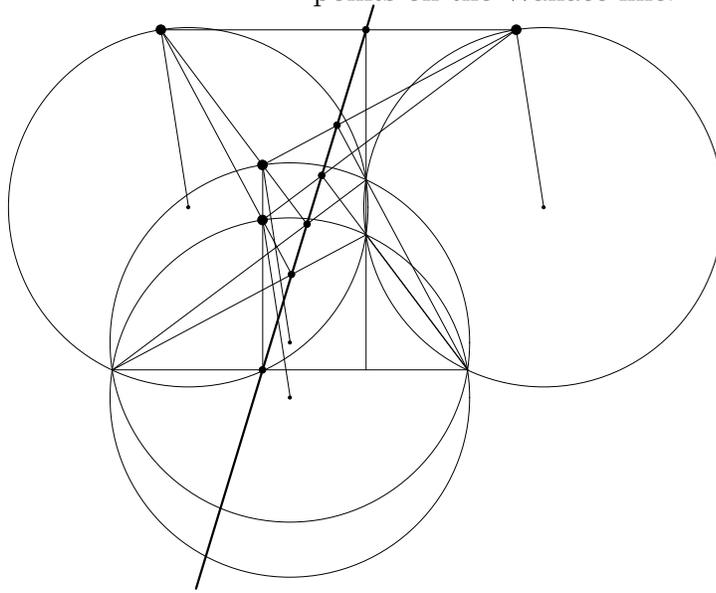
 
%
\caption{Quadration of the Simson-Wallace line theorem.}
\label{quadwall}
\end{figure}

\begin{figure}[ht]
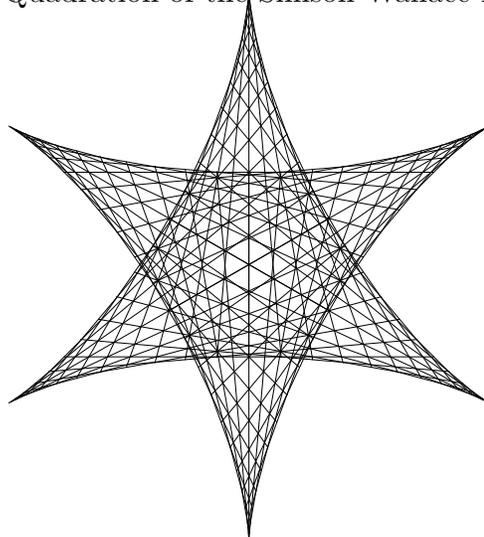
  
%
\caption{These parallels envelop a double-deltoid.}
\label{doubdelt}
\end{figure}

\newpage

There are 3 occasions when the parallels
coincide, with 12 feet on the same Wallace line.

\begin{figure}[h] 
%
\caption{Three-Symmetry}
\label{3dir}
\end{figure}

\vspace*{-30pt}


\begin{center}

\ldots and three occasions when the parallels are
perpendicular to these; the edges of the
BEAT twins, forming the Steiner Star of David.

\end{center}

\begin{figure}[h]
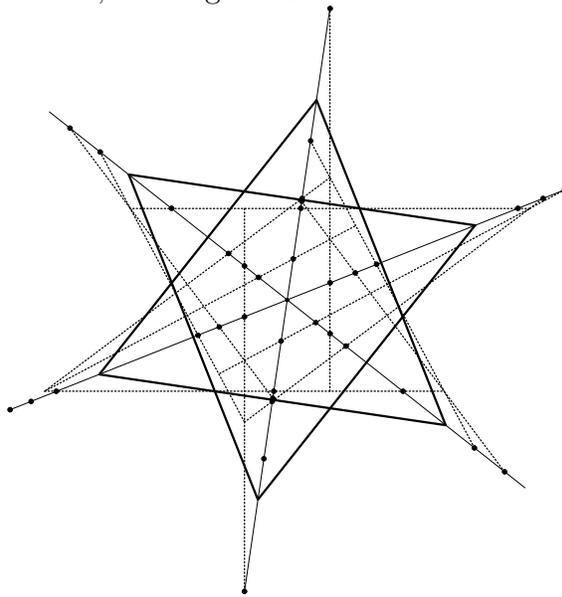
 
%
\caption{The Steiner Star of David}
\label{starofd}
\end{figure}

\newpage

The axes of symmetry are the trisectors
of the angles between the diameters through the midpoints
and diagonal points, nearer to the former.

\begin{figure}[h]
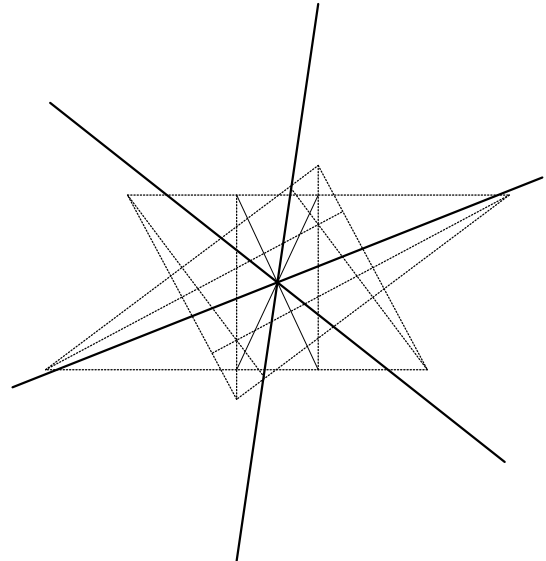
 
%
\caption{Trisector axes of symmetry}
\label{trisect}
\end{figure}

If the angles of the original triangle are $A$, $B$, $C$,
then the angles between the edges of the Steiner Star
and those of the triangle and its twin are
$$|B-C|/3, \qquad |C-A|/3, \qquad |A-B|/3$$

These edges are parallel to those of the 
Morley triangles.

\begin{figure}[h]
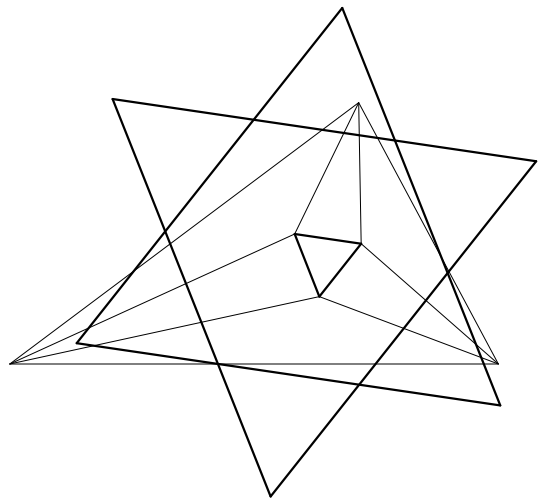
  
%
\caption{Morley's theorem and the Star of David.}
\label{morley}
\end{figure}

\begin{center}

\ldots of which there are $18\times8=144$.  Here are 18 : --

\end{center}

\newpage

\begin{figure}[h]
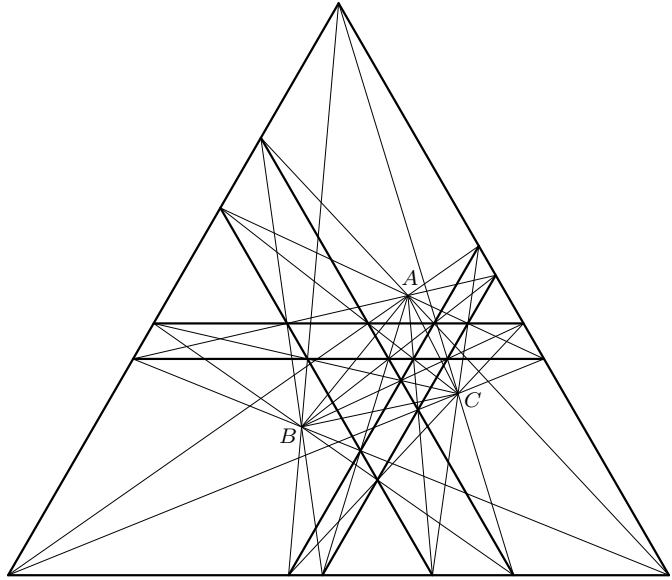
  
%
\caption{Eighteen Morley triangles.}
\label{18morley}S
\end{figure}

\begin{center}

In fact there are 27 equilateral triangles. \\
\bigskip
Eighteen are genuine Morley triangles. \\
\bigskip
Conway has given the name
Guy Faux triangle to the other nine. \\
\bigskip
For details see [{\it AMM}, 114(2007) 97--141] \\
\bigskip
 which
describes, amongst other things, the Lighthouse Theorem for $n=3$.

\end{center}

\chapter{Notes on the Droz-Farny Theorem}

\section{Preliminary Remarks}

The Droz-Farny theorem certainly deserves to reach a
wider audience.  So also should its several
generalizations.  Many of the references \cite{Bra4,Eh,Goor,KB} give
generalizations (and perhaps others \cite{Poh,PT}).

\section{Statement of theorem}

If a pair of perpendicular lines through $H$, the orthocentre
of a triangle $ABC$, intersect the edge $BC$ (respectively
$CA$ and $AB$) in $X_1$, $X_2$ (respectively $Y_1$, $Y_2$ and
$Z_1$, $Z_2$), then the midpoints of $X_1X_2$, $Y_1Y_2$,
$Z_1$, $Z_2$ are collinear.

Additionally, though it may not have been in Droz-Farny's
original statement:

The pair of perpendicular lines, the three edges of the triangle,
and the Droz-Farny line are all tangent to a parabola.

Note that if two tangents to a parabola are perpendicular,
then they intersect on the directrix; moreover, the join
of the points of contact passes through the focus.

\section{Generalizations}

1. Beginnings of a projective version:

The pair of perpendicular lines may be replaced by any
rectangular hyperbola having the orthocentre as centre.

To see this, take the pair of perpendicular lines,
or asymptotes to the hyperbola, as coordinate axes.
Then the equation to any edge of a triangle may be
written in intercept form,
$$\frac{x}{a}\ +\ \frac{y}{b}\ =\ 1$$
and the equation to the rectangular hyperbola as
$xy\ =\ c^2$. Then the midpoint of the segment between
the asymptotes (axes) is $(a/2,b/2)$ and the
midpoint of the chord of the hyperbola is given by
$$\frac{x}{a}\ +\ \frac{c^2}{bx}\ =\ 1$$
$$bx^2\ +\ ac^2\ =\ abx$$
the sum of whose roots is $ab/b$, with average $a/2$
so that the midpoint of the chord is (also) $(a/2,b/2)$
and the result is independent of the signs of $a$, $b$
and $c^2$.

See (also?) \cite{Eh}

2. Semi-affine version:

In place of the midpoints of the three segments one can take
any other ratio. \cite{Lam2} See also \cite{KB}.

It is of interest to investigate the family of pseudo
D-F lines that one gets by considering all ratios. In fact
they evidently envelop a parabola inscribed in the triangle
(the edges of the triangle and the pair of perpendicular
lines are all tangents to the parabola).  Then one can ask
about the family of such parabolas generated by varying the
pair of perpendicular lines through the orthocentre.
Will every parabola inscribed in the triangle correspond
to a suitable pair of perpendicular lines?  Probably ``yes''.
I note in passing that a pair of perpendicular tangents to
a parabola meet on the directrix, so the directrix passes
through the orthocentre.  Where are the focus, axis and vertex?

See (empty!) \S9.7 and Appendix (\S9.12.)  [this all written in
bits and cobbled together, so there's repetition and possibly
contradiction!]

What is the locus of the foci of the parabolas?  The envelope
of the directrices is the orthocentre.  What is the locus of the
vertices?  The envelope of the axes?

\section{Envelope of D-F line}

If one varies the pair of perpendicular lines, the resulting
D-F lines envelope a conic.  This degenerates to a point
(the right angle) in the case of a right-angled triangle.
Otherwise it is an ellipse or a hyperbola, according as the
triangle is acute-angled or obtuse-angled (othocentre inside
or outside the circumcircle).  Foci: orthocentre, $H$, and
circumcentre, $O$.  [That, \&/or something later, is wrong!!]
Centre: nine-point centre.  Major axis:
Euler line.  Vertices: intersections of Euler line with
nine-point circle.  Length of major axis: $OH$.  Length
of conjugate axis: $\sqrt{R^2-OH^2}$, where $R$ is the
circumradius.  Asymptotes:  diameters of nine-point circle
through the points of tangency of the tangents to the
nine-point circle from the orthocentre (imaginary if triangle
is acute, when orthocentre is inside nine-point circle).

What relations are there, if any, between this conic and the
several other conics associated with a triangle?

\section{Miscellaneous ramblings}

I want to write something about the
Droz-Farny theorem, and could make use of someone to
try things out on.  It seems to me that Theorems 131
and 133 in Durell's Projective Geometry are Feuerbach's
theorem and its dual, which is a projective generalization
of the Droz-Farny theorem.

               Here's a start.  Feuerbach's theorem is
the `11-pt conic' theorem.  Given a line and a 4-point
pencil of conics, then the locus of the pole of the line
is a conic.  This passes through 11 points:

               The 3 diagonal points of the quadrangle. \\
               The 6 harmonic conjugates with respect to
pairs of points of the quadrangle of the intersection
of the line with the edge. \\
               The 2 double points of the involution cut
on the line by the conics.

               Special cases: \\
               (a) Line is line at infinity. Conjugates
are the midpoints of the edges. \\
               (b) Quadrangle is orthocentric.  Conics
are rectangular hyperbolas. \\
               (c) (a) \& (b) together.  11-pt conic is
9-pt circle + circular points at infinity.

               The dual of Feuerbach's theorem concerns
a point and a range of conics touching 4 fixed lines.
The envelope of the polar of the point with respect to
the conics is the 11-line conic.  The 11 lines are

               The 6 harmonic conjugates with respect
to pairs of edges of the quadrilateral of the join of
the point to the intersection of the pair of edges. \\
               The 3 edges of the diagonal line triangle
of the quadrilateral. \\
               The tangents at the point to the two
conics which pass through the point.

               What's the connexion with Droz-Farny? I
hear you cry.  Take one of the four lines to be the
line at infinity.  The conics are now parabolas.  The
two parabolas which pass through the given point have
their foci on the circumcircle of the triangle formed
by the three lines not at infinity, and their
directrices passing through the orthocentre.  This
last is Steiner's theorem.  Durell {\it A Concise
Geometrical Conics}, p.22.  Euler line ???  Simson-
Wallace line.

               We're now getting close to the Droz-Farny
theorem.  Compare Cosmin Pohoata's Theorem B:  Let $P$
be a point in the plane of a given triangle $ABC$.  If
$A'$, $B'$, $C'$ are the points where the reflexions of lines
$PA$, $PB$, $PC$ in a given line through $P$ meet the edges
$BC$, $CA$, $AB$, then $A'$, $B'$, $C'$ are collinear.

\newpage

Let's try again, with a picture.

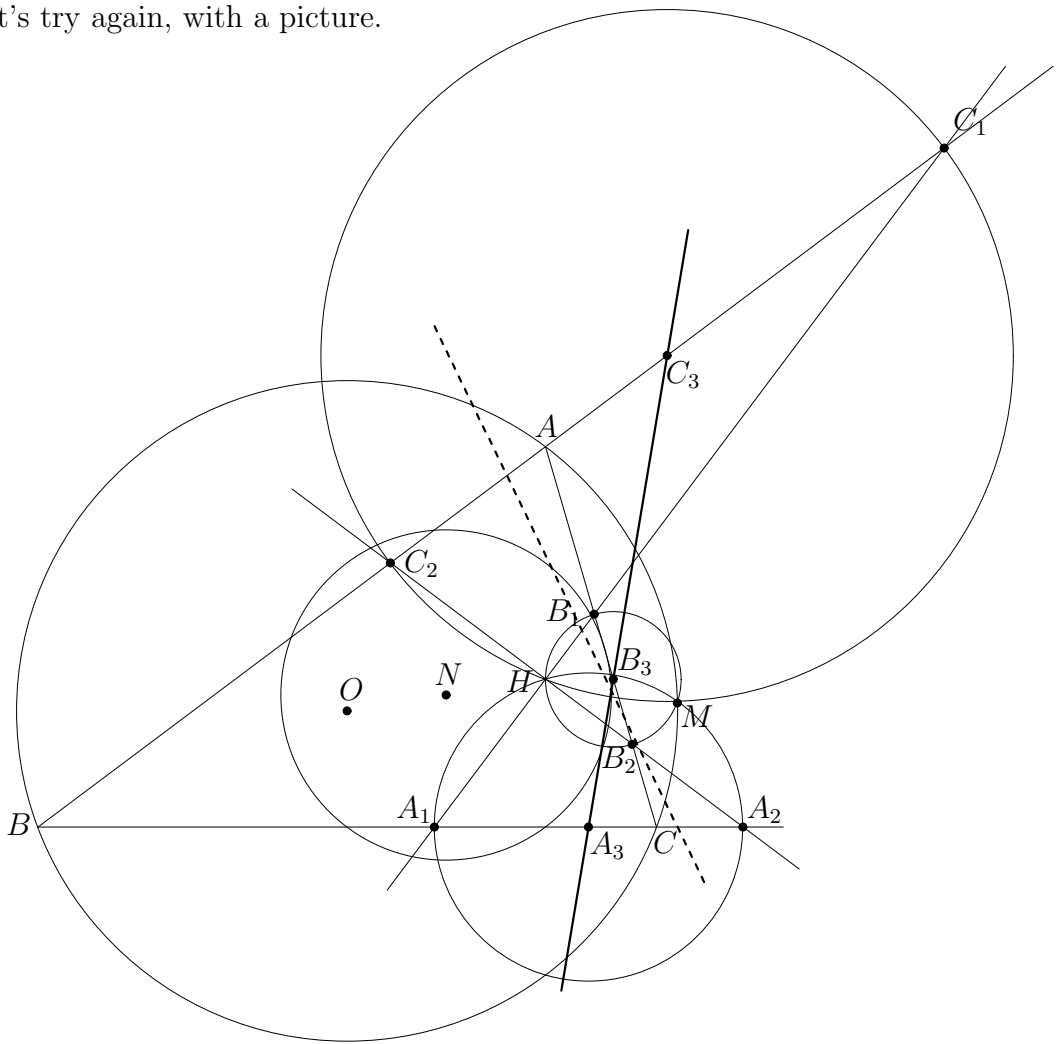
\begin{figure}[h]  
\begin{picture}(400,360)(-220,-70)
\setlength{\unitlength}{6pt}
\drawline(-32,0)(32,48)     
\put(-0.7,24.6){$A$}
\put(-34,-0.5){$B$}
\put(6.8,-1.6){$C$}
\drawline(-32,0)(15,0)     
\drawline(0,24)(7,0)       
\put(-7,0){\circle*{0.5}}     
\put(-9.4,0.6){$A_1$}
\put(25.1428,42.8571){\circle*{0.5}}  
\put(25.7,44){$C_1$}
\put(3.08,13.44){\circle*{0.5}}   
\put(0,13){$B_1$}
\put(12.4444,0){\circle*{0.5}}   
\put(12.7,0.6){$A_2$}
\put(-9.7778,16.6667){\circle*{0.5}}  
\put(-8.9,16.1){$C_2$}
\put(5.4756,5.2267){\circle*{0.5}}  
\put(3.5,3.7){$B_2$}
\put(2.7222,0){\circle*{0.5}}  
\put(2.8,-1.7){$A_3$}
\put(7.6825,29.7619){\circle*{0.5}}   
\put(7.6,28){$C_3$}
\put(4.2778,9.3333){\circle*{0.5}}  
\put(4.5,9.9){$B_3$}

\put(-2.5,8.5){$H$}   
\put(-12.5,7.3333){\circle{41.6667}}   
\put(-12.5,7.3333){\circle*{0.5}}      
\put(-13,8){$O$}
\put(-6.25,8.3333){\circle{20.8333}}   
\put(-6.25,8.3333){\circle*{0.5}}      
\put(-7,9){$N$}
\put(2.7222,0){\circle{19.4444}}       
\put(7.6825,29.7619){\circle{43.6507}} 
\put(4.2778,9.3333){\circle{8.5554}}   
\put(8.3243,7.8378){\circle*{0.5}}     
\put(8.5,6.3){$M$}
\drawline(-10,-4)(29,48)           
\drawline(16,-2.6667)(-16,21.3333) 
\thicklines
\drawline(1,-10.3333)(9,37.6667)    
\dashline[80]{0.5}(10,-3.4565)(-7,31.61056) 

\end{picture}
\caption{Droz-Farny (thickline) and Wallace-Simson (dashline) lines.}
\label{DFWS}
\end{figure}

\vspace*{-5pt}

{\bf Query:}  What are the relative rates of rotation
for the \\

\vspace*{-15pt}

1.  Pair of perpendicular lines through $H$,

2.  Point $M$ on the circumcircle,

3.  Droz-Farny line,

4.  Wallace-Simson line of $M$ ??

In connexion with the first of these four, note that
the pair of perpendicular lines may be replaced by a
rectangular hyperbola having the orthocentre as centre.
Now consider rectangular hyperbolas which are close to
the perpendiculars (which are their asymptotes).  They
may lie in the first and third quadrants formed by the
perpendicular lines or in the second and fourth.

\newpage

One way of describing Fig. \ref{DFWS} is as follows:

Choose a point $A_3$ on edge $BC$ and draw the circle
centre $A_3$ passing through the orthocentre $H$.  If
this cuts $BC$ in $A_1$ and $A_2$, then $A_3$ is the
midpoint of $A_1A_2$ and $\angle A_1HA_2$ is a right
angle (angle in a semicircle).  Let $M$ be the point
of intersection of this circle with the circumcircle
of $ABC$.  There are two choices for this --- take the one
on the same side of edge $BC$ as $A$.  [Investigate
what happens if you make the other choice.]

Let the line $A_1H$ meet the edges $CA$, $AB$ in $B_1$
and $C_1$, and the line $A_2H$ meet the edges $CA$, $AB$
in $B_2$ and $C_2$.

If you can see that circles $B_1HB_2$ and $C_1HC_2$
also pass through $M$, then we have a proof of the
Droz-Farny theorem, since the three circles each pass
through both $H$ and $M$, and so form a coaxal system
with $HM$ as radical axis and the line of centres is
the perpendicular bisector of $HM$, passing through
the midpoints $A_3$, $B_3$, $C_3$ of $A_1A_2$,
$B_1B_2$, $C_1C_2$ --- i.e., the Droz-Farny line.

Note that the Wallace line of $M$ also passes through
the midpoint of $HM$.

It is now easy to find the envelope of the Droz-Farny
line as the pair of perpendiculars through $H$ rotate.

\smallskip

$O$ is the circumcentre. Let the circumradius $OM$
cut the Droz-Farny line at $P$, so that $HP=PM$ and
$OP+PH=OP+PM=R$, the circumradius.  Moreover the
Droz-Farny line is equally inclined to $HP$ and $OP$
so that it is a tangent to the conic with foci
$O$ and $H$ and axis of length $R$.

The conic is an ellipse or a hyperbola according as
$H$ is internal or external to triangle $ABC$, i.e.,
according as $ABC$ is acute or obtuse; or, according
as the D-F line is the external or internal angle-bisector
of $\angle OPH$.

\newpage

Let's try yet again, with another figure (\ref{DFWS2})

\begin{figure}[h]  
\begin{picture}(400,430)(-220,-170)
\setlength{\unitlength}{0.7pt}
\drawline(-36,-77)(-36,103)(-132,-77)(204,-77)(-36,103)(22.6667,213) 
\put(-49,104){$A$}
\put(-145,-88){$B$}
\put(206,-92){$C$}

\put(-36,51){\circle*{4}}   
\put(-35,52){$H$}   
\put(36,-51){\circle{340}}     
\put(36,-51){\circle*{4}}      
\put(39,-56){$O$}   

\put(138,85){\circle*{4}}   
\put(140,87){$M$}
\put(51,68){\circle*{4}}    
\put(53,71){$Q$}    
\drawline(-36,51)(66.3465,-10.53799)   
\drawline(36,-51)(138,85)   

\put(66.3465,-10.53799){\circle*{4}}   
\put(68,-16){$P$}
\put(0,0){\circle*{4}}      

\put(79.333,-77){\circle*{4}}   
\put(79.333,-77){\circle{344.5912}}   
\put(80,-91){$A_3$}
\put(57.9259,32.5556){\circle*{4}}   
\put(57.9259,32.5556){\circle{191.4395}}   
\put(59,34){$B_3$}
\put(22.6667,213){\circle*{4}}   
\put(22.6667,213){\circle{344.5912}}  
\put(251.6289,-77){\circle*{4}}   
\put(24,214){$C_3$}
\put(253,-88){$A_1$}
\drawline(204,-77)(251.6289,-77)   

\put(3,-4){$N$}     

\thicklines
\drawline(22.6667,213)(79.333,-77)   

\end{picture}
\caption{A converse approach.}
\label{DFWS2}
\end{figure}
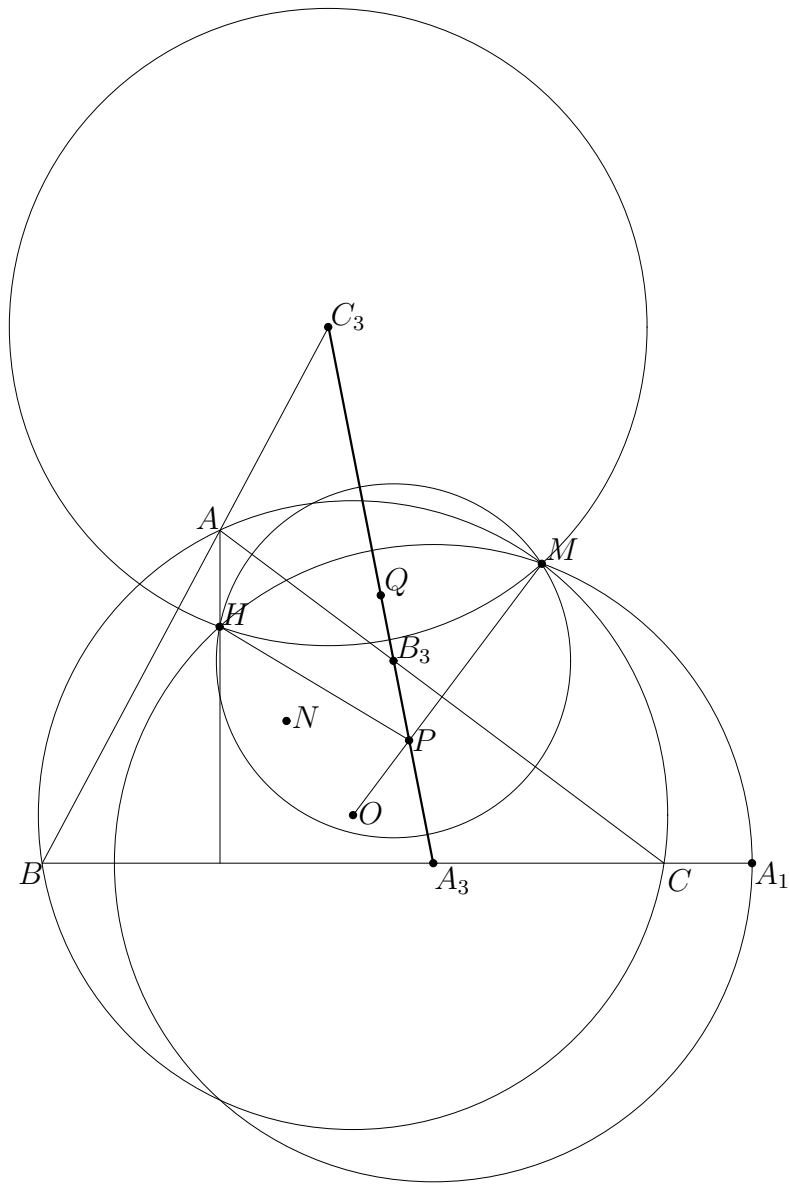

In Fig.\ref{DFWS2} $M$ is an arbitrary point on the circumcircle
of triangle $ABC$ and $H$ is the orthocentre.  The perpendicular
bisector of the segment $MH$ cuts the edges $BC$, $CA$, $AB$ of
the triangle in $A_3$, $B_3$, $C_3$ respectively.  This is a
putative Droz-Farny line.

The circles centres $A_3$, $B_3$, $C_3$ passing through $H$,
and therefore also through $M$, cut the respective edges
$BC$, $CA$, $AB$ in the pairs of points $A_1$, $A_2$;
$B_1$, $B_2$; $C_1$, $C_2$.  If we can prove that $A_1$, $B_1$,
$C_1$ are collinear with $H$ and that $A_2$, $B_2$, $C_2$ are
also collinear with $H$, then we have proved the Droz-Farny
theorem.  Note that $A_3$, $B_3$, $C_3$ are the midpoints of
$A_1A_2$, $B_1B_2$, $C_1C_2$ respectively, since the latter
are diameters of the three circles.  Moreover
angles $A_1HA_2$, $B_1HB_2$, $C_1HC_2$ are angles in semicircles
and therefore right-angles, so that it suffices to show, for
example, that $A_1HC_2$ is a right angle.

\smallskip

Third time lucky\,?? Yet another figure\,! (\ref{DFWS3})

\begin{figure}[h]
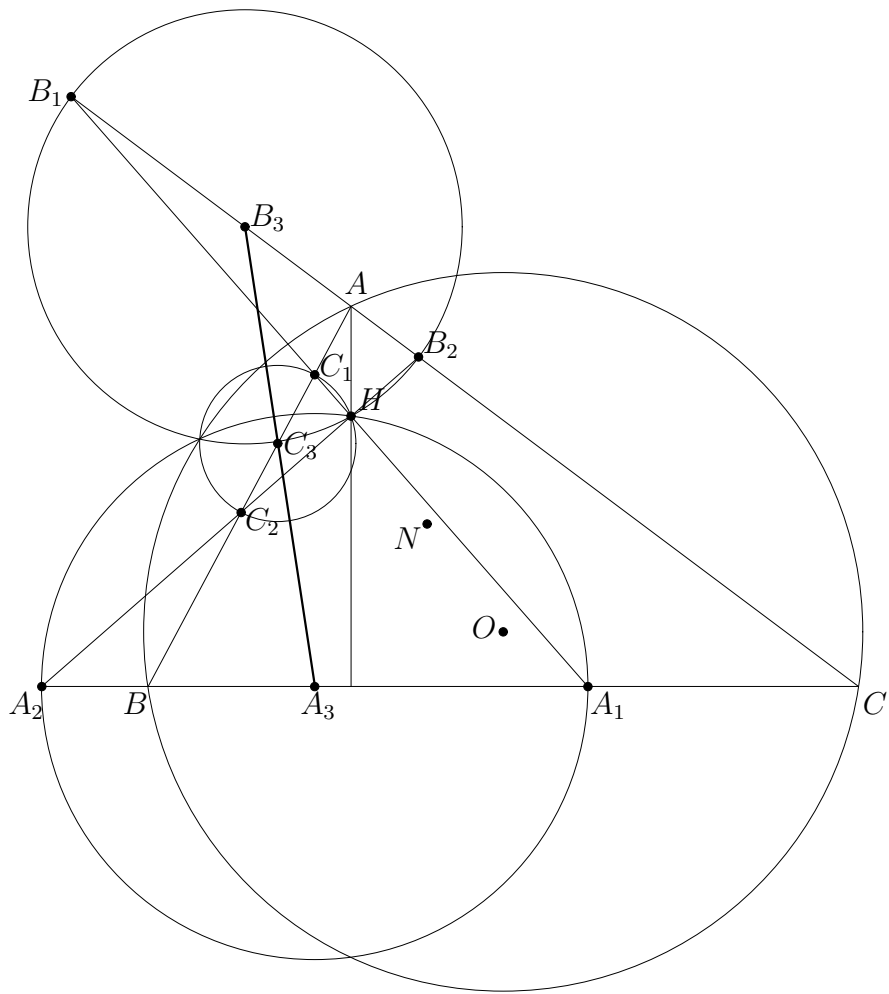
  

\caption{Another approach.}
\label{DFWS3}
\end{figure}

Choose a point $A_3$ arbitrarily on $BC$.  The circle
centre $A_3$ passing through $H$ cuts $BC$ in $A_1$, $A_2$.
The lines $A_1H$, $A_2H$ are perpendicular (angle in
a semicircle).  They meet the edge $CA$ in $B_1$ and $B_2$
and the edge $AB$ in $C_1$ and $C_2$.  Are the midpoints
$B_3$ and $C_3$ of $B_1B_2$ and $C_1C_2$ collinear with
$A_3$\,?

Do the 4 circles, centres $O$, $A_3$, $B_3$, $C_3$, all
concur\,??

\newpage

\section{There are two tangents from a point to a conic.}

In our converse approach (Fig.\ref{DFWS2}) there were,
in fact, two choices for our point $M$ where the circle
centre $A_3$ intersects the circumcircle.  Call them
$M_1$ and $M_2$.

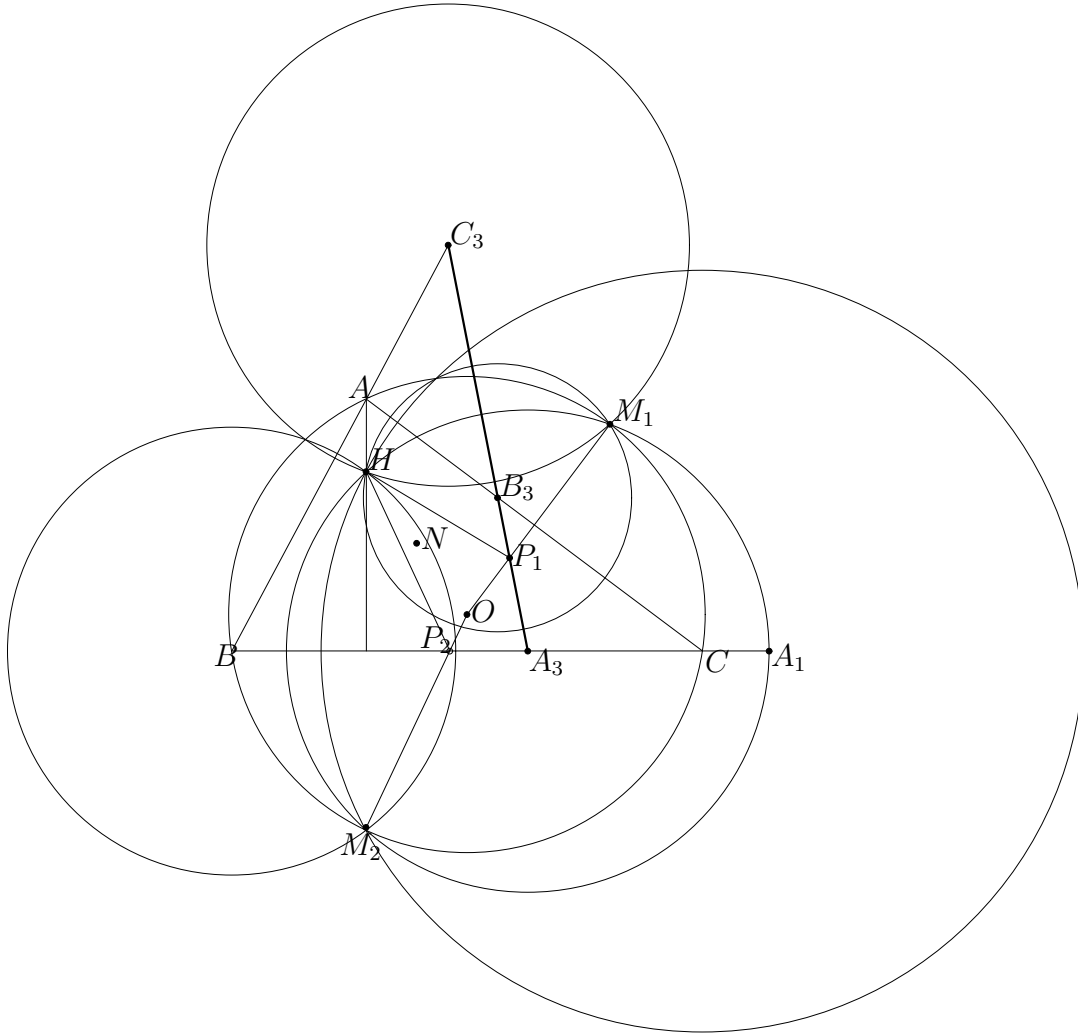
\begin{figure}[h]  
\begin{picture}(400,370)(-200,-175)
\setlength{\unitlength}{0.53pt}
\drawline(-36,-77)(-36,103)(-132,-77)(204,-77)(-36,103)(22.6667,213) 
\put(-49,104){$A$}
\put(-145,-88){$B$}
\put(206,-92){$C$}

\put(-36,51){\circle*{4}}   
\put(-35,52){$H$}   
\put(36,-51){\circle{340}}     
\put(36,-51){\circle*{4}}      
\put(39,-56){$O$}   

\put(138,85){\circle*{4}}   
\put(140,87){$M_1$}
\drawline(-36,51)(66.3465,-10.53799)   
\drawline(36,-51)(138,85)   

\put(-36,-203){\circle*{4}}   
\put(-55,-223){$M_2$}
\put(66.3465,-10.53799){\circle*{4}}   
\put(68,-16){$P_1$}

\put(-132,-77){\circle{320}}
\put(204,-77){\circle{544}}

\drawline(36,-51)(-36,-203)  

\put(23.68421,-77){\circle{4}}   
\put(2,-75){$P_2$}
\drawline(-36,51)(23.68421,-77)

\put(0,0){\circle*{4}}      

\put(79.333,-77){\circle*{4}}   
\put(79.333,-77){\circle{344.5912}}   
\put(80,-91){$A_3$}
\put(57.9259,32.5556){\circle*{4}}   
\put(57.9259,32.5556){\circle{191.4395}}   
\put(59,34){$B_3$}
\put(22.6667,213){\circle*{4}}   
\put(22.6667,213){\circle{344.5912}}  
\put(251.6289,-77){\circle*{4}}   
\put(24,214){$C_3$}
\put(253,-88){$A_1$}
\drawline(204,-77)(251.6289,-77)   

\put(3,-4){$N$}     

\thicklines
\drawline(22.6667,213)(79.333,-77)   

\end{picture}
\caption{Two tangents to the conic.}
\label{DFWS4}
\end{figure}

$M_2$ is the reflexion of $H$ in
the edge $BC$.  This point is well-known to lie
on the circumcircle.  Then $BC$
is the perpendicular bisector of $HM_2$ and is a
candidate for a Droz-Farny line through $A_3$.
Certainly the circles centres $B$ and $C$,
passing through $H$, also pass througn $M_2$.
$OM_2$ intersects $BC$ in $P_2$.
$HP_2+P_2O=M_2P_2+P_2O=R$, the circumradius.
$BC$ is the exterior angle-bisector of
$\angle HP_2O$, so that $BC$ is the tangent
at $P_2$ to the ellipse with foci $H$ and $O$
and major axis of length $R$.

$A_3BC$ and $A_3B_3C_3$ are the tangents at $P_2$
and $P_1$ to the conic.

\newpage

\section{Quadration and the Droz-Farny line.}

Let's find the Droz-Farny lines through $A_3$
for the two triangles $ABC$ and $HBC$.  For each
triangle, $BC$ is a D-F line.

\begin{figure}[h]
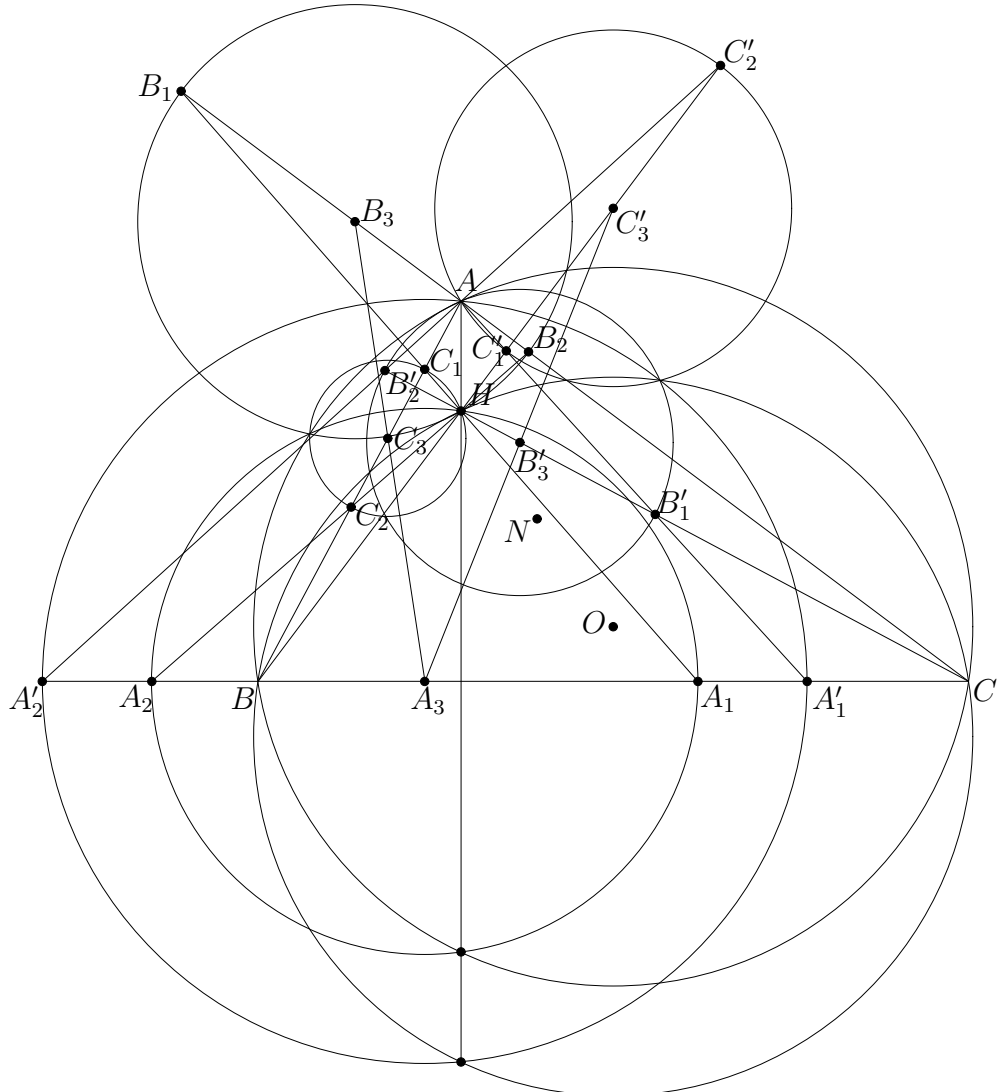
  

\caption{Droz-Farny lines for triangles $ABC$ and $HBC$.}
\label{DFWS5}
\end{figure}

\newpage

\begin{figure}[h]  
%
\caption{Quadration and the Droz-Farny lines.}
\label{DFWS6}
\end{figure}

\newpage

\begin{figure}[h]
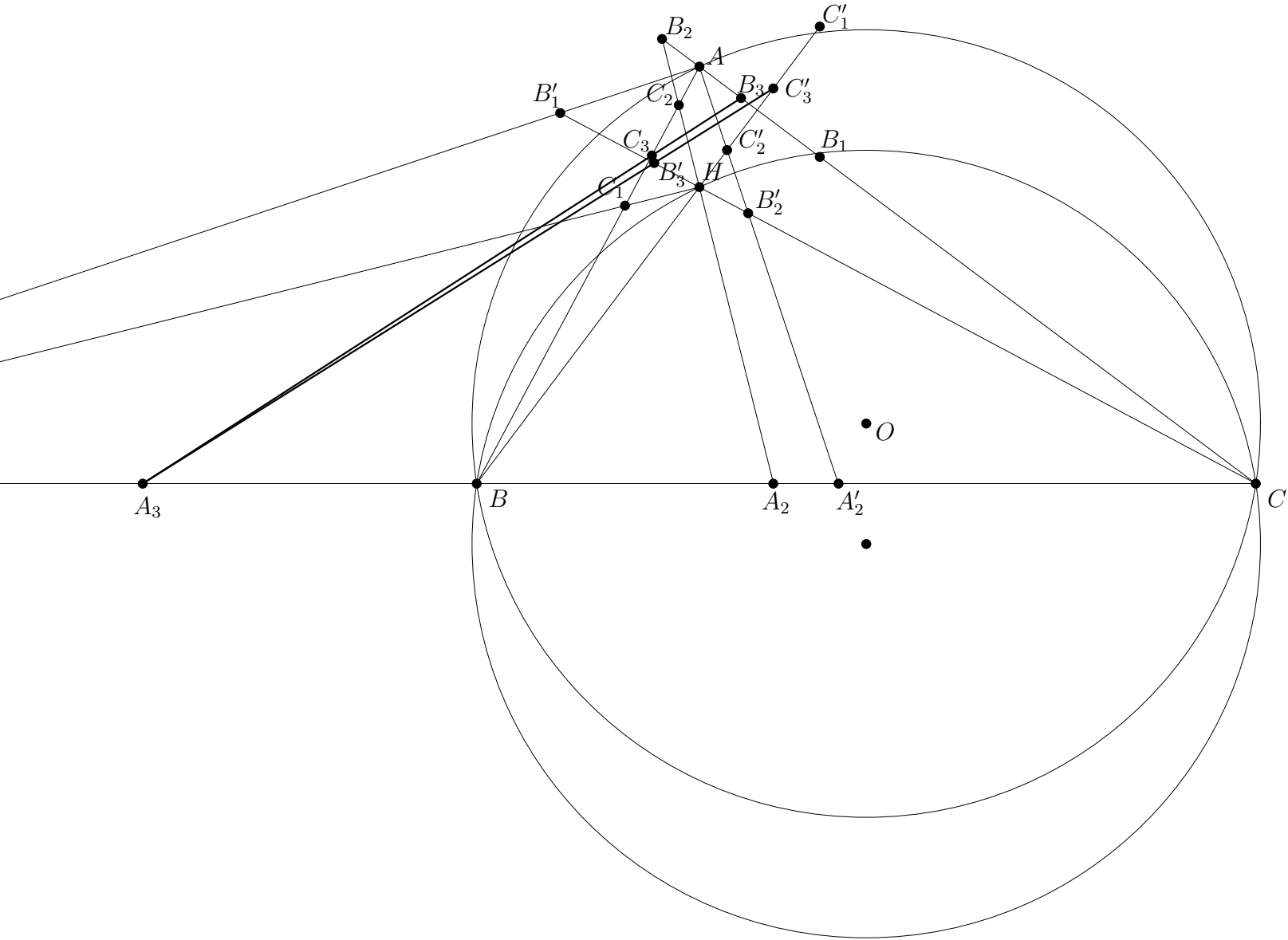
  
%
\caption{Enlargement of Fig.\ref{DFWS6}.}
\label{DFWS7}
\end{figure}

\newpage

\section{Some background theorems}

{\bf Carnot's theorem} \cite[p.101]{Car}.  The segment of an
altitude from the orthocentre to an edge equals its
extension from the edge to the circumcircle.

In fact the circumcircles of $BHC$, $CHA$, $AHB$ are the
reflexions of the circumcircle $ABC$ in the edges $BC$,
$CA$, $AB$.  Indeed

\bigskip

{\bf Theorem R}. \cite[p.99]{Lal}\cite[\S333]{J}

Let $\mathcal{L}$ be a line through the orthocentre $H$ of a
triangle $ABC$ which cuts $BC$, $CA$, $AB$ in $L_a$,
$L_b$, $L_c$.  The reflexions of $\mathcal{L}$ in the
edges of $ABC$, $H_aL_a$, $H_bL_b$, $H_cL_c$, concur at a
point, $P$, on the circumcircle.  Moreover, the Wallace line
of $P$ is parallel to $\mathcal{L}$.

The proof is an offshoot of that of the Wallace (Simson) line
theorem.

\bigskip

{\bf Miquel's theorem} \cite{Miq}\cite[\S184]{J}.  Also
known as the {\bf Pivot Theorem}.

$X$, $Y$, $Z$ are arbitrary points on the respective edges $BC$,
$CA$, $AB$ of a triangle $ABC$.  Then the circles $AYZ$, $BZX$,
$CXY$ pass through a common point, the {\bf Miquel point} of $XYZ$.

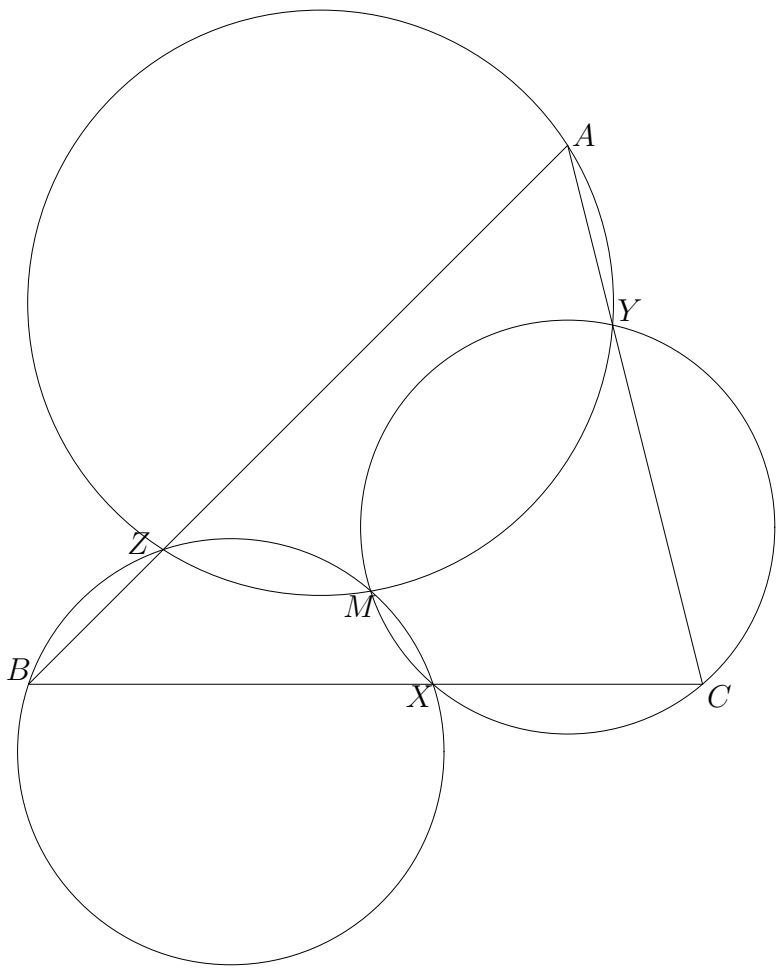
\begin{figure}[h]  
\begin{picture}(360,370)(-300,-110)
\setlength{\unitlength}{17pt}
\drawline(0,12)(-12,0)(3,0)(0,12)
\put(-7.5,-1.5){\circle{9.48683}}
\put(0,3.5){\circle{9.21954}}
\put(-5.5,8.5){\circle{13.0384}}
\put(0.1,12){$A$}
\put(-12.5,0.1){$B$}
\put(3.1,-0.5){$C$}
\put(-3.6,-0.5){$X$}
\put(1.1,8.1){$Y$}
\put(-9.8,2.9){$Z$}
\put(-5,1.5){$M$}
\end{picture}
\caption{Miquel's theorem}
\label{miquel}
\end{figure}


Miquel's theorem is easy to prove: Suppose that the circles $BZX$, $CXY$
meet again in $M$.  Then $\angle ZMX=\pi-B$, $\angle XMY=\pi-C$, so that
$\angle YMZ=\pi-A$ and $AZMY$ is cyclic.

\bigskip

{\bf The Wallace (aka Simson) Line theorem}.  This appears
elsewhere in this document; see especially \S8.7.

\bigskip

{\bf Ceva's theorem}.  $O$ is in the plane of triangle $ABC$.
$OA$ meets $BC$ in $A'$, $OB$ meets $CA$ in $B'$, $OC$ meets
$AB$ in $C'$. Then the product of the ratios
$$\frac{BA'}{A'C}\ \frac{CB'}{B'A}\ \frac{AC'}{C'B} \ = \ +1$$
and conversely, if the product is 1, the ``cevians'' concur.

\newpage

{\bf Menelaus's theorem}. If a line in the plane of triangle $ABC$
meets the edges $BC$, $CA$, $AB$ respectively in $A'$, $B'$, $C'$,
then the product of the ratios
$$\frac{BA'}{A'C}\ \frac{CB'}{B'A}\ \frac{AC'}{C'B} \ = \ -1$$
and conversely, if the product is -1, then the points
$A'$, $B'$, $C'$ are collinear.

\smallskip

Note that combination of Ceva's and Menelaus's theorems gives
a ``pole and polar'' relationship between points and lines
with respect to a triangle.  E.g. The polar of the centroid
is the line at infinity.

\newpage

{\bf Parabola theorems}.

{\bf Notation.} $S$ is the focus, $A$ the vertex,
$X$ the foot of the directrix.  $P$ any point on the curve.
$PN$ the ordinate from $P$ to the axis.  $PM$ perpendicular
from $P$ to the directrix.  Tangent at $P$ cuts directrix
at $R$ and axis at $T$.  Normal at $P$ cuts axis at $G$.
$SY$ is perpendicular from $S$ to $PT$,

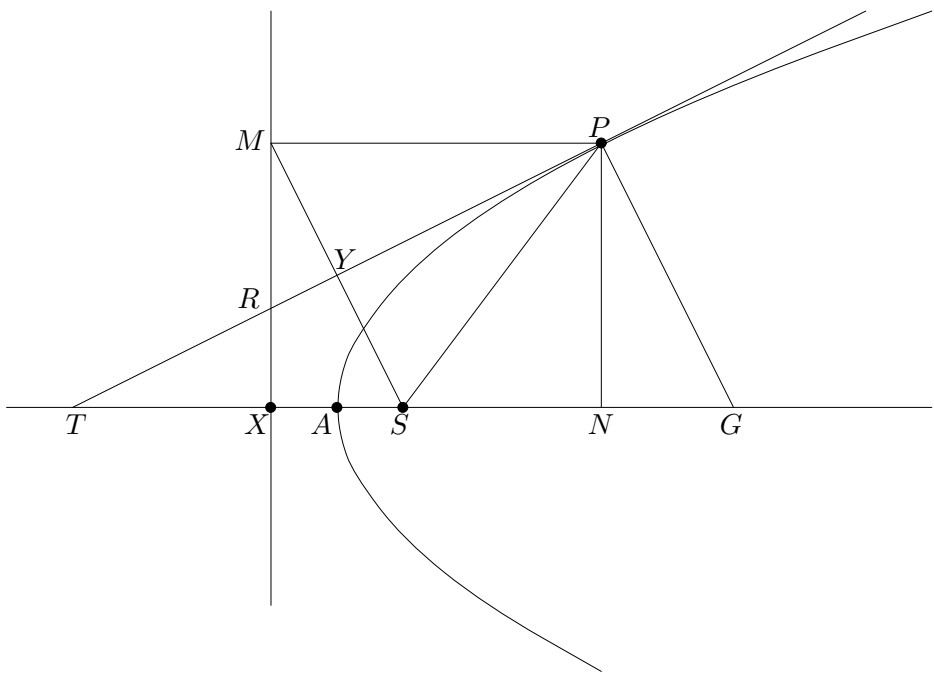
\begin{figure}[h]  
\begin{picture}(400,230)(-160,-90)
\setlength{\unitlength}{2.5pt}
\put(0,0){\circle*{1.5}}    
\put(10,0){\circle*{1.5}}   
\put(-10,0){\circle*{1,5}}
\drawline(-50,0)(90,0)      
\drawline(-10,60)(-10,-30)  
\spline(90,60)(62.5,50)(40,40)(22.5,30)(10,20)(2.5,10)(0.625,5)
(0,0)(0.625,-5)(2.5,-10)(10,-20)(22.5,-30)(40,-40) 
\put(40,40){\circle*{1.5}}   
\put(38,41){\small$P$}
\drawline(-40,0)(80,60)     
\drawline(-10,40)(40,40)(40,0)  
\put(-15.5,39){\small$M$}
\put(38,-4){\small$N$}
\drawline(60,0)(40,40)(10,0)(-10,40)   
\put(-41,-4){\small$T$}
\put(-4,-4){\small$A$}
\put(8,-4){\small$S$}
\put(-14,-4){\small$X$}
\put(58,-4){\small$G$}
\put(-15,15){\small$R$}
\put(-0.5,21){\small$Y$}
\end{picture}
\caption{The parabola}
\label{parabola}
\end{figure}

{\bf Theorem.} \cite[Theorem 9]{Du}. (i) $PT$ bisects
$\angle SPM$; $SP=ST=SG$. \\
(ii) $SM$, $PT$ bisect each other at right angles. \\
(iii) $Y$ lies on tangent at vertex. \\
(iv) $AN=AT$; $AY=\frac{1}{2}NP$; $NG=2AS$. \\
(v) $SY^2=SA\times SP$.

{\bf Theorem.} \cite[Theorem 10]{Du}. The tangents
at the ends of a focal chord $PSP^{\prime}$ meet at
$R$, on the directrix.  The line through $R$ parallel
to the axis bisects $PP^{\prime}$ at $V$.
$\angle PRP^{\prime}=90^{\circ}$. If $RV$ cuts the
parabola at $K$, then $RK=KV$ and $PP^{\prime}=4SK$.

{\bf Theorem.} \cite[Theorem 11]{Du}. If $OP$, $OQ$
are tangents to a parabola, then triangles $SOP$,
$SQO$ are similar; $SO^2=SP\cdot SQ$; the exterior angle
between the tangents is equal to the angle either subtends
at the focus; and $\angle SOQ$ equals the angle that
$OP$ makes with the axis.

{\bf Theorem.} \cite[Theorem 12]{Du}. The circumcircle
of a triangle circumscribing a parabola passes through
the focus and ({\bf Steiner's theorem}) the orthocentre
lies on the directrix.

\newpage

{\bf Jean-Louis Ayme's proof of the Droz-Farny theorem.}

$H$ is the orthocentre of triangle $ABC$.  Two
perpendicular lines through $H$, $\mathcal{L}$ and
$\mathcal{L^{\prime}}$, cut $BC$ in $X$ and
$X^{\prime}$, $CA$ in $Y$ and $Y^{\prime}$, and
$AB$ in $Z$ and $Z^{\prime}$.

Let $\mathcal{C}$, $\mathcal{C}_a$, $\mathcal{C}_b$,
$\mathcal{C}_c$, be the circumcircles of triangles
$ABC$, $HXX^{\prime}$, $HYY^{\prime}$, $HZZ^{\prime}$,
and their centres be $O$, $M_a$, $M_b$, $M_c$.  Let
$H_a$, $H_b$, $H_c$ be the reflexions of $H$ in
$BC$, $CA$ and $AB$, respectively.  By Carnot's
theorem, $H_a$, $H_b$, $H_c$ lie on $\mathcal{C}$.

Note that, as $XHX^{\prime}$ is a right angle,
the circle $\mathcal{C}_a$ has $XX^{\prime}$
as a diameter, so that $M_a$ lies on $BC$ and the
circle passes through $H_a$ which also lies on the
circle $\mathcal{C}$. (Carnot's theorem.) Similarly $H_b$ is an
intersection of $\mathcal{C}$ and $\mathcal{C}_b$
and the perpendicular to $CA$ through $H$.

\begin{figure}[h]
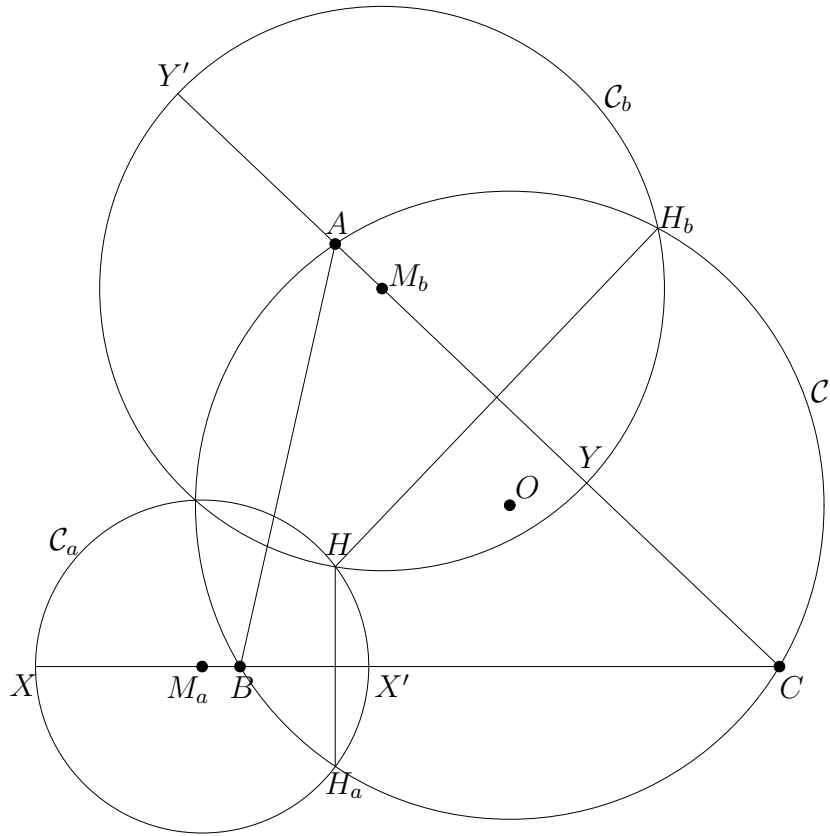
  

\caption{Start of Ayme's proof.}
\label{Ayme1}
\end{figure}

\newpage

\begin{figure}[h]
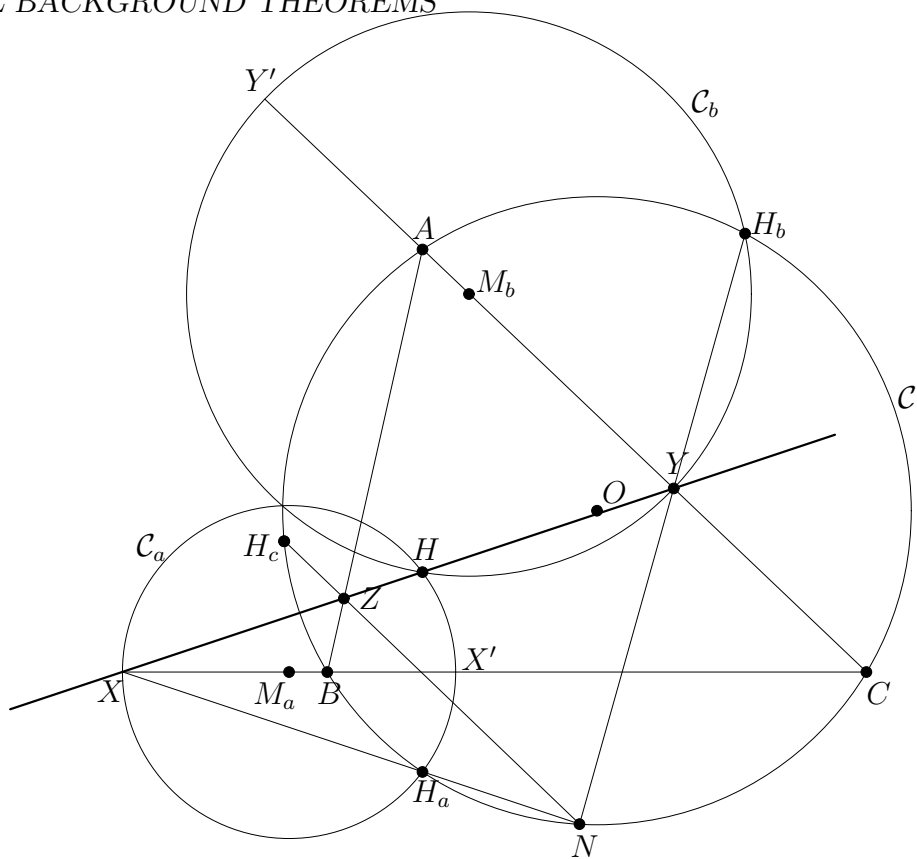
  
%
\caption{Conclusion of Ayme's proof.}
\label{Ayme2}
\end{figure}

By applying Theorem R to the line $XYZ$ through $H$,
we conclude that its reflexions $H_aX$, $H_bY$, $H_cZ$
in the respective edges $BC$, $CA$, $AB$
concur on the circumcircle of $ABC$ at a point $N$, say.
Now apply Miquel's pivot theorem to the triangle $YXN$
with the points $H_a$, $H_b$, $H$ on respective edges
$XN$, $NY$, $YX$, and we see that the circles
$HH_aX=\mathcal{C}_a$, $HH_bY=\mathcal{C}_b$, and
$NH_aH_b=\mathcal{C}$ concur in a point $M$, say.
Similarly it can be shown that the circles $\mathcal{C}$,
$\mathcal{C}_b$, $\mathcal{C}_c$ all pass through $M$.
So the circles $\mathcal{C}_a$, $\mathcal{C}_b$, and
$\mathcal{C}_c$ all pass through $H$ and $M$, and are
therefore coaxal, with their centres $M_a$, $M_b$, $M_c$
collinear on the Droz-Farny line.

Note that this is the perpendicular bisector of the
segment $HM$.  $HM$ is the axis of the parabola we're
seeking.  The D-F line is the tangent at the vertex.
The parallel to the D-F line through $H$ is the directrix.
Let this meet $BC$, $CA$, $AB$ in $D_a$, $D_b$, $D_c$.
Then the perpendiculars to $BC$, $CA$, $AB$ through
$D_a$, $D_b$, $D_c$ are also tangents to the parabola.

Is there a short proof that $H_aX^{\prime}$, $H_bY^{\prime}$,
$H_cZ^{\prime}$ concur at a point, $N^{\prime}$, say,
diametrically across the circumcircle from $N$\,?

In Figure\,\ref{Ayme2}, $O$ does not lie on the line $HXYZ$.
The points $MHH_cZZ^{\prime}$ are concyclic (the
circle $\mathcal{C}_c$). The lines $AHH_a$, $BHH_b$,
$CHH_c$ are the altitudes of triangle $ABC$.

\newpage

\section{The parabola associated with the D-F line}

The dual of Feuerbach's 11-point conic theorem is
the 11-line theorem \cite[Theorem 133]{Du1}:

$p$ is a variable line through a fixed point $L$;
$p^{\prime}$ is a line conjugate to $p$ w.r.t.\ a
range of conics touching four fixed lines $a$, $b$,
$c$, $d$; the envelope of $p^{\prime}$ is a conic
$\sigma$ which touches the polars of $L$ w.r.t.\
the conics of the range.

The conic $\sigma$ touches the following 11 lines:

If $h$ is the join of $L$ and $ab$ and $h^{\prime}$
is the harmonic conjugate of $h$ w.r.t.\ $a$, $b$,
then $\sigma$ touches $h^{\prime}$ and the corresponding
five lines through the other five vertices of the
quadrilateral $abcd$; $\sigma$ also touches the three
edges of the diagonal line triangle of $abcd$, and the
two double lines of the involution formed by the pairs
of tangents from $L$ to the range of conics.

\medskip

We are interested in the case where one of the four
lines, say $d$, is the line at infinity.  The conics
of the range are then parabolas.

\medskip

As we vary the pair of perpendicular lines, the locus
of the focus of the parabola is the circumcircle of $ABC$.
The directrix passes through the orthocentre $H$.  The
Droz-Farny line is the tangent at the vertex of the
parabola and we already know that this envelopes the
conic with foci $H$ and $O$, the circumcentre of $ABC$,
and major (real) axis of length $R$, the circumradius.
Centre is 50-point centre.  Axis is diameter of 50-point
circle (Euler line).

\newpage

Here's a description of the parabola associated with the
Droz-Farny line.  We'll use the notation of the last figure,
repeated here.

\begin{figure}[h]
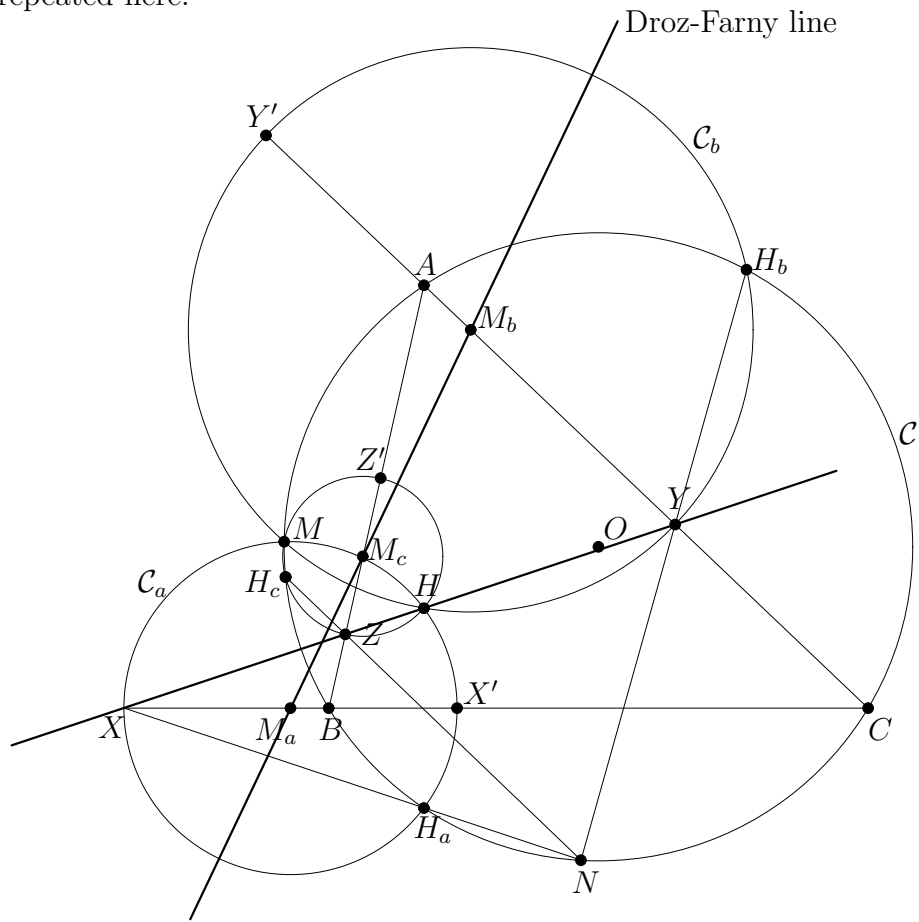
  

\caption{Parabola associated with Droz-Farny line (to appear\,!!)}
\label{parabola2}
\end{figure}

The Droz-Farny line is the tangent to the parabola at its
vertex, the midpoint of $MH$; where $H$ is the orthocentre
and $M$ is the point on the circumcircle where the circles
$\mathcal{C}_a$, $\mathcal{C}_b$ and $\mathcal{C}_c$ concur.

$M$ is the focus of the parabola, and $MH$ is its axis.  The
perpendicular to the axis through $H$ is the directrix,
parallel to the Droz-Farny line.

The edges $BC$, $CA$, $AB$ touch the parabola; where are the
points of contact\,?  The perpendicular lines $HXYZ$ and
$HX^{\prime}Y^{\prime}Z^{\prime}$ touch the parabola: where\,?
If the points of contact are $P$ and $P^{\prime}$, then the
line $PP^{\prime}$ passes through the focus, $M$.

The parabola also touches the line at infinity.

\newpage

Another drawing

\begin{figure}[h]
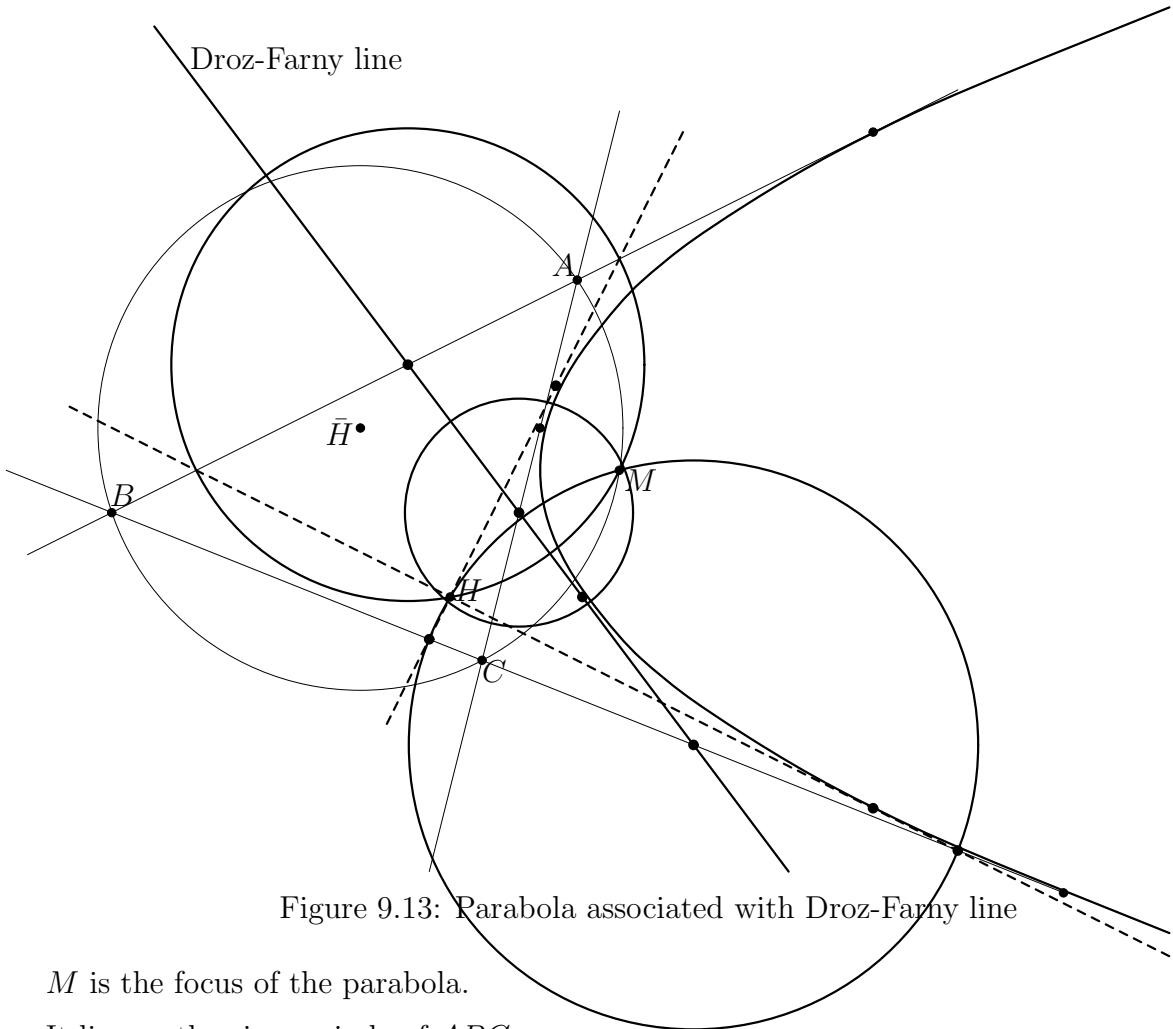
  

\caption{Parabola associated with Droz-Farny line}
\label{parabola3}
\end{figure}

$M$ is the focus of the parabola.

It lies on the circumcircle of $ABC$,

The D-F line is the perpendicular bisector of $MH$, where
$H$ is the orthocentre of $ABC$.  Let it cut $OM$ in $P$,
so that $OP+PH=OP+PM=OM=R$, the circumradius.  $MP$ and $HP$ are
equally inclined to the D-F line, so that it is a tangent to the
ellipse, foci $O$ and $H$ with major axis $R$.  In case $H$
is outside $ABC$ (now obtuse) and its circumcircle, then
$OP-PH=OM=R$ and the conic is a hyperbola.

D-F line cuts $BC$, $CA$, $AB$ in $X$, $Y$, $Z$.

Circles centres $X$; $Y$; $Z$ passing through $H$ cut
 $BC$; $CA$; $AB$ in $X_1$, $X_2$; $Y_1$, $Y_2$; $Z_1$, $Z_2$.

$H$, $X_1$, $Y_1$, $Z_1$ are collinear. Also $H$, $X_2$, $Y_2$, $Z_2$.
These (dashed) lines are perpendicular, and tangent to the parabola.
Synthetic proof wanted.  Compare next section (or not! according
to taste).

\newpage

\begin{figure}[h]
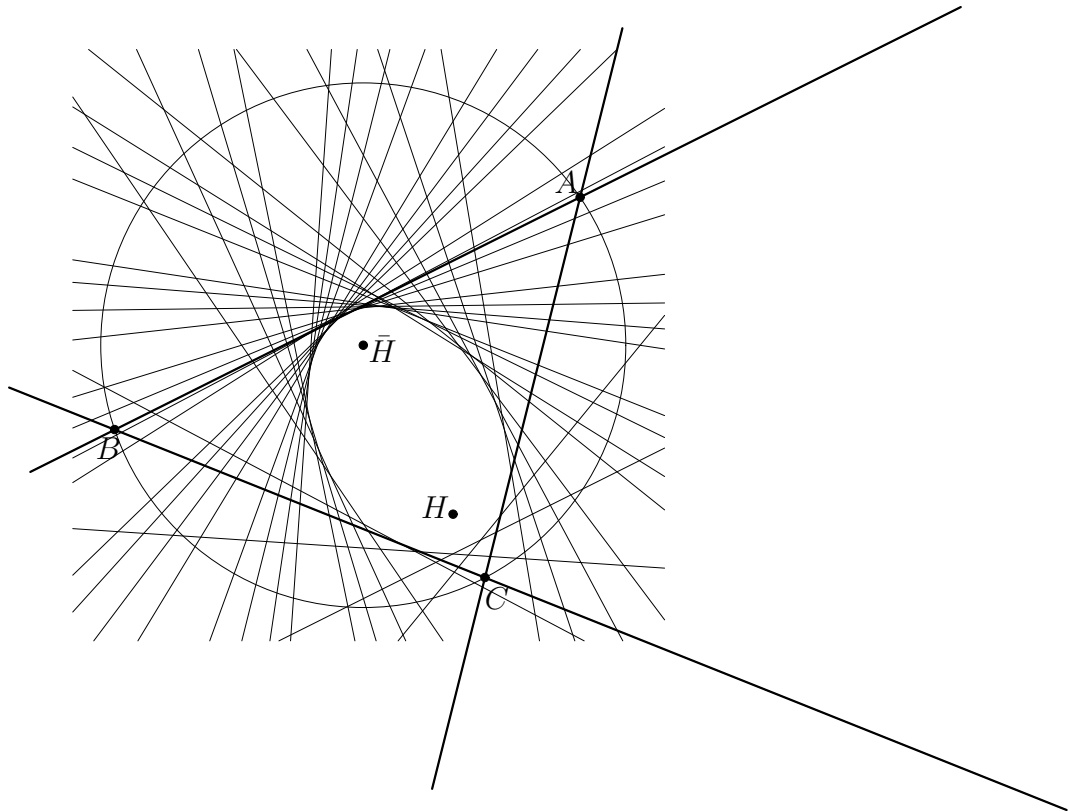
  

\caption{Envelope of Droz-Farny line}
\label{dfenvelope}
\end{figure}

$M$ is a variable point on the circumcircle of $ABC$.  The
perpendicular bisector of $MH$ is a Droz-Farny line, where
$H$ is the orthocentre of $ABC$.  In Fig. 9.14, if the
coordinates of $H$ are $(a,b)$ and those of $M$ are $(p,q)$,
then the D-F line has equation $$2x(p-a)+2y(q-b)=p^2-a^2+q^2-b^2$$.

Note that the reflexions of $H$ in the edges of $ABC$ lie on
the circumcircle, so that the edges are special cases of the
D-F line.  So also are the perpendicular bisectors of $AH$,
$BH$, $CH$, which form the twin triangle of $ABC$ (rotate
$ABC$ through an angle $\pi$ about its (one-\&-only, ``9-point'')
CENTRE, the midpoint of $OH$.  There is an involution
between $O$ and $H$.

\newpage

Let's look at the obtuse case.  Interchange $A$ and $H$; i.e.,
$A$ is the orthocentre of triangle $BHC$.  This is an example of
{\bf quadration.}

\begin{figure}[h]
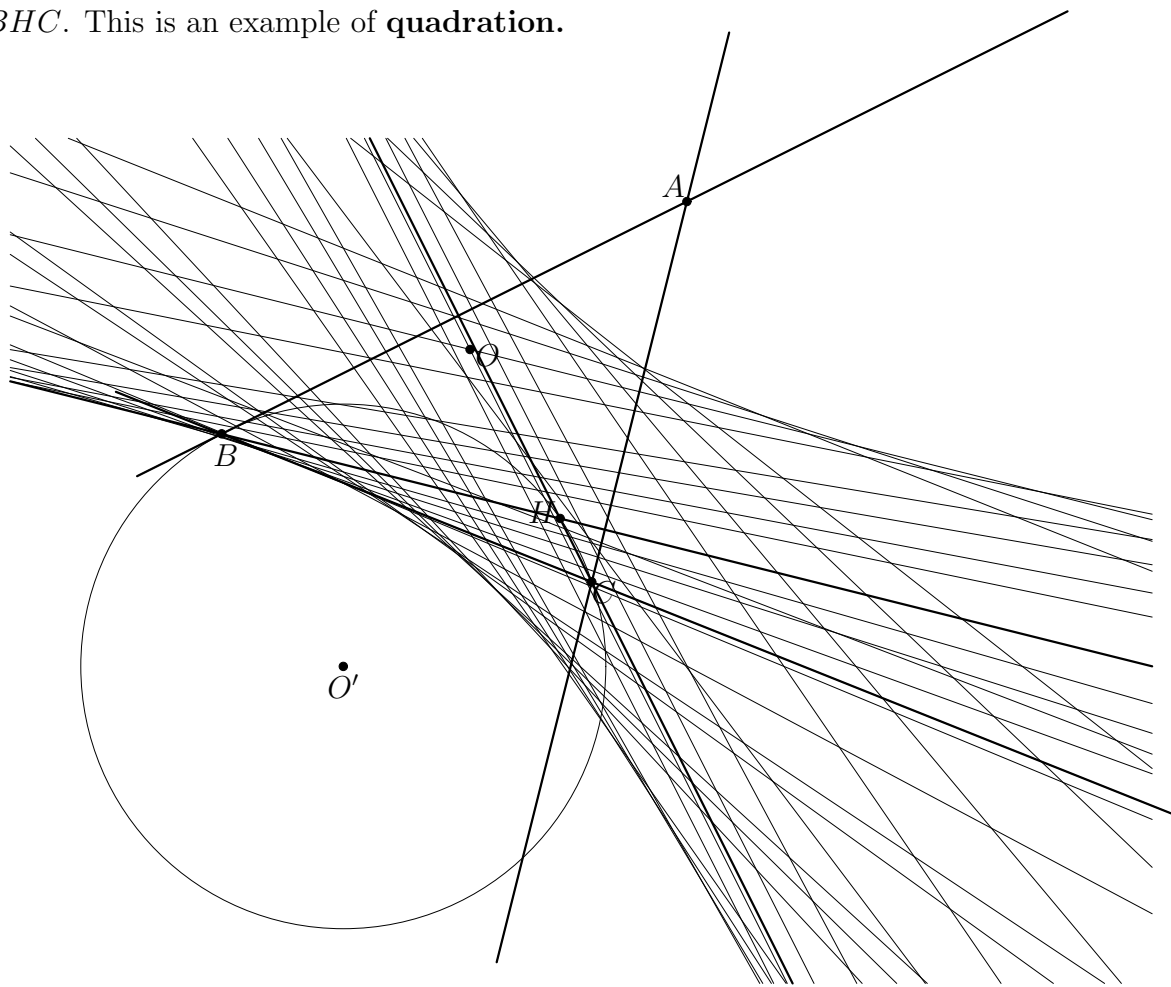
  

\caption{Envelope of Droz-Farny line in the obtuse case}
\label{dfobtuse}
\end{figure}

Envelope is a hyperbola, foci $A$, $O'$, real axis of
length $R$, the circumradius.  Note that the edges, $BC$, $CH$
and $HB$ are also tangents to the hyperbola, since the
reflexions of $A$ in these edges lie on the circumcircle
$BCH$.  Since $B$, $C$, $H$ are points on the circumcircle,
so also are the perpendicular bisectors of $AB$, $AC$, $AH$
tangents to the hyperbola.  Indeed these three are the edges
of the TWIN of $HBC$, its reflexion in the CENTRE of triangle
$HBC$, which is also the CENTRE of $ABC$ and the centre of
the hyperbola,which we may see as follows.  The asymptotes
of the hyperbola will be the perpendicular bisectors of
$AM_1$, $AM_2$, the tangents from $A$ to the circumcircle
$BHC$.  As tangent is perpendicular to radius, these two
bisectors will pass through the midpoints of $AM_1$, $AM_2$
and be parallel to the radii $O'M_1$, $O'M_2$, and hence
cut the third edge, $AO'$ of triangles $AM_1O'$, $AM_2O'$
in the midpoint of $AO'$, which is also the midpoint of
$OH$ and the CENTRE.

\newpage

\section{Interruption --- added later  --- can now throw away\,?}

Or earlier, depending on what edition of this work you are reading.
See \S 9.8

I outline what I hope will be a synthetic proof of
the Droz-Farny theorem (\&/or its converse?).  Perhaps
it's the same as J.-L.~Ayme's ?

\begin{figure}[h]  
\begin{picture}(425,360)(-220,-210)
\setlength{\unitlength}{0.4pt}
\drawline(0,240)(-217,-271)
\drawline(0,240)(0,-240)
\drawline(0,240)(-320,0)
\drawline(0,240)(90,-160)
\dashline{5}(-510,155)(512,-279)
\dottedline{7}(-510.48,-48)(141.8,320)
\dottedline{7}(428.77,-520)(-67.7,360)
\put(-145,0){\circle{560.8}}
\put(147,-124){\circle{785.1}}
\put(-205.6,25.7){\circle{593.9}}
\thicklines
\put(-133,-84){\circle{410}}
\drawline(-510,0)(550,0)
\drawline(-512,-43.2)(520,194)
\drawline(-216,360)(444,-520)
\put(-2,247){$H$}
\put(2,77){$A$}
\put(-320,-20){$B$}
\put(56,4){$C$}
\put(0,-260){$H_a$}
\put(-163,95){$H_b$}
\put(75,-87){$H_c$}
\put(-239,-294){$K$}
\put(-148,7){{\small$A_3$}}
\put(150,-122){{\small$B_3$}}
\put(-209,5){{\small$C_3$}}
\put(142,7){{\small$A_1$}}
\put(-457,7){{\small$A_2$}}
\put(392,-444){{\small$B_1$}}
\put(-96,205){{\small$B_2$}}
\put(48,98){{\small$C_1$}}
\put(-529,-34){{\small$C_2$}}
\end{picture}
\caption{}
\label{}
\end{figure}

\newpage

Depending on whether one is trying to prove
the theorem or its converse, select a point on an
edge of the triangle (say $A_3$ on $BC$) or a point
on the circumcircle of $ABC$ (say $K$).  These are
related as follows.  If you selected $A_3$ on $BC$,
then draw the circle, centre $A_3$, passing through
the orthocentre, $H$.  This will intersect the
circumcircle of $ABC$ in two points, $H_a$ and $K$.
$H_a$ is the reflexion of $H$ in the edge $BC$ ---
it's well known to those who well know it, that the
circumcircle of $HBC$ is the reflexion of the
circumcircle of $ABC$ in the edge $BC$.  Alternatively,
choose $K$ arbitrarily on the circumcircle of $ABC$
and draw the circle through $KHH_a$. By the symmetry
already mentioned, its centre will lie on $BC$; call
it $A_3$.  In either case the circle will intersect
$BC$ in $A_1$, $A_2$, say, and $A_1A_2$ will be a
diameter of the circle and $HA_1$, $HA_2$ will be our
perpendicular rays through the orthocentre (shown
dotted in Figure 1).

Next, draw the circles $KHH_b$ and $KHH_c$.  They form
a coaxal system with the circle $KHH_a$.  By the
symmetry already noted, their centres, $B_3$ and $C_3$
say, will lie on $CA$ and $AB$ respectively, and will
be collinear with $A_3$, forming the Droz-Farny line,
the perpendicular bisector of $HK$ (shown dashed in Figure 1).

If these two circles respectively cut the edges $CA$
and $AB$ in the points $B_1$, $B_2$ and $C_1$, $C_2$,
then it remains to be shown that $H$, $A_1$, $B_1$,
$C_1$ and $H$, $A_2$, $B_2$, $C_2$ are collinear,
and that the two lines are perpendicular.  The angles
$A_1HA_2$, $B_1HB_2$, $C_1HC_2$ are each (right) angles in a
semicircle, so it suffices to show, for example, that
$HA_1$ is perpendicular to $HB_2$.  Because we haven't
made it clear as to which of $A_1$, $A_2$, etc., is
which, we run into difficulty --- a difficulty experienced
both here and by the author in his proofs;
the arguments become diagram dependent, and angles are
liable to be confused with their supplements.
For example, the specification {\it internal} angle-bisector
may be ambiguous.  In fact, if a line through an angle,
$APB$, say, is reflected about {\it either} angle-bisector,
the result is the same.

\section{Three (or more) new(?) theorems}

We have also shown:

{\bf Theorem.}  If the pair of perpendicular rays through
the orthocentre is allowed to vary, then the locus of the
reflexion of the orthocentre in the Droz-Farny line is the
circumcircle of $ABC$.

{\bf Theorem.}  If the pair of perpendicular rays through
the orthocentre is allowed to vary, then the locus of the
foot of the perpendicular from the orthocentre onto the
Droz-Farny line is the nine-point circle of $ABC$.

\newpage

It is also now straightforward to prove:

{\bf Theorem.}  The envelope of the Droz-Farny line is a
conic.  The centre is the nine-point centre of the triangle
$ABC$, and the major axis is the segment of the Euler line
which is a diameter of the nine-point circle, and has length
$R$, the circumradius of $ABC$.  The conic is an ellipse or
a hyperbola, according as the triangle is acute
or obtuse.  The conjugate axis has length $\sqrt{R^2-h^2}$,
where $h=HO$, the distance of the orthocentre from the
circumcentre.  If $h>R$, the (length of the) conjugate axis
(of the hyperbola) is imaginary.  In this case the asymptotes
are the diameters of the nine-point circle through the
points of contact of the tangents to the nine-point circle
from the orthocentre.

These asymptotes are special cases of the Droz-Farny line.
Other examples are the edges of the triangle, and the
perpendicular bisectors of the segments $HA$, $HB$ and $HC$.
Also the tangents to the nine-point circle at its intersections
with the Euler line, i.e., the vertices of the conic.

\begin{figure}[h]
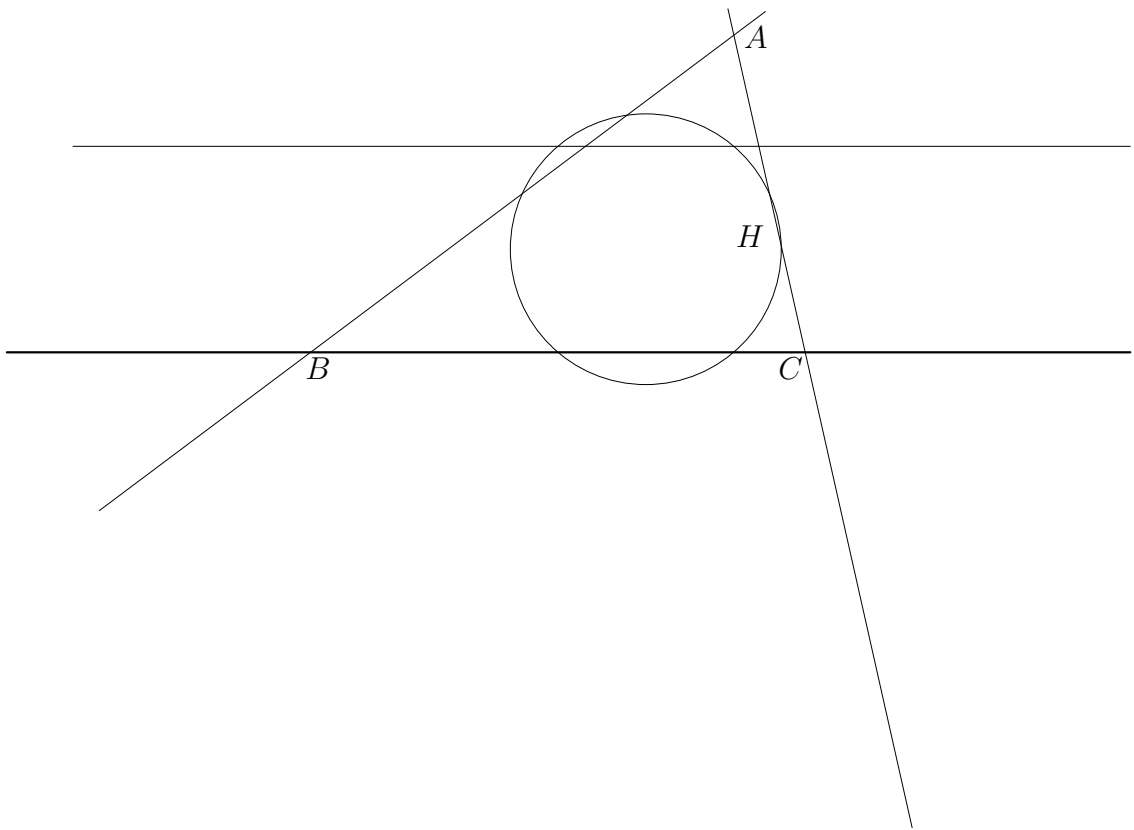
  

\caption{Droz-Farny lines touching an ellipse.  Incomplete. Add
perp bisectors of $AH$, $BH$, $CH$, tans to 9-pt circle where it
meets the Euler line, and further examples.}
\label{ellipse}
\end{figure}

\begin{figure}[h]
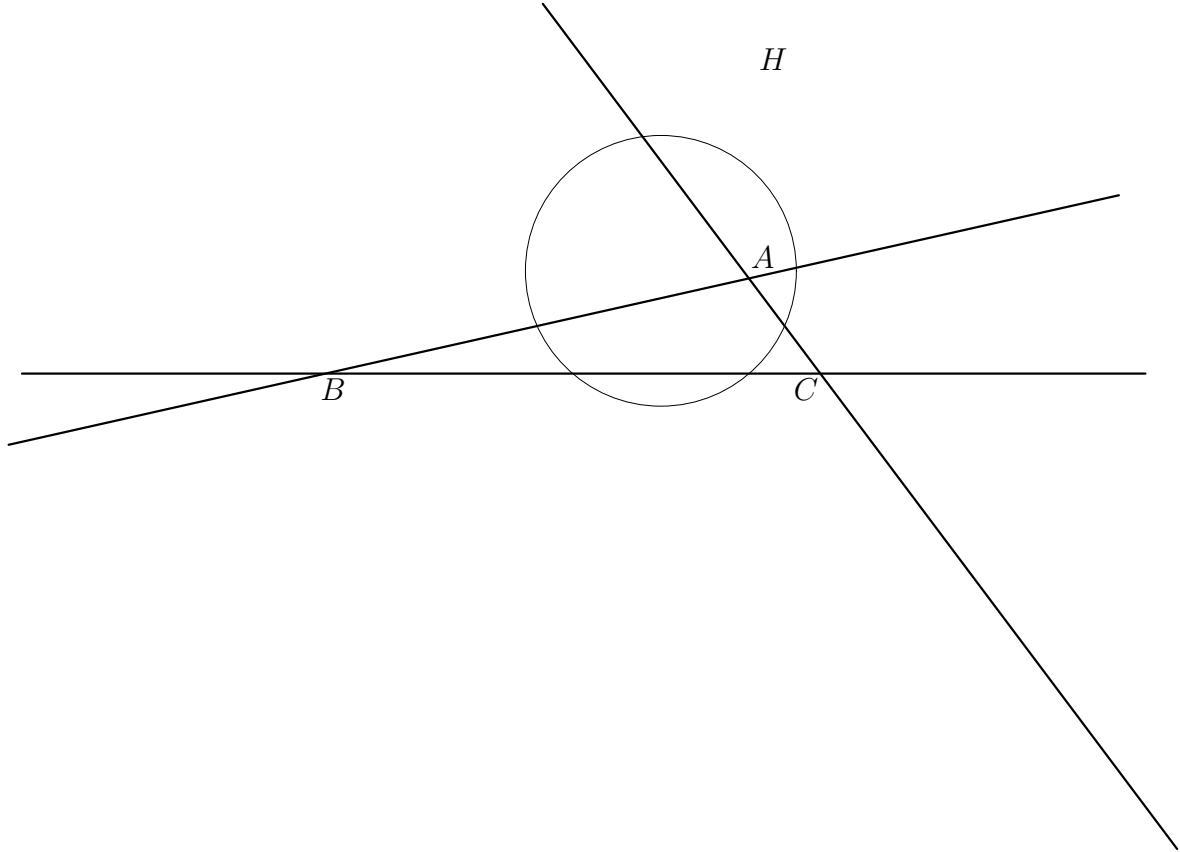
  
%
\caption{Droz-Farny lines touching a hyperbola.  Incomplete.
Add perp bisectors of $AH$, $BH$, $CH$ and tans to 9-pt circle
where it cuts the Euler line, the two asymptotes and some other
examples.}
\label{hyperb}
\end{figure}

For a right triangle the conic degenerates to a point, the
right angle itself, through which the Droz-Farny line
always passes.

The case of the equilateral triangle deserves special mention.
The orthocentre, circumcentre, nine-point centre and incentre
all coincide, the Euler line is indeterminate, and the nine-point
circle is also the incircle, and is the enveloped conic
mentioned in the last of the three theorems above.  We have
special cases of these earlier theorems:

{\bf Theorem.}  Let two perpendicular rays through the centre,
$N$, of an equilateral triangle $ABC$ meet the edges $BC$;
$CA$; $AB$ in $A_1,A_2$; \quad $B_1,B_2$; \quad $C_1,C_2$
respectively.  Then the midpoints of the segments $A_1A_2$,
$B_1B_2$, $C_1C_2$ lie on a line tangent to the incircle.

Conversely,

{\bf Theorem.}  If a tangent to the incircle of an
equilateral triangle $ABC$ meets the edges $BC$;
$CA$; $AB$ in $X$, $Y$, $Z$, respectively, then the circles,
centres $X$, $Y$ and $Z$, passing through the incentre $I$
of the triangle, concur again at a point on the circumcircle
of $ABC$.

\section{Clarifications(??) of Cosmin Pohoata's original Note}

p.3, l.--6.  $\angle BA'C'$ should be $\angle CA'C'$,
though I'm suspicious that the angles (or their supplements)
may depend on the position of $P$ on the circumcircle
(e.g., between $B$ and $C$ or $C$ and $A$ or $A$ and $B$).

p.4, l.3.  $\angle BQC = \angle BAC +\angle CPB$  nowhere
near in my figure!  I suspect that the figure depends on
the relationship between $P$ and the triangle --- e.g.,
whether it's inside the circumcircle or not; \&/or in which
of the seven regions deternined by the triangle edges $P$
lies.  For example, in the next line of the text,
$\angle C'AB'$ is not equal to $\angle BAC$ in my diagram,
but rather to its supplement.

\begin{figure}[h]  
\begin{picture}(425,340)(-275,-175)
\setlength{\unitlength}{0.5pt}
\dashline{7}(30,256)(-562.5,-60)  
\dashline{7}(-8,255)(320,-360)    
\drawline(0,260)(0,-350)          
\drawline(24,258)(-480,-120)      
\drawline(-4.5,260)(135,-360)     
\drawline(13.417,260)(-402.5,-360) 
\drawline(-7,260)(175.5,-271.5) 
\drawline(21.1,260)(-430,-167.55) 
\dashline{7}(-513.4,97.2)(279.5,-121.5)
\put(-161,0){\circle{578}}
\put(110.077,-74.77){\circle{666.92}}
\put(-232.432,19.702){\circle{640.486}}
\thicklines
\drawline(-550,0)(300,0) 
\drawline(324,-360)(-156,280) 
\drawline(-560,-54)(280,135) 
\put(8,230){$P$}
\put(2,80){$A$}
\put(-324,-20){$B$}
\put(34,-20){$C$}
\put(-162,-22){{\small$X$}}
\put(107,-68){{\small$Y$}}
\put(-246,26){{\small$Z$}}
\put(132,7){{\small$A'$}}
\put(-470,7){{\small$A^{\prime\prime}$}}
\put(293,-340){{\small$B'$}}
\put(-96,205){{\small$B^{\prime\prime}$}}
\put(82,100){{\small$C'$}}
\put(-570,-46){{\small$C^{\prime\prime}$}}
\end{picture}
\caption{}
\label{}
\end{figure}

\newpage

3. Concluding Remarks.  Simple enough, right?  No, too
complicated and wrong!

I attempt here to clarify the author's ``reformulation''.
Figure 5 is intended to illustrate that reformulation in
the special case of the Droz-Farny line.

$P$ is taken as the special case of the orthocentre of the
triangle $ABC$ and the lines through $P$ are shown dashed in
Figure 4.  They meet the edges of $ABC$ in $A'$, $B'$, $C'$
and $A^{\prime\prime}$, $B^{\prime\prime}$, $C^{\prime\prime}$.
If the lines through $P$ are perpendicular, then the
reflexion of $PA$, say, in the internal bisector of
$\angle A'PA^{\prime\prime}$ will be the same as the reflexion
in the external bisector --- in fact, even in the general case,
this is still so, so there seems to be no need to specify
``internal''.

Apologies that this is incomplete.  There's a good paper
somewhere, but this isn't it.

\section{Other details concerning C.~P.'s Note}

Abstract perhaps too long and rambling.

Introduction is same as Abstract.

Abstract, l.6 `overridden' (spelling), l.7 `inesthetic' (?)

The author does well to avoid the word ``side'' in connexion
with a triangle, but ``edge'' would be better than `sidelines'.
I don't like ``altitudines'' on p.4.  ``Altitudes'' is the usual
word.

The labels $B_3$ and $C_3$ in the author's Figure 1 should be
interchanged.

p.1, l.--2. ``intersection'' (omit ess)

p.3, l.--6. space in $P, A'$

p.4, l.1. `represents'

p.4, l.5. space in `in triangle'

p.4, l.--3. `are' for `ar'

\section{Appendix}

[Not included in submitted report.]

Note that the pair of perpendicular line through $H$, the orthocentre,
can be replaced by any rectangular hyperbola, centre $H$.

Note that the pair of perpendicular lines together with
the Droz-Farny line and the three edges of the triangle,
all touch the same parabola.  What can one say
about the family of parabolas produced in this way?

It appears (or does it? is this too wild a guess?
I think so.) that the focus of the parabola is the point
that I've called $K$, the reflexion of $H$ in the D-F line.
The directrix is the line through $H$ parallel to the
D-F line.  The axis is $HK$ and the vertex its
intersection with the D-F line.

\newpage

\section{Start again with converse. Generalizations.}

Our generalized (Quadrated \& Twinned) triangle comprises
eight triangles:
\vspace*{-12pt}
\begin{center}
\begin{tabular}{ccc}
triangle & orthocentre & circumcentre \\
{\bf 124} & {\bf 7} & $\bar{\bf 7}$ \\
{\bf 741} & {\bf 2} & $\bar{\bf 2}$ \\
{\bf 472} & {\bf 1} & $\bar{\bf 1}$ \\
{\bf 217} & {\bf 4} & $\bar{\bf 4}$ \\
$\bar{\bf 1}\bar{\bf 2}\bar{\bf 4}$ & $\bar{\bf 7}$ & {\bf 7}  \\
$\bar{\bf 7}\bar{\bf 4}\bar{\bf 1}$ & $\bar{\bf 2}$ & {\bf 2} \\
$\bar{\bf 4}\bar{\bf 7}\bar{\bf 2}$ & $\bar{\bf 1}$ & {\bf 1} \\
$\bar{\bf 2}\bar{\bf 1}\bar{\bf 7}$ & $\bar{\bf 4}$ & {\bf 4} \\
\end{tabular}
\end{center}
\vspace*{-12pt}
where, to avoid special cases, we assume that the triangles
are scalene, that is, neither right-angled nor isosceles,
and the labels are chosen so that {\bf 7} and $\bar{\bf 7}$
are vertices of obtuse angles.

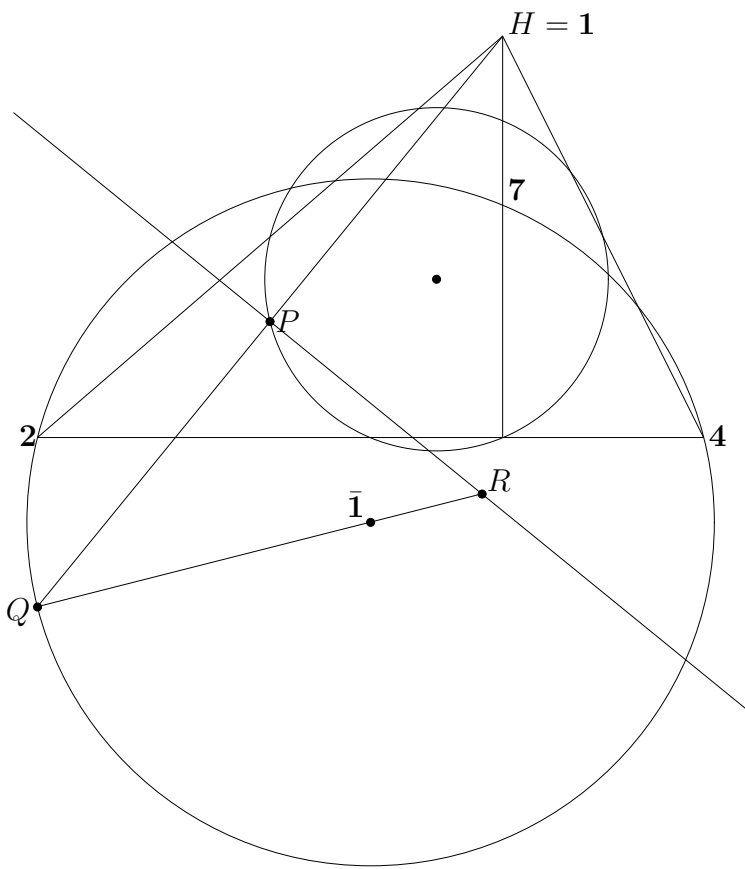
\begin{figure}[h]  
\begin{picture}(400,320)(-250,-230)
\setlength{\unitlength}{1.0pt}
\drawline(-151,-60)(101,-60)(25,92)(-151,-60)    
\drawline(25,92)(25,-60)
\put(0,0){\circle{130}}  
\put(-25,-92){\circle{260}}  
\put(-25,-92){\circle*{3}}
\put(0,0){\circle*{3}}
\put(-34,-90){$\bar{\bf1}$}
\put(27,93){$H={\bf1}$}
\put(27,30){\bf7}
\put(-158,-63){\bf2}
\put(103,-63){\bf4}
\put(-63,-16){\circle*{3}}
\put(-61,-20){$P$}
\put(-151,-124){\circle*{3}}
\put(-163,-129){$Q$}
\put(17.1386,-81.2981){\circle*{3}}
\put(19,-80){$R$}
\drawline(25,92)(-151,-124)            
\drawline(-160,63.037)(120,-165.111)  
\drawline(-151,-124)(17.1386,-81.2981)  
\end{picture}
\caption{Converse of Droz-Farny theorem}
\label{dfenvelope}
\end{figure}

Choose any one of the eight triangles with orthocentre $H$,
say ($H={\bf1}$ in Fig.~9.20), and join $H$ to an arbitrary
point $P$ on the Centre Circle.
Then the line through $P$ perpendicular to $HP$ is a suitable
candidate for a Droz-Farny Line in the following sense.

[Note that if $Q$ is on $HP$ produced, so that $HP=PQ$, then
$Q$ lies on the appropriate circumcircle.]

If the suggested Droz-Farny Line cuts the edges of the
selected triangle in $X$, $Y$, $Z$, and the circles centres
$X$, $Y$, $Z$ passing though $H$ cut the respective edges
in $X_1$, $X_2$ and $Y_1$, $Y_2$ and $Z_1$, $Z_2$, then,
if the subscripts have been chosen suitably, the points
$X_1$, $Y_1$, $Z_1$ are collinear, as are $X_2$, $Y_2$, $Z_2$
and the two lines are perpendicular and pass through $H$.

Wanted: a perspicuous synthetic proof.  It is immediate
that angles $X_1HX_2$, $Y_1HY_2$, $Z_1HZ_2$ are right
angles, since they are angles in semicircles.

It is easy to see\footnote{For example, with orthocentre {\bf1}
of triangle {\bf724} and a point $P$ on the Central Circle,
and corresponding point $Q$ on the circumcircle {\bf724}, the
Droz-Farny line is the perpendicular bisector of $Q${\bf1},
say $PR$, where $R$ is the intersection of this perpendicular
bisector with $Q{\bf\bar1}$ and ${\bf\bar1}$ the circumcentre
of triangle {\bf724}.  Then
${\bf\bar1}R\pm{\bf1}R={\bf\bar1}R\pm QR={\bf\bar1}Q$,
the circumradius of triangle {\bf724}, and diameter of
the Central Circle, so that $R$ lies on a conic, foci
${\bf\bar1}$ and ${\bf1}$ and major (transverse) axis
a diameter of the Central Circle.  Moreover, since
angle ${\bf1}RP$ = angle $QRP$, the Droz-Farny line
$PR$ is a tangent to the conic.}
that the envelope of the Droz-Farny lines
for any of the eight triangles is a conic with centre at the
Centre of the triangle, foci $H$ and $\bar H$, and vertices
at the intersections of the appropriate Euler line with the
Central Circle.  Twins generate the same conic.  Acute
triangles generate an ellipse.  Obtuse triangles generate
three hyperbolas whose asymptotes pass through the Centre
of the triangle and intersect the Central Circle in six
further noteworthy points.

Six of the eight vertices of the triangle lie outside the
Central Circle.  The points of contact of the six pairs of
tangents from these vertices with the Central Circle define
six diameters of the Central Circle which are three pairs
of asymptotes of the three hyperbolas enveloped by Droz-Farny
lines.

\newpage

\begin{figure}[h]
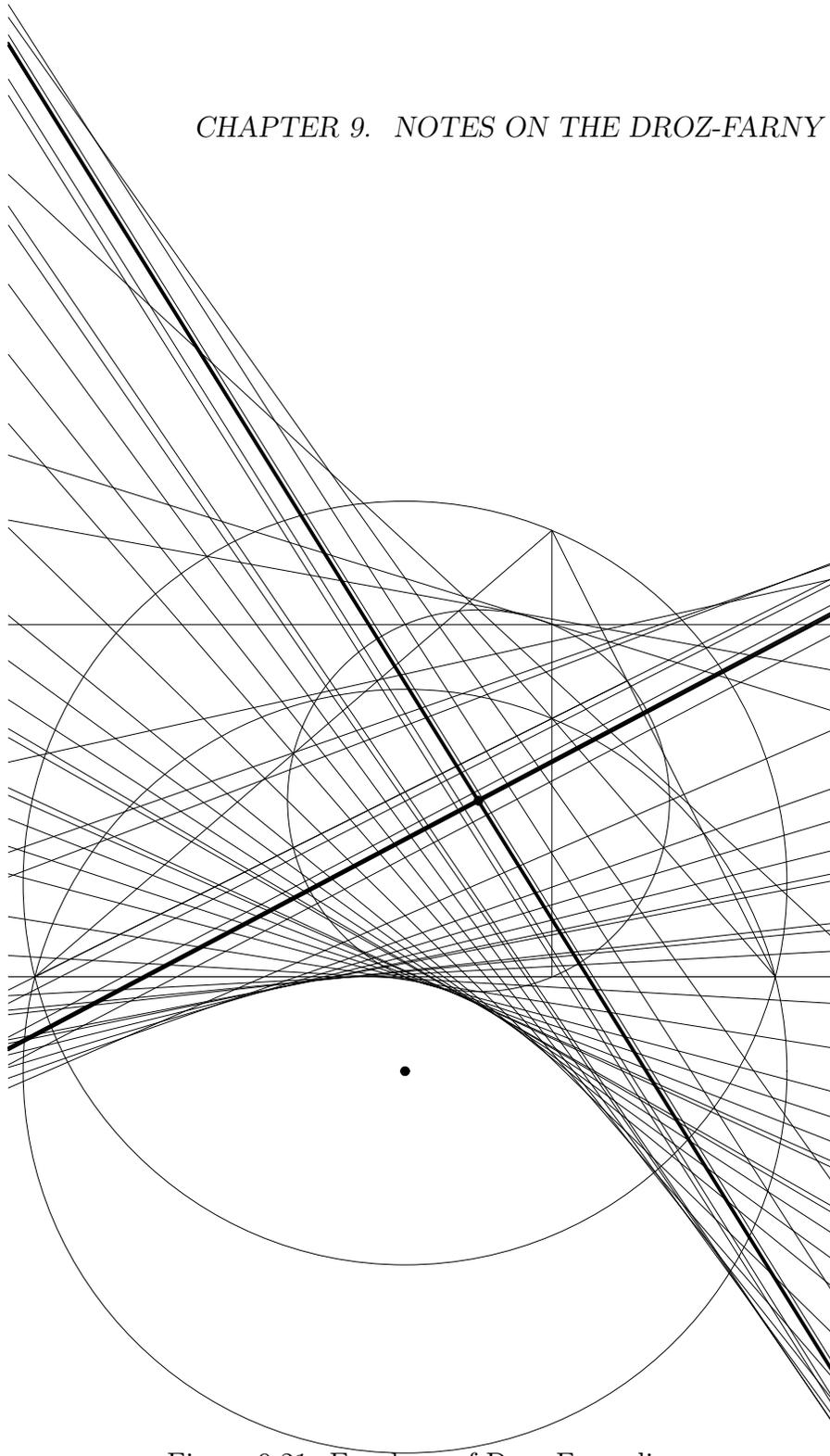
  

\caption{Envelope of Droz-Farny line}
\label{dfenvelope2}
\end{figure}

\end{document}